\documentclass[11pt]{amsart}
\usepackage[T1]{fontenc}
\usepackage{lmodern}
\usepackage{amsmath,amssymb,amsthm,mathtools,mathrsfs,microtype,booktabs,array}
\usepackage[totalwidth=430pt,totalheight=640pt]{geometry}
\usepackage{xcolor}
\definecolor{linkteal}{RGB}{0,104,112}
\usepackage[colorlinks=true,linkcolor=linkteal,citecolor=linkteal,
 urlcolor=linkteal,pdfborder={0 0 0},linktoc=all]{hyperref}
\usepackage{enumitem}
\hypersetup{pdftitle={Asymptotics of Hitchin metrics and modularity of exceptional ALG gravitational instantons}, pdfauthor={Shiyu Cao}, pdfsubject={Hitchin metric asymptotics and exceptional ALG modularity}}
\allowdisplaybreaks
\newcommand{\C}{\mathbb C}
\newcommand{\R}{\mathbb R}
\newcommand{\Z}{\mathbb Z}
\newcommand{\Q}{\mathbb Q}
\newcommand{\Pp}{\mathbb P}
\newcommand{\M}{\mathcal M}
\newcommand{\Acal}{\mathscr A}
\newcommand{\Id}{\operatorname{Id}}
\newcommand{\tr}{\operatorname{tr}}
\newcommand{\diag}{\operatorname{diag}}
\newcommand{\Pic}{\operatorname{Pic}}

\newcommand{\End}{\operatorname{End}}
\newcommand{\Hom}{\operatorname{Hom}}
\newcommand{\ParEnd}{\operatorname{ParEnd}}
\newcommand{\SParEnd}{\operatorname{SParEnd}}
\newcommand{\Herm}{\operatorname{Herm}}
\newcommand{\ad}{\operatorname{ad}}

\newcommand{\pardeg}{\operatorname{pardeg}}
\newcommand{\rank}{\operatorname{rank}}

\newcommand{\Vol}{\operatorname{Vol}}
\newcommand{\Rm}{\operatorname{Rm}}
\newcommand{\ii}{\mathrm i}
\newcommand{\dd}{\mathrm{d}}
\newcommand{\sfm}{\mathrm{sf}}
\newcommand{\app}{\mathrm{app}}
\newcommand{\dec}{\mathrm{dec}}
\newcommand{\ex}{\mathrm{ex}}

\newtheorem{theorem}{Theorem}[section]
\newtheorem{proposition}[theorem]{Proposition}
\newtheorem{lemma}[theorem]{Lemma}
\newtheorem{corollary}[theorem]{Corollary}
\theoremstyle{definition}

\theoremstyle{remark}

\numberwithin{equation}{section}

\title[Hitchin metrics and exceptional ALG instantons]{Asymptotics of Hitchin metrics and modularity of exceptional ALG gravitational instantons}
\author{Shiyu Cao}
\address{Chern Institute of Mathematics and LPMC, Nankai University}
\email{shiyucao@126.com}
\date{}

\begin{document}
\begin{abstract}
We prove asymptotic estimates for the Hitchin metrics on four-dimensional moduli spaces of parabolic $\mathrm{SL}(n,\C)$-Higgs bundles on the three-punctured sphere with flag types corresponding to the Dynkin types $E_6$, $E_7$ and $E_8$. For fixed weights in the alcove interiors and fixed scalar graded residues, the Hitchin metric converges in $C^\infty$ to the semiflat metric along the end. The decay is polynomial in general and exponential when all graded residues vanish. For every admissible parameter with weights in the closed alcoves, the moduli space with its Hitchin metric is a complete nonflat ALG gravitational instanton of order two. Extending the period calculation of Dimakis--Fredrickson--Yae to the boundary and allowing compactly supported changes of marking, we obtain every point of the normalized period domain. It follows from the ALG Torelli theorem that every ALG gravitational instanton with cone parameter $1/3$, $1/4$ or $1/6$ is isometric, up to homothety and hyperk\"ahler rotation, to one of these moduli spaces with its Hitchin metric.
\end{abstract}
\maketitle
\tableofcontents
\clearpage

\section{Introduction}

On the three-punctured sphere, parabolic $\mathrm{SL}(n,\C)$-Higgs bundles of ranks $n=3,4,6$ with the flag types of $E_6,E_7,E_8$ have four-dimensional moduli spaces. These flag types include partial flags. We study the Hitchin metric at infinity in the Hitchin base with the parabolic weights and scalar graded residues fixed. For weights in the interiors of the alcoves, our local construction at the punctures gives a comparison with the semiflat metric even when graded residues coincide or vanish. When semistability equals stability, the moduli spaces are complete ALG gravitational instantons of order two. The construction extends to every admissible parameter in the closed weight alcoves and realizes every exceptional ALG gravitational instanton up to homothety and hyperk\"ahler rotation.

On compact Riemann surfaces, Hitchin constructed the hyperk\"ahler metric on the moduli space of solutions to his equations. Hitchin, Donaldson, Corlette and Simpson established the correspondences between solutions of these equations, polystable Higgs bundles and reductive representations~\cite{Hitchin1987,Donaldson1987,Corlette1988,Simpson1992}. Simpson extended the correspondence to noncompact curves by incorporating filtrations and growth conditions at the punctures~\cite{Simpson1990}. Biquard developed the logarithmic correspondence in the presence of a smooth divisor~\cite{BiquardLogarithmic}, and Biquard--Boalch constructed the wild correspondence and hyperk\"ahler moduli spaces for prescribed irregular data~\cite{BiquardBoalch}.

Over the locus of smooth spectral curves, the abelian fibers of the Hitchin fibration and the special K\"ahler geometry of the base determine the semiflat metric~\cite{FreedSpecialKahler}. Gaiotto--Moore--Neitzke proposed a detailed description of the corrections to this metric in terms of wall-crossing and Bogomolny--Prasad--Sommerfield (BPS) invariants~\cite{GMN}. They predicted exponential decay along suitable directions escaping to infinity away from the discriminant. Establishing this prediction requires estimates for harmonic metrics and for the infinitesimal deformations whose $L^2$ pairings define the metric on the moduli space.

Mazzeo--Swoboda--Weiss--Witt developed a gluing construction for rank-two Higgs bundles with simple zeros of the determinant~\cite{MSWWEnds}. A limiting configuration solves the decoupled equations away from those zeros, and rotationally symmetric local solutions resolve its singularities. Estimates for the linearized equations then give genuine solutions close to the approximate ones. Mazzeo--Swoboda--Weiss--Witt subsequently analyzed infinitesimal deformations and obtained polynomial asymptotic expansions for the Hitchin metric relative to the semiflat metric~\cite{MSWWMetric}.

Dumas--Neitzke proved exponential decay of the difference between the Hitchin and semiflat metrics on the rank-two Hitchin section by expressing the local contribution near a simple zero through a variation induced by a holomorphic flow~\cite{DumasNeitzke}. The resulting cancellation explains why estimates for the individual terms in the metric pairing can give a weaker rate of decay than an estimate for their sum. Fredrickson constructed higher-rank local models and approximate solutions on a generic spectral locus~\cite{FredricksonEnds}, and combined gluing with this cancellation to prove exponential convergence of the moduli-space metric~\cite{FredricksonMetric}. In higher rank, Fredrickson's hypotheses restrict the comparison to simple eigenvalue crossings.

Another approach develops quantitative estimates directly for harmonic bundles and their deformation complexes. Mochizuki proved exponential decoupling for generically regular semisimple Higgs fields away from the spectral discriminant, together with convergence to limiting configurations in rank two~\cite{MochizukiAsymptotics}. Mochizuki--Szab\'o subsequently established exponential convergence in higher rank when the Higgs bundle comes from a line bundle on the normalization of the spectral curve. Their estimates are uniform for families satisfying their hypotheses~\cite{MochizukiSzabo}. Mochizuki also proved exponential decay of the difference between the Hitchin and semiflat metrics along scaling rays throughout the smooth spectral-curve locus for stable degree-zero Higgs bundles on a compact curve~\cite{MochizukiMetric}. In this last comparison, the harmonic representatives of tangent classes and their pairings are estimated for arbitrary ramification of a smooth spectral curve.

For singular spectral curves, He--Horn--Li constructed limiting configurations and proved exponential convergence under a local fiducial condition at spectral singularities. Their metric results concern the restrictions to specified Hitchin subintegrable systems and selected subsets over quadratic-differential strata~\cite{HeHornLi}. Sagman--Smillie extended local estimates and asymptotic decoupling to arbitrary Higgs bundles, with applications to the Wentzel--Kramers--Brillouin (WKB) problem and high-energy harmonic maps~\cite{SagmanSmillie}. Their analysis described the geometry away from the critical set even when the Higgs field has a persistent nilpotent part.

For parabolic rank-two Higgs bundles, Fredrickson--Mazzeo--Swoboda--Weiss adapted the gluing and tangent-space analysis to simple poles~\cite{FMSW}. Their local constructions treated simple branch points, strongly parabolic points, and weakly parabolic points, and yielded exponential convergence to the semiflat metric along generic scaling rays, with the tangent variations specified in~\cite[Theorems~7.1 and~7.2]{FMSW}. On the four-punctured sphere, they also established the ALG geometry of the strongly parabolic moduli space and the predicted optimal exponential rate. Scaling preserves the moduli problem in the strongly parabolic case. In the weakly parabolic case the comparison in~\cite[Theorem~7.2]{FMSW} takes place in the moduli spaces for the equations with scaling parameter $t$. Rewriting these as Hitchin equations replaces the Higgs field by $t\varphi$ and scales its residue eigenvalues.

Chen--Li studied generic curves tending to infinity inside moduli spaces of rank-two irregular Higgs bundles over the projective line with prescribed singularity data~\cite{ChenLi}. The branch points move along these curves, and the local rescalings depend on their changing positions. Near a weakly parabolic point they used a further model obtained from a harmonic bundle on the projective line with one weakly parabolic pole and one twisted irregular pole. Under their assumption of trivial parabolic structure at every weakly parabolic point, the relative metric error decays faster than every inverse power of the curve parameter. In the absence of weakly parabolic points they obtained a stretched exponential bound. For the higher-rank families considered here, we also need uniform estimates for the local models and their deformations when several residue eigenvalues coalesce.

In real dimension four, these metric questions connect Hitchin moduli spaces with gravitational instantons. Tian--Yau constructed complete Ricci-flat K\"ahler metrics on open varieties~\cite{TianYau}. Hein applied their method to complements of fibers in rational elliptic surfaces, obtaining several asymptotic geometries, including ALG ends when the removed singular fiber has finite monodromy~\cite{Hein}. Chen--Chen classified the possible asymptotic geometries under the assumption that curvature decays faster than quadratically and established compactifications by elliptic surfaces in the ALG case~\cite{ChenChenI}. Sun--Zhang's analysis of collapsing hyperk\"ahler four-manifolds gave the corresponding classification for metrics with square-integrable curvature, including the ALG$^*$ and ALH$^*$ geometries~\cite{SunZhang}.

The gauge-theoretic constructions of Fredrickson--Mazzeo--Swoboda--Weiss and Chen--Li gave explicit connections between Hitchin fibrations and these asymptotic models~\cite{FMSW,ChenLi}. The semiflat metric makes this connection accessible through the periods and monodromy of the spectral curves. Once an ALG structure has been established with the required decay, global classification results relate its hyperk\"ahler forms to cohomological data. In particular, Chen--Viaclovsky--Zhang proved a Torelli theorem for ALG gravitational instantons of order two relative to a common asymptotic coordinate system, with the period conditions expressed using integral classes of self-intersection $-2$~\cite{ChenViaclovskyZhang}.

Boalch posed the modularity question in Problem~1.6 of the 2015 AIM workshop \emph{Spectral data for Higgs bundles}~\cite{Boalch2015}. The metric formulation for ALG and ALG$^*$ gravitational instantons appears in~\cite[Conjecture~1.2]{FMSWTorelli}. The exceptional ranks and pole data considered here occur in~\cite[Section~4.1]{BoalchH3}. Fredrickson--Mazzeo--Swoboda--Weiss computed the Torelli parameters for rank-two parabolic Hitchin moduli spaces on the four-punctured sphere and proved that the image of the period map is open and dense~\cite[Theorem~9.2 and Remark~9.28]{FMSWTorelli}. Dimakis--Fredrickson--Yae computed the integral hyperk\"ahler periods for the exceptional families with weights in the interiors of the alcoves and parameters in their generic locus. They realized the ALG period domain outside a subset of real codimension at least two~\cite[Theorems~6.9 and~6.11]{DFY}. We extend their calculation to admissible boundary weights and prove completeness and quantitative asymptotics at every admissible parameter in the closed alcoves. With compactly supported changes of marking we then prove that the marked period map is surjective and deduce modularity.

\subsection{The moduli spaces and their asymptotic geometry}
\label{sec:parabolic-data}
Throughout, $X=\Pp^1$ and $D=\{0,1,\infty\}$. We consider the ranks and flag types
\begin{equation}\label{eq:types}
 (n;e_0,e_1,e_\infty)=(3;3,3,3),\quad (4;4,4,2),\quad (6;6,3,2).
\end{equation}
The corresponding root systems have ranks $\ell=6,7,8$, respectively.
A parabolic bundle consists of a rank-$n$ holomorphic vector bundle $E$ with weighted flags $E_p=F_p^0\supset\cdots\supset F_p^{e_p}=0$ at $p\in D$. Each graded quotient has dimension $m_p=n/e_p$. We fix an isomorphism $\det E\simeq\mathcal O_X$ and weights satisfying
\begin{equation*}
 \sum_{j=0}^{e_p-1}\alpha_{p,j}=0,\qquad
 \alpha_{p,0}<\cdots<\alpha_{p,e_p-1}<\alpha_{p,0}+1.
\end{equation*}
The space of these weights is $\Acal$. Here $K_X$ denotes the canonical line bundle. We consider trace-free logarithmic Higgs fields $\Phi:E\to E\otimes K_X(D)$ whose residues preserve the flags. On each graded quotient the residue acts by
\begin{equation}\label{eq:graded-residue}
 \operatorname{Gr}_j\operatorname{Res}_p\Phi
   =\sigma_{p,j}\Id,
 \qquad \sum_{j=0}^{e_p-1}\sigma_{p,j}=0.
\end{equation}
The complex numbers $\sigma_{p,j}$ may coincide. The condition $\sigma=0$ gives strongly parabolic Higgs fields. For the flag types in~\eqref{eq:types}, the spaces of residue parameters have complex dimensions $6$, $7$ and $8$, respectively.

For a saturated subbundle $F\subset E$, $k_{p,j}$ denotes the dimension of the induced graded quotient. The parabolic degree of the subbundle is
\begin{equation*}
 \pardeg F=\deg F+\sum_{p,j}k_{p,j}\alpha_{p,j}.
\end{equation*}
The pair is stable if this degree is negative for every nonzero proper $\Phi$-invariant saturated subbundle, and semistable if these degrees are at most zero. For suitable residues, semistability equals stability even at weights lying on walls for strongly parabolic Higgs bundles.

We write $\mathscr B$ for the vector space of centered scalar graded residues in~\eqref{eq:graded-residue}, and $\M_{\alpha,\sigma}$ for the coarse moduli space of semistable objects with these data.

For each fixed residue parameter, the Hitchin base is an affine line with coordinate $t$. The fibers of the Hitchin morphism over $|t|$ sufficiently large are the degree-$n$ Picard varieties of smooth normalized spectral curves. They consist of stable objects for every allowed weight. We call the inverse image of this region the end of the Hitchin fibration. We use the normalization because the spectral curve can be singular over punctures with partial flags.

\begin{theorem}\label{thm:main-geometry}
Here $\M_{\alpha,\sigma}$ denotes the moduli space for any of the three rank and flag types in Section~\ref{sec:parabolic-data}. The weights are centered, strictly increasing and of width less than one, and the scalar graded residues satisfy the trace condition~\eqref{eq:graded-residue}. We assume that every semistable object with these data is stable.

Then $\M_{\alpha,\sigma}$ is a nonempty smooth connected rational surface. Its Hitchin morphism is a projective surjection onto an affine line and has connected fibers. For $t\in U_{\sigma}$, with $U_{\sigma}$ as in Proposition~\ref{prop:normalized-spectral-family}, the Hitchin fiber is $\Pic^n(Y_{t,\sigma})$, the Picard variety of the normalized genus-one spectral curve $Y_{t,\sigma}$. The harmonic metrics inducing the unit metric on the trivialized determinant line bundle define the complete Hitchin hyperk\"ahler metric on this surface.

There is a smooth projective rational elliptic surface $f:S\to\Pp^1$ with a section, relatively minimal over $\Pp^1$, for which
\begin{equation*}
 \M_{\alpha,\sigma}
       \simeq S\setminus F_\infty,
 \qquad F_\infty=f^{-1}(\infty),
 \qquad K_S\simeq\mathcal O_S(-F_\infty).
\end{equation*}
This is an algebraic isomorphism over the affine Hitchin base. The fiber $F_\infty$ has Kodaira type $\mathrm{IV}$, $\mathrm{III}$ or $\mathrm{II}$ for $n=3$, $4$ or $6$, respectively. The graded residues may coincide or vanish throughout this statement.
\end{theorem}

We prove this theorem at the end of Section~\ref{subsec:metric-completeness}.

The limiting ALG model depends only on the rank and flag type, as recorded in Table~\ref{tab:main-ALG-parameters}. We use the metric normalization in~\eqref{eq:tangent-trial-norm}, defined by the matrix trace.

\begin{table}[htbp]
\centering
\begin{tabular}{c|c|c|c}
 $n$ & $(e_0,e_1,e_\infty)$ & $\tau_\infty$ & Cone angle\\
 \hline
 $3$ & $(3,3,3)$ & $\exp(2\pi\ii/3)$ & $2\pi/3$\\
 $4$ & $(4,4,2)$ & $\ii$ & $\pi/2$\\
 $6$ & $(6,3,2)$ & $\exp(2\pi\ii/3)$ & $\pi/3$
\end{tabular}
\caption{The ALG models. Each limiting torus has area $4\pi^2$.}
\label{tab:main-ALG-parameters}
\end{table}

We use the notation $b=\operatorname{Im}\tau_\infty$, $\zeta=\exp(2\pi\ii/n)$ and
$\Lambda_\infty=(2\pi/\sqrt b)(\Z+\tau_\infty\Z)$. The model is
\begin{equation}\label{eq:main-ALG-model}
 \bigl(\{|U|>R_0\}\times\C/\Lambda_\infty\bigr)
 \big/\bigl((U,V)\sim(\zeta U,\zeta^{-1}V)\bigr),
 \qquad R=|U|,
\end{equation}
with
\begin{equation*}
 \begin{gathered}
 g_0=|\dd U|^2+|\dd V|^2,\qquad
 \omega_0=\frac{\ii}{2}
             (\dd U\wedge \dd\bar U+\dd V\wedge \dd\bar V),\\
 \Omega_0=\dd U\wedge \dd V,\qquad
 (\omega_{1,0},\omega_{2,0},\omega_{3,0})
       =(\omega_0,\operatorname{Re}\Omega_0,\operatorname{Im}\Omega_0).
 \end{gathered}
\end{equation*}
The smooth end coordinates are specified in~\eqref{eq:semiflat-ALG-coordinates}. On the cover $t=\lambda^n$, we write $T=|\lambda|$. The radius $R$ is a fixed positive multiple of $T$, and the angular coordinates $\theta_A,\theta_B$ of the unitary characters have period $2\pi$.

\begin{theorem}\label{thm:main-asymptotics}\label{thm:zero-residue-exponential}
For each rank and flag type in Table~\ref{tab:main-ALG-parameters} with parameters in $\Acal\times\mathscr B$, we denote the Hitchin metric on the end of the fibration by $g_H$ and write $(\omega_1,\omega_2,\omega_3)=(\omega_I,\omega_J,\omega_K)$. We use the smooth end coordinates~\eqref{eq:semiflat-ALG-coordinates} and the model~\eqref{eq:main-ALG-model}.

\begin{enumerate}[label=\textup{(\roman*)}]
\item For every nonempty compact set $K\subset\Acal\times\mathscr B$ and every $0<\nu<1/4$, there exist $R_0=R_0(K,\nu)>0$ and constants $C_{k,K,\nu}$ for all integers $k\geq0$ such that
\begin{equation}\label{eq:main-ALG-rate}
 \bigl|(\nabla^{g_0})^k(g_H-g_0)\bigr|_{g_0}
 +\sum_{j=1}^3
       \bigl|(\nabla^{g_0})^k(\omega_j-\omega_{j,0})\bigr|_{g_0}
       \leq C_{k,K,\nu}R^{-\nu-k},
 \qquad R\geq R_0.
\end{equation}
Thus these Hitchin metrics have ALG ends of every order $0<\nu<1/4$.
\item If $\sigma=0$, then for every nonempty compact set $K\subset\Acal$, there exist $c=c(K)>0$, $R_0=R_0(K)>0$ and constants $C_{k,K}$ for all integers $k\geq0$ such that
\begin{equation}\label{eq:main-zero-residue-rate}
 \bigl|(\nabla^{g_0})^k(g_H-g_0)\bigr|_{g_0}
 +\bigl|(\nabla^{g_0})^k(\omega_I-\omega_0)\bigr|_{g_0}
       \leq C_{k,K}e^{-cR},
 \qquad R\geq R_0.
\end{equation}
Moreover,
\begin{equation}\label{eq:zero-residue-symplectic}
 \omega_J+\ii\omega_K=\Omega_0=\dd U\wedge \dd V.
\end{equation}
\end{enumerate}
In each part the estimates are uniform over $K$, the Picard fibers and all arguments of $\lambda$, where $t=\lambda^n$. The radius $R_0$ and, in \textup{(ii)}, the rate $c$ are independent of $k$. All derivatives are taken along the moduli space, with the weights and residues fixed.
\end{theorem}

The parameter sets in \textup{(i)} may contain coinciding graded residues. Both assertions hold on the stable end even when strictly semistable objects occur elsewhere in the moduli space. We prove the two assertions in Section~\ref{sec:metric-derivatives}.

\begin{theorem}\label{thm:main-global-ALG}\label{prop:order-two-ALG}
For each parameter $\xi=(\alpha,\sigma)\in\Acal\times\mathscr B$ for which semistability equals stability, the whole moduli space $\M_\xi$ with its Hitchin metric is a connected, nonflat, complete ALG gravitational instanton of order two, with the model~\eqref{eq:main-ALG-model}. There is a smooth chart $\Psi_\xi$ on the end, obtained by composing the chart defined by the spectral curves with a global diffeomorphism isotopic to the identity, such that, for every integer $k\geq0$,
\begin{equation}\label{eq:main-order-two-rate}
 \bigl|(\nabla^{g_0})^k(\Psi_\xi^*g_H-g_0)\bigr|_{g_0}
 +\sum_{j=1}^3
   \bigl|(\nabla^{g_0})^k
        (\Psi_\xi^*\omega_j-\omega_{j,0})\bigr|_{g_0}
       =O(R^{-2-k}).
\end{equation}
The corrected chart and the constants may depend on $\xi$.

\end{theorem}
We prove this theorem in Section~\ref{sec:ALG-order-two}.

\subsection{Boundary weights, periods and modularity}
\label{sec:main-modularity}

For modularity we also allow weights in the closed alcoves. Labels with equal real weights belong to the same graded piece. We prescribe the conjugacy class of the semisimple residue on each such piece. The numerical lattice~$\Lambda$ in~\eqref{eq:numeric-def} records ranks, degrees and labelled dimensions. Its primitive radical generator~$F$ corresponds to a smooth Hitchin fiber. A parameter is admissible if the formal period triple~\eqref{eq:formal-periods} is nonzero on every integral class of square $-2$. We denote the admissible set, including boundary weights, by~$\mathscr U$, as defined in Section~\ref{sec:closed-parameters}.

\begin{theorem}\label{thm:boundary-ALG}
For every $(\alpha,\sigma)\in\mathscr U$, semistability equals stability, and $\M_{\alpha,\sigma}$ is a nonempty smooth connected rational surface with the elliptic compactification of Theorem~\ref{thm:main-geometry}. Its Hitchin metric is a complete nonflat ALG gravitational instanton with the model~\eqref{eq:main-ALG-model}. In the spectral coordinates~\eqref{eq:semiflat-ALG-coordinates}, for every $0<\nu<1/4$ there are $R_0>0$ and constants $C_{k,\nu}$ for all integers $k\geq0$ such that
\begin{equation}\label{eq:boundary-ALG-all-orders}
 \bigl|(\nabla^{g_0})^k(g_H-g_0)\bigr|_{g_0}
 +\sum_{j=1}^3\bigl|(\nabla^{g_0})^k
                    (\omega_j-\omega_{j,0})\bigr|_{g_0}
 \leq C_{k,\nu}R^{-\nu-k},\qquad R\geq R_0.
\end{equation}
The constants may depend on $(\alpha,\sigma)$ and are uniform over the Picard fibers and all arguments of~$\lambda$. The radius~$R_0$ is independent of~$k$. A global diffeomorphism isotopic to the identity corrects this spectral chart to an ALG chart of order two.
\end{theorem}

The boundary theorem is proved in Section~\ref{sec:boundary}. By Proposition~\ref{prop:hyperkahler-periods}, the integral period formula of Dimakis--Fredrickson--Yae holds for every parameter in~$\mathscr U$. The normalized period domain~$\mathcal D$ consists of homomorphisms $P:\Lambda\to\R^3$ with $P(F)=(4\pi^2,0,0)$ and $P(C)\ne0$ for every $C^2=-2$.

We choose a generic reference surface~$\M_*$ with vanishing graded residues. Its compact fiber has type $\mathrm{IV}^*$, $\mathrm{III}^*$, or $\mathrm{II}^*$. With the marking of second homology by~$\Lambda$ in~\eqref{eq:reference-integral-marking}, we take $G$ to be the group of compactly supported diffeomorphisms generated by generalized Dehn twists about the components of this fiber.

\begin{theorem}\label{thm:full-period-realization}
The family over~$\mathscr U$ admits integral markings compatible with a single ALG chart of decay order two on~$\M_*$. Its marked period map satisfies
\begin{equation}\label{eq:surjective-marked-periods}
 \bigcup_{g\in G}g^*P(\mathscr U)=\mathcal D,
 \qquad (g^*P)(C)=P(g_*C).
\end{equation}
\end{theorem}

The compact support of the changes of marking preserves the coordinates needed for the ALG Torelli theorem. We obtain the following modularity theorem.

\begin{theorem}\label{thm:modularity}
Every ALG gravitational instanton with cone parameter $1/3$, $1/4$, or $1/6$ is, up to homothety and hyperk\"ahler rotation, isometric to a member of the rank-$3$, rank-$4$, or rank-$6$ family, respectively, equipped with its Hitchin $L^2$ metric. With the stated homothety and rotation, the isometry preserves the ordered hyperk\"ahler triple.
\end{theorem}

Both theorems are proved in Section~\ref{sec:periods}.

\subsection{Outline of the proof}

The normalized spectral curve describes the local Higgs bundles near the punctures and determines the semiflat metric at infinity. At a puncture of flag length $e$, the change of variables $\xi=\lambda^e x$ gives a direct sum of models whose spectral equations are $P(w)=c_\nu\xi$, where $P(w)=\prod_j(w-\sigma_j)$ and the constants $c_\nu$ are distinct. Each normalized component is a smooth $w$-line even when roots of $P$ coincide. The products $\prod_{i<j}(w-\sigma_i)$ recover the flags with their graded multiplicities, without any separation condition on the branch points of the local model.

By the existence theorem for good filtered Higgs bundles, each model admits a harmonic metric. Its asymptotics at infinity are uniform on compact parameter sets, so we can match it to the flat metric on the spectral line bundle, fixing the relative positive constants of the blocks when the flag is partial. On the shrinking annuli used for gluing, the perturbation of the inner model and the difference between the diagonal entries produce an error in Hitchin's equation that decays polynomially. Coercivity and the logarithmic trace inequality bound the ratio of the approximate metric to the harmonic metric with the same unit determinant in terms of this error. This estimate is uniform even though the initial comparison constant for the two metrics may depend on the individual Higgs bundle.

For horizontal and vertical tangent classes, we construct Dolbeault representatives whose $L^2$ norms approach the corresponding Hodge norms on the spectral curve. Orthogonal projection onto the harmonic representatives gives the upper bound for $g_H$. By the trace pairing, $g_H$ and $g_{\sfm}$ have the same holomorphic symplectic form and hence the same volume form. In complex dimension two the reciprocal lower bound follows from the upper bound. When $\sigma=0$, the local harmonic metrics approach their diagonal limits exponentially by cyclic symmetry. The algebraic terms in the tangent pairings cancel by a local Green identity. Harmonic projection and the volume identity then yield exponential decay of the difference between the two metrics.

Both metrics are Ricci-flat K\"ahler metrics in the same complex structure, so interior estimates for their $C^0$ difference give bounds for its covariant derivatives. To obtain algebraic decay, we rescale holomorphic covering charts of radii comparable to $R$. On returning to the original coordinates, each derivative contributes a factor $R^{-1}$. For exponential decay, charts of fixed radius give the same exponential rate for every derivative order. The derivatives are taken on the moduli space at fixed weights and residues, and the bounds are uniform on compact parameter sets.

Properness of the Hitchin morphism and the description of its normalized spectral curves give a connected elliptic compactification. From the initial metric comparison we obtain completeness, quadratic volume growth and square-integrable curvature. The section and the monodromy at infinity imply that the fundamental group is finite. With a compact elliptic fiber the metric therefore cannot be flat. By the theorem on improvement of the ALG decay rate, for each parameter there is an ALG chart in which the decay has order two. The original estimates in spectral coordinates are uniform on compact parameter sets.

At an admissible boundary parameter, labels with equal weights modulo $\Z$ have distinct complex residues. The eigenspaces of the graded residue determine the filtration, and the joint normal operator has a positive spectral gap on the complement of the constant centralizer. We use this bound in the deformation theory at the boundary. Choosing local identifications that match the graded residues in the polynomial models gives the estimate and the horizontal representatives needed to compare the metrics there. The integral periods extend from $\mathscr U^\circ$ to $\mathscr U$ by smooth variation on compact subsets of the moduli space. Generalized Dehn twists about compact spheres act by affine root reflections. After these changes of marking, reduction to the closed alcoves realizes every normalized period. We apply the ALG Torelli theorem in compatible coordinates at infinity.

In Sections~\ref{sec:preliminaries}--\ref{sec:approximate-metrics} we construct the spectral curves and harmonic metrics used in Section~\ref{sec:L2-geometry} to compare the two metrics. In Sections~\ref{sec:elliptic-ALG-geometry}--\ref{sec:periods} we establish the complete ALG geometry, extend the metrics to admissible boundary weights, and prove surjectivity of the marked period map and modularity. The appendices contain the coefficient calculations, the algebraic moduli construction, and the local weighted analysis.

\subsection*{Declaration of AI-Use}
The mathematical idea and blueprint of this work belongs to the author. In preparing the manuscript, the author used GPT-6 astra to assist with literature searches, language editing and  routine computations. The author  takes full responsibility for the paper in its entirety.

\section{Preliminaries}
\label{sec:preliminaries}

We describe the filtered lattices and deformation complex for the parabolic data of Section~\ref{sec:parabolic-data}, then specify the normalizations of the Hitchin and semiflat metrics. The Hitchin metric uses harmonic representatives of tangent classes on the punctured curve. For the semiflat metric we use the Hodge pairings on the normalized spectral curves. In Section~\ref{sec:hitchin-base} we describe $h^{-1}(U_{\sigma})$ as a relative Picard variety and compute the semiflat metric near infinity.

\subsection{Parabolic Higgs bundles}
\label{subsec:parabolic-preliminaries}

With the parabolic data and stability condition of Section~\ref{sec:parabolic-data}, we write $X^*=X\setminus D$. Locally at $p$, with coordinate $x$, the inverse images of the flags under $E\to E/xE$ form a chain
\begin{equation}\label{eq:periodic-chain}
 E=E^0\supset E^1\supset\cdots\supset E^{e_p}=xE,\qquad
 E^{j+e_p}=xE^j.
\end{equation}
The successive quotients have length $m_p$. Associating the weight $\alpha_{p,j}+k$ to $E^{j+ke_p}$ gives the equivalent description by filtered lattices in the meromorphic bundle $E(*D)$.

A logarithmic Higgs field is weakly parabolic when each residue preserves the flag, and strongly parabolic when $\operatorname{Res}_p\Phi(F_p^j)\subset F_p^{j+1}$. With the scalar graded residues~\eqref{eq:graded-residue}, the latter condition is equivalent to $\sigma=0$.

\subsection{The moduli space}\label{sec:moduli-deformations}

Here $\ParEnd_0E$ denotes the sheaf of trace-free endomorphisms preserving the flags, and $\SParEnd_0E$ the subsheaf of those sending $F_p^j$ into $F_p^{j+1}$. With the weights, determinant and graded residues fixed, the deformation complex is
\begin{equation}\label{eq:deformation-complex-weak}
 \mathcal C^\bullet:\quad
 \ParEnd_0E\xrightarrow{[\Phi,\,\cdot\,]}
 \SParEnd_0E\otimes K_X(D).
\end{equation}
Its Higgs variations have zero graded residue. Scalar residue actions ensure that the commutator has this property even when some prescribed residues coincide. 

\begin{proposition}\label{prop:local-deformations}
At every stable pair with the prescribed determinant and scalar graded residues, the complex~\eqref{eq:deformation-complex-weak} is self-dual under Serre duality and satisfies
\begin{equation*}
 \mathbb H^0(\mathcal C^\bullet)=\mathbb H^2(\mathcal C^\bullet)=0,
 \qquad \dim_\C\mathbb H^1(\mathcal C^\bullet)=2.
\end{equation*}
The automorphisms of the parabolic Higgs bundle preserving the determinant trivialization form $\mu_n$. The deformation functor of the parabolic Higgs bundle on local Artin $\C$-algebras, with a specified identification of the special fiber and with determinant trivialization and scalar graded residues fixed, is represented by the smooth formal disc $\operatorname{Spf}\C[[s_1,s_2]]$. The group $\mu_n$ acts trivially on this functor.
\end{proposition}

\begin{proof} We give the deformation-theoretic calculation in Appendix~\ref{proof:prop:local-deformations}. \end{proof}

\begin{proposition}\label{prop:joint-walls}
Every strictly semistable pair has a nonzero proper invariant saturated subbundle of rank $r$ and degree $d$ for which
\begin{equation}\label{eq:joint-walls}
 d+\sum_{p,j}k_{p,j}\alpha_{p,j}=0,
 \qquad \sum_{p,j}k_{p,j}\sigma_{p,j}=0,
 \qquad |d|<\frac{3r(n-r)}{n}.
\end{equation}
Here $0\leq k_{p,j}\leq m_p$ and $\sum_jk_{p,j}=r$ at each puncture. The finitely many nonempty loci defined by these equations have real codimension three in the joint parameter space. Their complement is open and dense, and semistability equals stability there.

For any fixed $\sigma$ and any $\alpha\in\Acal$ for which semistability equals stability, there are arbitrarily close rational weights $\alpha'\in\Acal$ defining the same stable and semistable objects with those residues. The equality also holds for the fiberwise stability conditions in families.
\end{proposition}

\begin{proof}
A strictly semistable pair has an invariant saturated subbundle $F$ of parabolic degree zero. The trace of its induced Higgs field is a meromorphic one-form on $X$ whose residue at $p$ is $\sum_jk_{p,j}\sigma_{p,j}$. The second equation in~\eqref{eq:joint-walls} follows by summing these residues.

At each puncture we write $S_p=\sum_jk_{p,j}\alpha_{p,j}$. The selected $r$ weights have mean $S_p/r$, and the remaining $n-r$ weights have mean $-S_p/(n-r)$. Since both means lie in an interval of length less than one, $|S_p|<r(n-r)/n$, and the bound on the degree follows from $d=-\sum_pS_p$. Only finitely many integers $r,d,k_{p,j}$ can occur. For each such array, the linear functional $\sum_{p,j}k_{p,j}\alpha_{p,j}$ is nonzero on the centered weight space because the flag at zero is full and its coefficient vector has $r$ entries equal to one and $n-r$ equal to zero. Where the weight equation meets $\Acal$ it cuts out a real hyperplane. The independent complex residue equation imposes two more real conditions. The union of these common zero sets over the finitely many arrays is closed and has empty interior in the joint parameter space.

For the last assertion we choose a small neighborhood of $\alpha$ contained in $\Acal$. Only finitely many expressions $d+\sum_{p,j}k_{p,j}\alpha_{p,j}$ can change sign there, because the weight sums are bounded and $d$ is an integer. We can therefore take a rational centered $\alpha'$ sufficiently close that every expression nonzero at $\alpha$ has the same sign there. An $\alpha$-stable object is then $\alpha'$-stable. Conversely, an $\alpha'$-semistable object is $\alpha$-semistable, since a positive $\alpha$-slope would remain positive. The hypothesis makes it $\alpha$-stable. These implications prove equality of all four classes of objects, and they apply to each geometric fiber of a family.
\end{proof}

The characteristic coefficients take values in $\mathcal H_0=\bigoplus_{i=2}^n H^0(X,K_X(D)^{\otimes i})$.

\begin{proposition}\label{prop:algebraic-moduli}
If every semistable object with the fixed data $(\alpha,\sigma)$ is stable, the parabolic Higgs bundles with these data have a quasi-projective coarse moduli scheme $\M_{\alpha,\sigma}$. This scheme is smooth, and each nonempty component has complex dimension two. The Hitchin morphism
\begin{equation*}
 h:\M_{\alpha,\sigma}\longrightarrow\mathcal H_0
\end{equation*}
is projective. The same coarse scheme is obtained by specifying an isomorphism $\det E\simeq\mathcal O_X$ or by fixing the isomorphism class of the determinant line bundle. The scheme may be empty.
\end{proposition}

\begin{proof} The proof is given in Appendix~\ref{proof:prop:algebraic-moduli}. \end{proof}

\subsection{Harmonic metrics and the Hitchin metric}

A harmonic metric is a Hermitian metric $h$ on $E|_{X^*}$ satisfying
\begin{equation}\label{eq:prelim-hitchin-equations}
 F_{D_h}+[\Phi,\Phi^{\dagger_h}]=0,\qquad
 \bar\partial_E\Phi=0,
\end{equation}
where $D_h$ is its Chern connection and $\dagger_h$ denotes Hermitian adjoint, including conjugation of the form type. The connection $D_h+\Phi+\Phi^{\dagger_h}$ is then flat. We use the tame nonabelian Hodge correspondence in the form of~\cite[Theorem~7.1, p.~54]{MochizukiGood}. As verified in Proposition~\ref{prop:harmonic-metric-existence}, every stable parabolic Higgs bundle with our data has a unique harmonic metric that induces its filtration and the unit metric on its trivialized determinant line bundle.

The filtration records the growth of $h$ at the punctures. In a holomorphic frame adapted to the flag, a vector of weight $\alpha_{p,j}$ has norm comparable to $|x|^{\alpha_{p,j}}$. For scalar graded residues, this comparison follows from Lemma~\ref{lem:tame-power-comparison}. For the limits of the metric blocks and the estimates for their derivatives, we use Proposition~\ref{prop:harmonic-metric-existence}.

A Dolbeault deformation is a pair $v=(b,\varphi)$, where $b$ is an $\End_0E$-valued $(0,1)$-form and $\varphi$ is an $\End_0E$-valued $(1,0)$-form. The differentials of the deformation complex on $X^*$ are
\begin{equation}\label{eq:tangent-differentials}
 \dd_1\eta=(\bar\partial_E\eta,[\Phi,\eta]),\qquad
 \dd_2(b,\varphi)=\bar\partial_E\varphi+[\Phi,b].
\end{equation}
The second commutator includes the wedge product. A harmonic representative satisfies $\dd_2v=0$ and $\dd_1^*v=0$, where the adjoint uses the unweighted $L^2$ pairing. In Propositions~\ref{prop:weighted-Dolbeault-comparison} and~\ref{prop:ordinary-L2-projection} we specify the weighted Sobolev domains of the differentials and prove that the first cohomology of the weighted complex is isomorphic to the space of algebraic deformation classes.

All traces in the metric and symplectic pairings are matrix traces. We set
\begin{equation}\label{eq:tangent-trial-norm}
 \|v\|_h^2
 =2\ii\int_{X^*}\tr\bigl(
       \varphi\wedge\varphi^{\dagger_h}
       -b\wedge b^{\dagger_h}\bigr),
 \qquad v=(b,\varphi).
\end{equation}
This is twice the unweighted $L^2$ norm squared of the pair. Its real polarization on harmonic representatives defines the Hitchin metric $g_H$. The induced complex structures are
\begin{equation}\label{eq:hk-quaternionic-operators}
 \begin{aligned}
 I(b,\varphi)&=(\ii b,\ii\varphi),\\
 J(b,\varphi)&=(\ii\varphi^{\dagger_h},-\ii b^{\dagger_h}),\qquad
 K(b,\varphi)=(-\varphi^{\dagger_h},b^{\dagger_h}).
 \end{aligned}
\end{equation}
Here $I$ agrees with the complex structure of the Higgs moduli space. Writing $\omega_Q(v,w)=g_H(Qv,w)$, we use $\Omega=\omega_J+\ii\omega_K$ with
\begin{equation}\label{eq:semiflat-trace-form}
 \Omega\bigl((b_1,\varphi_1),(b_2,\varphi_2)\bigr)
 =2\int_{X^*}\tr(\varphi_1\wedge b_2-\varphi_2\wedge b_1).
\end{equation}
For algebraic Dolbeault representatives the same pairing is written over $X$. In Propositions~\ref{prop:hyperkahler-structure} and~\ref{prop:regular-end-hyperkahler} we prove that the pairing descends, establish the smooth hyperk\"ahler structure and give its holomorphic symplectic form on the end of the Hitchin fibration.

\subsection{Spectral data and the semiflat metric}

The Hitchin morphism records the characteristic coefficients:
\begin{equation*}
 h(E,\Phi)=(a_2,\ldots,a_n)\in
 \mathcal H_0,\qquad
 \det(\eta\Id-\Phi)=\eta^n+\sum_{i=2}^n a_i\eta^{n-i}.
\end{equation*}
Here $\eta$ is the tautological section on $\operatorname{Tot}K_X(D)$. The vanishing of this polynomial defines the spectral curve. The residue conditions restrict its coefficients to the affine line with coordinate $t$ computed in Proposition~\ref{prop:coefficient-line}. We write $C_{t,\sigma}$ for the spectral curve and $Y_{t,\sigma}$ for its normalization, with induced projection $\pi:Y_{t,\sigma}\to X$.

In Proposition~\ref{prop:normalized-spectral-family} we construct a smooth proper family of connected genus-one curves over $U_{\sigma}=\C\setminus\Delta_{\sigma}$, with $\Delta_{\sigma}$ finite and $\pi$ finite flat of degree $n$. The normalization separates the branches above punctures with partial flags, and the construction extends to parameters with coincident graded residues.

A spectral line bundle $L$ on $Y_{t,\sigma}$ determines $E=\pi_*L$ and a Higgs field by multiplication by $\eta$. Near $p$, we write $\eta=v\,\dd x/x$ and set $b_j=\sigma_{p,j}$. Under the spectral correspondence of Proposition~\ref{prop:filtered-spectral-correspondence}, the periodic lattices are
\begin{equation}\label{eq:newton-spectral-lattices}
 f_j(v)=\prod_{i=0}^{j-1}(v-b_i),\qquad
 E^j=\pi_*\bigl(f_j(v)L\bigr),\qquad f_0=1.
\end{equation}
These lattices recover the flags in families even when residue values coincide. The Hitchin fiber is $\Pic^n(Y_{t,\sigma})$. The line bundle $L_0=\pi^*\mathcal O_X(1)$ provides its origin through the map $L\mapsto L\otimes L_0^{-1}$ to $\Pic^0(Y_{t,\sigma})$.

With $\sigma$ fixed we abbreviate $Y_{t,\sigma}$ to $Y_t$. The differential induced by $\eta$ on $Y_t$ can have poles above $D$. For a locally transported closed cycle $\Gamma_t$ avoiding $\pi^{-1}(D)$, the variation of its period at fixed residues is represented by a nowhere-zero holomorphic differential $\nu_t$:
\begin{equation*}
 \frac{\dd}{\dd t}\int_{\Gamma_t}\eta=\int_{\Gamma_t}\nu_t.
\end{equation*}
We compute $\nu_t$ and prove this identity in Lemma~\ref{lem:fixed-residue-differential}. Changing the representative of a class in $H_1(Y_t;\Z)$ adds only a constant residue period.

On a simply connected subset of $U_{\sigma}$, we choose an integral basis $A,B$ with $A\cdot B=1$ and representatives avoiding the poles. The periods and their derivative ratio are
\begin{equation}\label{eq:semiflat-special-coordinates}
 a(t)=\int_A\eta,\qquad
 a_D(t)=\int_B\eta,\qquad
 \tau(t)=\frac{a_D'(t)}{a'(t)}.
\end{equation}
Then $a'(t)\ne0$, $\operatorname{Im}\tau>0$ and $\dd a_D=\tau\,\dd a$. A point of $\Pic^0(Y_t)$ has a unitary flat representative $\dd+\ii\chi$, with real periods $\theta_A=\int_A\chi$ and $\theta_B=\int_B\chi$ in $\R/2\pi\Z$. We define the Gauss--Manin horizontal distribution by holding these periods constant~\cite[Remark~7.8]{FMSW}.

The semiflat metric makes horizontal and vertical tangent vectors orthogonal. The squared norm in a horizontal direction is $2\ii\int_{Y_t}(\dot t\nu_t)\wedge\overline{\dot t\nu_t}$. For a vertical harmonic $(0,1)$-form $\upsilon$ the squared norm is $-2\ii\int_{Y_t}\upsilon\wedge\overline\upsilon$. We use the normalization of~\cite[Proposition~7.6]{FMSW} with scaling parameter equal to one. In these coordinates, we obtain from Proposition~\ref{prop:semiflat-period-formulas}
\begin{equation}\label{eq:semiflat-period-formulas}
 \begin{aligned}
 g_{\sfm}
   &=4\operatorname{Im}\tau\,|\dd a|^2
       +\frac{|\dd\theta_B-\tau\,\dd\theta_A|^2}
                    {\operatorname{Im}\tau},\\
 \Omega
   &=2\ii\,\dd a\wedge(\dd\theta_B-\tau\,\dd\theta_A).
 \end{aligned}
\end{equation}
The second formula agrees with~\eqref{eq:semiflat-trace-form}. By~\cite[Theorem~3.8]{FreedSpecialKahler}, $g_{\sfm}$ is hyperk\"ahler. In this normalization each Picard fiber has area $4\pi^2$. In Section~\ref{sec:hitchin-base} we prove independence of the choice of cycles and compute the asymptotics for large $|t|$. We then compare $g_H$ and $g_{\sfm}$ with $\alpha$ and $\sigma$ fixed.

\subsection{Weight alcoves and numerical classes}
\label{sec:closed-parameters}

For the modularity statement we also allow weights on the boundary of the alcove. We retain the labelled periodic lattices and prescribe the semisimple residue on each graded piece of the filtration by distinct weights. The numerical condition below ensures that these residues recover every labelled step; we prove this in Proposition~\ref{prop:canonical-refinement}.

The product of closed weight alcoves is
\begin{equation}\label{eq:closed-parameters}
 \overline{\Acal}=
 \left\{\alpha:\sum_j\alpha_{p,j}=0,\quad
 \alpha_{p,0}\leq\cdots\leq\alpha_{p,e_p-1}
                  \leq\alpha_{p,0}+1\text{ for each }p\right\}.
\end{equation}
Its interior is the space~$\Acal$ used above. In the labelled periodic lattices~\eqref{eq:periodic-chain} the successive quotients have dimensions~$m_p$ and residue actions~$\sigma_{p,j}\Id$, with periodic labels given by $\alpha_{p,j+e_p}=\alpha_{p,j}+1$ and $\sigma_{p,j+e_p}=\sigma_{p,j}$. Omitting intermediate steps with equal weights gives the associated real filtration. On each graded piece we prescribe a semisimple residue with the indicated eigenvalues and multiplicities, up to conjugacy. Together with the trivial filtered determinant line bundle, these data define a regular filtered Higgs bundle~\cite[Section~1.1.3]{MochizukiGood}. A logarithmic extension then gives the logarithmic parabolic data of~\cite[Sections~1 and~8]{BiquardLogarithmic}, with the parabolic degree and stability condition of Section~\ref{sec:parabolic-data}.

Ranks, degrees and labelled dimensions define the lattice
\begin{equation}\label{eq:numeric-def}
 \begin{gathered}
 \Lambda=
 \frac{\{(r,d,(k_{p,j}))\in\Z^{2+\sum_pe_p}:\sum_jk_{p,j}=r\text{ for each }p\}}
      {\Z\delta_0},
 \qquad \delta_0=(n,0,(m_p)_{p,j}),\\
 B(C_{r,d,k},C_{r',d',k'})=rr'-\sum_{p,j}k_{p,j}k'_{p,j},
 \qquad F=C_{0,-1,0}.
 \end{gathered}
\end{equation}
Here $C_{r,d,k}$ denotes the class of~$(r,d,k)$, and $C^2=B(C,C)$. The identities $e_pm_p=n$ and $\sum_pm_p=n$ show that~$B$ descends to the quotient. For centered weights and residues, the affine formulas
\begin{equation}\label{eq:formal-periods}
 \begin{aligned}
 P_I(C_{r,d,k})&=-4\pi^2\left(d+\sum_{p,j}k_{p,j}\alpha_{p,j}\right),\\
 P_{\C}(C_{r,d,k})&=8\pi^2\sum_{p,j}k_{p,j}\sigma_{p,j}
 \end{aligned}
\end{equation}
also descend to~$\Lambda$. We set $P_J=\operatorname{Re}P_{\C}$ and $P_K=\operatorname{Im}P_{\C}$, and write $P_{\alpha,\sigma}=(P_I,P_J,P_K)$. Their interpretation as metric periods follows in Proposition~\ref{prop:hyperkahler-periods} from the calculation of Dimakis--Fredrickson--Yae for $\alpha\in\Acal$ and continuity at the boundary.

We call $(\alpha,\sigma)\in\overline{\Acal}\times\mathscr B$ admissible if
\begin{equation}\label{eq:root-condition}
 C\in\Lambda,\quad C^2=-2
 \quad\Longrightarrow\quad P_{\alpha,\sigma}(C)\ne0.
\end{equation}
The admissible set is denoted by~$\mathscr U$, and $\mathscr U^\circ=\mathscr U\cap(\Acal\times\mathscr B)$. The entries of a class in~\eqref{eq:root-condition} are arbitrary integers satisfying~\eqref{eq:numeric-def}. With $A_F=4\pi^2$, the normalized period domain is
\begin{equation}\label{eq:full-period-domain}
 \mathcal D=\left\{P\in\Hom(\Lambda,\R^3):
 P(F)=(A_F,0,0),\quad P(C)\ne0\text{ whenever }C^2=-2\right\}.
\end{equation}

\section{Spectral geometry and the semiflat metric}
\label{sec:hitchin-base}

The characteristic polynomial and the scalar actions on the graded quotients determine the spectral curve. Over a puncture with a partial flag, the singularities persist as $|t|$ grows. With the local notation $\Phi=H\,\dd x/x$ and $P(v)=\prod_j(v-\sigma_{p,j})$, adjoining the endomorphism $P(H)/x$ resolves these singularities and recovers the flags from line bundles on the normalization, even when graded residues coincide. The periods of the resulting genus-one family determine $g_{\sfm}$ and the scale of its ALG model.

\subsection{Local coefficient conditions}

\begin{lemma}\label{lem:local-coefficients}
Our coefficient ring is $R=\C\{x\}$. For integers $e,m\geq1$, with $n=em$ and fixed $\sigma_0,\ldots,\sigma_{e-1}\in\C$, we use the polynomial
\begin{equation*}
 P(v)=\prod_{j=0}^{e-1}(v-\sigma_j).
\end{equation*}
For a monic polynomial $F(x,v)\in R[v]$ of degree $n$, the following conditions are equivalent:
\begin{enumerate}
\item There is a free $R$-module $E$ of rank $n$, an endomorphism $H$ with characteristic polynomial $F$, and a flag
\begin{equation*}
 E/xE=E_0^0\supset E_0^1\supset\cdots\supset E_0^e=0
\end{equation*}
whose graded quotients have dimension $m$ and on which $H(0)$ induces $\sigma_j\Id$ on $E_0^j/E_0^{j+1}$.
\item In the expansion $F(x,v)=\sum_{\ell\geq0}x^\ell F_\ell(v)$, one has
\begin{equation}\label{eq:local-divisibility}
 F_0=P^m,\qquad P^{m-\ell}\mid F_\ell\quad(1\leq\ell<m).
\end{equation}
\item There are unique polynomials $c_\ell(v,x)\in R[v]$, with $\deg_v c_\ell<e$, such that
\begin{equation}\label{eq:local-monic-expansion}
 F=P^m+\sum_{\ell=1}^m x^\ell c_\ell(v,x)P^{m-\ell}.
\end{equation}
\end{enumerate}
The realization in the first condition can be chosen holomorphically in the polynomial coefficients and the prescribed scalars, also at points where two or more of the scalars coincide. If the coefficient of $v^{n-1}$ vanishes, the resulting endomorphism is trace-free.
\end{lemma}

\begin{proof}
Starting with $H$ and the flag, we write $H_0=H(0)$. The scalar actions on the graded quotients imply
\begin{equation*}
 (H_0-\sigma_j\Id)E_0^j\subset E_0^{j+1},\qquad
 P(H_0)=0,\qquad \det(v\Id-H_0)=P(v)^m.
\end{equation*}
Polynomial division gives a matrix polynomial $T(v)$ of degree at most $e-1$ satisfying
\begin{equation*}
 P(v)\Id=(v\Id-H_0)T(v).
\end{equation*}
We decompose $H$ as $H_0+xC(x)$. Expanding the determinant in $\C(v)[[x]]$ gives
\begin{equation*}
 \det(v\Id-H)
 =\sum_{k=0}^n(-x)^kP(v)^{m-k}
     \operatorname{tr}\!\left(\bigwedge\nolimits^k(T(v)C(x))\right).
\end{equation*}
For $\ell<m$, only terms with $k\leq\ell$ contribute to the coefficient of $x^\ell$, and each is divisible by $P^{m-\ell}$. This proves~\eqref{eq:local-divisibility}.

Successive division by the monic polynomial $P$ gives a unique expansion
\begin{equation*}
 F=P^m+\sum_{\ell=1}^m r_\ell(v,x)P^{m-\ell},
 \qquad \deg_v r_\ell<e.
\end{equation*}
For each $k<m$, divisibility of $F_k$ by $P^{m-k}$ says that its remainders in lower powers of $P$ vanish. In particular, the coefficient of $x^k$ in $r_\ell$ is zero whenever $k<\ell$, which gives $r_\ell=x^\ell c_\ell$ and proves the third condition. Conversely,~\eqref{eq:local-monic-expansion} implies~\eqref{eq:local-divisibility}.

For the construction of $H$ we use the companion matrix $B(v,x)$ with characteristic polynomial
\begin{equation*}
 w^m+c_1(v,x)w^{m-1}+\cdots+c_m(v,x).
\end{equation*}
Its subdiagonal entries are one, its last column is $(-c_m,\ldots,-c_1)^{\mathsf T}$, and its remaining entries are zero. For $m=1$ it is the matrix $(-c_1)$. Then
\begin{equation*}
 \det\bigl(P(v)\Id_m-xB(v,x)\bigr)=F(x,v).
\end{equation*}
The matrix polynomial $P(v)\Id_m-xB(v,x)$ is monic of degree $e$ in $v$. Division by this matrix polynomial gives a unique vector remainder of degree less than $e$. Hence
\begin{equation*}
 E=R[v]^m\big/\bigl(P(v)\Id_m-xB(v,x)\bigr)R[v]^m
\end{equation*}
is free over $R$, with basis $v^a\epsilon_i$ for $0\leq a<e$ and $1\leq i\leq m$, where $\epsilon_1,\ldots,\epsilon_m$ is the standard basis of $R[v]^m$. Multiplication by $v$ defines $H$. In this basis $H$ is the block companion matrix of $P(v)\Id_m-xB(v,x)$. Eliminating the blocks of successive degrees in $v\Id-H$ gives
\begin{equation*}
 \det(v\Id-H)=\det\bigl(P(v)\Id_m-xB(v,x)\bigr)=F(x,v).
\end{equation*}
Modulo $x$, the module is $(\C[v]/P)^m$. The monic polynomials $f_0=1$ and $f_a(v)=\prod_{j<a}(v-\sigma_j)$ give a triangular change of basis from the powers of $v$. We take $E_0^a$ to be the span of $f_b\epsilon_i$ with $b\geq a$. The identity
\begin{equation*}
 vf_a=\sigma_a f_a+f_{a+1},\qquad f_e=P=0\quad\text{in }\C[v]/P,
\end{equation*}
gives the required scalar action on each graded quotient. All divisions used above are by monic polynomials or by powers of $x$ with the stated vanishing order. The construction therefore depends holomorphically on the allowed parameters. Finally, the coefficient of $v^{n-1}$ is $-\operatorname{tr}H$.
\end{proof}

\subsection{Global coefficient formulas}

On the affine $z$-chart of $X$ the section
\begin{equation*}
 s=\frac{\dd z}{z(z-1)}
\end{equation*}
of $K_X(D)$ determines a fiber coordinate by $\eta=ys$. The section has no zeros on this chart and a simple zero at infinity. Since $K_X(D)\simeq\mathcal O_X(1)$, a trace-free characteristic polynomial has the form
\begin{equation}\label{eq:global-characteristic-polynomial}
 F(z,y)=y^n+\sum_{k=2}^n a_k(z)y^{n-k},\qquad \deg a_k\leq k.
\end{equation}
The residue polynomials are
\begin{equation}\label{eq:residue-polynomials}
 \begin{aligned}
 P_0(y)&=\prod_{j=0}^{e_0-1}(y+\sigma_{0,j}),&
 P_1(y)&=\prod_{j=0}^{e_1-1}(y-\sigma_{1,j}),\\
 P_\infty(v)&=\prod_{j=0}^{e_\infty-1}(v+\sigma_{\infty,j}).
 \end{aligned}
\end{equation}
The signs come from the expressions $s=(z-1)^{-1}\dd z/z$ at $0$ and $s=z^{-1}\dd z/(z-1)$ at $1$ for the logarithmic frame. At infinity we use $u=1/z$ and $v=uy$. Then $s=-\dd u/(1-u)$, so the section $s/u=-(1-u)^{-1}\dd u/u$ is a logarithmic frame with residue $-1$. The local polynomial of Lemma~\ref{lem:local-coefficients} at infinity is therefore
\begin{equation}\label{eq:infinity-polynomial}
 F^\infty(u,v)=u^nF(1/u,v/u).
\end{equation}
At the three punctures we impose~\eqref{eq:local-divisibility} with the polynomials in~\eqref{eq:residue-polynomials}, the local variables $z,z-1,u$, and the corresponding multiplicities $m_p$.

\begin{proposition}\label{prop:coefficient-line}
For every choice of the scalar graded residues satisfying~\eqref{eq:graded-residue}, the coefficients in~\eqref{eq:global-characteristic-polynomial} satisfying these local conditions form a nonempty affine line in $\mathcal H_0$. It has a coordinate $t$ for which
\begin{equation}\label{eq:coefficient-line}
 \begin{gathered}
 F_{t,\sigma}(z,y)=F_{0,\sigma}(z,y)-t\,r_n(z),\\
 r_3(z)=r_4(z)=z(z-1),\qquad r_6(z)=z(z-1)^2.
 \end{gathered}
\end{equation}
The origin can be chosen so that each coefficient $a_k(z;0,\sigma)$ is a polynomial homogeneous of degree $k$ in the residues. Every point of the line admits the local matrix realizations constructed in Lemma~\ref{lem:local-coefficients}.
\end{proposition}

\begin{proof} The coefficient calculation is given in Appendix~\ref{app:coefficients}. It uses only monic division, and its formulas are polynomial in the prescribed residues. \end{proof}

In particular, the formulas give $a_k(z;0,0)=0$ for $2\leq k\leq n$, and
\begin{equation}\label{eq:zero-residue-polynomial}
 F_{t,0}(z,y)=y^n-t\,r_n(z).
\end{equation}
Homogeneity gives, for $\lambda\in\C^*$,
\begin{equation}\label{eq:coefficient-scaling}
 \begin{aligned}
 F_{\lambda^nt,\lambda\sigma}(z,\lambda y)
       &=\lambda^nF_{t,\sigma}(z,y),\\
 \lambda^{-n}F_{\lambda^n,\sigma}(z,\lambda y)
       &=F_{1,\sigma/\lambda}(z,y).
 \end{aligned}
\end{equation}
The second equality describes the rescaled polynomial on the cover $t=\lambda^n$ of the Hitchin base with the residues fixed. The residues after rescaling are $\sigma/\lambda$.

The Hitchin morphism therefore factors through this affine line. We use the spectral correspondence below to realize these coefficients globally.

\subsection{Normalized spectral curves}
\label{subsec:normalized-spectral-curves}

The spectral curve $C_{t,\sigma}\subset\operatorname{Tot}K_X(D)$ is defined by~\eqref{eq:coefficient-line}, and its family $\mathcal C_{\sigma}\to\C_t$ is parametrized by $t$. The coefficient conditions allow singularities above punctures with a partial flag, even for large $|t|$. We construct a finite modification that separates the local branches and varies holomorphically when graded residues coincide. For strongly parabolic Higgs bundles, the local module construction and the resulting spectral correspondence are developed in~\cite[Theorem~4.3, p.~10, and Theorem~4.12, p.~15]{SuWenWang}. The construction below includes the prescribed scalar actions on the graded quotients.

\begin{proposition}\label{prop:normalized-spectral-family}
For every $\sigma$, there is a finite set $\Delta_{\sigma}\subset\C$ such that, over $U_{\sigma}:=\C\setminus\Delta_{\sigma}$, the spectral curves admit a finite simultaneous normalization
\begin{equation*}
 \mathcal Y_{\sigma}\longrightarrow\mathcal C_{\sigma}|_{U_{\sigma}}.
\end{equation*}
The family $\mathcal Y_{\sigma}\to U_{\sigma}$ is smooth and proper, with connected fibers of genus one, and its projection $\pi:\mathcal Y_{\sigma}\to X\times U_{\sigma}$ is finite flat of degree $n$.

There is also a neighborhood $S$ of zero in the residue parameter space for which the normalizations of $C_{1,\boldsymbol\gamma}$ form a smooth proper holomorphic family over $\boldsymbol\gamma\in S$, with the same genus and covering degree. For each $\sigma$, its pullback by $\boldsymbol\gamma=\sigma/\lambda$ gives the normalized curves over $t=\lambda^n$ for all sufficiently large $|\lambda|$.
\end{proposition}

\begin{proof}
We first construct the family near $(t,\sigma)=(1,0)$. In a coordinate $x$ vanishing at the puncture, the spectral differential has the form $v\,\dd x/x$. Omitting the puncture index, we have $e$ graded quotients of dimension $m$ and $n=em$. For $\boldsymbol\gamma$ near zero we write $P_{\boldsymbol\gamma}(v)=\prod_{j=0}^{e-1}(v-\gamma_{p,j})$. For the local polynomial in this frame, the expansion in Lemma~\ref{lem:local-coefficients} reads
\begin{equation*}
 F(x,v,\boldsymbol\gamma)
 =P_{\boldsymbol\gamma}(v)^m+
   \sum_{\ell=1}^m x^\ell c_\ell(v,x,\boldsymbol\gamma)
                  P_{\boldsymbol\gamma}(v)^{m-\ell},
 \qquad \deg_v c_\ell<e.
\end{equation*}
The coefficients are holomorphic in $x$ and $\boldsymbol\gamma$. We adjoin a coordinate $u$ satisfying
\begin{equation}\label{eq:spectral-finite-modification}
 \begin{aligned}
 P_{\boldsymbol\gamma}(v)&=xu,\\
 Q(u,v,x,\boldsymbol\gamma)
   &:=u^m+\sum_{\ell=1}^m
             c_\ell(v,x,\boldsymbol\gamma)u^{m-\ell}=0.
 \end{aligned}
\end{equation}
These equations imply $F=0$. Since $F$ is monic in $v$ and $Q$ is monic in $u$, the new local algebra is finite over both the spectral algebra and the $x$-disc times the parameter space. On $x\ne0$, the coordinate $u=P_{\boldsymbol\gamma}(v)/x$ is uniquely determined, so the finite modification is an isomorphism there.

At $\boldsymbol\gamma=0$, equation~\eqref{eq:zero-residue-polynomial} has the local form
\begin{equation*}
 F(x,v,0)=v^{em}-x^m h(x),\qquad h(0)\ne0.
\end{equation*}
Thus $Q(u,0,0,0)=u^m-h(0)$ has $m$ distinct nonzero roots. By the holomorphic implicit function theorem, the solutions near these roots form $m$ disjoint branches
\begin{equation}\label{eq:local-normalized-branches}
 u=U_\nu(x,v,\boldsymbol\gamma),\qquad
 P_{\boldsymbol\gamma}(v)=xU_\nu(x,v,\boldsymbol\gamma),
 \qquad 1\leq\nu\leq m,
\end{equation}
where every $U_\nu$ is a unit. The derivative of the second equation with respect to $x$ is nonzero at the reference point. After shrinking the neighborhoods, each branch therefore has the form $x=X_\nu(v,\boldsymbol\gamma)$ and is smooth over $S$. These branches cover the finite modification near the puncture. Indeed, monicity bounds both $v$ and $u$ as $x$ and $\boldsymbol\gamma$ tend to zero, and every limiting point has $v=0$ and $u^m=h(0)$. By Weierstrass preparation in $v$, the second equation in~\eqref{eq:local-normalized-branches} has degree $e$ up to a holomorphic unit. Each branch is consequently finite flat of degree $e$ over the $x$-disc times $S$.

The same local graph gives the ramification index above a repeated residue. If $b$ is a root of $P_{\boldsymbol\gamma}$ of multiplicity $r$, then
\begin{equation*}
 X_\nu(v,\boldsymbol\gamma)
 =\frac{P_{\boldsymbol\gamma}(v)}
        {U_\nu(X_\nu(v,\boldsymbol\gamma),v,\boldsymbol\gamma)}
\end{equation*}
vanishes to order $r$ at $v=b$. The normalized branch is smooth, and its projection has ramification index $r$. The finite modification is an isomorphism for $x\ne0$ and every branch meets that set. Its smooth fibers are therefore the normalizations of the local spectral curves.

We glue these modifications at the three punctures to the original spectral curve on $X\setminus D$. The gluing maps are specified by $u=P_{\boldsymbol\gamma}(v)/x$ on each punctured overlap. At $\boldsymbol\gamma=0$, the spectral projection is \'{e}tale on the complement of small puncture discs. By compactness of this complement, the projection remains \'{e}tale there for all $\boldsymbol\gamma$ in a sufficiently small $S$. The glued map to the original spectral family is finite, since finiteness holds on every spectral chart. The original family is finite over $X\times S$, so the modified family is proper over $S$. Its charts are smooth over $S$, and their projections to $X\times S$ are finite flat of total degree $n$.

The reference polynomial $y^n-r_n(z)$ is irreducible by Eisenstein's criterion at $z=0$, because $r_n$ has a simple zero there. The normalization of the corresponding spectral curve is connected. Above a puncture $p$, the reference curve has $m_p$ points, each of ramification index $e_p$. There is no ramification away from $D$, so by the Riemann--Hurwitz formula
\begin{equation}\label{eq:spectral-genus-one}
 2g-2=-2n+\sum_{p\in D}m_p(e_p-1)=0.
\end{equation}
A proper smooth submersion is locally differentiably trivial by Ehresmann's fibration theorem. For connected $S$, all the fibers are therefore connected curves of genus one. Each fiber is smooth and connected, hence irreducible. Since the finite modification is generically an isomorphism, each associated characteristic polynomial is irreducible over $\C(z)$.

For arbitrary fixed $\sigma$, equation~\eqref{eq:coefficient-scaling} sends the curve at $(t,\sigma)=(\lambda^n,\sigma)$ to the curve at $(1,\sigma/\lambda)$. For all sufficiently large $|\lambda|$, the latter parameter lies in $S$. The same scaling extends to~\eqref{eq:spectral-finite-modification}, so the curves for all sufficiently large $|t|$ admit the same normalization construction, with finite flat projection.

For the algebraic family over a Zariski-open subset of the $t$-line, we cover $X$ by the three affine open sets obtained by omitting two punctures and fix a logarithmic trivialization on each. Monic division in Lemma~\ref{lem:local-coefficients} and division by the prescribed powers of $x$ give coefficients $c_\ell$ regular on the chart and polynomial in $t$. Thus~\eqref{eq:spectral-finite-modification} defines a finite algebraic modification on each chart. The overlaps avoid $D$, where $u=P(v)/x$ gives the isomorphism to the original spectral curve. These algebras glue to a finite modification $Z_{\sigma}\to\mathcal C_{\sigma}$ over the entire $t$-line, proper over that line.

The locus where $Z_{\sigma}\to\C_t$ is not smooth is closed, and its image is closed by properness. The image contains no sufficiently large $|t|$ and is therefore finite. The locus where the finite coherent module $\pi_*\mathcal O_{Z_{\sigma}}$ is not locally free of rank $n$ is closed in $X\times\C_t$. Its projection to $\C_t$ is finite for the same reason. At each puncture we also remove the images of
\begin{equation*}
 \{x=0,\ u=0\},\qquad
 \{x=0,\ \partial_uQ=0\}.
\end{equation*}
These sets are finite over $\C_t$, since the modification is finite over $X\times\C_t$, and they contain no points with sufficiently large $|t|$. Their images are therefore finite. The union of these finite sets is denoted by $\Delta_{\sigma}$, and we take $\mathcal Y_{\sigma}=Z_{\sigma}|_{U_{\sigma}}$.

The family $\mathcal Y_{\sigma}\to U_{\sigma}$ is smooth and proper, and $\pi$ is finite flat of degree $n$. The connected base $U_{\sigma}$ contains points with arbitrarily large $|t|$. Local differentiable triviality therefore gives connected fibers of genus one throughout $U_{\sigma}$. The modification is an isomorphism away from $D$, and its restriction to each fiber is the normalization because the modified fiber is smooth. This proves the asserted simultaneous normalization.
\end{proof}

For sufficiently large $|t|$, the family descends under the cover $t=\lambda^n$ as follows. We may take $S$ invariant under multiplication by $n$th roots of unity. For $\zeta^n=1$, the change
\begin{equation}\label{eq:spectral-cyclic-descent}
 (\lambda,\boldsymbol\gamma,y)
 \longmapsto(\zeta\lambda,\zeta^{-1}\boldsymbol\gamma,
                         \zeta^{-1}y)
\end{equation}
preserves the equations of the rescaled family by~\eqref{eq:coefficient-scaling}. In a local logarithmic frame it acts by $v\mapsto\zeta^{-1}v$ and $u\mapsto\zeta^{-e}u$, as follows from $P_{\zeta^{-1}\boldsymbol\gamma}(\zeta^{-1}v)=\zeta^{-e}P_{\boldsymbol\gamma}(v)$. This gives a $\mu_n$-action on the normalized family over the $\lambda$-cover. The spectral differential $\lambda ys$ is invariant, so the quotient is the normalized family with fixed residues over all sufficiently large $|t|$.

\subsection{The spectral correspondence}
\label{subsec:spectral-correspondence}

Here $Y_{t,\sigma}$ is the fiber of $\mathcal Y_{\sigma}$ over $t$. The action of $u$ determines which spectral modules extend to the normalization. The scalar action on each graded quotient makes $P(H)/x$ holomorphic, where $H$ is the Higgs matrix in a logarithmic frame. This endomorphism extends the spectral module to the finite modification, and the ordered factors of $P$ recover the flag even when some residues coincide.

\begin{proposition}\label{prop:filtered-spectral-correspondence}
\begin{enumerate}
\item For each $t\in U_{\sigma}$, direct image under $\pi:Y_{t,\sigma}\to X$, with Higgs field induced by multiplication by the spectral differential, gives a bijection between $\Pic^n(Y_{t,\sigma})$ and the isomorphism classes of degree-zero parabolic Higgs bundles with graded residues $\sigma$ and characteristic polynomial $F_{t,\sigma}$. Every such Higgs bundle is trace-free, has trivial determinant line bundle, and is stable for every $\alpha\in\Acal$.
\item The construction commutes with base change by locally Noetherian analytic parameter spaces, including nonreduced ones. The sheafification of the relative functor of isomorphism classes in the analytic topology is $\Pic^n(\mathcal Y_{\sigma}/U_{\sigma})$, both with and without a specified determinant trivialization. The section $L_0=\pi^*\mathcal O_X(1)$ gives an isomorphism with $\Pic^0(\mathcal Y_{\sigma}/U_{\sigma})$ by tensoring with $L_0^{-1}$.
\end{enumerate}
\end{proposition}

\begin{proof}
We start with a Higgs bundle $(E,\Phi)$ having the prescribed characteristic polynomial. At a puncture with local coordinate $x$ and $\Phi=H\,\dd x/x$, we write $b_j=\sigma_{p,j}$ and $P(v)=\prod_{j=0}^{e-1}(v-b_j)$. The flag on $E/xE$ lifts to lattices
\begin{equation*}
 E=E^0\supset E^1\supset\cdots\supset E^e=xE,
 \qquad \operatorname{length}(E^j/E^{j+1})=m.
\end{equation*}
The graded residue condition gives $(H-b_j\Id)E^j\subset E^{j+1}$. We extend the indices by $E^{j+e}=xE^j$ and $b_{j+e}=b_j$. Starting with any lattice $E^j$, the successive factors $H-b_j\Id,\ldots,H-b_{j+e-1}\Id$ send it into $E^{j+e}=xE^j$. All factors commute, so
\begin{equation}\label{eq:normalized-algebra-action}
 P(H)E^j\subset xE^j,\qquad B:=x^{-1}P(H)\in\End(E),
 \qquad B(E^j)\subset E^j.
\end{equation}
The endomorphism $B$ commutes with $H$. By Cayley--Hamilton, after substituting $P(H)=xB$ into~\eqref{eq:local-monic-expansion} we obtain
\begin{equation*}
 x^m\left(B^m+\sum_{\ell=1}^m
                    c_\ell(H,x)B^{m-\ell}\right)=0.
\end{equation*}
Multiplication by $x$ is injective on $E$, hence $Q(B,H,x)=0$. Thus the local algebra in~\eqref{eq:spectral-finite-modification} acts on $E$, with $v$ acting by $H$ and $u$ by $B$.

Away from $D$, the Higgs field already gives a module over the original spectral algebra. Its restriction agrees with the extended action because $u=P(v)/x$ wherever $x\ne0$. We therefore obtain a coherent module on $Y_{t,\sigma}$. Over the generic \'{e}tale locus every spectral root occurs once, so this module has rank one. Any torsion on the smooth curve would be supported at finitely many points, and its finite direct image would be torsion in the locally free sheaf $E$. There is consequently no torsion. The local rings of a smooth curve are discrete valuation rings, over which a finite torsion-free module of rank one is free. The spectral module is a line bundle $L$.

We recover the flag by proving $f_j(H)E=E^j$, as asserted in~\eqref{eq:newton-spectral-lattices}. The filtered relations first give the inclusion $f_j(H)E\subset E^j$. Over a sufficiently small puncture neighborhood, $u$ is a unit by the construction of $U_{\sigma}$. Setting $v=b_i$ in the normalized algebra gives $xu=0$, hence $x=0$, and leaves the degree-$m$ monic equation $Q(u,b_i,0)=0$. The divisor $v=b_i$ on the normalized curve therefore has length $m$. Additivity of divisor lengths, with multiplicity, gives
\begin{equation*}
 \operatorname{length}(E/f_j(H)E)=jm
             =\operatorname{length}(E/E^j).
\end{equation*}
The inclusion is thus an equality. Repeated values among the $b_i$ contribute repeated factors to the divisor and satisfy the same length calculation. Near $(t,\sigma)=(1,0)$ we write $L_\nu$ for the restriction of $L$ to the $\nu$th branch in~\eqref{eq:local-normalized-branches} and $\pi_\nu$ for the projection of this branch. Then $E_\nu^j=(\pi_\nu)_*(f_j(v)L_\nu)$, with each successive quotient of length one.

Conversely, a line bundle $L$ on $Y_{t,\sigma}$ has locally free direct image $E=\pi_*L$ of rank $n$, since $\pi$ is finite flat. Multiplication by the spectral differential gives its Higgs field. Generically over $X$, this is multiplication by the spectral coordinate on a rank-one module over the function field of the spectral curve, a degree-$n$ extension of $\C(z)$, so its characteristic polynomial is $F_{t,\sigma}$. At a puncture, formula~\eqref{eq:newton-spectral-lattices} defines a lattice chain with
\begin{equation*}
 E^e=\pi_*\bigl(P(v)L\bigr)=\pi_*(xuL)=xE.
\end{equation*}
Each successive quotient has dimension $m$, and multiplication by $v$ acts on it as $b_j$. These are precisely the prescribed flag and graded residue conditions. The two constructions are inverse on modules and morphisms.

The degrees follow from finite direct image and the Riemann--Roch formula
\begin{equation}\label{eq:spectral-line-degree}
 \deg E+n=\chi(X,E)=\chi(Y_{t,\sigma},L)=\deg L.
\end{equation}
Here $X$ has genus zero and $Y_{t,\sigma}$ has genus one. Thus $\deg E=0$ is equivalent to $\deg L=n$. The characteristic polynomial has vanishing coefficient of degree $n-1$, so $\operatorname{tr}\Phi=0$. Since it is irreducible over $\C(z)$, the generic fiber $E\otimes\C(z)$ has no nonzero proper invariant subspace. There is therefore no proper invariant subbundle, and stability holds for every allowed weight. Every degree-zero determinant line bundle on $X=\Pp^1$ is trivial. The line bundle $L_0=\pi^*\mathcal O_X(1)$ has degree $n$, so it gives a Higgs bundle in every fiber and the asserted origin for the relative Picard sheaf.

For a locally Noetherian analytic base $S'$ mapping to $U_{\sigma}$ or to the residue neighborhood in Proposition~\ref{prop:normalized-spectral-family}, we write $\mathcal Y'$ for the normalized family over $S'$. For a family $(E,\Phi)$ of parabolic Higgs bundles over $S'$ with the prescribed characteristic polynomials, local splittings of the flags show that the lattices $E^j$ are locally free on $X\times S'$, since the graded quotients are locally free. The relative coordinate $x$ acts injectively on these lattices, even when $S'$ is nonreduced. Thus the construction in~\eqref{eq:normalized-algebra-action} gives the unique endomorphism $B$ satisfying $xB=P(H)$, with $Q(B,H,x)=0$. Pullback preserves the defining identity, and injectivity of $x$ after base change preserves uniqueness. The normalized algebra action therefore commutes with base change.

The associated coherent spectral module $M$ on $\mathcal Y'$ is flat over $S'$. Indeed, its finite direct image $E$ is flat over $S'$, and localization at a point of the finite algebra preserves this flatness. On every fiber, $M$ is a line bundle by the argument above. At a point $q\in\mathcal Y'$ over $s\in S'$ the relevant Noetherian local rings are $R=\mathcal O_{S',s}$ and $A=\mathcal O_{\mathcal Y',q}$. The finite $A$-module $M_q$ is $R$-flat, and $M_q/\mathfrak m_RM_q$ is free of rank one over $A/\mathfrak m_RA$. By the local freeness criterion~\cite[Tag~00MH]{Stacks}, the module $M_q$ is free over $A$. Reduction modulo $\mathfrak m_R$ shows that it has rank one. Coherence then makes $M$ a line bundle.

The cokernel of $E\xrightarrow{f_j(H)}E^j$ vanishes on every fiber, so~\eqref{eq:newton-spectral-lattices} holds over $S'$ by Nakayama's lemma. In the converse construction, the divisor $v=b_j$ is finite flat of degree $m$ over $S'$ because its algebra is defined by the monic equation $Q(u,b_j,0)=0$, with $u$ invertible. Near this divisor, $\partial_uQ\ne0$ permits us to solve for $u$. The $x$ derivative of $P(v)-xu$ at $x=0$ is then $-u\ne0$. Thus $v-b_j$ is a relative local coordinate. These coordinates persist after nonreduced base change, so $v=b_j$ remains a relative effective Cartier divisor. Multiplication by $f_j$, a product of non-zero-divisors, gives
\begin{equation*}
 f_jM/f_{j+1}M\simeq M/(v-b_j)M.
\end{equation*}
The direct image of $M/(v-b_j)M$ is locally free of rank $m$ over $S'$. Hence the product lattices define flags with locally free graded quotients, and the lattices and finite direct images commute with base change. Away from $D$, multiplication by $v$ on a line bundle over the monic algebra $\mathcal O[v]/(F)$ has characteristic polynomial $F$; injectivity of $x$ extends this identity across each puncture.

The determinant argument in Proposition~\ref{prop:algebraic-moduli} applies over $S'$, so specifying a determinant trivialization leaves the relative moduli sheaf unchanged. The relative moduli sheaf is therefore $\Pic^n$, with line bundles taken up to tensoring by pullbacks from $S'$. The section $L_0$ and tensoring with $L_0^{-1}$ commute with base change.
\end{proof}

If semistability equals stability for $(\alpha,\sigma)$, then $\M_{\alpha,\sigma}$ is a nonempty smooth complex surface and its Hitchin morphism is surjective onto the affine line of Proposition~\ref{prop:coefficient-line}. Indeed, the smooth coarse moduli scheme and projective Hitchin morphism follow from Proposition~\ref{prop:algebraic-moduli}. The degree-$n$ line bundles above realize every point of $U_{\sigma}$ by stable Higgs bundles. Properness makes the image closed, hence equal to the entire affine line. Nonemptiness and the dimension assertion follow.

The spectral correspondence is compatible with~\eqref{eq:spectral-cyclic-descent}. Each polynomial $f_j$ is multiplied by the nonzero constant $\zeta^{-j}$, so the ordered lattices are preserved. The action fixes the projection to $X$ and therefore preserves $L_0$. Transport of spectral line bundles gives an order-$n$ action on the relative Picard family over the $\lambda$-cover, and its quotient describes the restriction of the fixed-residue Hitchin fibration to sufficiently large $|t|$.

\subsection{The semiflat metric}
\label{subsec:semiflat-geometry}

For fixed graded residues $\sigma$ we write $Y_{t,\sigma}$ as $Y_t$ for $t\in U_{\sigma}$. Over $U_{\sigma}$ the Hitchin fibration is the relative Picard variety of Proposition~\ref{prop:filtered-spectral-correspondence}, namely $\Pic^n(\mathcal Y_{\sigma}/U_{\sigma})$, with origin $L_0=\pi^*\mathcal O_X(1)$. The special coordinates on the base are periods of the spectral differential.

\begin{lemma}\label{lem:fixed-residue-differential}
For the spectral differential $\eta=y\,\dd z/[z(z-1)]$, the expression
\begin{equation}\label{eq:fixed-residue-differential}
 \nu_t=\frac{r_n(z)}{\partial_yF_{t,\sigma}(z,y)}
                   \frac{\dd z}{z(z-1)}
\end{equation}
extends to a nowhere-zero holomorphic differential on $Y_t$. If $\Gamma_t$ is a locally transported closed cycle avoiding $\pi^{-1}(D)$, then
\begin{equation}\label{eq:fixed-residue-period-variation}
 \frac{\dd}{\dd t}\int_{\Gamma_t}\eta=\int_{\Gamma_t}\nu_t.
\end{equation}
The derivative depends only on the class of $\Gamma_t$ in $H_1(Y_t;\Z)$. Different lifts to $H_1(Y_t\setminus\pi^{-1}(D);\Z)$ change the period of $\eta$ by a constant determined by the fixed residues.
\end{lemma}

\begin{proof}
Differentiating~\eqref{eq:coefficient-line} at fixed $z$ gives $\partial_yF_{t,\sigma}\,\partial_t y=r_n$, hence~\eqref{eq:fixed-residue-differential} on the \'{e}tale locus. Away from $D$, the spectral curve is smooth. At a ramification point, $\partial_zF\ne0$ and $\dd z/\partial_yF=-\dd y/\partial_zF$, which proves holomorphic extension there.

At a puncture we work in the logarithmic frame $\eta=v\,\dd x/x$ of Proposition~\ref{prop:normalized-spectral-family}. On each branch of the normalization, equations~\eqref{eq:spectral-finite-modification} give
\begin{equation*}
 P(v)=xU(x,v,t),\qquad
 P(v)=\prod_{j=0}^{e-1}(v-\sigma_{p,j}),
\end{equation*}
where $U$ is a holomorphic unit. On the graph $x=X(v,t)$ in the normalization we write $\widehat U(v,t)=U(X(v,t),v,t)$, so that $X=P/\widehat U$. Holding $x$ fixed and then restricting to a fiber gives
\begin{equation}\label{eq:puncture-period-differential}
 \nu_t=-\frac{\partial_tX}{X}\,\dd v
       =\partial_t\log\widehat U\,\dd v.
\end{equation}
Even at a repeated root of $P$, the last expression is holomorphic because $P$ is independent of $t$ and $\widehat U$ is a unit. Thus $\nu_t$ is holomorphic on the compact curve. It is nonzero on the generic \'{e}tale locus, and a nonzero holomorphic differential on a genus-one curve has no zeros.

On the total normalized family, exterior differentiation gives $\dd\eta=\dd t\wedge\nu_t$. If $b$ is a root of $P$ of multiplicity $r$, then $X(v,t)$ has a zero of order $r$ at $v=b$, and the residue of $\eta$ there is $rb$, independent of $t$. Near this point we use the lift $V$ of $\partial_t$ that holds $v$ fixed. Its contraction with $\eta$ is $v\,\partial_tX/X=-v\,\partial_t\log\widehat U$, which is holomorphic at $v=b$. For a smooth lift on the family extending these local lifts, Cartan's formula on each fiber reads
\begin{equation*}
 (\mathcal L_V\eta)|_{Y_t}=\nu_t+\dd_{Y_t}(\iota_V\eta).
\end{equation*}
Integration over a transported closed cycle proves~\eqref{eq:fixed-residue-period-variation}. Two lifts of a class in $H_1(Y_t;\Z)$ differ by puncture loops. Their $\eta$ periods differ by integral multiples of $2\pi\ii$ times the fixed residues $rb$, proving the final assertion.
\end{proof}

For the cycles, special coordinates and angles of the unitary characters in Section~\ref{sec:preliminaries}, we obtain from Lemma~\ref{lem:fixed-residue-differential} and the Riemann bilinear relations that $a'(t)\ne0$ and $\operatorname{Im}\tau>0$. Thus $a$ is a local coordinate. Even when graded residues coincide, the differential extends holomorphically across the points over the punctures. Its Hodge norm is therefore defined on the whole normalized spectral curve. By the trace duality in Proposition~\ref{prop:local-deformations}, the pairing~\eqref{eq:semiflat-trace-form} is independent of the algebraic Dolbeault representatives.

\begin{proposition}\label{prop:semiflat-period-formulas}
The Hodge norms on the spectral curve and the trace pairing are given by~\eqref{eq:semiflat-period-formulas}.

The metric is hyperk\"ahler, and each Picard fiber has area $4\pi^2$. These tensors are independent of the choices of cycles and their representatives avoiding the poles, and are compatible with the cyclic descent~\eqref{eq:spectral-cyclic-descent}.
\end{proposition}

\begin{proof}
We write $b=\operatorname{Im}\tau$ and $\vartheta=\nu_t/a'(t)$. The normalized differential $\vartheta$ has periods $1,\tau$, and $\ii\int_{Y_t}\vartheta\wedge\overline\vartheta=2b$. For a real harmonic one-form $\dot\chi$ with periods $p_A,p_B$, its induced Dolbeault variation is
\begin{equation*}
 \upsilon=(\ii\dot\chi)^{0,1}
          =-\frac{p_B-\tau p_A}{2b}\,\overline\vartheta.
\end{equation*}
Substitution in the Hodge norms on the spectral curve gives the two terms of $g_{\sfm}$ in~\eqref{eq:semiflat-period-formulas}. The complex structure on the relative Picard variety, as described by the exponential sequence, gives local holomorphic coordinates $(a,w)$ with $w=\theta_B-\tau\theta_A$ and fiber lattice $2\pi\Z+2\pi\tau\Z$.

We compute the symplectic form from the trace pairing~\eqref{eq:semiflat-trace-form}. A vertical tangent vector has a representative $\upsilon$. By the local Dolbeault lemma, smooth cutoffs give a cohomologous form $\upsilon'=\upsilon-\bar\partial f$ that vanishes near every ramification point and every point over $D$. On the remaining \'{e}tale locus it pushes forward to a diagonal variation $b_V$ of the holomorphic structure, with $\varphi_V=0$, and extends by zero over the omitted neighborhoods. Subtracting $(\tr b_V/n)\Id$ imposes the determinant condition. Its scalar part is $\bar\partial$-exact because $H^1(X,\mathcal O_X)=0$, and this correction has zero trace pairing with a trace-free Higgs variation.

For a horizontal tangent vector, the diagonal entries of its Higgs variation on the \'{e}tale locus are the eigenvalue derivatives $\dot t\nu_t$. Therefore the mixed pairing is
\begin{equation*}
 2\int_X\tr(\varphi_H\wedge b_V)
   =2\int_{Y_t}\dot t\nu_t\wedge\upsilon'
   =2\int_{Y_t}\dot t\nu_t\wedge\upsilon
   =2\ii\,\dd a(\dot t)(p_B-\tau p_A).
\end{equation*}
The middle equality follows from Stokes' theorem on the compact normalized curve, since $\nu_t$ is holomorphic. The localized representative contributes no puncture boundary term. Since a holomorphic two-form on this complex surface is determined by its pairing of base and fiber directions, the second formula in~\eqref{eq:semiflat-period-formulas} follows.

The hyperk\"ahler property follows from the semiflat construction~\cite[Theorem~3.8, p.~18 of the author version]{FreedSpecialKahler}, applied with the principal polarization and with $\Omega/(4\pi^2)$; its metric is $g_{\sfm}/(4\pi^2)$. In our normalization the K\"ahler form is
\begin{equation*}
 \omega_{\sfm}=2\ii b\,\dd a\wedge \dd\bar a
                       +\dd\theta_A\wedge \dd\theta_B.
\end{equation*}
Restriction to a fiber gives $\dd\theta_A\wedge \dd\theta_B$, whose integral is $4\pi^2$. The Hodge pairings on the spectral curves, the section $L_0$ and the differential with fixed residues are all intrinsic to the normalized family. They are preserved by changes of integral homology basis and by~\eqref{eq:spectral-cyclic-descent}. Changing the representatives of $A$ and $B$ adds constants to $a$ and $a_D$, so it leaves the tensors unchanged.
\end{proof}

\subsection{The semiflat end}
\label{subsec:semiflat-end}

The periods of the smooth spectral family in Proposition~\ref{prop:normalized-spectral-family} determine the end of $g_{\sfm}$. We write $Y_{1,\boldsymbol\gamma}$ as $Y_{\boldsymbol\gamma}$ and denote its spectral differential by $\eta_{\boldsymbol\gamma}$, for $\boldsymbol\gamma$ in a small ball $S$ centered at zero. The ball $S$ can be taken invariant under multiplication by $n$th roots of unity. By~\eqref{eq:coefficient-scaling}, the curve with fixed residues $\sigma$ over $t=\lambda^n$ is obtained from $Y_{\sigma/\lambda}$, with differential $\lambda\eta_{\sigma/\lambda}$.

We first fix a primitive integral marking of the reference curve. When all graded residues vanish, equations~\eqref{eq:zero-residue-polynomial} and~\eqref{eq:fixed-residue-differential} give $\eta_0=n\nu_1$, so $\eta_0$ is a nowhere-zero holomorphic differential. Its periods form a lattice $\mathcal L$. For $\zeta=\exp(2\pi\ii/n)$, the automorphism $\iota:(z,y)\mapsto(z,\zeta y)$ satisfies $\iota^*\eta_0=\zeta\eta_0$. Thus $\mathcal L$ is preserved by multiplication by
\begin{equation}\label{eq:limiting-period-ratio}
 \tau_\infty=
 \begin{cases}
  \ii,&n=4,\\
  \exp(2\pi\ii/3),&n=3,6.
 \end{cases}
\end{equation}
For any shortest nonzero period $A_0$ we have $\mathcal L=A_0(\Z+\tau_\infty\Z)$. Indeed, any additional period could be reduced modulo this sublattice to $A_0(x+\tau_\infty y)$ with $|x|,|y|\leq1/2$; its norm would be at most $|A_0|/\sqrt2$ in the square case and $\sqrt3|A_0|/2$ in the hexagonal case, contradicting minimality. The period map therefore gives a primitive integral basis $A,B$ with periods $A_0,\tau_\infty A_0$. By the Riemann bilinear relations and $\operatorname{Im}\tau_\infty>0$, we have $A\cdot B=1$.

\begin{lemma}\label{lem:exterior-period-expansion}
After shrinking $S$, we can represent the classes $A,B$ by cycles on the unramified coverings of the complement in $X$ of three fixed puncture discs, using a smooth identification of these coverings. The integral of $\eta_{\boldsymbol\gamma}$ over such a cycle depends only on its class in $H_1(Y_{\boldsymbol\gamma};\Z)$. The resulting period functions are holomorphic on $S$ and satisfy
\begin{equation}\label{eq:small-residue-periods}
 \begin{aligned}
 f_A(\boldsymbol\gamma)
       &:=\int_{A_{\boldsymbol\gamma}}\eta_{\boldsymbol\gamma}
         =A_0+O(|\boldsymbol\gamma|^2),\\
 f_B(\boldsymbol\gamma)
       &:=\int_{B_{\boldsymbol\gamma}}\eta_{\boldsymbol\gamma}
         =\tau_\infty A_0+O(|\boldsymbol\gamma|^2).
 \end{aligned}
\end{equation}
In particular, the special coordinates pulled back to the $\lambda$-cover have convergent Laurent expansions
\begin{equation}\label{eq:fixed-residue-Laurent-periods}
 \begin{pmatrix}a(\lambda)\\a_D(\lambda)\end{pmatrix}
 =\lambda\begin{pmatrix}A_0\\\tau_\infty A_0\end{pmatrix}
   +\sum_{d\geq2}\lambda^{1-d}
       \begin{pmatrix}f_{A,d}(\sigma)\\
                       f_{B,d}(\sigma)\end{pmatrix},
\end{equation}
where $f_{A,d}$ and $f_{B,d}$ are the homogeneous degree-$d$ Taylor terms of $f_A$ and $f_B$. There are no constant or logarithmic terms.
\end{lemma}

\begin{proof}
We take pairwise disjoint closed coordinate discs about the punctures. At $\boldsymbol\gamma=0$, each inverse image is a disjoint union of discs, one for each normalized local branch, and the projection is unramified outside the puncture discs and on annular neighborhoods of the boundary circles. By compactness the covering remains unramified there for small $\boldsymbol\gamma$. Both the normalized family and its restriction to the boundary circles are proper smooth submersions. By Ehresmann's fibration theorem for manifolds with boundary~\cite[Lemma~3.2]{CataneseDifferentiable}, the family admits a smooth trivialization preserving these disc components. Over the complement of the discs, the finite covering has constant monodromy on the contractible parameter ball. We may consequently lift fixed projected paths to mark its homology throughout $S$.

Near $p$ the spectral differential has the local expression $\eta_{\boldsymbol\gamma}=v\,\dd x/x$. In this fiber coordinate each disc $\Delta_{p,\nu}(\boldsymbol\gamma)$ is the graph $x=P_{p,\boldsymbol\gamma}(v)/\widehat U(v,\boldsymbol\gamma)$, with $\widehat U$ a unit, as in~\eqref{eq:local-normalized-branches}. At a root $b$ of $P_{p,\boldsymbol\gamma}$ of multiplicity $r_b$, the spectral differential has residue $r_bb$ by the local calculation in Lemma~\ref{lem:fixed-residue-differential}. The roots in each disc, counted with multiplicity, are $\gamma_{p,0},\ldots,\gamma_{p,e_p-1}$. Therefore
\begin{equation}\label{eq:whole-disc-residue-cancellation}
 \int_{\partial\Delta_{p,\nu}(\boldsymbol\gamma)}
                 \eta_{\boldsymbol\gamma}
 =2\pi\ii\sum_b r_bb
 =2\pi\ii\sum_{j=0}^{e_p-1}\gamma_{p,j}=0.
\end{equation}
The boundary circles generate the kernel of the map from the first homology of the covering of the complement of the three puncture discs to $H_1(Y_{\boldsymbol\gamma};\Z)$. Integration over cycles in this cover therefore factors through $H_1(Y_{\boldsymbol\gamma};\Z)$.

We represent $A,B$ by cycles outside the puncture discs and continue their lifts over fixed paths in the base curve. The spectral roots on these paths depend holomorphically on $\boldsymbol\gamma$, so both periods are holomorphic on $S$. By homogeneity of the coefficients (Proposition~\ref{prop:coefficient-line}), we have
\begin{equation*}
 F_{1,\boldsymbol\gamma}=F_{1,0}+O(|\boldsymbol\gamma|^2)
\end{equation*}
outside the puncture discs. Since $\partial_yF_{1,0}$ is bounded away from zero there, implicit differentiation shows that the first derivatives of each continued root with respect to the residues vanish at zero. Since the projected paths are fixed, the first derivatives of the periods with respect to the residues also vanish at zero. This proves~\eqref{eq:small-residue-periods} on $S$.

The scaling of the differential gives $a(\lambda)=\lambda f_A(\sigma/\lambda)$ and $a_D(\lambda)=\lambda f_B(\sigma/\lambda)$. Substitution of the convergent Taylor series proves~\eqref{eq:fixed-residue-Laurent-periods}.
\end{proof}

By the differentiation formula in Lemma~\ref{lem:fixed-residue-differential}, the period ratio for the family with fixed residues is $\tau(\lambda)=\partial_\lambda a_D/\partial_\lambda a$. Consequently,
\begin{equation}\label{eq:semiflat-period-derivative-rates}
 \partial_\lambda a=A_0+O(|\lambda|^{-2}),\qquad
 \tau(\lambda)=\tau_\infty+O(|\lambda|^{-2}).
\end{equation}
Differentiating the convergent series gives $O(|\lambda|^{-2-k})$ after $k$ derivatives of either error. On a common smaller ball $S$, the same bounds hold for derivatives of any fixed order in the real residue coordinates, uniformly for $\sigma$ in any fixed compact set, by Cauchy's estimates.

We use the ordered basis $(A,B)$ and express the action of $\iota_*$ as
\begin{equation}\label{eq:reference-cycle-action}
 \begin{gathered}
 \iota_*\begin{pmatrix}A\\B\end{pmatrix}
       =M_n\begin{pmatrix}A\\B\end{pmatrix},\\
 M_3=\begin{pmatrix}0&1\\-1&-1\end{pmatrix},\quad
 M_4=\begin{pmatrix}0&1\\-1&0\end{pmatrix},\quad
 M_6=\begin{pmatrix}1&1\\-1&0\end{pmatrix}.
 \end{gathered}
\end{equation}
These matrices follow directly from multiplication of the primitive periods by $\zeta$: for $n=3,4$, $\zeta=\tau_\infty$, while for $n=6$, $\zeta=1+\tau_\infty$. In each case $M_n(1,\tau_\infty)^{\mathsf T}=\zeta(1,\tau_\infty)^{\mathsf T}$.

The descent map from~\eqref{eq:spectral-cyclic-descent} is
\begin{equation*}
 D_{\boldsymbol\gamma}:Y_{\boldsymbol\gamma}
       \longrightarrow Y_{\zeta^{-1}\boldsymbol\gamma},
 \qquad y\longmapsto\zeta^{-1}y.
\end{equation*}
At zero it equals $\iota^{-1}$. Its matrix in the marking of $H_1(Y_{\boldsymbol\gamma};\Z)$ is therefore $M_n^{-1}$ everywhere on the connected ball, since this matrix is integral and locally constant. It preserves the part of the covering outside the puncture discs. Changing representatives adds only boundary periods, which vanish by~\eqref{eq:whole-disc-residue-cancellation}. Moreover, $D_{\boldsymbol\gamma}^*\eta_{\zeta^{-1}\boldsymbol\gamma}=\zeta^{-1}\eta_{\boldsymbol\gamma}$. Writing $f=(f_A,f_B)^{\mathsf T}$, we obtain
\begin{equation}\label{eq:period-deck-action}
 M_n^{-1}f(\zeta^{-1}\boldsymbol\gamma)
       =\zeta^{-1}f(\boldsymbol\gamma),\qquad
 \begin{pmatrix}a(\zeta\lambda)\\a_D(\zeta\lambda)\end{pmatrix}
       =M_n\begin{pmatrix}a(\lambda)\\a_D(\lambda)\end{pmatrix}.
\end{equation}

The Picard descent sends $L$ to $(D_{\boldsymbol\gamma}^{-1})^*L$. It preserves $L_0=\pi^*\mathcal O_X(1)$ and hence acts as a group isomorphism on the relative Jacobian. Evaluation of the pulled-back unitary characters on the marked cycles gives
\begin{equation}\label{eq:angle-deck-action}
 \begin{pmatrix}\theta_A\\\theta_B\end{pmatrix}
       \longmapsto M_n\begin{pmatrix}\theta_A\\\theta_B\end{pmatrix}
       \quad\text{in }(\R/2\pi\Z)^2.
\end{equation}

With $b=\operatorname{Im}\tau_\infty$, smooth coordinates on the end are given by
\begin{equation}\label{eq:semiflat-ALG-coordinates}
 \begin{gathered}
 U=2\ii\sqrt b\,A_0\lambda,\qquad
 V=\frac{\theta_B-\tau_\infty\theta_A}{\sqrt b}
       \pmod{\Lambda_\infty},\\
 \Lambda_\infty=\frac{2\pi}{\sqrt b}(\Z+\tau_\infty\Z).
 \end{gathered}
\end{equation}
The angular coordinates use the marking transported by the Gauss--Manin connection. For each matrix in~\eqref{eq:reference-cycle-action}, equation~\eqref{eq:angle-deck-action} gives the deck action
\begin{equation}\label{eq:semiflat-ALG-deck-action}
 (U,V)\longmapsto(\zeta U,\zeta^{-1}V).
\end{equation}
In these coordinates the flat model and its radial coordinate are
\begin{equation}\label{eq:semiflat-ALG-model}
 \begin{aligned}
 g_0&=|\dd U|^2+|\dd V|^2,\\
 \omega_0&=\frac{\ii}{2}
            (\dd U\wedge \dd\bar U+\dd V\wedge \dd\bar V),\\
 \Omega_0&=\dd U\wedge \dd V,\\
 R&=|U|=2\sqrt b\,|A_0|\,|t|^{1/n}.
 \end{aligned}
\end{equation}
These tensors descend to the quotient by~\eqref{eq:semiflat-ALG-deck-action}. This is the standard ALG model of~\cite[Definition~2.5]{CV}, with $\beta=1/n$, $\tau=\tau_\infty$, and $L=2\pi/\sqrt b$.

\begin{proposition}\label{prop:semiflat-ALG-end}
For every fixed centered residue parameter $\sigma$, coordinates~\eqref{eq:semiflat-ALG-coordinates} give a smooth identification of the end of the Hitchin fibration with
\begin{equation*}
 \bigl(\{|U|>R_0\}\times\C/\Lambda_\infty\bigr)
       \big/\bigl((U,V)\sim(\zeta U,\zeta^{-1}V)\bigr).
\end{equation*}
In these coordinates, for every integer $k\geq0$,
\begin{equation}\label{eq:semiflat-ALG-rate}
 \begin{aligned}
 &\bigl|(\nabla^{g_0})^k(g_{\sfm}-g_0)\bigr|_{g_0}
  +\bigl|(\nabla^{g_0})^k(\omega_{\sfm}-\omega_0)\bigr|_{g_0}\\
 &\hspace{35mm}
  +\bigl|(\nabla^{g_0})^k(\Omega-\Omega_0)\bigr|_{g_0}
       \leq C_kR^{-2-k}.
 \end{aligned}
\end{equation}
For every compact set $K$ of centered residue parameters, the radius and constants can be chosen uniformly on $K$. For any fixed derivative $\partial_{\sigma}^{\ell}$ in the real coordinates of the residues, the differentiated tensor differences in this common chart satisfy the same estimate with constant $C_{K,k,\ell}$.

The model has cone angle $2\pi/n$ and fiber area $4\pi^2$. At $\sigma=0$, all three tensor differences vanish identically.
\end{proposition}

\begin{proof}
The relative Picard description and the invariant origin give smooth angular coordinates defined by this marking over the contractible ball $S$. Pulling back by $\boldsymbol\gamma=\sigma/\lambda$ and using~\eqref{eq:angle-deck-action} identifies the end with the indicated smooth quotient for all sufficiently large $|t|$.

After substituting $\partial_\lambda a=A_0$ and $\tau=\tau_\infty$ in~\eqref{eq:semiflat-period-formulas} and in the K\"ahler form of Proposition~\ref{prop:semiflat-period-formulas}, we obtain the tensors in~\eqref{eq:semiflat-ALG-model}. The fiber and symplectic normalizations follow from
\begin{equation*}
 \begin{aligned}
 \frac{\ii}{2}\dd V\wedge \dd\bar V&=\dd\theta_A\wedge \dd\theta_B,\\
 \dd U\wedge \dd V
       &=2\ii A_0\,\dd\lambda\wedge
                     (\dd\theta_B-\tau_\infty \dd\theta_A).
 \end{aligned}
\end{equation*}
For large $|\lambda|$, $\operatorname{Im}\tau$ has a positive lower bound. Each coefficient of a tensor difference is a smooth rational expression in $\partial_\lambda a$, $\tau$, and their conjugates, vanishing at their constant limits. By~\eqref{eq:semiflat-period-derivative-rates} and its derivatives, the bound is $O(|\lambda|^{-2-k})$ after any derivative of order $k$ in the Cartesian coordinates of the base.

The coefficients are independent of $\theta_A,\theta_B$, and the model metric has constant coefficients in the coframe of Cartesian differentials in $\lambda$ and angular differentials on the torus. The same bounds consequently hold for all $g_0$-covariant derivatives. The fixed conversion between $R$ and $|\lambda|$ in~\eqref{eq:semiflat-ALG-model} proves~\eqref{eq:semiflat-ALG-rate}. The bounds remain uniform on compact sets of residues after differentiation in the residue parameters, by the period estimates following Lemma~\ref{lem:exterior-period-expansion}. The finite deck action is an isometry preserving all three model tensors, so the estimates descend to the quotient.

The lattice $\Lambda_\infty$ has Euclidean covolume $4\pi^2$, and the base rotation has order $n$, giving the stated area and cone angle. For $\sigma=0$, $f_A$ and $f_B$ are evaluated only at zero; thus $a=A_0\lambda$ and $\tau=\tau_\infty$ exactly, and every tensor difference vanishes.
\end{proof}

\section{Local harmonic metrics}
\label{sec:local-models}\label{sec:harmonic-analysis}

At a puncture with flag length $e$, we use the rescaling $\xi=\lambda^e x$, where $t=\lambda^n$. In this coordinate, the Higgs bundle approaches a direct sum of rank-$e$ models. Each model has a Higgs matrix with diagonal entries $\sigma_j$, subdiagonal entries $1$, and upper-right entry $c\xi$, multiplied by $\dd\xi/\xi$. We derive the matrices for all three ranks from the spectral curve. The filtrations prescribe the growth of the metrics at zero and at the irregular singularity at infinity. Each summand corresponds to a branch of the normalization, and comparison with the metric on the spectral line bundle fixes the relative constants.

The model exists for all centered weights satisfying $\alpha_0<\cdots<\alpha_{e-1}<\alpha_0+1$ and all centered scalar graded residues. At infinity its eigenline bundles become asymptotically orthogonal, and the diagonal entries of the harmonic metric relative to the flat metric on the eigenline bundles have expansions in inverse powers of the radial variable, with coefficients depending on the angular variable. The comparison estimate below gives uniform bounds on compact parameter sets. When all graded residues vanish, cyclic symmetry gives exponential decay of the difference between these two metrics, so the metrics can be glued on fixed discs.

\subsection{The scale near a puncture}

At a fixed puncture $p$ we write $e=e_p$ and $m=m_p$, omitting $p$ from the local weights and residues. Thus $n=em$, $e\in\{2,3,4,6\}$, and
\begin{equation*}
 P(w)=\prod_{j=0}^{e-1}(w-\sigma_j),\qquad
 \sum_{j=0}^{e-1}\sigma_j=0.
\end{equation*}
In a coordinate $x$ vanishing at $p$ and the logarithmic frame $\dd x/x$, the spectral differential before rescaling is $\eta=w\,\dd x/x$. On the cover of the Hitchin base we use the variables
\begin{equation}\label{eq:local-model-scale}
 t=\lambda^n,\qquad T=|\lambda|,\qquad
 v=\lambda^{-1}w,\qquad
 \boldsymbol\gamma=\lambda^{-1}\sigma,\qquad
 \xi=\lambda^e x,\qquad \varepsilon=\lambda^{-e}.
\end{equation}
Here $\boldsymbol\gamma$ is the residue parameter of the rescaled spectral family in Proposition~\ref{prop:normalized-spectral-family}; the prescribed residues $\sigma$ remain fixed as $T$ tends to infinity.

In this frame the monic expansion of Lemma~\ref{lem:local-coefficients} and the homogeneity of the family~\eqref{eq:coefficient-line} give the local characteristic polynomial
\begin{equation}\label{eq:local-unscaled-pencil}
 F^{(p)}_{t,\sigma}(x,w)
 =P(w)^m+\sum_{\ell=1}^m
       x^\ell A_\ell(w,x,\sigma)P(w)^{m-\ell}
       -t x^m h_p(x).
\end{equation}
The function $h_p$ is a holomorphic unit. Each $A_\ell$ is holomorphic in $x$, polynomial in the residues, of degree less than $e$ in $w$, and homogeneous of total degree $e\ell$ in $(w,\sigma)$. Indeed, the coefficient of $t$ has vanishing order $m$ at $p$, by the formula for $r_n$ in~\eqref{eq:coefficient-line}. The remaining assertions follow from the monic $P$-adic expansion~\eqref{eq:local-monic-expansion} and~\eqref{eq:coefficient-scaling}.

\begin{proposition}\label{prop:local-polynomial-models}
In the coordinates~\eqref{eq:local-model-scale}, the finite modification of the spectral curve has equations
\begin{equation}\label{eq:local-rescaled-modification}
 \begin{aligned}
 P(w)&=\xi u,\\
 0&=u^m-h_p(\varepsilon\xi)
      +\sum_{\ell=1}^m\varepsilon^\ell
       A_\ell(w,\varepsilon\xi,\sigma)u^{m-\ell}.
 \end{aligned}
\end{equation}
The constants $c_1,\ldots,c_m$ are the distinct roots of $c^m=h_p(0)$. When $(\xi,w,\sigma)$ ranges over a compact set and $T$ is sufficiently large, the second equation has branches $u=W_\nu(\xi,w,\varepsilon,\sigma)$ satisfying
\begin{equation*}
 W_\nu=c_\nu+O(T^{-e}),\qquad 1\leq\nu\leq m,
\end{equation*}
with every fixed number of derivatives in $\xi,w$. The bounds are uniform on the chosen compact subsets. The limiting normalized curves, spectral differential, and ordered lattices on each branch are
\begin{equation}\label{eq:polynomial-spectral-model}
 \xi=\frac{P(w)}{c_\nu},\qquad
 \eta=w\frac{\dd\xi}{\xi},\qquad
 f_j(w)L,\quad f_j(w)=\prod_{i=0}^{j-1}(w-\sigma_i),
 \quad 0\leq j\leq e.
\end{equation}
Here $L$ is a spectral line bundle on the corresponding normalized component, and $f_j(w)L$ are its ordered submodules. Each curve is a smooth affine $w$-line. A root of $P$ of multiplicity $r$ has ramification index $r$ over $\xi=0$. A root $a$ of $P'$ of multiplicity $d$ has ramification index $d+1$, with branch value $P(a)/c_\nu$.
\end{proposition}

\begin{proof}
After substituting $x=\varepsilon\xi$ into~\eqref{eq:local-unscaled-pencil} we write $P(w)=\xi u$. Since $t\varepsilon^m=1$, division by $\xi^m$ on $\xi\ne0$ gives the second equation in~\eqref{eq:local-rescaled-modification}. This equation extends across $\xi=0$ and is the finite modification already constructed in~\eqref{eq:spectral-finite-modification}, with coordinate
\begin{equation*}
 \frac{P_{\boldsymbol\gamma}(v)}{x}
       =\frac{P(w)}{\xi}=u.
\end{equation*}
At $\varepsilon=0$, the $u$-equation is $u^m=h_p(0)$, whose roots are distinct and nonzero. The branches $W_\nu$ therefore exist by the holomorphic implicit function theorem, with the stated bounds on compact subsets. The equation $P(w)=\xi W_\nu$ has nonzero $\xi$ derivative at $\varepsilon=0$. It therefore defines a graph $\xi=\Xi_\nu(w,\varepsilon,\sigma)$ converging to $P(w)/c_\nu$.

With $\lambda$ fixed, $\dd x/x=\dd\xi/\xi$. The products for the rescaled residues satisfy $\prod_{i<j}(v-\gamma_i)=\lambda^{-j}f_j(w)$. Multiplication by this nonzero constant leaves the lattices in~\eqref{eq:newton-spectral-lattices} unchanged. This proves both the differential and filtration assertions in~\eqref{eq:polynomial-spectral-model}. The ramification indices follow by differentiating $\xi=P(w)/c_\nu$.
\end{proof}

We compute the ramification from the roots of $P'$, with arbitrary multiplicities. As $a$ runs through the distinct roots of $P'$, with respective multiplicities $d_a$, the branch divisor of the $\nu$th model is the zero divisor of the polynomial
\begin{equation}\label{eq:polynomial-branch-divisor}
 \operatorname{Res}_w(P(w)-c_\nu\xi,P'(w))
       =e^e\prod_a(P(a)-c_\nu\xi)^{d_a}.
\end{equation}
If $P$ has $s$ distinct roots, the total ramification multiplicity over zero is $e-s$ and that away from zero is $s-1$. Two branch values in the same or different blocks coincide precisely when the corresponding values of $P(a)/c_\nu$ agree. Distinct blocks have no common eigenvalue for $\xi\ne0$ since $c_\nu\ne c_\kappa$. Over zero their coincident eigenvalues lie on separate components of the normalization.

At finite $T$, the graph equation also records any splitting of a higher-order critical point. Its derivative is
\begin{equation}\label{eq:finite-scale-critical-equation}
 \partial_w\Xi_\nu
   =\frac{P'(w)-\Xi_\nu\partial_wW_\nu}
          {W_\nu+\Xi_\nu\partial_\xi W_\nu},
\end{equation}
with $W_\nu$ evaluated on the graph. The denominator stays nonzero on each fixed compact set for large $T$. At an $r$-fold root $b$ of $P$, the function $\Xi_\nu=P/W_\nu$ vanishes to order $r$. Near any other critical point of multiplicity $d$, all nearby critical points, counted with multiplicity, are the roots of the Weierstrass polynomial of the numerator in~\eqref{eq:finite-scale-critical-equation}. The estimates below require no lower bound on the distances between these critical points.

\begin{proposition}\label{prop:local-spectral-matching}
For every compact set $B$ of residue parameters there exist $R,\delta,c_*>0$ and $T_0$ such that, for $\sigma\in B$ and $T\geq T_0$, all $n$ exact eigenvalues in the puncture disc $|x|\leq\delta$ are distinct on
\begin{equation}\label{eq:outer-spectral-annulus}
 R\leq|\xi|\leq\delta T^e,
\end{equation}
and their pairwise distances are at least $c_*|\xi|^{1/e}$. Every ramification point over $|x|\leq\delta$ lies over $|x|<RT^{-e}$.

For fixed $0<a_0<a_1$ and $0<\theta<1$ the following estimates hold on the annulus
\begin{equation}\label{eq:spectral-matching-annulus}
 a_0T^{-e\theta}\leq|x|\leq a_1T^{-e\theta},
\end{equation}
with a local choice of sheets and a root $\chi^e=\xi$. If $w^{\ex}$ is an exact eigenvalue, $w^{\mathrm{mod}}$ is the associated root of $P(w)=c_\nu\xi$, and $\kappa_{\nu,j}^e=c_\nu$ labels their common leading term, then
\begin{equation}\label{eq:spectral-matching-errors}
 \frac{w^{\ex}-w^{\mathrm{mod}}}{\chi}
       =O_B(T^{-e\theta}),\qquad
 \frac{w^{\mathrm{mod}}}{\chi}-\kappa_{\nu,j}
       =O_B(T^{-2(1-\theta)}).
\end{equation}
Both estimates hold with every fixed number of derivatives in the rescaled annulus coordinate $T^{e\theta}x$, on local annular sectors. The constants may depend on the derivative order, $B$, $a_0$, $a_1$, and $\theta$.
\end{proposition}

\begin{proof}
The polynomial has the expansion $P(w)=w^e+\sum_{k=2}^e p_k(\sigma)w^{e-k}$. The coefficients $p_k$ are bounded on $B$. Under $w=\chi z$, the model equation becomes
\begin{equation*}
 z^e+\sum_{k=2}^e p_k(\sigma)\chi^{-k}z^{e-k}
       =c_\nu.
\end{equation*}
Its roots are $z=\kappa_{\nu,j}+O_B(|\xi|^{-2/e})$. Taken over all blocks, the $\kappa_{\nu,j}$ are the $n$ distinct roots of $z^n=h_p(0)$.

The degree bound in~\eqref{eq:local-unscaled-pencil} gives $|A_\ell(w,x,\sigma)|\leq C_B(1+|w|)^{e-1}$ for $|x|\leq\delta$. After eliminating $u$ from~\eqref{eq:local-rescaled-modification} and dividing by $\xi^m$, we substitute $w=\chi z$. On bounded sets of $z$, the resulting polynomial differs from $z^n-h_p(0)$ by at most
\begin{equation*}
 C_B\left(|\xi|^{-2/e}+T^{-e}|\xi|
                 +T^{-e}(1+|\xi|^{(e-1)/e})\right).
\end{equation*}
We enclose the roots $\kappa_{\nu,j}$ in pairwise disjoint circles of fixed radii. Taking $R$ large, then $\delta$ small, and then $T_0$ large makes this perturbation smaller than $|z^n-h_p(0)|$ on every circle throughout~\eqref{eq:outer-spectral-annulus}. Each circle therefore encloses exactly one root, counted with multiplicity, by Rouch\'e's theorem. These are all $n$ roots, so they are simple and uniformly separated in $z$. Rescaling by $\chi$ proves the claimed eigenvalue separation and locates all ramification in the inner disc.

On these roots, the equation for the branch in the variable $u$ gives $W_\nu-c_\nu=O_B(T^{-e}|\xi|)$. Hence
\begin{equation*}
 \frac{P(w^{\ex})-c_\nu\xi}{\xi}
             =O_B(T^{-e}|\xi|).
\end{equation*}
Since the rescaled model roots are uniformly simple inside these circles, $(w^{\ex}-w^{\mathrm{mod}})/\chi=O_B(T^{-e}|\xi|)$ follows by the implicit function theorem. On~\eqref{eq:spectral-matching-annulus}, comparability of $|\xi|$ and $T^{e(1-\theta)}$ gives the first estimate in~\eqref{eq:spectral-matching-errors}, and the expansion of the model roots gives the second. The roots are distinct and admit holomorphic labels on slightly larger annular sectors, where the derivative bounds follow from Cauchy's estimates in $T^{e\theta}x$.
\end{proof}

\subsection{Model filtered Higgs bundles}\label{subsec:filtered-polynomial-models}

For a nonzero constant $c$ and local data $(\alpha_j,\sigma_j)_{j=0}^{e-1}$, we take on $\Pp^1_\xi$ the bundle $\mathcal E_c=\mathcal O_{\Pp^1}^{\oplus e}$ with frame $f_0,\ldots,f_{e-1}$. Over the affine $\xi$-line, the Newton polynomials $f_j(w)=\prod_{i<j}(w-\sigma_i)$, $0\leq j<e$, form a $\C[\xi]$-basis of $\C[\xi,w]/(P(w)-c\xi)$ and define the frame $(f_0,\ldots,f_{e-1})$. Multiplication by $w$ defines the Higgs field
\begin{equation}\label{eq:polynomial-model-matrix}
 \begin{aligned}
 A_c&=\diag(\sigma_0,\ldots,\sigma_{e-1})
       +\sum_{j=0}^{e-2}E_{j+1,j}+c\xi E_{0,e-1},\\
 \Phi_c&=A_c\frac{\dd\xi}{\xi}.
 \end{aligned}
\end{equation}
Here $E_{ab}$ sends $f_b$ to $f_a$. Indeed, $wf_j=\sigma_jf_j+f_{j+1}$ for $j<e-1$, while $wf_{e-1}=\sigma_{e-1}f_{e-1}+c\xi f_0$. In particular,
\begin{equation*}
 \det(q\Id-A_c)=P(q)-c\xi.
\end{equation*}
In the following four matrices the superscript denotes the size of the block.
\begin{equation}\label{eq:local-model-small-matrices}
 A_c^{(2)}=
 \begin{pmatrix}
  \sigma_0&c\xi\\
  1&\sigma_1
 \end{pmatrix},\qquad
 A_c^{(3)}=
 \begin{pmatrix}
  \sigma_0&0&c\xi\\
  1&\sigma_1&0\\
  0&1&\sigma_2
 \end{pmatrix}.
\end{equation}
\begin{equation}\label{eq:local-model-large-matrices}
 A_c^{(4)}=
 \begin{pmatrix}
  \sigma_0&0&0&c\xi\\
  1&\sigma_1&0&0\\
  0&1&\sigma_2&0\\
  0&0&1&\sigma_3
 \end{pmatrix},\qquad
 A_c^{(6)}=
 \begin{pmatrix}
  \sigma_0&0&0&0&0&c\xi\\
  1&\sigma_1&0&0&0&0\\
  0&1&\sigma_2&0&0&0\\
  0&0&1&\sigma_3&0&0\\
  0&0&0&1&\sigma_4&0\\
  0&0&0&0&1&\sigma_5
 \end{pmatrix}.
\end{equation}

For the coordinates $x=z$, $z-1$ and $1/z$ at $p=0,1,\infty$, respectively, the equality $\eta=y\,\dd z/[z(z-1)]=w\,\dd x/x$ gives
\begin{equation*}
 w=\frac{y}{x-1}\quad(p=0),\qquad
 w=\frac{y}{1+x}\quad(p=1),\qquad
 w=-\frac{xy}{1-x}\quad(p=\infty).
\end{equation*}
Substitution into~\eqref{eq:coefficient-line} gives the values of the holomorphic units in~\eqref{eq:local-unscaled-pencil}. The triples $(h_p(0))_{p=0,1,\infty}$ are $(1,1,-1)$, $(-1,1,1)$ and $(1,1,1)$ for $n=3,4,6$, respectively. The roots $c_\nu^{m_p}=h_p(0)$ give the following limiting matrices at the three punctures.
\begin{equation}\label{eq:local-model-block-table}
 \begin{array}{c|ccc}
  n&p=0&p=1&p=\infty\\ \hline
  3&A_1^{(3)}&A_1^{(3)}&A_{-1}^{(3)}\\[2pt]
  4&A_{-1}^{(4)}&A_1^{(4)}&A_1^{(2)}\oplus A_{-1}^{(2)}\\[2pt]
  6&A_1^{(6)}&A_1^{(3)}\oplus A_{-1}^{(3)}
       &A_1^{(2)}\oplus A_\omega^{(2)}\oplus A_{\omega^2}^{(2)}
 \end{array}
 \qquad \omega=e^{2\pi\ii/3}.
\end{equation}
Each entry uses the residues at its puncture and is multiplied by $\dd\xi/\xi$, with $\xi=\lambda^{e_p}x$. The direct sum has $m_p$ blocks of size $e_p$ and characteristic polynomial $P(q)^{m_p}-h_p(0)\xi^{m_p}$.

At zero, the residue of a block is $A_c(0)=\diag(\sigma_j)+\sum_{j=0}^{e-2}E_{j+1,j}$. For the descending flag $F^j=\langle f_j,\ldots,f_{e-1}\rangle$ with weight $\alpha_j$, the induced residue on $F^j/F^{j+1}$ is $\sigma_j$. The direct sum of these flags is the given partial flag with graded quotients of dimension $m_p$. The Newton basis is related to the monomial basis by a triangular matrix with determinant one, so these formulas remain valid when residues coincide. The Higgs field is logarithmic at zero. At infinity, with $u=\xi^{-1}$, we prescribe the growth exponents
\begin{equation}\label{eq:model-infinity-weights}
 \delta_j=\frac{e-1-2j}{2e},\qquad 0\leq j<e.
\end{equation}

In a local coordinate $z$, the frame $f_j$ and growth exponents $d_j$ determine the filtered lattices
\begin{equation}\label{eq:model-filtered-lattices}
 \mathcal P_a\mathcal E_c
       =\sum_j\mathcal O\,z^{-\lfloor a+d_j\rfloor}f_j,
 \qquad a\in\R.
\end{equation}
This is the filtration induced by the diagonal reference metric $|f_j|=|z|^{d_j}$. Membership in $\mathcal P_a$ means growth $O(|z|^{-a-\epsilon})$ for every $\epsilon>0$. We use $d_j=\alpha_j$ at zero and $d_j=\delta_j$ at infinity. An adapted metric is one inducing the filtration~\eqref{eq:model-filtered-lattices}, as in~\cite[Definition~2.11]{MochizukiGood}.

\begin{proposition}\label{prop:polynomial-harmonic-model}
The data in~\eqref{eq:polynomial-model-matrix} and~\eqref{eq:model-filtered-lattices} define a stable good filtered Higgs bundle of degree zero, with trace-free Higgs field and trivial filtered determinant line bundle. It admits a unique adapted harmonic metric $h_c$ satisfying $\det h_c=1$ in the frame $f_0\wedge\cdots\wedge f_{e-1}$.

On the cover $u=r^e$ at infinity, the irregular values are $e\kappa_i/r$, where $\kappa_i$ runs through the $e$ distinct roots of $c$; after subtracting their differentials, the eigenvalues of $\Phi_c$ are holomorphic one-forms. These assertions hold for all coincidences and zeros among the $\sigma_j$.

For the distinct constants $c_1,\ldots,c_m$ in Proposition~\ref{prop:local-polynomial-models}, the direct sum of these filtered Higgs bundles is polystable, with filtered Higgs automorphism group $(\C^*)^m$. Its adapted harmonic metrics are precisely
\begin{equation}\label{eq:model-block-constants}
 \bigoplus_{\nu=1}^m s_\nu h_{c_\nu},\qquad s_\nu>0.
\end{equation}
For the determinant metric of the direct sum to equal one, the constants must satisfy $\prod_{\nu=1}^m s_\nu^e=1$, leaving $m-1$ positive real parameters. Requiring the determinant metric of each block to equal one fixes all these constants.
\end{proposition}

\begin{proof}
At zero, the matrix $A_c$ preserves every lattice in~\eqref{eq:model-filtered-lattices}. For the subdiagonal terms this follows from $\lfloor a+\alpha_j\rfloor\leq\lfloor a+\alpha_{j+1}\rfloor$; for the term $c\xi f_0$ it follows from $\lfloor a+\alpha_{e-1}\rfloor-\lfloor a+\alpha_0\rfloor\leq1$. Thus the Higgs field is logarithmic on every lattice, and its graded residues indexed by the flag are $\sigma_j$.

At infinity we write $g_j=r^jf_j$. By~\eqref{eq:model-infinity-weights}, every $g_j$ has growth exponent $(e-1)/2$ in $r$. The spectral roots have the form $w_i=v_i(r)/r$, where
\begin{equation}\label{eq:model-infinity-roots}
 \prod_{j=0}^{e-1}(v_i-r\sigma_j)=c,\qquad
 v_i(0)=\kappa_i,\qquad \kappa_i^e=c.
\end{equation}
Since the roots $\kappa_i$ are distinct, each $v_i$ is holomorphic, and differentiation at zero using $\sum_j\sigma_j=0$ gives $v_i'(0)=0$. Evaluation at the roots takes $g_j$ to $\prod_{a<j}(v_i-r\sigma_a)$. The resulting holomorphic matrix has invertible value $(\kappa_i^j)$ at $r=0$, so the filtration splits as a sum of filtered eigenline bundles of common growth exponent $(e-1)/2$. The eigenvalues of $\Phi_c$ satisfy
\begin{equation}\label{eq:model-irregular-values}
 w_i\frac{\dd\xi}{\xi}
       =-e\kappa_i\frac{\dd r}{r^2}+O(1)\,\dd r
       =\dd\left(\frac{e\kappa_i}{r}\right)+O(1)\,\dd r.
\end{equation}
The differences of distinct irregular values have simple poles with nonzero leading coefficients, and the remainder on each eigenline bundle is holomorphic, as required by the goodness conditions of~\cite[Section~2.4.2]{MochizukiGood}. The deck action $r\mapsto\omega r$, $\omega^e=1$, permutes the eigenline bundles and multiplies $g_j$ by the scalar $\omega^j$. The common weight makes the filtration invariant, and it descends to~\eqref{eq:model-infinity-weights}.

Since both $\sum_j\alpha_j$ and $\sum_j\delta_j$ vanish, the filtered determinant line bundle is trivial with zero growth exponents and the filtered degree is zero. The trace of the Higgs field is $\sum_j\sigma_j=0$. The characteristic polynomial $P(w)-c\xi$ is irreducible over $\C(\xi)$ by Gauss' lemma, since it is monic and its quotient ring over $\C[\xi]$ is $\C[w]$. The generic spectral algebra is thus a field acting on itself and has no nonzero proper invariant subspace over $\C(\xi)$. The absence of proper invariant filtered subbundles proves stability.

By the wild Kobayashi--Hitchin correspondence recalled in~\cite[Theorem~1.2]{LiMochizuki}, there is a harmonic metric inducing the prescribed filtration. Any two adapted harmonic metrics admit a common orthogonal Higgs decomposition and differ by a positive scalar on each summand. The orthogonal projections are bounded and holomorphic, so they preserve the growth filtration. Generic irreducibility makes this decomposition trivial. The metric is therefore unique up to one positive constant. Its determinant is a harmonic metric on the trivial filtered line bundle with zero Higgs field. In rank one such a metric is unique up to a positive constant, so a constant rescaling gives the unit metric on the trivialized determinant line bundle.

For distinct $c_\nu$, the characteristic polynomials $P(w)-c_\nu\xi$ are pairwise coprime over $\C(\xi)$, so all Higgs intertwiners between different factors vanish. A filtered endomorphism of one factor has characteristic coefficients holomorphic on $\Pp^1$, hence constant. Subtracting a constant eigenvalue gives an endomorphism with zero determinant, which vanishes by generic irreducibility. Thus every filtered endomorphism is scalar. Applying the cited uniqueness theorem to the direct sum yields~\eqref{eq:model-block-constants}, and taking determinants gives the stated constraint.
\end{proof}

\subsection{The harmonic metric at infinity in the local model}

For the model with data $(e,c,\alpha,\sigma)$ in Proposition~\ref{prop:polynomial-harmonic-model}, we work on the cover $\xi=r^{-e}$. In terms of the roots $w_i=v_i(r)/r$ from~\eqref{eq:model-infinity-roots}, the evaluation matrices are
\begin{equation}\label{eq:model-evaluation-matrices}
 W_{ij}(r)=f_j(w_i),\qquad
 V_{ij}(r)=r^jf_j(w_i)
          =\prod_{a<j}(v_i-r\sigma_a).
\end{equation}
The matrix $V$ is holomorphic and invertible at zero. Near infinity, where the spectral roots are distinct, the columns of $W^{-1}$ form a holomorphic frame $\ell_i$ of eigenvectors. In the frame $(f_0,\ldots,f_{e-1})$ we define the metric by the matrix
\begin{equation}\label{eq:model-flat-eigenline-metric}
 \begin{aligned}
 h_0&=|\det W|^{-2/e}W^*W,\\
 h_0(\ell_i,\ell_j)&=d(r)\delta_{ij},\qquad
 d(r)=|\det W|^{-2/e}=|r|^{e-1}|\det V(r)|^{-2/e}.
 \end{aligned}
\end{equation}
Thus the eigenline bundles are orthogonal, and the common squared norm $d(r)$ makes the determinant equal to one in this frame.

\begin{proposition}\label{prop:polynomial-model-tail}
The metric $h_0$ is invariant under the deck action, flat on each eigenline bundle, and has determinant one. It induces the prescribed filtration~\eqref{eq:model-infinity-weights} at infinity. In the coordinates $\zeta=r^{-1}$ and $\rho=|\zeta|$, for every integer $k\geq0$ there are constants $C_k,R_k$ depending on the fixed model such that the harmonic metric $h_c$ with determinant one satisfies
\begin{equation}\label{eq:polynomial-model-tail}
 \bigl|\nabla_0^k(h_0^{-1}h_c-\Id)\bigr|_{h_0,|\dd\zeta|^2}
       \leq C_k\rho^{-1-k},\qquad \rho\geq R_k.
\end{equation}
Here $\nabla_0$ combines the Chern connection of $h_0$ on $\End\mathcal E_c$ with the Euclidean connection in the $\zeta$ coordinate. More precisely, the matrix
\begin{equation}\label{eq:model-normalized-metric-matrix}
 G_{ij}(r)=\frac{h_c(\ell_i,\ell_j)}{d(r)}
\end{equation}
extends smoothly to $r=0$, with $G(0)=\Id$. Its off-diagonal entries vanish to every order there, and every Taylor polynomial of $\log G_{ii}$ is harmonic.

For the sum of $m$ blocks in~\eqref{eq:model-block-constants}, the corresponding orthogonal sum of flat metrics is $\bigoplus_\nu s_\nu h_{0,\nu}$, with the same blockwise estimates and the determinant constraint in Proposition~\ref{prop:polynomial-harmonic-model}.
\end{proposition}

\begin{proof}
The identity $\det W=r^{-e(e-1)/2}\det V$ gives the second expression for $d(r)$ in~\eqref{eq:model-flat-eigenline-metric}. Since $\det W$ is holomorphic and nonzero near infinity, $\log d$ is harmonic there. Thus the metric on each eigenline bundle is flat, and their orthogonal sum is harmonic for the diagonal Higgs field. In the frame $g_j=r^jf_j$, the matrix is $dV^*V$, uniformly comparable to $|r|^{e-1}\Id$. Hence $h_0$ induces the prescribed filtration at infinity already checked in Proposition~\ref{prop:polynomial-harmonic-model}. Under $r\mapsto\omega r$, $\omega^e=1$, the rows of $W$ are permuted and the columns of $V$ acquire factors $\omega^j$. Both $d$ and $h_0$ are invariant.

The metrics $h_c$ and $h_0$ define good wild harmonic bundles with the same filtration~\eqref{eq:model-filtered-lattices} near infinity. By~\cite[Theorem~2.21]{MochizukiGood}, there are positive constants $C,C'$ such that $Ch_0\leq h_c\leq C'h_0$ in a neighborhood of infinity. In particular, the functions $q_i=\log G_{ii}$ are bounded there.

We apply local decay estimates on discs whose radius is comparable to $\rho$. By~\eqref{eq:model-irregular-values}, the eigenvalues of $\Phi_c$ in the $\zeta$ coordinate are
\begin{equation*}
 \bigl(e\kappa_i+O(\zeta^{-2})\bigr)\,\dd\zeta.
\end{equation*}
For $|\zeta_0|$ large, we use $\zeta=\zeta_0+\frac14|\zeta_0|z$ on $|z|<1$. The coefficients in this coordinate lie within $O(|\zeta_0|^{-1})$ of the constants $\frac14|\zeta_0|e\kappa_i$. Their minimum separation is a positive constant times $|\zeta_0|$, and their moduli are bounded by a fixed multiple of this separation. Thus Assumption~2.2 of~\cite{MochizukiAsymptotics} holds for all sufficiently large $|\zeta_0|$, with constants depending on the fixed model.

On smaller discs the eigenline bundles are exponentially close to orthogonal by Proposition~2.3 and Corollary~2.6 of~\cite{MochizukiAsymptotics}. The curvatures of their induced metrics satisfy the bounds in Proposition~2.12. In the holomorphic frames of Section~2.2.5, these bounds hold for derivatives of every fixed order by Corollary~2.15. Since $\log d$ is harmonic, the curvature coefficient of the metric on the $i$th eigenline bundle is, up to sign, $\partial_r\partial_{\bar r}q_i$. Changing variables from $z$ to $r$ introduces only powers of $|\zeta_0|$. Thus, for every mixed derivative $D^k$ in $r,\bar r$, there are constants $C_k,N_k,\epsilon_k>0$ such that
\begin{equation*}
 \bigl|D^k(\partial_r\partial_{\bar r}q_i)\bigr|
       \leq C_k|r|^{-N_k}\exp(-\epsilon_k/|r|).
\end{equation*}
The right-hand side is smaller than every positive power of $|r|$. Therefore $\partial_r\partial_{\bar r}q_i$ extends smoothly across zero, with all derivatives zero.

The resulting Poisson equation admits a smooth local solution $u_i$ on the full disc. The difference $q_i-u_i$ is bounded and harmonic on the punctured disc, so it extends harmonically over zero. Hence $q_i$ is smooth at zero. Since its Laplacian vanishes to every order there, each Taylor polynomial of $q_i$ is harmonic.

We also transfer the off-diagonal derivative estimates to the frames $\ell_i$. On each rescaled disc there is a nonvanishing holomorphic function $b$ with $|b|^2=d$ because $\log d$ is harmonic. Then $\ell_i/b$ is an $h_0$-unit holomorphic eigenframe. The $h_c$ norm of $\ell_i/b$ is bounded above and below by positive constants independent of the disc center. The eigenframes used in Section~2.2.5 of~\cite{MochizukiAsymptotics} satisfy the same bounds on smaller discs. The corresponding frame vectors differ by nonvanishing holomorphic factors that are uniformly bounded along with their inverses. On smaller discs the derivatives of these factors are uniformly bounded by Cauchy's estimates. The exponential bounds in Corollary~2.15 therefore hold for all derivatives of the off-diagonal entries of $G$ in these frames. Changing variables to $r$ again introduces only polynomial factors, so the off-diagonal entries extend smoothly and vanish to every order at zero.

A primitive deck transformation permutes the eigenline bundles in an $e$-cycle. Deck invariance therefore makes all numbers $q_i(0)$ equal. Since $\det G=\det h_c/\det h_0=1$ and the off-diagonal entries vanish at zero, $\sum_iq_i(0)=0$. Hence $q_i(0)=0$ for every $i$. The smooth extension of the metric matrix then gives $G=\Id+O(|r|)$, with the derivative estimates obtained above.

In the $\ell_i$ frame the Chern connection of $h_0$ is the common scalar $\partial\log d$, so its induced connection on endomorphisms is entrywise differentiation. Applying the chain rule to the smooth matrix $G-\Id$, which vanishes at $r=0$, gives~\eqref{eq:polynomial-model-tail} under $r=\zeta^{-1}$. Indeed, a derivative of order $k\geq1$ involving $j\geq1$ derivatives in $r,\bar r$ has coefficient $O(\rho^{-k-j})$, while the undifferentiated matrix is $O(\rho^{-1})$. The assertion for a direct sum follows by applying this argument to each factor in~\eqref{eq:model-block-constants}.
\end{proof}

\subsection{Existence and behavior of harmonic metrics at the punctures}
\label{sec:harmonic-metrics}

We fix the weights and scalar graded residues in~\eqref{eq:graded-residue}. The local estimates below allow repeated residues and partial flags. They also specify the asymptotic behavior of the harmonic metric used to correct the approximate solution.

Near each puncture $p$ we use a coordinate $z$ vanishing at $p$ and a holomorphic frame $(e_1,\ldots,e_n)$ compatible with the flag and the fixed determinant trivialization. Repeating the weights and residues according to their multiplicities, we write
\begin{equation}\label{eq:harmonic-local-data}
 A=\diag(\alpha_1,\ldots,\alpha_n),\qquad
 B=\diag(\sigma_1,\ldots,\sigma_n),\qquad
 \Phi=R(z)\frac{\dd z}{z}.
\end{equation}
Thus $R(0)$ preserves the flag and induces $\sigma_{p,j}\Id$ on its $j$th graded quotient. In the coordinates $r=|z|$ and $\tau=-\log r$ we use
\begin{equation}\label{eq:harmonic-power-model}
 h_0=\diag(r^{2\alpha_1},\ldots,r^{2\alpha_n}),\qquad
 f_a=r^{-\alpha_a}e_a.
\end{equation}
The frame $(f_a)$ is unitary for $h_0$. In this frame the Chern connection of $h_0$ and the Higgs field are
\begin{equation}\label{eq:harmonic-cylinder-model}
 D_0=\dd+\ii A\,\dd\theta,\qquad
 \Phi=(B+q)(-\dd\tau+\ii\,\dd\theta),\qquad
 q=r^A R(z)r^{-A}-B.
\end{equation}
All cylindrical norms below use $\dd\tau^2+\dd\theta^2$, with $\theta$ of period $2\pi$. The positive number used in the decay estimates is
\begin{equation}\label{eq:harmonic-weight-frequency}
 \varepsilon_p=
 \min\bigl\{|k+\alpha_a-\alpha_b|:
 k\in\Z,\ 1\leq a,b\leq n,\ k+\alpha_a-\alpha_b\ne0\bigr\}>0.
\end{equation}
The width of the weight interval is less than one. Consequently $\varepsilon_p$ is the minimum of $1$, the positive weight differences, and one minus the positive weight differences. Flag preservation gives $R_{ab}(0)=0$ when $\alpha_a<\alpha_b$. On an equal-weight block, $R(0)=B$. Taylor expansion therefore gives, for every $k\geq0$,
\begin{equation}\label{eq:harmonic-taylor-decay}
 \sum_{i+j\leq k}
 |\partial_\tau^i\partial_\theta^j q|
 \leq C_k e^{-\varepsilon_p\tau}.
\end{equation}

\subsubsection{Local comparison and decay}

\begin{lemma}\label{lem:tame-power-comparison}
For any tame harmonic metric $h$ inducing the prescribed filtered bundle on a punctured disc with local data~\eqref{eq:harmonic-local-data}, the metrics $h$ and $h_0$ are mutually bounded on a smaller punctured disc. This conclusion holds for every holomorphic frame compatible with the flag.
\end{lemma}

\begin{proof}
For this Higgs bundle we use Mochizuki's norm estimate~\cite[Proposition~8.1.1, equation~(157), and Section~8.1.2]{MochizukiWild}. The irregular value is zero, and the nilpotent part of the residue on each parabolic graded quotient is zero. Hence the associated nilpotent weight filtration is concentrated in degree zero. In the notation of that estimate, $a(e_a)=-\alpha_a$ and $k(e_a)=0$, so the metric in the norm estimate is the diagonal metric~\eqref{eq:harmonic-power-model}. The local decomposition by distinct residue eigenvalues required there is obtained from the holomorphic spectral projections of $R(z)$ after shrinking the disc. A repeated eigenvalue remains in one summand.

Because the nilpotent weight filtration is concentrated in degree zero, every holomorphic frame compatible with the flag satisfies the hypotheses of Mochizuki's estimate, so the comparison holds in each such frame.
\end{proof}

The following comparison identity applies both on a punctured disc and on $X\setminus D$.

\begin{lemma}\label{lem:harmonic-trace-comparison}
For two harmonic metrics $h_1,h_2$ on the same Higgs bundle, the function $\tr(h_1^{-1}h_2)$ is subharmonic. If the two metrics have determinant one, are mutually bounded near the puncture, and agree on the outer boundary of a disc, then they agree throughout the punctured disc.
\end{lemma}

\begin{proof}
With $S=h_1^{-1}h_2$, we work at a given point in a holomorphic frame normal for $h_1$ in which $S$ is diagonal with eigenvalues $\ell_a>0$. Expressing $\Phi$ as $\phi\,\dd z$ in this frame, we obtain from Simpson's metric comparison formula~\cite[Lemma~3.1(c) and proof of Lemma~3.1(d), p.~876]{Simpson1988}
\begin{equation}\label{eq:harmonic-trace-identity}
 \partial_{\bar z}\partial_z\tr S
 =\tr\bigl((\partial_z S)^\dagger S^{-1}\partial_z S\bigr)
  +\sum_{a,b}\frac{(\ell_a-\ell_b)^2}{\ell_b}|\phi_{ab}|^2
 \geq0.
\end{equation}
In a general frame the first derivative is the covariant derivative for $h_1$. Mutual boundedness makes $\tr S$ bounded near the puncture, so it extends subharmonically across that point. Its boundary value is $n$. By the maximum principle and the inequality $\tr S\geq n$, which follows from $\det S=1$, we have $S=\Id$.
\end{proof}

We use $H^s$ for the unweighted $L^2$ Sobolev space. For nonnegative integers $k$ we have $H^k=W^{k,2}$, and on $S^1$ we use Fourier series to define the spaces of fractional order.

We first solve a boundary problem on the half-cylinder. Here $V$ is the real vector space of trace-free Hermitian matrices, and $K\subset V$ is the subspace of constant matrices that preserve the equal-weight subspaces. For a partial flag, each Hermitian block in $K$ is arbitrary, subject to the total trace condition. For $\nu>0$ the weighted spaces and norms are
\begin{equation*}
 H^j_\nu=\{v:e^{\nu x}v\in H^j([0,\infty)\times S^1;V)\},
 \qquad \|v\|_{H^j_\nu}=\|e^{\nu x}v\|_{H^j}.
\end{equation*}

\begin{lemma}\label{lem:harmonic-half-cylinder}
For each $0<\nu<\varepsilon_p$, if $\tau_0$ is sufficiently large and $b\in H^{7/2}(S^1;V)$ sufficiently small, there is a harmonic metric of determinant one on $\tau\geq\tau_0$ whose matrix in the frame $(f_a)$ is
\begin{equation}\label{eq:harmonic-half-cylinder-metric}
 \widetilde H=e^{c+v},\qquad
 c\in K,\quad v\in H^4_\nu,\quad
 \widetilde H(\tau_0,\theta)=e^{b(\theta)},
 \qquad x=\tau-\tau_0.
\end{equation}
Here $c$ and $v$ are small in the indicated spaces. Every cylindrical derivative of $\widetilde H-e^c$ is $O(e^{-\nu x})$ for $x\geq1$. The smallness threshold and $\tau_0$ can be chosen uniformly for constant changes of frame in a compact set of positive block diagonal matrices of determinant one.
\end{lemma}

\begin{proof} We solve the boundary problem in Appendix~\ref{proof:lem:harmonic-half-cylinder}. \end{proof}

\begin{lemma}\label{lem:local-harmonic-exponential-decay}
For a harmonic metric $h$ of determinant one on a punctured disc with data~\eqref{eq:harmonic-local-data} satisfying $C^{-1}h_0\leq h\leq Ch_0$ for some $C\geq1$, we denote its matrix in the frame $(f_a)$ by $H$. There is a unique positive Hermitian block diagonal matrix $H_\infty$, with determinant one, such that, for every $0<\nu<\varepsilon_p$ and every $k\geq0$,
\begin{equation}\label{eq:harmonic-local-exponential-decay}
 \sum_{i+j\leq k}
 |\partial_\tau^i\partial_\theta^j(H-H_\infty)|
 \leq C_{k,\nu}e^{-\nu\tau}
\end{equation}
on a sufficiently small punctured disc. The blocks are the entire equal-weight subspaces, including when residues coincide. The constants may depend on the metric and the fixed local data.
\end{lemma}

\begin{proof} We prove the decay estimate in Appendix~\ref{proof:lem:local-harmonic-exponential-decay}. \end{proof}

At $\xi=0$ we can apply Lemmas~\ref{lem:tame-power-comparison} and~\ref{lem:local-harmonic-exponential-decay} to the harmonic metrics on the local models of Proposition~\ref{prop:polynomial-harmonic-model}. The graded quotients of these models are one-dimensional, so the limiting metric matrices are positive diagonal matrices. The constants in this application depend on the fixed model parameters.

\subsubsection{Existence and the asymptotic form}

\begin{proposition}\label{prop:harmonic-metric-existence}
Every stable parabolic Higgs bundle with the data~\eqref{eq:graded-residue} and a prescribed determinant trivialization admits a unique tame harmonic metric $h$ adapted to its filtration and inducing the unit metric on the determinant line bundle. To each local holomorphic frame compatible with the flag and the determinant trivialization we associate the frame $(f_a)$, which is unitary for $h_0$ by~\eqref{eq:harmonic-power-model}. The matrix $H=(h(f_a,f_b))_{a,b}$ has a positive Hermitian block diagonal limit $H_p$ of determinant one and satisfies~\eqref{eq:harmonic-local-exponential-decay} for every $0<\nu<\varepsilon_p$ and every derivative order. Each block of $H_p$ acts on the whole corresponding parabolic graded quotient.
\end{proposition}

\begin{proof}
The logarithmic Higgs field defines a regular filtered Higgs bundle, hence a good filtered Higgs bundle with irregular value zero. Its parabolic degree is zero by the determinant trivialization and the centered weights. On a curve this is the vanishing of the filtered first Chern class, and the condition on the second Chern number is vacuous. The filtered determinant line bundle has weight zero and zero Higgs field, so its prescribed unit metric is flat. By Theorem~7.1 of~\cite[p.~54]{MochizukiGood} with $\lambda=0$, there is a unique adapted harmonic metric inducing this metric on the determinant line bundle.

Because the scalar graded residues have zero nilpotent part, we may apply Lemma~\ref{lem:tame-power-comparison} to conclude that the metric and $h_0$ are mutually bounded in every such local frame. By Lemma~\ref{lem:local-harmonic-exponential-decay} its matrix has a positive Hermitian block diagonal limit of determinant one, with each block acting on an entire parabolic graded quotient. The difference from this limit satisfies~\eqref{eq:harmonic-local-exponential-decay} for every cylindrical derivative.
\end{proof}

For the weighted estimates we take a smooth reference metric $h_{\mathrm{ref}}$ of determinant one on $X\setminus D$ equal to $h_0$ near each puncture. For $k\geq0$ and $d\in\R$ we equip $W^{k,2}_d$ with the unweighted Sobolev norm on a fixed compact subset and the squared norm
\begin{equation}\label{eq:harmonic-weighted-norm}
 \sum_{p\in D}\sum_{j=0}^k
 \int_{\tau_p\geq\tau_*}e^{d\tau_p}
       |\nabla_0^j v|^2\,\dd\tau_p\,\dd\theta_p
\end{equation}
on the ends. We set $L^2_d=W^{0,2}_d$. Here $\nabla_0$ is induced by the model connection in~\eqref{eq:harmonic-cylinder-model}. On each end, this norm is equivalent to $\|e^{d\tau/2}v\|_{W^{k,2}}$; thus its exponential weight agrees with that of $H^k_{d/2}$ used in Lemma~\ref{lem:harmonic-half-cylinder}.

\begin{proposition}\label{prop:harmonic-weighted-form}
For the metric $h=h_{\mathrm{ref}}e^s$ of Proposition~\ref{prop:harmonic-metric-existence} and cutoff functions $\chi_p$ equal to one near $p$, we have
\begin{equation}\label{eq:harmonic-Levi-plus-decay}
 s=\sum_{p\in D}\chi_p c_p+v,\qquad
 c_p=\log H_p,\qquad
 v\in W^{k,2}_d
\end{equation}
for every $k\geq0$ and every $0<d<\min_p\varepsilon_p$. Each $c_p$ is a trace-free Hermitian matrix preserving the equal-weight subspaces. In an $h$-unitary frame obtained from~\eqref{eq:harmonic-power-model}, the Chern connection and Higgs field have the form
\begin{equation}\label{eq:harmonic-unitary-end}
 D_h=\dd+\ii A\,\dd\theta+a_h,\qquad
 \Phi=(B+b_h)(-\dd\tau+\ii\,\dd\theta),
\end{equation}
where all cylindrical derivatives of $a_h$ and $b_h$ are $O(e^{-\nu\tau})$ for every $0<\nu<\varepsilon_p$.
\end{proposition}

\begin{proof}
Taking the logarithm preserves~\eqref{eq:harmonic-local-exponential-decay}, because $H$ tends to a positive matrix. For the weighted estimate we take $\nu$ with $d/2<\nu<\min_p\varepsilon_p$. The integral of $e^{(d-2\nu)\tau_p}$ on each end is finite. On the compact regions where the cutoffs vary, the smooth remainder has finite Sobolev norm. Applying the change of frame $H^{-1/2}$ to the connection and Higgs field gives~\eqref{eq:harmonic-unitary-end}, since the constant limit commutes with $A$ and $B$ and all derivatives of the remainder have the stated decay.
\end{proof}

\subsection{Coercivity and comparison of Hermitian metrics}
\label{sec:metric-correction}

The estimates for the form domain, coercivity and metric comparison below apply on $X=\Pp^1$ with any finite puncture set $D$. The applications to our original Higgs bundles use $D=\{0,1,\infty\}$, while the local models use $D=\{0,\infty\}$. A fixed smooth conformal metric $g$ on $X$ has K\"ahler form $\omega$. We write $M=X\setminus D$. Given a Hermitian metric $h_0$ on $E|_M$, we denote both its Chern connection and the induced connection on endomorphisms by $D_{h_0}$. The quadratic form on trace-free $h_0$-Hermitian endomorphisms is
\begin{equation}\label{eq:correction-energy}
 Q_{h_0}(u)=\int_M\bigl(|D_{h_0}u|^2+2|[\Phi,u]|^2\bigr)\,\dd V_g.
\end{equation}
All norms in this formula use $g$ and $h_0$. The completion $V_{h_0}$ of $C_c^\infty(M,\Herm_0(E,h_0))$ for the norm $(\|u\|_2^2+Q_{h_0}(u))^{1/2}$ is the form domain of the associated nonnegative Friedrichs extension.

\subsubsection{The form domain and a uniform coercivity estimate}

\begin{lemma}\label{lem:correction-form-domain}
For $u\in V_{h_0}$, the function $|u|$ belongs to $H^1(X)$ and satisfies
\begin{equation}\label{eq:correction-kato-sobolev}
 \|\dd|u|\|_2\leq\|D_{h_0}u\|_2,\qquad
 \|u\|_q\leq C_q\bigl(\|u\|_2+Q_{h_0}(u)^{1/2}\bigr),
 \qquad 2\leq q<\infty.
\end{equation}
The constants $C_q$ depend only on $(X,g)$ and $q$.
\end{lemma}

\begin{proof}
For integration by parts we use the following cutoffs. In a coordinate disc about a puncture the cutoff equals zero for $|z|\leq a^2$ and one for $|z|\geq a$ and is linear in $\log|z|$ between these circles. Smoothing its corners preserves the bound $C/|\log a|$ for its Dirichlet integral. Combining these functions at the finitely many punctures, with $a\to0$, gives $\chi_j\in C_c^\infty(M)$ satisfying
\begin{equation}\label{eq:correction-logarithmic-cutoffs}
 0\leq\chi_j\leq1,\qquad
 \chi_j\longrightarrow1\text{ pointwise on }M,\qquad
 \|\dd\chi_j\|_2\longrightarrow0.
\end{equation}
Multiplication by $\chi_j$ approximates every smooth function on $X$ in $H^1$. Thus the scalar minimal form domain on $M$ is $H^1(X)$. Kato's inequality $|\dd|u||\leq|D_{h_0}u|$ for compactly supported smooth sections passes to $V_{h_0}$ by approximation and weak lower semicontinuity. Applying the scalar Sobolev inequality on the fixed compact surface to $|u|$ proves~\eqref{eq:correction-kato-sobolev}.
\end{proof}

\begin{lemma}\label{lem:energy-bound-under-convergence}
For a sequence $(E_j,\Phi_j,h_j)$ of Hermitian Higgs bundles of the same rank on $M$, the quadratic forms $Q_{h_j}$ and their domains $V_{h_j}$ are defined by~\eqref{eq:correction-energy} using $\Phi_j$. We assume smooth convergence, under unitary bundle isomorphisms on compact subsets of $M$, of the Chern connections to $D_{A_\infty}$ and of $b_j\Phi_j$ to $\Phi_\infty$, where $0<b_j\leq1$. If $\Phi_\infty$ has distinct eigenvalues outside a finite subset of $M$ and the covering parametrizing its one-dimensional eigenspaces over this locus is connected, then there are $c>0$ and $j_0$ such that
\begin{equation*}
 Q_{h_j}(u)\geq c\|u\|_2^2
 \qquad(j\geq j_0,\ u\in V_{h_j}).
\end{equation*}
\end{lemma}

\begin{proof}
If the assertion fails, we pass to a subsequence and choose $u_j\in V_{h_j}$ with $\|u_j\|_2=1$ and $Q_{h_j}(u_j)\to0$. Applying Lemma~\ref{lem:correction-form-domain}, Kato's inequality and the scalar Poincar\'e inequality on $X$, we obtain
\begin{equation*}
 \bigl\||u_j|-a_j\bigr\|_2\longrightarrow0,
 \qquad a_j\longrightarrow\Vol_g(X)^{-1/2},
\end{equation*}
where $a_j$ is the average of $|u_j|$. In the given unitary identifications, convergence of the connections gives local $H^1$ bounds for $u_j$. By Rellich compactness and a diagonal argument, a subsequence converges to a trace-free Hermitian limit $u_\infty$, weakly in $H^1_{\mathrm{loc}}$ and strongly in $L^2_{\mathrm{loc}}$, with
\begin{equation}\label{eq:gap-limit-section}
 |u_\infty|=\Vol_g(X)^{-1/2},\qquad
 D_{A_\infty}u_\infty=0,\qquad
 [\Phi_\infty,u_\infty]=0.
\end{equation}
The commutator vanishes because, on each compact subset,
\begin{equation*}
 \|[b_j\Phi_j,u_j]\|_2
 \leq b_j Q_{h_j}(u_j)^{1/2}/\sqrt2\longrightarrow0.
\end{equation*}
The parallel section $u_\infty$ is smooth. On the locus of distinct eigenvalues of $\Phi_\infty$, its diagonal entries in a local frame of eigenvectors are locally constant by $D_{A_\infty}u_\infty=0$. They define a locally constant real function on the connected covering of one-dimensional eigenspaces and are therefore all equal. The trace-free condition makes this value zero, contrary to~\eqref{eq:gap-limit-section}.
\end{proof}

\subsubsection{Metric comparison in terms of the error in Hitchin's equation}

For a Hermitian metric $h_0$ the error in Hitchin's equation is
\begin{equation}\label{eq:correction-moment}
 \mu(h_0)=\ii\Lambda_\omega
       \bigl(F_{h_0}+[\Phi,\Phi^{\dagger_{h_0}}]\bigr).
\end{equation}
In the proof we use a bound on the ratio of the two metrics to justify integration by parts for each fixed pair. The constants in the conclusion depend only on the parameters indicated below and on the Riemannian surface $(X,g)$.

\begin{proposition}\label{prop:small-defect-metric-comparison}
We consider smooth Hermitian metrics $h_0$ and $h$ on a holomorphic Higgs bundle of rank $n\geq2$ over $M$ with trace-free Higgs field and a fixed holomorphic determinant trivialization. Both metrics induce the unit metric in this trivialization. We assume that $h$ is harmonic and that
\begin{equation}\label{eq:small-defect-assumptions}
 s=\log(h_0^{-1}h),\qquad
 \|s\|_\infty<\infty,\qquad
 Q_{h_0}(u)\geq c\|u\|_2^2\quad(u\in V_{h_0})
\end{equation}
for some $c>0$. For fixed $1<p<\infty$ we write $\epsilon=\|\mu(h_0)\|_p$. There are $\epsilon_0,C>0$, depending only on $(X,g)$, $n$, $p$, and $c$, such that $\epsilon\leq\epsilon_0$ implies
\begin{equation}\label{eq:small-defect-conclusion}
 s\in V_{h_0},\qquad
 Q_{h_0}(s)^{1/2}\leq C\epsilon,\qquad
 \|s\|_\infty\leq C\sqrt\epsilon.
\end{equation}
In particular,
\begin{equation}\label{eq:small-defect-form-comparison}
 e^{-C\sqrt\epsilon}h_0\leq h\leq e^{C\sqrt\epsilon}h_0.
\end{equation}
The constants are independent of the initial finite value of $\|s\|_\infty$, the Chern connection, and the size of $\Phi$.
\end{proposition}

\begin{proof}
We use $r_0=\mu(h_0)$, $M_s=\|s\|_\infty$, and the path of metrics $h_u=h_0e^{us}$ for $0\leq u\leq1$. The endomorphism $s$ is Hermitian for each $h_u$, has the same pointwise norm for these metrics, and is trace-free. The operators on endomorphism-valued forms are
\begin{equation*}
 D''=\bar\partial_E+[\Phi,\cdot],\qquad
 D'_{h_u}=\partial_{h_u}+[\Phi^{\dagger_{h_u}},\cdot].
\end{equation*}
By Simpson's metric variation identity and the K\"ahler identities~\cite[Lemma~3.1(a)--(c) and proof of Proposition~5.1, pp.~876 and 883]{Simpson1988}, we have locally
\begin{equation}\label{eq:correction-local-variation}
 \frac{\dd}{\dd u}\mu(h_u)=\ii\Lambda_\omega D''D'_{h_u}s,
 \qquad
 2\|D''s\|_{h_0}^2=Q_{h_0}(s),
\end{equation}
where the second equality is used whenever either side is finite. With the norms on endomorphism-valued forms induced by $g$ and $h_u$, integration of the first identity against $s$ gives the quadratic term $\|D''s\|_{h_u}^2$ with coefficient one.

Using the cutoffs in~\eqref{eq:correction-logarithmic-cutoffs}, we test the metric variation against $\chi_j^2s$ to justify integration by parts. Integration over $0\leq u\leq1$ gives the term
\begin{equation*}
 I_j=\int_0^1\int_M\chi_j^2|D''s|_{h_u}^2\,\dd V_g\,\dd u.
\end{equation*}
The contribution from $u=0$ and $u=1$ is $-\int_M\chi_j^2\tr(sr_0)\,\dd V_g$, since $\mu(h)=0$. The term in which a derivative falls on the cutoff is bounded by $CM_s\|\dd\chi_j\|_2 I_j^{1/2}$. Since $s$ is $h_u$-Hermitian, $\partial_{h_u}s=(\bar\partial_Es)^{\dagger_{h_u}}$; a scalar cutoff commutes with $\Phi$. Thus
\begin{equation}\label{eq:correction-cutoff-energy}
 I_j\leq2M_s\|r_0\|_1+CM_s^2\|\dd\chi_j\|_2^2.
\end{equation}
It follows by Fatou's lemma that $I:=\int_0^1\|D''s\|_{h_u}^2\,\dd u<\infty$. The eigenvalues of $s$ lie in $[-M_s,M_s]$, so the norms of endomorphism-valued forms for $h_u$ and $h_0$ differ by bounded factors for this individual pair of metrics. By~\eqref{eq:correction-local-variation} we have $Q_{h_0}(s)<\infty$. The convergence $\chi_js\to s$ holds in the form norm because the term obtained by differentiating the cutoff is bounded by $M_s\|\dd\chi_j\|_2$ and the other derivative and commutator terms converge by dominated convergence. This proves $s\in V_{h_0}$. As the cutoff error tends to zero, we obtain the identity
\begin{equation}\label{eq:correction-integrated-variation}
 -\int_M\tr(sr_0)\,\dd V_g
     =\int_0^1\|D''s\|_{h_u}^2\,\dd u=I.
\end{equation}

For a bound independent of $M_s$, we work in a frame diagonalizing $s$ that is unitary for $h_0$, and denote the eigenvalues by $a_i$. The squared norm of the entry in column $j$ and row $i$ is multiplied by $e^{u(a_i-a_j)}$. For every real $x$,
\begin{equation*}
 \int_0^1 e^{ux}\,\dd u\geq\frac1{1+|x|}.
\end{equation*}
For $x\geq0$ the left side is at least one; for $x=-y<0$, the inequality follows from $e^y\geq1+y$. Since $|a_i-a_j|\leq2M_s$, equation~\eqref{eq:correction-integrated-variation} yields
\begin{equation}\label{eq:correction-convex-energy-bound}
 \frac{Q_{h_0}(s)}{2(1+2M_s)}
 \leq-\int_M\tr(sr_0)\,\dd V_g
 \leq\epsilon\|s\|_{p'},\qquad \frac1p+\frac1{p'}=1.
\end{equation}
By Lemma~\ref{lem:correction-form-domain}, the coercivity estimate, and finiteness of the volume, we have $\|s\|_{p'}\leq C_p Q_{h_0}(s)^{1/2}$. When the quadratic form vanishes, $s=0$. Otherwise we may divide to obtain
\begin{equation}\label{eq:correction-energy-initial-bound}
 Q_{h_0}(s)^{1/2}\leq C\epsilon(1+M_s).
\end{equation}

For the supremum estimate we use the function
\begin{equation}\label{eq:correction-normalized-trace}
 f=\log\frac{\tr(e^s)}n.
\end{equation}
Since the endomorphism is trace-free, $0\leq f\leq|s|$. By Simpson's logarithmic trace inequality~\cite[Lemma~3.1(d), pp.~876--877]{Simpson1988} for the nonnegative scalar Laplacian $\Delta_X=\dd^*\dd$,
\begin{equation}\label{eq:correction-log-trace-inequality}
 \Delta_X f\leq2|r_0|\quad\text{on }M.
\end{equation}
Cyclicity of the trace gives
\begin{equation*}
 \dd f=\frac{\tr(e^sD_{h_0}s)}{\tr(e^s)},\qquad
 |\dd f|\leq|D_{h_0}s|.
\end{equation*}
Thus $f$ has finite Dirichlet integral. Its boundedness and the cutoffs~\eqref{eq:correction-logarithmic-cutoffs} extend it to $H^1(X)$. Testing~\eqref{eq:correction-log-trace-inequality} against $\chi_j\varphi$, for a smooth nonnegative function $\varphi$ on $X$, introduces a term bounded by $\|\varphi\|_\infty\|\dd f\|_2\|\dd\chi_j\|_2$. Letting $j\to\infty$ proves the same inequality distributionally on $X$.

We use the scalar estimate in the proof of Simpson's Proposition~2.1~\cite[pp.~873--874]{Simpson1988} on the compact surface $X$ with $p>1$, applied to $e^{-a\Delta_X}f$. The heat semigroup commutes with $\Delta_X$, is positive, and is contractive on $L^1$ and $L^p$, so the required bounds hold uniformly as $a\downarrow0$. Passing to the limit and using~\eqref{eq:correction-energy-initial-bound} and~\eqref{eq:correction-normalized-trace} gives
\begin{equation}\label{eq:correction-scalar-supremum}
 \|f\|_\infty
 \leq C\|f\|_1+C_p\epsilon
 \leq C\epsilon(1+M_s).
\end{equation}

For a trace-free Hermitian matrix, the normalization in~\eqref{eq:correction-normalized-trace} gives
\begin{equation}\label{eq:correction-matrix-trace-bound}
 |s|\leq C_n(f+\sqrt f).
\end{equation}
Indeed, if $a=\max_i a_i\geq0$, then $|s|\leq\sqrt{n(n-1)}\,a$, while convexity of the exponential on the other $n-1$ eigenvalues gives
\begin{equation*}
 f\geq\log\left(\frac{e^a+(n-1)e^{-a/(n-1)}}n\right).
\end{equation*}
The function on the right is comparable to $a^2$ near zero, bounded below by a positive multiple of $a$ for large $a$, and positive for $a>0$. This proves~\eqref{eq:correction-matrix-trace-bound}. Combining it with~\eqref{eq:correction-scalar-supremum} yields
\begin{equation*}
 M_s\leq C\bigl[\epsilon(1+M_s)+\sqrt{\epsilon(1+M_s)}\bigr].
\end{equation*}
After division by $1+M_s$,
\begin{equation*}
 \frac{M_s}{1+M_s}\leq C(\epsilon+\sqrt\epsilon).
\end{equation*}
For $\epsilon_0\leq1$ sufficiently small, the right side is at most $1/2$. Then $M_s\leq C\sqrt\epsilon$, and substitution into~\eqref{eq:correction-energy-initial-bound} gives $Q_{h_0}(s)^{1/2}\leq C\epsilon$. This also covers $\epsilon=0$. The eigenvalue bound for $s$ proves~\eqref{eq:small-defect-form-comparison}.
\end{proof}

\subsection{Uniform estimates for the local models}
\label{sec:compact-model-estimates}

For fixed $e\in\{2,3,4,6\}$ and $c\in\C^*$, the parameters for the model of Proposition~\ref{prop:polynomial-harmonic-model} form the space
\begin{equation}\label{eq:compact-model-parameter-space}
 \begin{split}
 \mathscr P_e=\bigl\{(\alpha,\sigma)
       \in\R^e\times\C^e:\;&
       \alpha_0<\cdots<\alpha_{e-1}<\alpha_0+1,\\
       &\textstyle\sum_j\alpha_j=\sum_j\sigma_j=0\bigr\}.
 \end{split}
\end{equation}
The polynomial model with parameters $\mathbf p=(\alpha,\sigma)$ has Higgs field $\Phi_{\mathbf p}$ and harmonic metric $h_{\mathbf p}$ of determinant one. We compare their matrices in the frames $f_0,\ldots,f_{e-1}$ of~\eqref{eq:polynomial-model-matrix}. The metric $h_{0,\sigma}$ defined in~\eqref{eq:model-flat-eigenline-metric} is independent of $\alpha$.

\begin{proposition}\label{prop:compact-polynomial-models}
For each compact set $K\subset\mathscr P_e$, there are $a_0>0$ and $C\geq1$ such that, for $\mathbf p\in K$ and $0<|\xi|\leq a_0$,
\begin{equation}\label{eq:compact-model-power-bound}
 C^{-1}D_{\alpha}^2
       \leq h_{\mathbf p}\leq CD_{\alpha}^2,
 \qquad
 D_{\alpha}
       =\diag\bigl(|\xi|^{\alpha_0},\ldots,
                         |\xi|^{\alpha_{e-1}}\bigr).
\end{equation}
On each compact subset of $\C^*$, all derivatives on the punctured curve of the matrices of $h_{\mathbf p}$ and $h_{\mathbf p}^{-1}$ are bounded uniformly for $\mathbf p\in K$. The map $\mathbf p\mapsto h_{\mathbf p}$ is continuous in the local $C^\infty$ topology.

At infinity we use $\xi=r^{-e}$, $\zeta=r^{-1}$, and $\rho=|\zeta|$, as in Proposition~\ref{prop:polynomial-model-tail}. For every integer $k\geq0$, there are $C_k,R_k$, depending on $K,e,c,k$, such that
\begin{equation}\label{eq:compact-model-tail}
 \bigl|\nabla_{0,\sigma}^k
       (h_{0,\sigma}^{-1}h_{\mathbf p}-\Id)
       \bigr|_{h_{0,\sigma},|\dd\zeta|^2}
       \leq C_k\rho^{-1-k}
       \qquad(\mathbf p\in K,\ \rho\geq R_k).
\end{equation}
Here $\nabla_{0,\sigma}$ combines the Chern connection of $h_{0,\sigma}$ on endomorphisms with Euclidean differentiation in $\zeta$. The normalized matrices $G_{\mathbf p}$ of~\eqref{eq:model-normalized-metric-matrix} extend smoothly to a common disc about $r=0$, with uniform bounds on derivatives of every order in the disc coordinate. Their off-diagonal entries vanish to every order at zero, uniformly for $\mathbf p\in K$, and $\log G_{\mathbf p,ii}(0)=0$.

For a fixed finite collection of blocks, these conclusions hold blockwise for the metrics $\bigoplus_\nu s_\nu h_{\mathbf p_\nu}$ of~\eqref{eq:model-block-constants}, when the parameters and positive factors range over compact subsets of their respective parameter spaces and $\prod_\nu s_\nu^e=1$. On the exterior region the flat metric is $\bigoplus_\nu s_\nu h_{0,\sigma_\nu}$.
\end{proposition}

The compact parameter sets may contain points with coincident graded residues. The derivatives of the harmonic metrics are taken on the punctured curve, and the metrics depend continuously on the parameters. The constants may depend on the distance to the boundary of the weight alcove in~\eqref{eq:compact-model-parameter-space}.

\begin{proof}
For a fixed $\mathbf p_*=(\alpha_*,\sigma_*)$, we write $h_*=h_{\mathbf p_*}$. We first work in a sufficiently small convex neighborhood of $\mathbf p_*$. Throughout the proof, the moment map~\eqref{eq:correction-moment}, the quadratic form~\eqref{eq:correction-energy}, and all global $L^s$ norms use one smooth conformal metric $g$ on the compact surface $\Pp^1$, with punctures $0$ and $\infty$. We write $\mu(k,\Phi)$ when the dependence on the Higgs field matters.

We construct a reference metric $k_{\mathbf p}$ from the harmonic metric $h_*$. Near zero we write
\begin{equation}\label{eq:compact-model-reference-zero}
 H_*=D_{\alpha_*}^{-1}h_*
                    D_{\alpha_*}^{-1},
 \qquad
 k_{\mathbf p}=D_{\alpha}H_*
                    D_{\alpha}.
\end{equation}
Applying Lemma~\ref{lem:local-harmonic-exponential-decay} to the local model as explained after its proof, we obtain a positive diagonal matrix $H_{*,\infty}$ of determinant one. The difference $H_*-H_{*,\infty}$ and all its derivatives in the cylindrical coordinates are $O(|\xi|^\nu)$ for some $\nu>0$. In particular, $H_{*,\infty}$ commutes with every diagonal weight and residue matrix in the neighborhood.

We take a closed neighborhood on which $\kappa_0>0$ is a lower bound for the successive weight differences and for $1+\alpha_0-\alpha_{e-1}$, and fix
\begin{equation*}
 0<\nu_0<\min(\nu,\kappa_0,1).
\end{equation*}
We write $a=|\xi|$. In the frame $(|\xi|^{-\alpha_j}f_j)_{j=0}^{e-1}$, the coefficient of the Higgs field relative to $\dd\xi/\xi$ is the diagonal matrix $\diag(\sigma_j)$ plus entries of orders
\begin{equation*}
 a^{\alpha_{j+1}-\alpha_j}
       \quad(0\leq j<e-1),
 \qquad a^{1+\alpha_0-\alpha_{e-1}}.
\end{equation*}
Their derivatives in the cylindrical coordinates satisfy the same bounds. In this frame, the Chern connection of $D_{\alpha}^2$ is $D_{0,\alpha}=\dd+\ii\diag(\alpha_j)\,\dd\theta$, where $\theta=\arg\xi$. We use the same notation for its induced connection on endomorphisms. Since $D_{0,\alpha}$ is flat, the curvature and adjoint for the relative metric $H_*$ are
\begin{equation*}
 \begin{aligned}
 F_{k_{\mathbf p}}
   &=D_{0,\alpha}^{0,1}
       \bigl(H_*^{-1}D_{0,\alpha}^{1,0}H_*\bigr),\\
 \Phi_{\mathbf p}^{\dagger_{k_{\mathbf p}}}
   &=H_*^{-1}\Phi_{\mathbf p}^{\dagger_{D_{\alpha}^2}}H_*\,.
 \end{aligned}
\end{equation*}
At $H_*=H_{*,\infty}$ and with the off-diagonal Higgs entries set to zero, the curvature and Higgs commutator vanish because the weight, residue and metric matrices are diagonal. Thus each first derivative of the error in Hitchin's equation with respect to the parameters contains $H_*-H_{*,\infty}$, a derivative of this difference in the cylindrical coordinates, or an off-diagonal Higgs entry. Derivatives of the weight powers contribute $\log a$, so these terms are bounded by
\begin{equation*}
 C(a^\nu+a^{\kappa_0})(1+|\log a|).
\end{equation*}
We integrate these parameter derivatives along the segment from $\mathbf p_*$ to $\mathbf p$. Since the reference metric solves Hitchin's equation at $\mathbf p_*$, the error at $\mathbf p$ has norm at most $C|\mathbf p-\mathbf p_*|a^{\nu_0}$ with respect to the cylindrical metric. After conversion to the fixed smooth metric $g$, we obtain
\begin{equation}\label{eq:compact-model-zero-defect}
 |\mu(k_{\mathbf p},\Phi_{\mathbf p})|_{k_{\mathbf p},g}
       \leq C|\mathbf p-\mathbf p_*|a^{\nu_0-2}.
\end{equation}
We take
\begin{equation}\label{eq:compact-model-integrability}
 1<s<\frac{2}{2-\nu_0}.
\end{equation}
The right-hand side of~\eqref{eq:compact-model-zero-defect} has $L^s$ norm $O(|\mathbf p-\mathbf p_*|)$ near zero.

To extend this family of metrics near infinity, we use the evaluation matrices $W_{\sigma}$ and $V_{\sigma}$ in~\eqref{eq:model-evaluation-matrices}, matching roots by their limits $\kappa_i$ with $\kappa_i^e=c$. On a common exterior region we define
\begin{equation}\label{eq:compact-model-holomorphic-alignment}
 C_{\sigma}
       =W_{\sigma}^{-1}W_{\sigma_*},
 \qquad
 B_{\sigma}
       =(\det C_{\sigma})^{-1/e}
                         C_{\sigma}.
\end{equation}
The two evaluation matrices undergo the same row permutation under the deck action. Thus $C_{\sigma}$ is single-valued and holomorphic on the exterior in $\xi$. With
$S=\diag(1,\zeta,\ldots,\zeta^{e-1})$, one has $W_{\sigma}=V_{\sigma}S$. The root expansion~\eqref{eq:model-infinity-roots} gives
\begin{equation}\label{eq:compact-model-alignment-estimate}
 \begin{split}
 V_{\sigma}-V_{\sigma_*}
       &=O(|\sigma-\sigma_*|\rho^{-1}),\\
 S(C_{\sigma}-\Id)S^{-1}
       &=O(|\sigma-\sigma_*|\rho^{-1}).
 \end{split}
\end{equation}
After $k$ derivatives in $\zeta,\bar\zeta$, each right-hand side is $O(|\sigma-\sigma_*|\rho^{-1-k})$. Hence $\det C_{\sigma}\to1$ uniformly, and the branch near one defines the root in~\eqref{eq:compact-model-holomorphic-alignment}. In particular, $\det B_{\sigma}=1$.

On this region we define the reference metric for each parameter by
\begin{equation}\label{eq:compact-model-reference-infinity}
 k_{\mathbf p}
       =(B_{\sigma}^{-1})^*
                         h_*B_{\sigma}^{-1}.
\end{equation}
The Higgs fields $\Phi_{\mathbf p}$ and
$B_{\sigma}\Phi_{\mathbf p_*}B_{\sigma}^{-1}$
have the same eigenline bundles. In the coordinate on the cover, the derivatives of order $j$ of the eigenvalue differences are bounded by
$C_j|\sigma-\sigma_*|\rho^{N_j}$ for some $N_j$.
Holomorphic transport preserves the harmonic equation for the second Higgs field and the metric~\eqref{eq:compact-model-reference-infinity}. For a diagonal matrix $Z$ and a positive Hermitian matrix $H$, we denote the off-diagonal part of $H$ by $H_{\mathrm{off}}$. The identity
\begin{equation*}
 [Z,H^{-1}Z^*H]=[Z,H^{-1}[Z^*,H_{\mathrm{off}}]]
\end{equation*}
gives, for two diagonal matrices $Z,Z'$,
\begin{equation}\label{eq:diagonal-commutator-difference}
 \left|[Z,H^{-1}Z^*H]-[Z',H^{-1}Z'^*H]\right|
 \leq C|Z-Z'|(|Z|+|Z'|)|H_{\mathrm{off}}|\,,
\end{equation}
where $C$ depends only on a bound for $H^{-1}$. In the frame diagonalizing both Higgs fields, the normalized metric matrix is $G_*$. By Proposition~\ref{prop:polynomial-model-tail}, its off-diagonal entries vanish to every order in $r$. The factor measuring the change in the parameters comes from the eigenvalue difference, so for every $N$,
\begin{equation}\label{eq:compact-model-infinity-defect}
 |\mu(k_{\mathbf p},\Phi_{\mathbf p})|_{k_{\mathbf p},g}
       \leq C_N|\mathbf p-\mathbf p_*|\rho^{-N}.
\end{equation}
The conversion to $g$ introduces only a fixed power of $\rho$, absorbed by the arbitrary choice of $N$.

The change of frame with determinant one in~\eqref{eq:compact-model-holomorphic-alignment} takes the flat metric for the reference parameter to the flat metric for the new parameter. With
$d_{\sigma}=|\det W_{\sigma}|^{-2/e}$, we obtain
\begin{equation}\label{eq:compact-model-exact-normalization}
 \begin{split}
 (B_{\sigma}^{-1})^*
       h_{0,\sigma_*}B_{\sigma}^{-1}
       &=h_{0,\sigma},\\
 k_{\mathbf p}
       &=W_{\sigma}^*
                      d_{\sigma}G_*W_{\sigma}.
 \end{split}
\end{equation}
The second equality uses $h_*=W_{\sigma_*}^*
d_{\sigma_*}G_*W_{\sigma_*}$ in matched eigenframes. Thus $k_{\mathbf p}$ has the fixed normalized matrix $G_*$ at infinity.

On two fixed disjoint annuli, we interpolate the changes in~\eqref{eq:compact-model-reference-zero} and~\eqref{eq:compact-model-reference-infinity} to the identity. More explicitly, near zero the change is
$D_{\alpha}D_{\alpha_*}^{-1}$,
acting on both sides of $h_*$, and near infinity it is
$B_{\sigma}^{-1}$, acting by Hermitian pullback. On the annuli these matrices of determinant one are close to the identity. Cutting off their logarithms defines a smooth metric $k_{\mathbf p}$ on $\C^*$ with determinant one, equal to $h_*$ at $\mathbf p_*$. On the annuli its error in Hitchin's equation is $O(|\mathbf p-\mathbf p_*|)$. Together with~\eqref{eq:compact-model-zero-defect}--\eqref{eq:compact-model-infinity-defect}, this proves
\begin{equation}\label{eq:compact-model-reference-defect}
 \|\mu(k_{\mathbf p},\Phi_{\mathbf p})\|_{L^s(k_{\mathbf p},g)}
       \leq C|\mathbf p-\mathbf p_*|.
\end{equation}

On a smaller parameter neighborhood the quadratic forms of the reference metrics satisfy the uniform lower bound
\begin{equation}\label{eq:compact-model-reference-energy}
 Q_{k_{\mathbf p}}(u)\geq c_*\|u\|_2^2,
 \qquad u\in V_{k_{\mathbf p}},
\end{equation}
with $c_*>0$. Otherwise there are $\mathbf p_j\to\mathbf p_*$ and sections of unit $L^2$ norm on which the quadratic forms tend to zero. The reference metrics and Higgs fields converge smoothly on compact subsets of $\C^*$. Positive square roots of their matrices give unitary bundle isomorphisms with the required connection convergence. The spectral cover of the limiting Higgs field is $P_{\sigma_*}(w)=c\xi$ with the preimages of zero and branch values removed. It is a complex $w$-plane minus finitely many points and is connected, even when residues or branch values coincide. This contradicts Lemma~\ref{lem:energy-bound-under-convergence} with $b_j=1$.

For each fixed $\mathbf p$, the ratio $k_{\mathbf p}^{-1}h_{\mathbf p}$ and its inverse are bounded. The bounds follow from Lemma~\ref{lem:local-harmonic-exponential-decay} and~\eqref{eq:compact-model-reference-zero} at zero and from Proposition~\ref{prop:polynomial-model-tail} and~\eqref{eq:compact-model-exact-normalization} at infinity. We apply Proposition~\ref{prop:small-defect-metric-comparison} using~\eqref{eq:compact-model-integrability}, \eqref{eq:compact-model-reference-defect}, and~\eqref{eq:compact-model-reference-energy} to obtain
\begin{equation}\label{eq:compact-model-reference-comparison}
 \begin{split}
 u_{\mathbf p}&=\log(k_{\mathbf p}^{-1}h_{\mathbf p}),\\
 \|u_{\mathbf p}\|_\infty
       &\leq C|\mathbf p-\mathbf p_*|^{1/2},\qquad
 Q_{k_{\mathbf p}}(u_{\mathbf p})^{1/2}
       \leq C|\mathbf p-\mathbf p_*|.
 \end{split}
\end{equation}
The constants are independent of the initial bounds for the individual metric ratios. Equations~\eqref{eq:compact-model-reference-zero} and~\eqref{eq:compact-model-reference-comparison} prove~\eqref{eq:compact-model-power-bound} locally in parameter space. Equation~\eqref{eq:compact-model-exact-normalization} also gives uniform mutual bounds for $h_{\mathbf p}$ and $h_{0,\sigma}$ near infinity.

On compact subsets of $\C^*$, equation~\eqref{eq:compact-model-reference-comparison} uniformly bounds the matrices of $h_{\mathbf p}$ and their inverses, and the coefficients of $\Phi_{\mathbf p}$ have uniformly bounded derivatives of every order in the curve coordinate. By Lemma~\ref{lem:harmonic-local-compactness}(i), the metric matrices have uniformly bounded derivatives of every order on smaller discs. Every sequence $\mathbf p\to\mathbf p_*$ thus has a subsequence converging in $C^\infty_{\mathrm{loc}}$. Its limit is $h_*$ by~\eqref{eq:compact-model-reference-comparison}, which proves continuity.

For the estimates at infinity, we check that the constants used in the proof of Proposition~\ref{prop:polynomial-model-tail} can be chosen uniformly. On the discs
\begin{equation*}
 \zeta=\zeta_0+\tfrac14|\zeta_0|z,\qquad |z|<1,
\end{equation*}
the leading eigenvalue coefficients are $\tfrac14|\zeta_0|e\kappa_i$, with minimum separation $d$ a fixed positive multiple of $|\zeta_0|$ and moduli bounded by a fixed multiple of $d$. The remaining coefficients are uniformly $O(|\zeta_0|^{-1})$, hence at most $d/100$ for large $|\zeta_0|$. On fixed smaller discs, the eigenline bundles are therefore exponentially close to orthogonal, uniformly in the parameters, by the estimates of Mochizuki~\cite[Assumption~2.2, Proposition~2.3, Corollary~2.6, Proposition~2.12, Section~2.2.5, and Corollary~2.15, pp.~10--16]{MochizukiAsymptotics}. The curvature bounds for the induced metrics on the eigenline bundles and the derivative estimates in the disc coordinate are uniform on the same discs.

Equations~\eqref{eq:compact-model-exact-normalization} and~\eqref{eq:compact-model-reference-comparison} uniformly bound the $h_{\mathbf p}$ lengths of the $h_{0,\sigma}$-unit holomorphic frames of the eigenline bundles above and below. Their transitions to Mochizuki's frames, and the inverse transitions, are uniformly bounded. On smaller discs the derivatives of the transition functions and their inverses are uniformly bounded by Cauchy's estimates. We can now apply the remaining argument of Proposition~\ref{prop:polynomial-model-tail} with uniform constants. Its scalar Poisson equations for $\log G_{\mathbf p,ii}$ have right-hand sides flat at $r=0$ with uniform derivative bounds, while the metric comparison bounds $\log G_{\mathbf p,ii}$ in $C^0$. Repeated interior Schauder estimates on nested discs~\cite[Theorem~6.2, p.~90]{GilbargTrudinger} give uniform bounds for derivatives of every fixed order on a common smaller disc. The deck action and determinant normalization give $G_{\mathbf p}(0)=\Id$. Thus $G_{\mathbf p}-\Id=O(|r|)$ with uniform Taylor remainders, and the change from $r$ to $\zeta$ gives~\eqref{eq:compact-model-tail}.

A finite cover of $K$ gives common constants and domains. For finitely many blocks, we apply these estimates at each value of $c$, using the same $s_\nu$ on both sides of each comparison. Compactness of the positive factors preserves the bounds and determinant constraint.
\end{proof}

\subsection{Strongly parabolic models}
\label{subsec:strongly-parabolic-models}

We now take $\sigma=0$, so $P(w)=w^e$ in~\eqref{eq:polynomial-spectral-model}. The estimates below hold for each nonempty compact set $K\subset\Acal$. The matrices~\eqref{eq:local-model-small-matrices}--\eqref{eq:local-model-large-matrices} now have zero diagonal. In the coordinate $q=c\xi$, each block is the cyclic Higgs field with parameters $r=e$ and $m=1$ in Mochizuki's notation~\cite[Section~3.1.1 and Theorem~3.32]{MochizukiToda}; his filtration weights at zero are $(-\alpha_j)_{j=0}^{e-1}$.

On the $\lambda$-cover, the spectral family is the fixed curve $Y_1$ with Higgs differential $\lambda\eta_0$, where
\begin{equation*}
 y^n=r_n(z),\qquad \eta_0=\frac{y\,\dd z}{z(z-1)}.
\end{equation*}
The differential $\eta_0$ is nowhere zero. Near a point of $Y_1$ above $p$ with ramification index $e=e_p$, we now use a coordinate $w$ obtained by integration, with $w=0$ at that point and $\eta_0=e\,\dd w$. The cyclic stabilizer multiplies $w$ by an $e$th root of unity, so $x=w^e$ is a coordinate on the base. This choice replaces the coordinates used in~\eqref{eq:local-model-block-table}. The other branches above $p$ change the differential by constant factors of modulus one, which we absorb into $\lambda$. For a tangent variation $\dot\lambda$, multiplying $\dot\lambda$ by the same factor preserves $T=|\lambda|$ and $\dot\lambda/\lambda$.

For a local generator $\ell$ of the spectral line bundle, the direct image has the base frame $e_j=w^j\ell$, $0\leq j<e$. Multiplication by $\lambda w$ gives
\begin{equation}\label{eq:zero-residue-base-field}
 \Phi=\lambda A_1^{(e)}(x)\frac{\dd x}{x}
      =\lambda\left(\sum_{j=0}^{e-2}E_{j+1,j}
                         +xE_{0,e-1}\right)\frac{\dd x}{x},
 \qquad \sigma=0.
\end{equation}
For $e=2$, the unscaled field $\lambda^{-1}\Phi=\bigl(\begin{smallmatrix}0&1\\x^{-1}&0\end{smallmatrix}\bigr)\dd x$ is the strongly parabolic Higgs field in~\cite[Proposition~3.9]{FMSW}. On the punctured ramified cover, the frame $f_j=w^{-j}e_j$ makes the cyclic coefficient constant:
\begin{equation}\label{eq:zero-residue-cyclic-field}
 \Phi=e\lambda P\,\dd w,\qquad
 P=\sum_{j=0}^{e-1}E_{j+1,j},\qquad \rho=eT|w|,
\end{equation}
Here $P$ denotes the cyclic shift matrix, and indices on matrices and on $u_j$ below are taken modulo $e$. We define
\begin{equation}\label{eq:zero-residue-model-weights}
 \begin{gathered}
 b_j=\frac{j-(e-1)/2}{e},\qquad
 \gamma_j=e(\alpha_j-b_j),\\
 \delta_j=\alpha_{j+1}-\alpha_j\quad(0\leq j<e-1),\qquad
 \delta_{e-1}=1+\alpha_0-\alpha_{e-1},\\
 \kappa=\min_{\alpha\in K,\ p\in D,\ j}
                         e_p\delta_{p,j}>0.
 \end{gathered}
\end{equation}
Thus $\sum_j\gamma_j=0$, the identity $\gamma_{j+1}-\gamma_j=e\delta_j-1$ holds for every index modulo the block size, and $\kappa\leq1$.

\begin{lemma}\label{lem:zero-residue-cyclic-model}
We normalize the determinant metric of each block to one in the determinant frame $e_0\wedge\cdots\wedge e_{e-1}$. The harmonic metric on the block is diagonal in the base frame $(e_j)$, and its matrices in the base and ramified frames are
\begin{equation}\label{eq:zero-residue-model-metric}
 \begin{aligned}
 h(e_j,e_k)&=\delta_{jk}|x|^{2b_j}e^{2u_j(\rho)},\\
 h&=|w|^{-(e-1)}\diag(e^{2u_j(\rho)})
       \quad\text{in the frame }(f_j),\qquad \sum_j u_j=0.
 \end{aligned}
\end{equation}
The radial functions satisfy
\begin{equation}\label{eq:zero-residue-Toda}
 u_j''+\rho^{-1}u_j'
 =2\bigl(e^{2(u_j-u_{j-1})}-e^{2(u_{j+1}-u_j)}\bigr).
\end{equation}
There are constants $c_{\alpha,j}$, uniformly bounded over the indicated weights, such that for every integer $a\geq0$,
\begin{equation}\label{eq:zero-residue-model-zero}
 \left|(\rho\partial_\rho)^a
   \bigl(u_j-\gamma_j\log\rho-c_{\alpha,j}\bigr)\right|
       \leq C_a\rho^{2\kappa},\qquad 0<\rho\leq1.
\end{equation}
There is also a constant $c_1>0$, independent of $a$, such that
\begin{equation}\label{eq:zero-residue-model-infinity}
 |(\rho\partial_\rho)^a u_j|
       \leq C_a e^{-c_1\rho},\qquad \rho\geq1.
\end{equation}
All constants are uniform for the block weights induced by $K$.
\end{lemma}

At zero, the growth of the metric on the base frame is
\begin{equation}\label{eq:zero-residue-weight-growth}
 h(e_j,e_j)
   =(eT)^{2\gamma_j}e^{2c_{\alpha,j}}
       |x|^{2\alpha_j}\bigl(1+O(\rho^{2\kappa})\bigr)
       \qquad(\rho\to0).
\end{equation}
Indeed, $e_j=w^jf_j$, $x=w^e$ and $b_j+\gamma_j/e=\alpha_j$. Thus the weights change the growth on each graded quotient even though the strongly parabolic spectral curve is fixed. 

\label{subsec:cyclic-model-estimates}
\begin{proof}[Proof of Lemma~\ref{lem:zero-residue-cyclic-model}]
We apply uniqueness in Proposition~\ref{prop:polynomial-harmonic-model} with vanishing graded residues. Conjugation by the diagonal cyclic grading multiplies the Higgs field by a unit scalar and preserves its filtrations and determinant. The metric is therefore diagonal. Rotation of the base, accompanied by a diagonal change of frame, preserves the same data and makes its entries radial; these are also the symmetries in~\cite[Proposition~3.27 and Lemma~3.30]{MochizukiToda}. Scalar corrections to these frame changes normalize their determinants and leave their action on the metric unchanged. Since $e_j=w^jf_j$ and $x=w^e$, the two expressions in~\eqref{eq:zero-residue-model-metric} agree. The curvature equation gives~\eqref{eq:zero-residue-Toda}, since the squared norm of $E_{j+1,j}$ is $e^{2(u_{j+1}-u_j)}$ and $\Delta_w u_j=e^2T^2(u_j''+\rho^{-1}u_j')$.

After dilation and determinant normalization, the diagonal comparison in Proposition~\ref{prop:compact-polynomial-models}, with the weights~\eqref{eq:zero-residue-model-weights}, gives
\begin{equation*}
 u_j=\gamma_j\log\rho+O(1)
\end{equation*}
uniformly in $K$. We write $v_j=u_j-\gamma_j\log\rho$. The relations in~\eqref{eq:zero-residue-model-weights} and the equation~\eqref{eq:zero-residue-Toda} imply
\begin{equation*}
 (\rho v_j')'=O(\rho^{2\kappa-1}).
\end{equation*}
Since $v_j$ is bounded, the integration constant in $\rho v_j'$ vanishes. Two integrations give a finite limit $c_{\alpha,j}$ and an $O(\rho^{2\kappa})$ remainder. The resulting constants are uniformly bounded because the values at $\rho=1$ are uniformly bounded. Differentiating the equation with $\rho\partial_\rho$ proves~\eqref{eq:zero-residue-model-zero}; the estimates on fixed annuli in Proposition~\ref{prop:compact-polynomial-models} control the values and derivatives at $\rho=1$.

Near infinity we remove the common scalar in~\eqref{eq:zero-residue-model-metric} and write $H=\operatorname{diag}(e^{2u_j})$. We use a constant unitary Fourier matrix $F$ that diagonalizes $P$. Every diagonal entry of $G=F^*HF$ is
\begin{equation*}
 m=\frac1e\sum_j e^{2u_j}.
\end{equation*}
On a disc $w=w_0+a|w_0|z$, $|z|<1$, with a fixed small $a>0$, the coefficient eigenvalues of~\eqref{eq:zero-residue-cyclic-field} are $a|w_0|e\lambda$ times the distinct $e$th roots of unity. Their separation is a fixed positive multiple of $\rho_0=eT|w_0|$, and the ratio of their maximum modulus to their minimum separation is fixed. Above a common threshold these coefficients satisfy Assumption~2.2 in Mochizuki's paper. By his estimate for the relative orthogonality of the eigenline bundles~\cite[Assumption~2.2 and Corollary~2.6, pp.~10--11]{MochizukiAsymptotics},
\begin{equation}\label{eq:zero-residue-relative-orthogonality}
 |G_{ij}|\leq C e^{-c\rho}m\quad(i\ne j).
\end{equation}
The constants depend only on the fixed disc radii and $e$. We decompose $G$ as $m(\Id+E)$, where $E$ has zero diagonal and $|E|=O(e^{-c\rho})$. The determinant normalization gives
\begin{equation*}
 m=\det(\Id+E)^{-1/e}=1+O(e^{-c\rho}),
 \qquad H-\Id=O(e^{-c\rho}).
\end{equation*}
Consequently $u_j=O(e^{-c\rho})$ uniformly. On intervals of fixed length at large $\rho$, the radial equation has bounded coefficients and a smooth nonlinearity vanishing at $u=0$. The same exponential bound holds for $u'$ by the interior derivative estimate and for each successive derivative with respect to $\rho$ by differentiation of the radial equation. The inequality $\rho^a e^{-c\rho}\leq C_a e^{-c\rho/2}$ gives one rate for all powers of $\rho\partial_\rho$. The compact annulus estimates extend the result to $\rho\geq1$, proving~\eqref{eq:zero-residue-model-infinity}.
\end{proof}

\section{Approximate and harmonic metrics}
\label{sec:approximate-metrics}

We join the local harmonic metrics near the punctures to the flat metric induced by the spectral line bundle. At a puncture of flag length $e_p$ we match the eigenline bundles of the model to those of the given Higgs field on an annulus of radius comparable to $T^{-e_p/2}$. In a holomorphic eigenframe orthonormal for the exterior metric, the diagonal entries of the logarithm of the metric ratio are bounded in absolute value by a constant times $T^{-1/2}$. The off-diagonal entries of the logarithm are $O(T^{-N})$ for every $N$. These estimates bound the curvature and the Higgs commutator of the interpolated metric.

The estimates below hold on every nonempty compact set $K\subset\Acal\times\mathscr B$, where $\mathscr B$ is the residue space in~\eqref{eq:graded-residue}. The graded residues may coincide or vanish. We use the coordinates and parameters of~\eqref{eq:local-model-scale} and fix a smooth conformal metric $g$ on $X=\Pp^1$. Norms of endomorphism-valued forms are taken with respect to $g$ and the indicated Hermitian metric. The weights determine the two constants
\begin{equation}\label{eq:uniform-weight-margins}
 \begin{aligned}
 a_K&=\max_{(\alpha,\sigma)\in K,\ p\in D}
             (\alpha_{p,e_p-1}-\alpha_{p,0})<1,\\
 \kappa_K&=\min_{\substack{(\alpha,\sigma)\in K\\p\in D}}
       \min\bigl\{1+\alpha_{p,0}-\alpha_{p,e_p-1},\\
       &\hspace{7em}
          \alpha_{p,j+1}-\alpha_{p,j}\ (0\leq j<e_p-1)\bigr\}>0.
 \end{aligned}
\end{equation}
We choose
\begin{equation}\label{eq:approximate-exponent-range}
 0<\eta<1-a_K,\qquad
 1<p_0<\min\left\{\frac{12}{11},\frac{2}{1+a_K+\eta}\right\}.
\end{equation}
For this integrability exponent $p_0$ the decay exponent is
\begin{equation}\label{eq:approximate-defect-rate}
 \delta=\frac12-6\left(1-\frac1{p_0}\right)>0.
\end{equation}
For a metric $h$ we write $\mathscr E(h)=F_h+[\Phi,\Phi^{\dagger_h}]$.

\begin{proposition}\label{prop:approximate-metric}
There are constants $C,T_0>0$ depending on $K$, $g$, $\eta$ and $p_0$, as well as the puncture coordinates, Picard charts, local frames and cutoffs used below. For every $(\alpha,\sigma)\in K$, $t=\lambda^n$ with $T=|\lambda|\geq T_0$, and $L\in\Pic^n(Y_{t,\sigma})$, the associated filtered Higgs bundle admits a smooth Hermitian metric $h^{\app}$ on $X\setminus D$. It induces the prescribed filtration, has determinant one in the fixed trivialization and satisfies
\begin{equation}\label{eq:approximate-Lp-defect}
 \|\mathscr E(h^{\app})\|_{L^{p_0}(g,h^{\app})}
       \leq C T^{-\delta}.
\end{equation}
On the complement of the discs $|x_p|\leq2T^{-e_p/2}$, the Chern connection of $h^{\app}$ is flat and the eigenline bundles of the given Higgs field are orthogonal.

On a common disc in $\xi=\lambda^{e_p}x_p$ about zero, the metric is an orthogonal sum of the local harmonic metrics of Proposition~\ref{prop:polynomial-harmonic-model} with the given weights and residues. Each block is multiplied by its own positive constant and by a common positive function with harmonic logarithm. After removal of the common scalar, the constants satisfy $\prod_{\nu=1}^{m_p}s_{p,\nu}^{e_p}=1$.
\end{proposition}

\subsection{The flat metric induced from the spectral line bundle}

The scaling in~\eqref{eq:coefficient-scaling} gives an isomorphism from $Y_{t,\sigma}$ to $Y_{1,\boldsymbol\gamma}$, where $\boldsymbol\gamma=\sigma/\lambda$. We choose concentric balls $S_0,S_1$ in the space of graded residues with $\overline S_0\subset S_1$ and $\overline S_1$ in the neighborhood of Proposition~\ref{prop:normalized-spectral-family}. Since the residues range over a compact set, we can increase $T_0$ once so that $\boldsymbol\gamma\in S_0$ throughout $K$. The exterior metric is independent of the weights and depends on $\sigma$ through $\boldsymbol\gamma$.

We use a homology marking transported over $S_1$ by the Gauss--Manin connection and the origin $L_0=\pi^*\mathcal O_X(1)$ to parametrize $L\otimes L_0^{-1}$ by unitary characters in
\begin{equation*}
 \operatorname{Hom}\bigl(H_1(Y_{1,\boldsymbol0},\Z),U(1)\bigr).
\end{equation*}
The relative Picard family over $\overline S_1$ and the unit circle are compact. We cover their product by finitely many relatively compact charts. On a neighborhood of each chart closure we choose local generators of the spectral line bundles and frames for the determinant line bundles of their direct images. We estimate the coefficients on these closures.

At each puncture $p$, the normalized spectral curve has $m_p$ disjoint coordinate discs with coordinate $v$, the eigenvalue of $\lambda^{-1}\Phi$ in the logarithmic frame $\dd x_p/x_p$. We choose these discs over a neighborhood of $\overline S_1$. The functions defining their graphs $x_p=X_{p,\nu}(v,\boldsymbol\gamma)$ and the holomorphic units in~\eqref{eq:local-normalized-branches} are uniformly bounded on smaller discs over $\overline S_1$. In each disc we mark the auxiliary section $v=0$. For nonzero $\boldsymbol\gamma$, its image may lie away from $p$.

\begin{lemma}\label{lem:exterior-spectral-metric}
There is a flat metric $k_{\boldsymbol\gamma,L}$ on $L$ away from these auxiliary sections, unique up to a positive constant, whose squared norm in a nonvanishing local holomorphic generator has the form
\begin{equation}\label{eq:exterior-line-metric}
 k_{\boldsymbol\gamma,L}(v)
       =|v|^{-(e_p-1)}\exp f_{p,\nu}(v,\boldsymbol\gamma,L).
\end{equation}
Here $f_{p,\nu}$ is real-valued, smooth across $v=0$, and harmonic in $v$. After a scalar normalization, these metrics and their regular factors depend smoothly on the normalized curve parameter and on local Picard parameters.

For $T$ sufficiently large, their orthogonal direct images define a metric $h^{\dec}$ with determinant one outside discs $|x_p|<C T^{-e_p}$. On that region,
\begin{equation}\label{eq:exterior-flat-equations}
 F_{h^{\dec}}=0,\qquad [\Phi,\Phi^{\dagger_{h^{\dec}}}]=0.
\end{equation}
On compact subsets of $X\setminus D$, the metrics depend smoothly on $\boldsymbol\gamma$ and the local Picard parameters down to $\boldsymbol\gamma=0$.
\end{lemma}

\begin{proof}
We prescribe norm growth exponent $-(e_p-1)/2$ at each marked section over $p$. In each of the three flag types,
\begin{equation*}
 \sum_{p\in D}m_p(e_p-1)=2n.
\end{equation*}
Since $\deg L=n$ by~\eqref{eq:spectral-line-degree}, the filtered degree is zero. In the Picard charts above we follow the rank-one construction of~\cite[Proposition~3.2]{MochizukiGood}, choosing smooth reference metrics and inserting the powers $|v|^{-(e_p-1)}$ with coordinate cutoffs. The resulting metrics $k_{\mathrm{ref}}$ have curvature forms that extend smoothly across the marked sections and have integral zero by the degree calculation.

We choose smooth volume forms $\dd V$ on the curve family. By the scalar Poisson theorem and elliptic estimates in H\"older spaces~\cite[Theorem~1 and Section~4.2]{DonaldsonGeometricAnalysis}, we obtain a unique smooth real function $u$ satisfying
\begin{equation*}
 F_{k_{\mathrm{ref}}}+\bar\partial\partial u=0,
 \qquad \int_{Y_{1,\boldsymbol\gamma}}u\,\dd V=0\,.
\end{equation*}
The metric $k_{\boldsymbol\gamma,L}=k_{\mathrm{ref}}e^u$ is flat. Its regular logarithmic factors are smooth across the marked sections and harmonic, which gives~\eqref{eq:exterior-line-metric}. The logarithmic ratio of two metrics of this form extends harmonically to the compact curve and is constant, proving uniqueness up to a positive factor.

In a smooth trivialization of the curve family, the scalar Laplacians and the mean-zero projections vary smoothly. The restricted Laplacians are isomorphisms from mean-zero functions in $C^{k+2,\varpi}(Y_{1,\boldsymbol\gamma})$ to mean-zero functions in $C^{k,\varpi}(Y_{1,\boldsymbol\gamma})$ for $k\geq0$ and $0<\varpi<1$. Since inversion is smooth on bounded isomorphisms, the inverse Laplacians also vary smoothly. By elliptic regularity the solutions and their derivatives along the spectral curve have uniform bounds of every order on the compact chart closures. Thus the metrics and their regular factors have the asserted parameter dependence.

The auxiliary marked sections project into $|x_p|\leq C T^{-e_p}$ with a common $C$ for parameters in $K$. Indeed, in the rescaled equation $P(w)=\xi W_\nu$, the section $v=0$ has $w=0$, and $W_\nu$ stays uniformly bounded away from zero on a fixed disc containing this section. By Proposition~\ref{prop:local-spectral-matching}, the discs around the branch points can be chosen with a common radius, and the roots outside these discs are uniformly separated. Outside the discs the spectral map is unramified. The orthogonal direct image $\pi_*k_{\boldsymbol\gamma,L}$ of the flat metric on the spectral line bundle is therefore flat there, and $\Phi$ is normal for this metric.

The fixed trivialization of $\det E$ determines a nowhere-vanishing holomorphic section $d_E$. With $Q=|d_E|^2_{\det(\pi_*k_{\boldsymbol\gamma,L})}$ we define
\begin{equation}\label{eq:exterior-determinant-normalization}
 h^{\dec}=Q^{-1/n}\pi_*k_{\boldsymbol\gamma,L}.
\end{equation}
The function $\log Q$ is harmonic on this exterior region. Thus the normalization preserves both equations in~\eqref{eq:exterior-flat-equations} and gives the unit metric on the trivialized determinant line bundle. Smooth dependence follows from the scalar construction and the unramified direct image.
\end{proof}

\subsection{The approximate metric near a puncture}

At a fixed $p$ we use the notation $e=e_p$, $m=m_p$, $\xi=\lambda^e x$, and $\varepsilon=\lambda^{-e}$. We choose a holomorphic generator of $L$ on each of the $m$ normalized spectral discs, smoothly in a local parameter chart. The polynomials $f_j(w)=\prod_{i<j}(w-\sigma_i)$ multiplied by this generator give the direct image frame in~\eqref{eq:newton-spectral-lattices}. In this frame the prescribed flag remains fixed when residues coincide.

On a fixed compact $\xi$-disc the monic polynomial of the $\nu$th block is
\begin{equation}\label{eq:actual-inner-polynomial}
 Q_{\nu,\varepsilon}(w,\xi)
     =P(w)-c_\nu\xi+\varepsilon\xi b_{\nu,\varepsilon}(w,\xi),
 \qquad \deg_w b_{\nu,\varepsilon}<e.
\end{equation}
Its coefficients and every fixed number of their derivatives are bounded as $\varepsilon\to0$, uniformly for parameters in $K$. We obtain this expression from the decomposition of the finite flat local algebra into the separate $u$-branches of~\eqref{eq:local-rescaled-modification}. Contour integrals around the entire block give holomorphic monic coefficients even when roots within the block coincide. At $\xi=0$ the polynomial is $P$ for every $\varepsilon$, and at $\varepsilon=0$ it is $P-c_\nu\xi$. Its difference from the latter polynomial is therefore divisible by $\varepsilon\xi$.

We express the correction in the Newton basis by
\begin{equation*}
 -\varepsilon b_{\nu,\varepsilon}(w,\xi)
       =\sum_{j=0}^{e-1}d_{\nu,j}(\xi,\varepsilon)f_j(w).
\end{equation*}
The coefficients $d_{\nu,j}$ are holomorphic and $O(|\varepsilon|)$ on each fixed compact $\xi$-disc, together with every fixed number of derivatives in $\xi$. In this basis the Higgs field has matrix
\begin{equation}\label{eq:exact-inner-matrix}
 \begin{aligned}
 A_{\nu,\varepsilon}
   &=A_{c_\nu}
       +\xi\begin{pmatrix}
          0&\cdots&0&d_{\nu,0}\\
          \vdots&&\vdots&\vdots\\
          0&\cdots&0&d_{\nu,e-1}
        \end{pmatrix},\\
 \Psi_{\nu,\varepsilon}&=A_{\nu,\varepsilon}\frac{\dd\xi}{\xi}.
 \end{aligned}
\end{equation}
Indeed, the first $e-1$ columns still satisfy $wf_j=\sigma_jf_j+f_{j+1}$, while~\eqref{eq:actual-inner-polynomial} gives $wf_{e-1}=\sigma_{e-1}f_{e-1}+\xi(c_\nu f_0+\sum_jd_{\nu,j}f_j)$. Thus the residue is unchanged, and the entire correction lies in the last column. With $\Psi_{\nu,0}=\Phi_{c_\nu}$, we obtain
\begin{equation}\label{eq:inner-Higgs-difference}
 \Psi_{\nu,\varepsilon}-\Psi_{\nu,0}
       =\sum_{j=0}^{e-1}d_{\nu,j}E_{j,e-1}\,\dd\xi
       =\varepsilon\widetilde b_{\nu,\varepsilon}(\xi)\,\dd\xi,
\end{equation}
where $\widetilde b_{\nu,\varepsilon}$ is holomorphic and bounded on each fixed compact disc.

\begin{lemma}\label{lem:transported-inner-metric}
For some common $R_0$ and all sufficiently large $T$, there are smooth bundle automorphisms $B_\nu$ of determinant one on $|\xi|<3T^{e/2}$ with the following properties. They equal the identity on $|\xi|\leq R_0/2$, are holomorphic for $|\xi|\geq R_0$, and on that region send each eigenline bundle of the model to the corresponding eigenline bundle of the given Higgs field. We define
\begin{equation}\label{eq:transported-model-metric}
 h'_{\nu}(u_1,u_2)
       =h_{c_\nu}(B_\nu^{-1}u_1,B_\nu^{-1}u_2).
\end{equation}
Then, with the Euclidean metric in the $\xi$ coordinate,
\begin{equation}\label{eq:inner-Lp-defect}
 \left\|F_{h'_\nu}
       +[\Psi_{\nu,\varepsilon},\Psi_{\nu,\varepsilon}^{\dagger_{h'_\nu}}]
       \right\|_{L^{p_0}(|\xi|<3T^{e/2})}
       \leq C|\varepsilon|.
\end{equation}
The radius $R_0$, the threshold for $T$, and the constant $C$ are uniform over $K$, the Picard fibers and all arguments of $\lambda$. They may depend on the choices in~\eqref{eq:approximate-exponent-range} and Proposition~\ref{prop:approximate-metric}.
\end{lemma}

\begin{proof}
The puncture parameters range over compact subsets of~\eqref{eq:compact-model-parameter-space}, and the coefficients $c_\nu$, with $c_\nu^m=h_p(0)$, form a fixed finite set. The comparison with the powers of the radial coordinate in~\eqref{eq:compact-model-power-bound}, the derivative bounds on compact annuli and the estimates at infinity for the normalized matrices are uniform by Proposition~\ref{prop:compact-polynomial-models}. Only finitely many derivatives enter the estimate for the error in Hitchin's equation.

The positive lower bounds for the successive weight differences and for $1+\alpha_{p,0}-\alpha_{p,e_p-1}$ in~\eqref{eq:uniform-weight-margins} control the subdiagonal terms and the term $c_\nu\xi f_0$ in~\eqref{eq:polynomial-model-matrix}. Equations~\eqref{eq:compact-model-power-bound} and~\eqref{eq:inner-Higgs-difference} therefore give, on a common disc about zero,
\begin{equation*}
 \begin{aligned}
 |\Psi_{\nu,0}|_{h_{c_\nu},|\dd\xi|^2}
       &\leq C|\xi|^{-1},\\
 |\Psi_{\nu,\varepsilon}-\Psi_{\nu,0}|_{h_{c_\nu},|\dd\xi|^2}
       &\leq C|\varepsilon||\xi|^{-a_K}.
 \end{aligned}
\end{equation*}
The curvature cancels against the model Higgs commutator. Expanding the remaining commutators gives, on a small fixed punctured disc,
\begin{equation}\label{eq:inner-core-defect}
 \left|F_{h_{c_\nu}}
       +[\Psi_{\nu,\varepsilon},\Psi_{\nu,\varepsilon}^{\dagger_{h_{c_\nu}}}]\right|
       \leq C|\varepsilon||\xi|^{-1-a_K-\eta}.
\end{equation}
The quadratic term satisfies the same bound because $a_K<1$. For the chosen exponents we have
\begin{equation}\label{eq:uniform-puncture-integral}
 \int_0^1 r^{\,1-p_0(1+a_K+\eta)}\,\dd r
       =\frac{1}{2-p_0(1+a_K+\eta)}<\infty.
\end{equation}
On compact sets away from zero, the derivative bounds in Proposition~\ref{prop:compact-polynomial-models} imply that the error in Hitchin's equation is $O(|\varepsilon|)$ uniformly.

We choose $R_0$ so that all model and exact branch points lie in $|\xi|<R_0/2$ and the estimates in Propositions~\ref{prop:local-spectral-matching} and~\ref{prop:compact-polynomial-models} apply outside $R_0$. On sectors of $R_0\leq|\xi|\leq3T^{e/2}$ we label the corresponding roots and put
\begin{equation*}
 (W_0)_{ij}=f_j(w_i^{\mathrm{mod}}),\qquad
 (W_\varepsilon)_{ij}=f_j(w_i^{\ex}),\qquad
 C_\nu=W_\varepsilon^{-1}W_0.
\end{equation*}
As in the construction of~\eqref{eq:compact-model-holomorphic-alignment}, $C_\nu$ is single-valued and holomorphic on the annulus and sends the frames of the model eigenline bundles to those of the corresponding eigenline bundles of the given Higgs field.

We take $\zeta$ with $\zeta^e=\xi$ and put $S=\diag(1,\zeta,\ldots,\zeta^{e-1})$. In the factorization $W=VS$, the normalized matrices $V$ and $V^{-1}$ are bounded. Together with Proposition~\ref{prop:local-spectral-matching}, these bounds yield
\begin{equation}\label{eq:eigenline-alignment-bound}
 \bigl|S(C_\nu-\Id)S^{-1}\bigr|
       \leq C|\varepsilon||\xi|,
\end{equation}
with every fixed number of derivatives in a rescaled annulus coordinate. Since $|\varepsilon||\xi|\leq3T^{-e/2}$, the determinant has a single-valued $e$th root near one. We set $B_\nu=(\det C_\nu)^{-1/e}C_\nu$ on this annulus. It has determinant one, identifies the corresponding eigenline bundles, and satisfies~\eqref{eq:eigenline-alignment-bound}. The same formula is defined on $R_0/2\leq|\xi|\leq R_0$, which contains no branch point, and equals $\Id+O(|\varepsilon|)$ there with every derivative in the frame $f_0,\ldots,f_{e-1}$. On this fixed annulus we replace it by $\exp(\chi_0\log B_\nu)$, where $\chi_0$ is radial, zero near $R_0/2$, and one near $R_0$. Extending by the identity inside gives the asserted smooth map. Its antiholomorphic derivative is supported on this fixed annulus, where it and all its derivatives are $O(|\varepsilon|)$.

Where $B_\nu$ is holomorphic, $h'_\nu$ solves the harmonic metric equation for $B_\nu\Psi_{\nu,0}B_\nu^{-1}$. This form and $\Psi_{\nu,\varepsilon}$ are diagonal in the same frame of the eigenline bundles, and
\begin{equation*}
 |w_i^{\ex}-w_i^{\mathrm{mod}}|
       \leq C|\varepsilon||\xi|^{1+1/e}.
\end{equation*}
After division by the common squared norm of the frames of the model eigenline bundles, the matrix of $h'_\nu$ is the matrix $G$ in~\eqref{eq:model-normalized-metric-matrix}, up to a common scalar. The diagonal entries of $G$ and their reciprocals are uniformly bounded, and its off-diagonal entries vanish to every order at infinity. We apply~\eqref{eq:diagonal-commutator-difference} in this normalized frame, with $Z=\diag(w_i^{\ex})/\xi$ and $Z'=\diag(w_i^{\mathrm{mod}})/\xi$. By Proposition~\ref{prop:compact-polynomial-models}, for every $N$, we obtain
\begin{equation}\label{eq:aligned-exterior-defect}
 \left|F_{h'_\nu}
       +[\Psi_{\nu,\varepsilon},\Psi_{\nu,\varepsilon}^{\dagger_{h'_\nu}}]\right|
       \leq C_N|\varepsilon||\xi|^{-N}
       \quad (R_0\leq|\xi|\leq3T^{e/2}).
\end{equation}
Where the change of frame is holomorphic, the curvature is the conjugate of the original curvature. The additional curvature terms on the fixed interpolation annulus are $O(|\varepsilon|)$. We integrate this bound together with~\eqref{eq:inner-core-defect} and~\eqref{eq:aligned-exterior-defect} to obtain~\eqref{eq:inner-Lp-defect}.
\end{proof}

We now normalize the determinant and the relative sizes of the blocks. Here $\tau_v$ denotes the wedge of the frames $\prod_{i<j}(v-\gamma_i)$ in all $m$ blocks, including the chosen generators of the spectral line bundles. There is a factorization $\tau_v=f_{\boldsymbol\gamma}(x)d_E$, where $f_{\boldsymbol\gamma}$ is a holomorphic unit. This unit, its inverse and their derivatives along the base curve have common bounds on the chosen chart closures. The wedge $\tau_w$ of the frames $f_0,\ldots,f_{e-1}$ is
\begin{equation*}
 \tau_w=\lambda^{me(e-1)/2}\tau_v.
\end{equation*}
Consequently the positive scalar
\begin{equation}\label{eq:inner-determinant-factor}
 d_\lambda(x)=T^{e-1}|f_{\boldsymbol\gamma}(x)|^{2/n}
\end{equation}
satisfies $d_\lambda^n=|\tau_w/d_E|^2$ and has harmonic logarithm. For positive constants with $\prod_\nu s_\nu^e=1$, we define
\begin{equation}\label{eq:inner-metric}
 h^{\mathrm{in}}=d_\lambda(x)\bigoplus_{\nu=1}^m s_\nu h'_\nu.
\end{equation}
This metric has determinant one in $d_E$. Multiplication by $d_\lambda$ contributes no curvature, and constant factors on the holomorphic blocks preserve their Higgs adjoints and endomorphism norms. Thus~\eqref{eq:inner-Lp-defect} also holds for~\eqref{eq:inner-metric}, after summing the errors in Hitchin's equation over the blocks. Near zero the maps $B_\nu$ are the identity and $d_\lambda$ is smooth and positive for each fixed $\lambda$, so $h^{\mathrm{in}}$ induces the prescribed filtration.

We choose the relative block constants from the regular factors in~\eqref{eq:exterior-line-metric}, evaluated at $(v,\boldsymbol\gamma)=(0,\boldsymbol0)$ with the same unitary character and local generators of the spectral line bundle, by setting
\begin{equation}\label{eq:inner-block-constants}
 s_\nu=\exp\!\left(
 f_{p,\nu}(0,\boldsymbol0,L)
 -\frac1m\sum_{\kappa=1}^m f_{p,\kappa}(0,\boldsymbol0,L)
 \right)\,.
\end{equation}
In corresponding frames of the eigenline bundles, the leading coefficients of the orthogonal direct-image metric and the model metric~\eqref{eq:model-flat-eigenline-metric} have ratio $eT^{e-1}\exp f_{p,\nu}(0,\boldsymbol0,L)$. After the determinant normalization~\eqref{eq:exterior-determinant-normalization} and removal of $d_\lambda$, this ratio is~\eqref{eq:inner-block-constants}. The constants thus satisfy $\prod_\nu s_\nu^e=1$, leaving $m-1$ independent ratios as in~\eqref{eq:model-block-constants}. Smooth dependence and common bounds for the regular factors place $(s_\nu)$ in a compact subset of $\{(s_\nu):s_\nu>0,\ \prod_\nu s_\nu^e=1\}$, uniformly over $K$, the Picard fibers and all arguments of $\lambda$.

\subsection{Interpolation on the annulus and the error in Hitchin's equation}

\begin{lemma}\label{lem:inner-outer-metric-matching}
On the annulus $T^{-e/2}\leq|x|\leq2T^{-e/2}$ we write $z=T^{e/2}x$. In a local holomorphic frame of eigenvectors orthonormal for $h^{\dec}$,
\begin{equation}\label{eq:metric-matching-logarithm}
 \log\bigl((h^{\dec})^{-1}h^{\mathrm{in}}\bigr)=D_T+K_T,
\end{equation}
where $D_T$ is real diagonal and $K_T$ is off-diagonal. For every pair of integers $k,N\geq0$,
\begin{equation}\label{eq:metric-matching-bounds}
 \|D_T\|_{C^k_z}\leq C_kT^{-1/2},\qquad
 \|K_T\|_{C^k_z}\leq C_{k,N}T^{-N}.
\end{equation}
The logarithm has trace zero. These estimates hold on slightly larger annuli as well, uniformly over $K$, the Picard fibers and all arguments of $\lambda$.
\end{lemma}

\begin{proof}
On the annulus the eigenline bundles of the given Higgs field carry flat metrics, so a local holomorphic generator of each can be normalized to have constant norm. In the resulting frame the spectral estimates on the annulus, with every fixed number of $z$-derivatives, take the form
\begin{equation*}
 \begin{gathered}
 |v_i|\asymp T^{-1/2},\qquad
 |\boldsymbol\gamma|/|v_i|=O(T^{-1/2}),\\
 |\xi|^{-1/e}=O(T^{-1/2}),\qquad
 |\varepsilon||\xi|=O(T^{-e/2}).
 \end{gathered}
\end{equation*}
The regular factor in~\eqref{eq:exterior-line-metric} changes from its value at $(v,\boldsymbol\gamma)=(0,0)$ by $O(|v|+|\boldsymbol\gamma|)$. Each factor in $\prod_{i<j}(v-\gamma_i)$ is $v-\gamma_i=v(1+O(|\boldsymbol\gamma|/|v|))$. By the expansions at the separated roots, the normalized evaluation matrices, their inverses and their determinants differ from their leading terms by the same order. The unit $f_{\boldsymbol\gamma}(x)$ in~\eqref{eq:inner-determinant-factor} changes from $f_{\boldsymbol0}(0)$ by $O(|x|+|\boldsymbol\gamma|)$. All these differences are $O(T^{-1/2})$ with their $z$-derivatives.

The constants $s_\nu$ in~\eqref{eq:inner-block-constants} match the leading terms. By Proposition~\ref{prop:compact-polynomial-models}, the difference between the model diagonal entries and their leading terms is $O(|\xi|^{-1/e})$, and the off-diagonal part vanishes to every order. The isomorphism between the eigenline bundles in Lemma~\ref{lem:transported-inner-metric} preserves the decay of the off-diagonal entries. By the direct image and determinant formulas, the differences of the diagonal entries are bounded by a constant times $T^{-1/2}$. The corresponding differences off the diagonal are bounded by a constant times $T^{-N}$ for every $N$. The same statements hold after every fixed number of $z$-derivatives. These metric matrices and their inverses lie in a fixed bounded set of positive matrices, so the matrix logarithm preserves both estimates. Its diagonal and off-diagonal parts satisfy~\eqref{eq:metric-matching-bounds}. Both metrics have determinant one, so the logarithm has trace zero.
\end{proof}

\begin{proof}[Proof of Proposition~\ref{prop:approximate-metric}]
We choose a smooth radial cutoff $\chi(z)$ equal to one on $|z|\leq1$ and zero on $|z|\geq2$. On each gluing annulus we define
\begin{equation}\label{eq:approximate-metric-cutoff}
 h^{\app}=h^{\dec}
     \exp\!\left(\chi(z)\log\bigl((h^{\dec})^{-1}h^{\mathrm{in}}\bigr)\right).
\end{equation}
We use $h^{\mathrm{in}}$ on the inner disc and $h^{\dec}$ outside the annulus. The exponential of the cutoff times the Hermitian logarithm is smooth and positive. The logarithm has trace zero, so the metric induces the unit metric on the trivialized determinant line bundle. Near the punctures this construction agrees with~\eqref{eq:inner-metric}, giving the filtration and the asserted local model description.

In a local holomorphic frame of eigenvectors that is orthonormal for $h^{\dec}$, interpolation of the diagonal part alone has curvature $\bar\partial\partial(\chi D_T)$ and zero Higgs commutator. Two derivatives in $x$ contribute $T^e$, while~\eqref{eq:metric-matching-bounds} contributes $T^{-1/2}$. The diagonal curvature is therefore bounded by $CT^{e-1/2}$. The off-diagonal part changes the curvature by a quantity of order $T^{-N}$ for every $N$, because its derivatives satisfy the second estimate in~\eqref{eq:metric-matching-bounds}. For the Higgs commutator we obtain the same bound since, on the annulus,
\begin{equation*}
 |\Phi|_{g,h^{\dec}}\leq C T^{(e+1)/2},
\end{equation*}
so its quadratic contribution is bounded by $CT^{e+1}$ times the absolute value of an off-diagonal coefficient that decays to every order. The interpolated metric is uniformly comparable to $h^{\dec}$. Hence
\begin{equation*}
 |\mathscr E(h^{\app})|_{g,h^{\app}}
       \leq CT^{e-1/2}
\end{equation*}
on the gluing annulus. Its area is $O(T^{-e})$, and its contribution to the $L^{p_0}$ norm is at most
\begin{equation}\label{eq:annular-Lp-contribution}
 CT^{e(1-1/p_0)-1/2}.
\end{equation}

By Lemma~\ref{lem:transported-inner-metric}, the error in Hitchin's equation on the inner disc is bounded by $CT^{-e}$ in $L^{p_0}(|\dd\xi|^2)$. Under $\xi=\lambda^e x$, the norm of a two-form acquires the factor $T^{2e}$ and the area element acquires $T^{-2e}$. Since the fixed metric $g$ is uniformly comparable to $|\dd x|^2$ on this coordinate disc, the inner contribution is at most
\begin{equation}\label{eq:inner-Lp-contribution}
 CT^{2e-2e/p_0-e}=CT^{-e(2/p_0-1)}.
\end{equation}
The approximate metric solves Hitchin's equation outside the gluing annuli and inner discs by~\eqref{eq:exterior-flat-equations}. For $2\leq e\leq6$, the exponent in~\eqref{eq:annular-Lp-contribution} is at most $-\delta$, and $e(2/p_0-1)>\delta$ throughout~\eqref{eq:approximate-exponent-range}. Summing over the three punctures proves~\eqref{eq:approximate-Lp-defect}. For each $(\alpha,\sigma)\in K$, the local model metrics are unchanged as $\lambda$ and $L$ vary. Their uniform estimates, together with the bounds on the closures of finitely many charts covering the product of the Picard family and the unit circle, give common constants $C,T_0$.
\end{proof}

\subsection{Comparison with the harmonic metric}

\begin{proposition}\label{prop:approximate-to-harmonic-metric}
For $K$ and the exponents in~\eqref{eq:uniform-weight-margins}--\eqref{eq:approximate-defect-rate}, there are $C,c,T_0>0$, uniform over $K$, the Picard fibers and all arguments of $\lambda$, such that the metrics of Proposition~\ref{prop:approximate-metric} satisfy
\begin{equation}\label{eq:approximate-energy-gap}
 Q_{h^{\app}}(u)\geq c\|u\|_2^2
 \qquad(T\geq T_0,\ u\in V_{h^{\app}}).
\end{equation}
For $T\geq T_0$, the harmonic metric $h$ of determinant one on the same filtered Higgs bundle satisfies
\begin{equation}\label{eq:approximate-to-harmonic-rates}
 s=\log\bigl((h^{\app})^{-1}h\bigr)\in V_{h^{\app}},\qquad
 Q_{h^{\app}}(s)^{1/2}\leq CT^{-\delta},\qquad
 \|s\|_\infty\leq CT^{-\delta/2}.
\end{equation}
In particular,
\begin{equation}\label{eq:approximate-to-harmonic-forms}
 e^{-CT^{-\delta/2}}h^{\app}\leq h
                   \leq e^{CT^{-\delta/2}}h^{\app}.
\end{equation}
The compact set $K$ may contain parameters with coincident graded residues and parameters admitting strictly semistable objects outside the end.
\end{proposition}

\begin{proof}
Every object on the end of the Hitchin fibration is stable by Proposition~\ref{prop:filtered-spectral-correspondence}, so the harmonic metric $h$ exists by Proposition~\ref{prop:harmonic-metric-existence}. For each fixed object and $\lambda$, the metric $h$ and the metrics of the polynomial models used in $h^{\app}$ are each uniformly equivalent near the punctures to their prescribed diagonal metrics by Lemma~\ref{lem:tame-power-comparison}. The substitution $\xi=\lambda^{e_p}x_p$, the flag-compatible frame changes, and the positive scalar and block factors in~\eqref{eq:inner-metric} change only finite comparison constants. Positivity on the compact complement therefore gives
\begin{equation}\label{eq:harmonic-approximate-bounded-ratio}
 C^{-1}h^{\app}\leq h\leq C h^{\app}
 \qquad\text{on }X\setminus D
\end{equation}
for a finite constant $C$ depending on the object and $\lambda$.

If the lower bound for the quadratic forms failed, there would be parameters in $K$, $T_j\to\infty$, and $u_j\in V_{h^{\app}_j}$ with $\|u_j\|_2=1$ and $Q_{h^{\app}_j}(u_j)\to0$. After passage to a subsequence, the parameters, the unitary characters of the spectral line bundles and the unit complex numbers $\lambda_j/T_j$ converge. Since $\boldsymbol\gamma_j=\sigma_j/\lambda_j\to0$, the spectral metrics in the normalized family converge smoothly over compact subsets of $X\setminus D$ by Lemma~\ref{lem:exterior-spectral-metric}. The covering maps are unramified there, so the direct image metrics $h^{\dec}_j$ also converge smoothly. These compact subsets lie outside the gluing discs for large $j$, where $h^{\app}_j=h^{\dec}_j$. The Chern connections and $\lambda_j^{-1}\Phi_j$ therefore converge with all derivatives along the curve. Multiplication by the convergent sequence of unit complex numbers then gives convergence of $T_j^{-1}\Phi_j$.

The family of covers away from the punctures identifies the bundles compatibly on a compact exhaustion. Positive square roots of the metric matrices make these identifications unitary, preserving smooth convergence, trace and self-adjointness. By~\eqref{eq:exterior-flat-equations}, the limiting connection is flat and the limiting Higgs field is normal with distinct eigenvalues on $X\setminus D$. Its spectral cover is connected by Proposition~\ref{prop:normalized-spectral-family}, since it is obtained from $Y_{1,\boldsymbol0}$ by removing the preimages of the punctures. Convergence of the quadratic forms to zero on the chosen sections now contradicts Lemma~\ref{lem:energy-bound-under-convergence} with $b_j=T_j^{-1}$, proving~\eqref{eq:approximate-energy-gap} with uniform constants.

Contraction by the fixed K\"ahler form is bounded, so~\eqref{eq:approximate-Lp-defect} gives $\|\mu(h^{\app})\|_{p_0}\leq CT^{-\delta}$. The hypotheses of Proposition~\ref{prop:small-defect-metric-comparison} hold with $h_0=h^{\app}$ and $p=p_0$ by~\eqref{eq:harmonic-approximate-bounded-ratio} and~\eqref{eq:approximate-energy-gap}. The constants in that proposition are independent of the initial bound on the metric ratio, so we can increase $T_0$ once to obtain~\eqref{eq:approximate-to-harmonic-rates} and~\eqref{eq:approximate-to-harmonic-forms} for all these parameters.
\end{proof}

\subsection{Vanishing graded residues}

We now take $\sigma=0$ and allow the weights to vary over a compact set $K\subset\Acal$. The local coordinates and cyclic metrics are those of Section~\ref{subsec:strongly-parabolic-models}. Their exponential decay allows the gluing annuli to remain fixed in the coordinate $x$.

We use the origin $L_0=\pi^*\mathcal O_X(1)$ of the Picard variety to normalize the flat metric obtained by direct image. The section $\dd\pi/\eta_0$ of $L_0^2=\pi^*T_X$ vanishes on the ramification divisor. We prescribe squared norm one for this section off the divisor. The induced metric on $L_0$ is flat there and, in the coordinate above, has local expression $|w|^{-(e-1)}$, after absorbing a nonvanishing holomorphic factor in the local frame. The orthogonal direct image of its tensor product with the unitary flat metric on $L\otimes L_0^{-1}$ is $h^{\dec}$. In the base frame $e_j$ the exponent in the determinant metric is $\sum_j2b_j=0$, so the determinant metric extends as a flat metric over $X$. It is constant and can be normalized to one.

The local flat scalar factors and positive relative block factors from this construction also multiply the harmonic model metrics. The logarithms of the scalar factors are harmonic and the block factors are constant, so these rescalings preserve~\eqref{eq:zero-residue-Toda} and the norms of endomorphisms within each block. These normalizations agree on finitely many Picard charts, as in Lemma~\ref{lem:exterior-spectral-metric}.

On fixed disjoint discs about the punctures we choose smooth radial cutoffs $\chi_p$ equal to one on smaller discs and zero outside the chosen discs. Inside the discs we use~\eqref{eq:zero-residue-model-metric} with $\chi_pu_j$ in place of $u_j$ and with the scalar and block factors just described. Outside we use $h^{\dec}$. The resulting metric $h^{\app}$ has determinant one because $\sum_j u_j=0$. It solves the harmonic metric equation on the smaller discs and on the exterior.

\begin{lemma}\label{lem:zero-residue-metric-correction}
There are $C,c_2,T_0>0$ such that this approximate metric and the harmonic metric satisfy
\begin{equation}\label{eq:zero-residue-metric-correction}
 \bigl\|\log((h^{\app})^{-1}h)\bigr\|_\infty
          \leq C e^{-c_2T},\qquad T\geq T_0,
\end{equation}
uniformly over $K$, the Picard fibers and all arguments of $\lambda$.
\end{lemma}

\begin{proof}
On every cutoff annulus, $\rho$ is comparable to $T$. Formula~\eqref{eq:zero-residue-model-infinity} therefore bounds the error in Hitchin's equation exponentially. For any fixed $1<p<\infty$,
\begin{equation}\label{eq:zero-residue-defect}
 \|\ii\Lambda_g(F_{h^{\app}}+[\Phi,\Phi^{\dagger_{h^{\app}}}])\|_{L^p}
          \leq C e^{-cT},
\end{equation}
where $g$ is the fixed smooth metric on $X$. Derivatives of the cutoff and the change of scale introduce polynomial factors in $T$, which are absorbed by decreasing $c$ once.

On every compact subset of $X\setminus D$, the exponential estimate~\eqref{eq:zero-residue-model-infinity} makes the Chern connections of $h^{\app}$ converge to those of the exterior metric, along any subsequence for which the unit complex factors and the unitary characters of the spectral line bundles converge. The rescaled fields $T^{-1}\Phi$ have the same limit as in the proof of Proposition~\ref{prop:approximate-to-harmonic-metric}. The limiting Higgs field has distinct eigenvalues and connected spectral cover, so the uniform coercivity estimate~\eqref{eq:approximate-energy-gap} follows from Lemma~\ref{lem:energy-bound-under-convergence}. For each object the ratio between $h$ and $h^{\app}$ is bounded by Lemma~\ref{lem:tame-power-comparison}. Applying Proposition~\ref{prop:small-defect-metric-comparison} to~\eqref{eq:zero-residue-defect} proves~\eqref{eq:zero-residue-metric-correction}, with one decrease of the exponential rate and uniform constants.
\end{proof}

\section{Comparison of the Hitchin and semiflat metrics}

\label{sec:L2-geometry}

The comparison of Hermitian metrics in Section~\ref{sec:approximate-metrics} bounds the norms of representatives of tangent classes. Orthogonal projection to the $L^2$ harmonic representatives does not increase these norms, so it suffices for the upper bound on the Hitchin metric to construct closed representatives whose norms approach their spectral Hodge norms. The lower bound follows from equality of the volume forms. When all graded residues vanish, the algebraically decaying terms in the pairings cancel by the boundary integral identities. Interior estimates for the Ricci-flat K\"ahler equations then control the derivatives of the metric difference.

\subsection{\texorpdfstring{$L^2$}{L2} representatives of tangent classes}
\label{sec:tangent-representatives}

We fix a stable pair $(E,\Phi)$ with the prescribed weights and scalar graded residues. Its harmonic metric $h$ exists by Proposition~\ref{prop:harmonic-metric-existence}. The algebraic tangent space is isomorphic to the quotient of weighted Dolbeault cocycles in Proposition~\ref{prop:weighted-Dolbeault-comparison}. Orthogonal projection in $L^2$ gives the representatives used below. The gauge parameters may have arbitrary constant limits in each Levi block.

We write $X^*=X\setminus D$ and use on $X^*$ a conformal metric equal to $\dd\tau^2+\dd\theta^2$ near each puncture. For the weighted spaces~\eqref{eq:harmonic-weighted-norm} we fix
\begin{equation}\label{eq:tangent-weight-window}
 \varepsilon_*:=\min_{p\in D}\varepsilon_p,
 \qquad 0<d<\varepsilon_*\,,
\end{equation}
where $\varepsilon_p$ is defined in~\eqref{eq:harmonic-weight-frequency}. All formal adjoints below use the unweighted $L^2$ inner product defined by $h$. The weighted Sobolev norms defined by $h$ and $h_{\mathrm{ref}}$ are equivalent at every fixed order. Indeed, by Proposition~\ref{prop:harmonic-weighted-form} the metric ratio and its inverse tend to invertible block diagonal constants and have bounded cylindrical derivatives of all orders.

At $p$ we use the Lie algebra
\begin{equation*}
 \mathfrak l_p=
 \left\{\bigoplus_{j=0}^{e_p-1}c_j:
       c_j\in\End(\operatorname{Gr}_jE_p),\quad
       \sum_j\tr c_j=0\right\}.
\end{equation*}
For each end we choose a cutoff $\chi_p$ equal to one sufficiently far out and supported on that end. In the unitary frames of~\eqref{eq:harmonic-unitary-end} we define
\begin{equation}\label{eq:tangent-weighted-spaces}
 \begin{aligned}
 \mathcal R^0_d
   &=W^{2,2}_d(\End_0E)
      \oplus\left\{\sum_{p\in D}\chi_pc_p:c_p\in\mathfrak l_p\right\},\\
 \mathcal V_d
   &=W^{1,2}_d\bigl(\End_0E\otimes
                         (\overline K\oplus K)\bigr).
 \end{aligned}
\end{equation}
In these extended Sobolev domains the allowed limits lie in the Levi algebra, as in~\cite[Proposition~3.3 and equation~(3.7)]{CollierFredricksonWentworth}. Here $K$ and $\overline K$ denote the complex line bundles of $(1,0)$-forms and $(0,1)$-forms on $X^*$. We equip $\mathcal R^0_d$ with the sum norm. Different cutoffs define equivalent norms on this same space. Passing between the frames unitary for the diagonal metric and the $h$-unitary frames also preserves the space, since the limiting change of frame conjugates each matrix block to itself and its remainder decays exponentially.

We use the differentials $\dd_1,\dd_2$ in~\eqref{eq:tangent-differentials}. Each $c_p\in\mathfrak l_p$ commutes with both model matrices $A,B$ in~\eqref{eq:harmonic-local-data}. Thus $\dd_1(\chi_pc_p)$ is the sum of a compactly supported term and a term that decays exponentially. Using~\eqref{eq:harmonic-unitary-end} with a rate $\nu>d/2$, we obtain bounded maps
\begin{equation*}
 \dd_1:\mathcal R^0_d\longrightarrow\mathcal V_d,
 \qquad
 \dd_2:\mathcal V_d\longrightarrow
       L^2_d(\End_0E\otimes K\otimes\overline K),
 \qquad \dd_2\dd_1=0.
\end{equation*}

\subsubsection{Comparison with the deformation complex}

\begin{proposition}\label{prop:weighted-Dolbeault-comparison}
Inclusion of the smooth Dolbeault resolution of~\eqref{eq:deformation-complex-weak} induces an isomorphism
\begin{equation}\label{eq:weighted-Dolbeault-comparison}
 \mathbb H^1(X,\mathcal C^\bullet)
   \simeq\ker(\dd_2:\mathcal V_d\to L^2_d)
                    \big/\dd_1\mathcal R^0_d.
\end{equation}
Every class has a smooth representative $(b,\varphi)$ for which $b$ vanishes near $D$ and $\varphi$ is holomorphic there as a section of $\SParEnd_0E\otimes K_X(D)$.
\end{proposition}

\begin{proof} The proof is given in Appendix~\ref{proof:prop:weighted-Dolbeault-comparison}. \end{proof}

\subsubsection{Orthogonal projection and the trace pairing}

We use the index calculation for zero Higgs field in~\cite[Proposition~3.24]{CollierFredricksonWentworth} and verify its extension to the scalar graded residues considered here.

\begin{proposition}\label{prop:ordinary-L2-projection}
The operator
\begin{equation}\label{eq:tangent-gauge-isomorphism}
 L_h:=\dd_1^*\dd_1:\mathcal R^0_d
                    \longrightarrow L^2_d(\End_0E)
\end{equation}
is a bounded isomorphism. For every $v\in\mathcal V_d$ with $\dd_2v=0$, there is a unique $\eta\in\mathcal R^0_d$ such that
\begin{equation}\label{eq:tangent-coclosed-projection}
 v^h=v+\dd_1\eta,\qquad \dd_2v^h=0,
 \qquad \dd_1^*v^h=0.
\end{equation}
For a constant depending on the fixed pair, metric and $d$,
\begin{equation}\label{eq:tangent-projection-bound}
 \|\eta\|_{\mathcal R^0_d}
       +\|v^h\|_{W^{1,2}_d}
       \leq C\|v\|_{W^{1,2}_d}.
\end{equation}
The representative $v^h$ has finite $L^2$ norm and is orthogonal to $\dd_1\mathcal R^0_d$ with respect to the $L^2$ inner product.
\end{proposition}

\begin{proof} The proof is given in Appendix~\ref{proof:prop:ordinary-L2-projection}. \end{proof}

By Propositions~\ref{prop:weighted-Dolbeault-comparison} and~\ref{prop:ordinary-L2-projection}, the algebraic tangent space is isomorphic to
\begin{equation}\label{eq:tangent-harmonic-space}
 \mathcal H_{h,d}:=
 \{v\in\mathcal V_d:\dd_2v=0,\ \dd_1^*v=0\},
 \qquad \dim_\C\mathcal H_{h,d}=2.
\end{equation}
The squared norm~\eqref{eq:tangent-trial-norm} is twice the unweighted $L^2$ norm squared of the pair. For every $\zeta\in\mathcal R^0_d$, orthogonality of the harmonic representative to infinitesimal gauge transformations (Proposition~\ref{prop:ordinary-L2-projection}) gives
\begin{equation}\label{eq:tangent-minimum-norm}
 \|v^h+\dd_1\zeta\|_h^2
    =\|v^h\|_h^2+\|\dd_1\zeta\|_h^2,
 \qquad
 \|v^h\|_h^2
    =\inf_{\eta\in\mathcal R^0_d}\|v+\dd_1\eta\|_h^2.
\end{equation}
With this normalization, the unweighted norm of any representative in the weighted deformation complex bounds the Hitchin norm of its tangent class.

The trace form~\eqref{eq:semiflat-trace-form} extends to the weighted representatives with the same factor $2$. Its integral converges by the Cauchy--Schwarz inequality. To see that the form descends, we take $\varphi\in W^{1,2}_d$ and $\eta\in\mathcal R^0_d$. The form $\dd\tr(\varphi\eta)$ is integrable because derivatives of $\varphi$ pair with the bounded $\eta$ using the integrable weight $e^{-d\tau}$, and the remaining term is a product of $L^2$ factors. For the cutoffs $\kappa_R$ defined in Appendix~\ref{proof:prop:ordinary-L2-projection}, applying Stokes' theorem to $\kappa_R\tr(\varphi\eta)$ and using~\eqref{eq:tangent-cutoff-bound} give
\begin{equation}\label{eq:tangent-trace-boundary}
 \int_{X^*}\dd\tr(\varphi\eta)
   =-\lim_{R\to\infty}\int_{X^*}
                       \dd\kappa_R\wedge\tr(\varphi\eta)=0.
\end{equation}
If the first argument is changed by $\dd_1\eta$ and the second argument is $v_2=(b_2,\varphi_2)$, the change in the trace form is $2\int_{X^*}\dd\tr(\varphi_2\eta)-2\int_{X^*}\tr(\eta\,\dd_2v_2)$. Both terms vanish by~\eqref{eq:tangent-trace-boundary} and the cocycle condition on $v_2$. On the smooth representatives from Proposition~\ref{prop:weighted-Dolbeault-comparison}, the resulting pairing is the Dolbeault expression for the trace cup product followed by integration, with the factor $2$ fixed in~\eqref{eq:semiflat-trace-form}. It is nondegenerate by Proposition~\ref{prop:local-deformations}. Both the pairing and the isomorphism with algebraic deformation cohomology are independent of the metric used for the weighted Sobolev norms and the choice of cutoffs.

\subsection{The hyperk\"ahler structure}
\label{sec:hyperkahler-structure}

The tangent representatives that are harmonic for the unweighted $L^2$ inner product of Section~\ref{sec:tangent-representatives} define the Hitchin metric. We now construct smooth local solution spaces and prove that this metric is hyperk\"ahler. Dimakis--Fredrickson--Yae established smoothness of the local quotient and smooth dependence of the harmonic tangent representatives for these flag and residue data~\cite[Propositions~4.6--4.7 and Section~4.3]{DFY}. We use this construction with the derivative estimates and normalizations recorded in Appendix~\ref{proof:prop:smooth-harmonic-quotient}. The unitary gauge transformations have limits in the groups $K_p$ defined below, the maximal compact subgroups of the corresponding Levi subgroups. Those limits commute with both model matrices because the graded residues are scalar, and hence preserve the model connections and Higgs fields near the punctures.

Throughout this subsection we take $(\alpha,\sigma)$ such that semistability equals stability. We write $\M=\M_{\alpha,\sigma}$ and fix $d$ as in~\eqref{eq:tangent-weight-window}. The surface metric is cylindrical near $D$, and all adjoints and contractions use this metric and the relevant harmonic metric.

\subsubsection{The smooth quotient}

For $j\geq0$, we extend the notation~\eqref{eq:tangent-weighted-spaces} by setting
\begin{equation}\label{eq:hk-higher-spaces}
 \begin{aligned}
 \mathcal R^j_d
   &=W^{j+2,2}_d(\End_0E)
       \oplus\left\{\sum_p\chi_pc_p:c_p\in\mathfrak l_p\right\},\\
 \mathcal V^j_d
   &=W^{j+1,2}_d\bigl(\End_0E\otimes(\overline K\oplus K)\bigr),\\
 \mathcal U^j_d
   &=\{\xi\in\mathcal R^j_d:\xi^\dagger=-\xi\},
 \qquad
 \mathcal S^j_d=\ii\mathcal U^j_d .
 \end{aligned}
\end{equation}
Thus $\mathcal V^0_d=\mathcal V_d$. For a fixed $j\geq2$ we omit this index from $\mathcal G_d$ and $\mathcal Z_d$ below. Relative to a fixed smooth reference pair with the models~\eqref{eq:harmonic-unitary-end}, the configuration space has real tangent space $\mathcal V^j_d$, where a vector $(b,\varphi)$ changes the unitary connection by $b-b^\dagger$ and the Higgs field by $\varphi$. The group $\mathcal G_d$ consists of unitary gauge transformations of determinant one whose limit at $p$ belongs to
\begin{equation*}
 K_p=S\left(\prod_{a=0}^{e_p-1}U(\operatorname{Gr}_aE_p)\right)
\end{equation*}
and whose remainder belongs to $W^{j+2,2}_d$. Its Lie algebra is $\mathcal U^j_d$. We write $\mathcal Z_d$ for the space of solutions of the Hitchin equations whose filtered extensions are stable.

The products needed for the gauge action and metric equation are taken in
\begin{equation}\label{eq:hk-ambient-algebra}
 W^{j+2,2}_d(\End E)
 \oplus\left\{\sum_p\chi_pc_p:
        c_p\in\bigoplus_a\End(\operatorname{Gr}_aE_p)\right\}.
\end{equation}
Sobolev multiplication on unit cylindrical strips and the compact support of $(\chi_pc_p)(\chi_pd_p)-\chi_pc_pd_p$ and $\Id-\sum_p\chi_p\Id$ make~\eqref{eq:hk-ambient-algebra} a unital Banach algebra. Hence exponentiation, inversion on invertible elements, and local logarithms are smooth. The trace-free Hermitian and skew-Hermitian spaces in~\eqref{eq:hk-higher-spaces} impose the determinant and unitarity constraints, with matrix multiplication taken in~\eqref{eq:hk-ambient-algebra}.

\begin{proposition}\label{prop:smooth-harmonic-quotient}
\begin{enumerate}
\item The space $\mathcal Z_d$ is a smooth Banach manifold near every solution, and
\begin{equation}\label{eq:hk-principal-bundle}
 \mathcal Z_d\longrightarrow\M
\end{equation}
is a locally trivial smooth principal $\mathcal G_d/\mu_n$ bundle. The quotient complex structure agrees with the algebraic complex structure of $\M$.
\item In local analytic families of parabolic Higgs bundles with weights $\alpha$ and residues $\sigma$, we use $h_{\mathrm{ref}}$ as in Proposition~\ref{prop:harmonic-weighted-form}. For every $j\geq0$, the relative metric logarithms $\log(h_{\mathrm{ref}}^{-1}h)$ depend smoothly on the parameters in $\mathcal R^j_d$, and the differences between local solution representatives and the reference pair depend smoothly in $W^{j+1,2}_d$.

Starting with a holomorphic frame compatible with the flag, we write $H$ for the metric matrix in the associated frame unitary for $h_0$~\eqref{eq:harmonic-power-model}. Its positive Hermitian block diagonal limit $H_p$ depends smoothly on the parameters. On compact subsets of a parameter chart, all mixed derivatives of finite order in the moduli parameters and cylindrical variables of $H-H_p$ and of the solution coefficients minus the corresponding model coefficients near the punctures are $O(e^{-\nu\tau})$ for every $0<\nu<d/2$. The transition gauges satisfy the same estimates after subtraction of their smoothly varying limits in the Levi subgroup.
\end{enumerate}
\end{proposition}

\begin{proof} The proof is given in Appendix~\ref{proof:prop:smooth-harmonic-quotient}. \end{proof}

\subsubsection{The metric and the three complex structures}

\begin{proposition}\label{prop:hyperkahler-structure}
The real polarization of~\eqref{eq:tangent-trial-norm} on the harmonic representatives in $\mathcal H_{h,d}$ defines a smooth hyperk\"ahler metric $g_H$ on $\M$. The induced complex structures are~\eqref{eq:hk-quaternionic-operators}, with $I$ equal to the algebraic complex structure. If $\omega_Q(v,w)=g_H(Qv,w)$, then
\begin{equation}\label{eq:hk-symplectic-normalization}
 \Omega=\omega_J+\ii\omega_K,\qquad
 \Omega(v,w)=2\int_{X^*}\tr(\varphi\wedge c-q\wedge b)
\end{equation}
for $v=(b,\varphi)$ and $w=(c,q)$. In particular,
\begin{equation}\label{eq:hk-volume-normalization}
 \dd\operatorname{vol}_{g_H}
     =\frac{\omega_I^2}{2}
     =\frac{\Omega\wedge\overline\Omega}{4}.
\end{equation}
\end{proposition}

\begin{proof}
The constant forms $\omega_a$ on the configuration space are closed and gauge invariant. Their restrictions to $\mathcal Z_d$ vanish on gauge directions by~\eqref{eq:hk-moment-derivative}, so they descend along the principal bundle~\eqref{eq:hk-principal-bundle} to smooth closed two-forms on $\M$. The horizontal space for the $L^2$ pairing is $\mathcal H_{h,d}$ by~\eqref{eq:hk-four-operator-inverse}. The three endomorphisms $I,J,K$ preserve it by~\eqref{eq:hk-adjoint-identities}. Their restrictions define orthogonal quaternionic structures with the descended forms as fundamental forms. Thus the metric is hyperk\"ahler by~\cite[Lemma~6.8, p.~91]{Hitchin1987}, and $I$ agrees with the algebraic complex structure by Proposition~\ref{prop:smooth-harmonic-quotient}.

For $v=(B\,\dd\bar z,P\,\dd z)$ and $w=(C\,\dd\bar z,Q\,\dd z)$, polarization of~\eqref{eq:tangent-trial-norm} gives
\begin{equation*}
 g_H(v,w)=4\operatorname{Re}\int
              \tr(BC^\dagger+PQ^\dagger)\,\dd x\,\dd y.
\end{equation*}
Substitution of $J,K$ from~\eqref{eq:hk-quaternionic-operators} gives the real and imaginary parts of
$-4\ii\int\tr(PC-QB)\,\dd x\,\dd y$. Since
$\dd z\wedge \dd\bar z=-2\ii\,\dd x\,\dd y$, this is~\eqref{eq:hk-symplectic-normalization}. By the calculation of the trace pairing at the end of Section~\ref{sec:tangent-representatives}, the form in~\eqref{eq:hk-symplectic-normalization} agrees with~\eqref{eq:semiflat-trace-form}. In a quaternionic orthonormal frame in real dimension four,
$\omega_I^2=\omega_J^2=\omega_K^2=2\,\dd\operatorname{vol}_{g_H}$ and the distinct fundamental forms have zero wedge product. This proves~\eqref{eq:hk-volume-normalization}.
\end{proof}

\subsection{The metric on the stable locus}
\label{sec:Hitchin-semiflat-comparison}

On the end of the Hitchin fibration we fix the weights and scalar graded residues. We use the harmonic metric $h$ of determinant one and an admissible Sobolev weight as in~\eqref{eq:tangent-weight-window}. For an algebraic tangent vector $Z$ with representative $v_Z^h$ from Proposition~\ref{prop:ordinary-L2-projection}, we set
\begin{equation}\label{eq:harmonic-tangent-form}
 g_H(Z,Z)=\|v_Z^h\|_h^2,
\end{equation}
with the normalization~\eqref{eq:tangent-trial-norm}.

\begin{proposition}\label{prop:regular-end-hyperkahler}
For every $\alpha\in\Acal$ and every centered scalar graded residue parameter $\sigma$, the form $g_H$ in~\eqref{eq:harmonic-tangent-form} is a smooth hyperk\"ahler metric on the end of the Hitchin fibration. Its complex structure is the relative Picard complex structure, and its holomorphic symplectic form is~\eqref{eq:semiflat-trace-form}, with the period expression~\eqref{eq:semiflat-period-formulas}. When semistability equals stability for these data, it agrees with the restriction of the metric in Proposition~\ref{prop:hyperkahler-structure}.
\end{proposition}

\begin{proof}
With the weights and residues fixed, the correspondence of Proposition~\ref{prop:filtered-spectral-correspondence} sends a local Poincar\'e line bundle on the spectral family over a small contractible analytic chart $S$ in the end of the Hitchin fibration to a holomorphic family of stable parabolic Higgs bundles. Its determinant line bundle is pulled back from $S$ and becomes trivial after shrinking $S$. Base change gives the Kodaira--Spencer isomorphism.

On the product of each puncture disc with $S$ we use holomorphic frames compatible with the flags and the determinant trivialization. In these frames the Higgs variations have zero graded residue. By Lemma~\ref{lem:weighted-local-dolbeault}(i) and the lattice orders~\eqref{eq:tangent-lattice-orders}, they therefore depend smoothly on the parameters in the weighted spaces of~\eqref{eq:hk-metric-equation}. For this stable family we can now apply Lemma~\ref{lem:harmonic-family-dependence}. The resulting map from $S$ to the solution slice has derivative given by the Kodaira--Spencer map followed by $L^2$ harmonic projection, and is a local diffeomorphism.

The quaternionic calculation and descent in the proof of Proposition~\ref{prop:hyperkahler-structure} apply in each such chart, with the same normalization of the trace pairing. This yields the metric and the stated complex and holomorphic symplectic structures locally. The harmonic metrics are unique, so the local structures agree on overlaps.
\end{proof}

We work over a fixed nonempty compact set $K\subset\Acal\times\mathscr B$, with common exponents as in~\eqref{eq:uniform-weight-margins}--\eqref{eq:approximate-defect-rate}. At parameters in $K$ the graded residues may coincide, and strictly semistable objects may occur outside the end. For each parameter we use the harmonic metric $h$ and the metric $h^{\app}$ of Proposition~\ref{prop:approximate-metric}, with $t=\lambda^n$ and $T=|\lambda|$. All tangent vectors and derivatives are taken with the weights and residues fixed. We take the Sobolev weight
\begin{equation}\label{eq:trial-Sobolev-weight}
 0<d<\kappa_K/4
\end{equation}
in~\eqref{eq:tangent-weighted-spaces}. 
\subsection{Closed representatives for general residues}
\label{sec:metric-comparison}

The Gauss--Manin connection in Section~\ref{subsec:semiflat-geometry} decomposes each tangent vector as $Z=H+V$, with $H$ horizontal and $V$ vertical.

\begin{proposition}\label{prop:tangent-trial-comparison}
For every parameter in $K$ and all sufficiently large $T$, each point of the end of the Hitchin fibration admits a complex-linear choice $Z\mapsto v_Z$ of closed trace-free representatives in $\mathcal V_d$ with the following properties. The representatives $v_H$ and $v_V$ are orthogonal for both $h^{\app}$ and $h$, and
\begin{equation}\label{eq:tangent-trial-comparison}
 \begin{aligned}
 \|v_Z\|_{h^{\app}}^2
   &\leq (1+CT^{-1})g_{\sfm}(Z,Z),\\
 \|v_Z\|_h^2
   &\leq \bigl(1+C(T^{-1}+T^{-\delta/2})\bigr)g_{\sfm}(Z,Z).
 \end{aligned}
\end{equation}
The constants and the lower threshold for $T$ are uniform over $K$, the Picard fibers and all arguments of $\lambda$. They may depend on the exponents and choices listed in Proposition~\ref{prop:approximate-metric}.
\end{proposition}

\begin{proof}
We write $r=T^{-1/2}$. Near a puncture $p$ we use the coordinates in~\eqref{eq:local-model-scale} with $e=e_p$, $m=m_p$, and $\rho=|\xi|^{1/e}$. On the base we choose radial cutoffs satisfying
\begin{equation}\label{eq:trial-cutoff-regions}
 \begin{array}{c|cc}
       &\text{value }0&\text{value }1\\ \hline
 \chi_H& |x|\leq3r^e&|x|\geq4r^e\\
 \chi_V& |x|\geq6r^e&|x|\leq5r^e.
 \end{array}
\end{equation}
They have uniformly bounded Dirichlet energy. By Proposition~\ref{prop:local-spectral-matching} and smoothness of the normalized family, all ramification lies over $|x_p|<CT^{-e_p}$. Hence for large $T$ the cover is unramified over both cutoff annuli, where $h^{\app}=h^{\dec}$ and the eigenline bundles of the Higgs field are orthogonal.

We first construct the representative of the horizontal vector $H_\lambda$ defined by $\dd\lambda(H_\lambda)=1$. Over a small contractible disc in the $\lambda$-plane, the smooth proper family in Proposition~\ref{prop:normalized-spectral-family} is differentiably trivial, so inclusion of a fiber induces an isomorphism on fundamental groups. The unitary character representing $M\in\Pic^0(Y_t)$ thus extends to the total space. Its flat connection defines a holomorphic line bundle there, and
\begin{equation*}
 L=\pi^*\mathcal O_X(1)\otimes M
\end{equation*}
is a holomorphic family in which the unitary character is constant under parallel transport. Its image under the correspondence of Proposition~\ref{prop:filtered-spectral-correspondence} is the required horizontal family of parabolic Higgs bundles. Near each branch of the normalized curve we choose a parallel section of $M$ and a fixed local frame of $\mathcal O_X(1)$. Multiplying these sections by $f_j(w)=\prod_{i<j}(w-\sigma_i)$ gives the frame $f_0,\ldots,f_{e-1}$ in~\eqref{eq:newton-spectral-lattices}.

The Higgs field matrix in~\eqref{eq:exact-inner-matrix} can be written as
\begin{equation}\label{eq:horizontal-Newton-matrix}
 \begin{aligned}
 A_{\nu,\varepsilon}
    &=B+\xi\sum_{j=0}^{e-1}
           (c_\nu\delta_{j0}+d_{\nu,j}(\xi,\varepsilon))E_{j,e-1},\\
 B&=\diag(\sigma_0,\ldots,\sigma_{e-1})
                     +\sum_{j=0}^{e-2}E_{j+1,j}.
 \end{aligned}
\end{equation}
The residue $B=A_{c_\nu}(0)$ is independent of $\nu$ and $\varepsilon$. By~\eqref{eq:actual-inner-polynomial}, the coefficients $d_{\nu,j}$ and their derivatives $\xi\partial_\xi d_{\nu,j}$ and $\varepsilon\partial_\varepsilon d_{\nu,j}$ are $O(|\varepsilon|)$ on each fixed $\xi$-disc.

To estimate these coefficients on growing annuli, we omit $\nu$ from the notation. For $\rho$ large and $|\varepsilon\xi|$ small,
\begin{equation}\label{eq:horizontal-coefficient-bound}
 \rho^j\bigl(|d_j|+|\xi\partial_\xi d_j|
                 +|\varepsilon\partial_\varepsilon d_j|\bigr)
                 \leq C|\varepsilon\xi|.
\end{equation}
Indeed, evaluating the last column at the roots $w_i$ of the polynomial for this block gives
\begin{equation*}
 \sum_jd_j f_j(w_i)=u(w_i)-c_\nu.
\end{equation*}
The right side is $O(|\varepsilon\xi|)$ by Proposition~\ref{prop:local-spectral-matching}. If $\zeta^e=\xi$, the normalized roots stay uniformly separated, so $(\zeta^{-j}f_j(w_i))$ and its inverse are uniformly bounded. Solving this linear system proves the bound for $d_j$. Holomorphic branches of the roots on slightly larger $\xi$-sectors and $\varepsilon$-discs satisfy the same estimates as long as $|\varepsilon\xi|$ remains small. The two derivative bounds follow by Cauchy's estimates on discs with radii comparable to $|\xi|$ and $|\varepsilon|$. The coefficients extend holomorphically to $\varepsilon=0$. The coefficients defining the finite flat algebras and the normalized evaluation matrices have common bounds on the larger spectral charts of Section~\ref{sec:approximate-metrics}, so the coefficient and Cauchy estimates are uniform over $K$, including coincident residues.

A dot denotes $\partial_\lambda$ with $x$, the weights and the residues fixed. Since $\xi=\lambda^e x$ and $\varepsilon=\lambda^{-e}$,
\begin{equation}\label{eq:horizontal-local-variation}
 \begin{aligned}
 \dot A_{\nu,\varepsilon}
   &=\frac e\lambda(\xi\partial_\xi-\varepsilon\partial_\varepsilon)
                                      A_{\nu,\varepsilon}\\
   &=\frac{e\xi}{\lambda}\sum_{j=0}^{e-1}
       \bigl(c_\nu\delta_{j0}+d_{\nu,j}
          +\xi\partial_\xi d_{\nu,j}
          -\varepsilon\partial_\varepsilon d_{\nu,j}\bigr)E_{j,e-1}.
 \end{aligned}
\end{equation}
The residue $B$ is fixed. The local representative is $(0,\varphi_{\mathrm{in}})$ with $\varphi_{\mathrm{in}}=\bigoplus_\nu\dot A_{\nu,\varepsilon}\,\dd x/x$. Its trace is zero because it is the derivative of the global trace-free family of Higgs fields, for which the sum of the block traces vanishes.

We estimate this representative in the polynomial metrics $h_{c_\nu}$. The estimate in Proposition~\ref{prop:compact-polynomial-models} is uniform over $K$ and gives
\begin{equation*}
 h_{c_\nu}\asymp\diag(|\xi|^{2\alpha_0},\ldots,
                              |\xi|^{2\alpha_{e-1}})
\end{equation*}
near zero. With $\kappa_p=1+\alpha_{p,0}-\alpha_{p,e_p-1}\geq\kappa_K$ we have $|\xi E_{j,e-1}|_{h_{c_\nu}}\leq C|\xi|^{1+\alpha_j-\alpha_{e-1}}$, so the contribution of a small fixed $\xi$-disc is bounded by $CT^{-2}\int_0^1s^{2\kappa_p}\,\dd s/s\leq CT^{-2}/(2\kappa_K)$. Near infinity we have $h_{c_\nu}\asymp\diag(\rho^{2j-(e-1)})$ by Proposition~\ref{prop:compact-polynomial-models}, and hence
\begin{equation*}
 |\xi E_{j,e-1}|_{h_{c_\nu}}\leq C\rho^{j+1},\qquad
 |\dot A_{\nu,\varepsilon}|_{h_{c_\nu}}
             \leq CT^{-1}\rho(1+|\varepsilon|\rho^e).
\end{equation*}
On $|x|\leq4r^e$ one has $\rho\leq CT^{1/2}$ and $|\varepsilon|\rho^e\leq CT^{-e/2}$. A fixed intermediate annulus contributes $O(T^{-2})$. Conformal invariance of the norm of one-forms and $\dd x/x=\dd\xi/\xi$ therefore give
\begin{equation}\label{eq:horizontal-inner-energy}
 \|\varphi_{\mathrm{in}}\|_{\oplus h_{c_\nu},\,|x|\leq4r^e}^2
 \leq CT^{-2}\left(1+
                  \int_1^{CT^{1/2}}\rho^2\frac{\dd\rho}{\rho}\right)
 \leq CT^{-1}.
\end{equation}
We use the normalization of~\eqref{eq:tangent-trial-norm} for the squared norm of each component.

On the exterior the unramified spectral cover over the product of a fixed base domain with the parameter disc is a product covering. Parallel sections of $M$ and the fixed $\mathcal O_X(1)$ frames give a representative $(0,\varphi_{\mathrm{ex}})$, whose diagonal entries are $\dot w_i\,\dd x/x$. On the overlap we write
\begin{equation*}
 W_{ij}=f_j(w_i),\qquad S=W^{-1},\qquad
 \Gamma=\dot S S^{-1}.
\end{equation*}
Differentiating $S^{-1}AS=\diag(w_i)$ and expressing the result in the frame $f_0,\ldots,f_{e-1}$ gives
\begin{equation}\label{eq:horizontal-transition-variation}
 \varphi_{\mathrm{ex}}=\varphi_{\mathrm{in}}+[\Phi,\Gamma],
 \qquad \bar\partial_E\Gamma=0.
\end{equation}
These formulas are invariant under sheet permutations. The root equations give $|\dot w_i|\leq C\rho/T$ on the transition annulus. Thus $|\dot W_{ij}|\leq C\rho^j/T$ and $\dot W W^{-1}=O(T^{-1})$. Formula~\eqref{eq:exterior-line-metric} shows that the squared norms of the local generators of the eigenline bundles within each block have uniformly bounded ratios there, since $|v_i|\asymp r$ and the regular factors are uniformly bounded above and below. It follows that $|\Gamma|_{h^{\app}}\leq C/T$.

The determinant line bundle of the family is pulled back from the parameter disc by Proposition~\ref{prop:filtered-spectral-correspondence}. A holomorphic trivialization of this line bundle allows us to normalize the local frames by scalar $n$th roots. The Higgs matrices remain unchanged, and the transition derivative becomes
\begin{equation*}
 \Gamma^0=\Gamma-\frac{\tr\Gamma}{n}\Id,
 \qquad |\Gamma^0|_{h^{\app}}\leq C/T.
\end{equation*}
The derivative is independent of the choice of root of unity. We use $(0,\varphi_{\mathrm{in}})+\dd_1(\chi_H\Gamma^0)$ on the puncture neighborhoods and $(0,\varphi_{\mathrm{ex}})$ outside to define the horizontal representative. In the local frame its components are
\begin{equation}\label{eq:horizontal-trial-representative}
 b_H=(\bar\partial\chi_H)\Gamma^0,\qquad
 \varphi_H=(1-\chi_H)\varphi_{\mathrm{in}}
                               +\chi_H\varphi_{\mathrm{ex}}.
\end{equation}
The transition relation~\eqref{eq:horizontal-transition-variation} shows that these formulas represent the horizontal family. They give a closed trace-free pair, holomorphic in this frame at ramification points, with $b_H=0$ and $\operatorname{Res}_p\varphi_H=0$ at every $p\in D$. The bound on the Dirichlet energy of the cutoff gives $\|b_H\|_{h^{\app}}^2\leq CT^{-2}$.

To use the estimate~\eqref{eq:horizontal-inner-energy} for the approximate metric, we compare the norms in that estimate with those defined by $h^{\app}$. In the chosen frames,~\eqref{eq:inner-metric} reads
\begin{equation*}
 d_\lambda(x)\bigoplus_\nu s_\nu
                (B_\nu^{-1})^*h_{c_\nu}B_\nu^{-1}.
\end{equation*}
The common scalar and the positive block constants cancel from norms of block diagonal endomorphisms. By~\eqref{eq:eigenline-alignment-bound} and Proposition~\ref{prop:compact-polynomial-models}, the endomorphisms $B_\nu$ and $B_\nu^{-1}$, which equal the identity on a fixed inner disc, are uniformly bounded in the model norms on the growing annulus. The interpolation preserves this comparison by~\eqref{eq:metric-matching-bounds}. On $2r^e\leq|x|\leq4r^e$, the comparison follows from~\eqref{eq:exterior-line-metric} and the normalized evaluation matrices in the proof of Lemma~\ref{lem:inner-outer-metric-matching}. The same bounds apply on this fixed enlargement of the annulus, so the required endomorphism norms are uniformly comparable throughout $|x|\leq4r^e$.

The contribution of the inner disc to $\|\varphi_H\|_{h^{\app}}^2$ is at most $CT^{-1}$ by~\eqref{eq:horizontal-inner-energy}. The exterior variation on the transition annulus satisfies $|\varphi_{\mathrm{ex}}|\leq CT^{-1}\rho\,|\dd x/x|$, giving the same bound for its contribution. Outside the discs $|x_p|\leq4r^{e_p}$, orthogonal direct image preserves the spectral Hodge norm. The derivative of the eigenvalue differential at fixed $x$ is
\begin{equation*}
 \nu_\lambda=n\lambda^{n-1}\nu_t.
\end{equation*}
Its squared $L^2$ norm over $Y_{t,\sigma}$ is $g_{\sfm}(H_\lambda,H_\lambda)$. Under the scaling to $Y_{1,\boldsymbol\gamma}$, the differential becomes $n\nu_{1,\boldsymbol\gamma}$, where $\nu_{1,\boldsymbol\gamma}$ is defined by~\eqref{eq:fixed-residue-differential} for the normalized family. It converges to a nonzero holomorphic differential as $\boldsymbol\gamma\to0$. Its squared norm is therefore bounded above and away from zero. Combining these bounds proves
\begin{equation}\label{eq:horizontal-trial-bound}
 \|v_H\|_{h^{\app}}^2
       \leq(1+CT^{-1})g_{\sfm}(H,H)
\end{equation}
for every horizontal vector $H$, by complex-linear extension from $H_\lambda$.

For the vertical representative we use the harmonic $(0,1)$-form $\upsilon$ on $Y_{1,\boldsymbol\gamma}$ representing $V$. The notation $\|\cdot\|_{\mathrm{sp}}$ refers to the spectral $L^2$ norm of Section~\ref{subsec:semiflat-geometry}. To compute the trace of this form under $\pi$, we use coordinates $x=z^k$ at a ramification point. The trace of a holomorphic form $\sum_{a\geq0}c_a z^a\,\dd z$ is $\sum_{\ell\geq0}c_{k\ell+k-1}x^\ell\,\dd x$, which extends across the branch value. The global trace vanishes since $\Pp^1$ has no holomorphic one-forms. By conjugation, $\operatorname{Tr}_\pi\upsilon=0$.

On the $\nu$th spectral coordinate disc we take the antiholomorphic primitive $f_\nu^0$ of $\upsilon$ with $f_\nu^0(0)=0$. Smooth dependence of harmonic $(0,1)$-forms in the normalized family, restricted to a fixed compact neighborhood of $\boldsymbol\gamma=0$ in its parameter space, gives
\begin{equation*}
 |\upsilon|\leq C\|\upsilon\|_{\mathrm{sp}},\qquad
 |f_\nu^0(v)|\leq C|v|\,\|\upsilon\|_{\mathrm{sp}}.
\end{equation*}
The function $\sum_\nu\operatorname{Tr}_\pi f_\nu^0$ extends antiholomorphically across the branch values. Its $\bar\partial$ derivative is zero, so it is a constant $C_p$. Taking the limit at $x=0$ gives
\begin{equation}\label{eq:vertical-primitive-trace}
 C_p=\sum_{\nu=1}^{m_p}\sum_{j=0}^{e_p-1}
                 f_\nu^0(\gamma_{p,j}),\qquad
 f_\nu=f_\nu^0-C_p/n,
 \qquad \sum_\nu\operatorname{Tr}_\pi f_\nu=0.
\end{equation}
Repeated values of $\gamma_{p,j}$ are counted with multiplicity in this formula. The estimates for the local primitives imply $|C_p|\leq C|\boldsymbol\gamma|\,\|\upsilon\|_{\mathrm{sp}}$ throughout the same ball in the space of graded residues.

The local graph~\eqref{eq:local-normalized-branches} has the form $x=P_{\boldsymbol\gamma}(v)/U_\nu(v,\boldsymbol\gamma)$ with $U_\nu$ a bounded holomorphic unit. Since $|\boldsymbol\gamma|=O(T^{-1})$, the inverse image of $|x|\leq6r^e$ lies in $|v|\leq Cr$, where $|\partial_vx|\leq Cr^{e-1}$. Hence the derivative of the pulled-back cutoff $\chi_V$ is $O(r^{-1})$, its support has area $O(r^2)$, and $|f_\nu|\leq Cr\|\upsilon\|_{\mathrm{sp}}$. We set $F=\chi_V f_\nu$ on these discs and sum the extensions by zero over the punctures. Then
\begin{equation}\label{eq:vertical-cutoff-form}
 \upsilon_r=\upsilon-\bar\partial F,
 \qquad
 \|\bar\partial F\|_{\mathrm{sp}}^2
                 \leq Cr^2\|\upsilon\|_{\mathrm{sp}}^2.
\end{equation}
The form $\upsilon_r$ vanishes over $|x_p|\leq5r^{e_p}$, including every ramification point. Its trace is zero, since the cutoff is pulled back from the base on every sheet and~\eqref{eq:vertical-primitive-trace} gives $\operatorname{Tr}_\pi F=0$. Harmonic orthogonality on the compact curve gives the exact identity
\begin{equation}\label{eq:vertical-trial-Gram}
 \|\upsilon_r\|_{\mathrm{sp}}^2
       =\|\upsilon\|_{\mathrm{sp}}^2
                           +\|\bar\partial F\|_{\mathrm{sp}}^2.
\end{equation}

The diagonal direct image of $\upsilon_r$ on the unramified locus extends by zero near the branch values to a smooth trace-free form $b_V$. This form commutes with $\Phi$, and the pair $v_V=(b_V,0)$ is closed. The line bundle deformation $\bar\partial_L+s\upsilon_r$, for a small complex parameter $s$, is constant near the punctures and ramification points, and its direct image is represented by $(b_V,0)$. Since $[\upsilon_r]=[\upsilon]$, this deformation represents $V$ under the relative spectral correspondence. On the support of $b_V$, the eigenline bundles of the given Higgs field are orthogonal for $h^{\app}$, so the norm of $v_V$ equals $\|\upsilon_r\|_{\mathrm{sp}}$. Equations~\eqref{eq:vertical-cutoff-form} and~\eqref{eq:vertical-trial-Gram} imply
\begin{equation}\label{eq:vertical-trial-bound}
 \|v_V\|_{h^{\app}}^2
        \leq(1+CT^{-1})g_{\sfm}(V,V).
\end{equation}
The primitives, trace correction and cutoff are complex-linear in $\upsilon$.

It remains to check that these representatives lie in the weighted Sobolev space defined using $h$. Every nonzero frequency $k+\alpha_a-\alpha_b$ in~\eqref{eq:harmonic-weight-frequency} has absolute value at least $\kappa_K$. In~\eqref{eq:tangent-normal-operator} the residue term is nonnegative, and the scalar graded residues preserve the constant Levi summands. These facts make the weight $d$ in~\eqref{eq:trial-Sobolev-weight} admissible throughout $K$. Near $D$, the horizontal representative is $(0,C(x)\,\dd x)$ with $C$ holomorphic in the flag-adapted frame, and the vertical representative vanishes. In a frame unitary for $h_0$, the logarithmic coefficients of the horizontal representative and all their cylindrical derivatives are $O(e^{-\kappa_K\tau})$. For each object, the positive Hermitian limits on the graded quotients imply that the weighted Sobolev norms defined by the diagonal metric and by $h$ are equivalent, by Proposition~\ref{prop:harmonic-weighted-form}. Thus both representatives belong to $\mathcal V_d$ and represent the required classes under~\eqref{eq:weighted-Dolbeault-comparison}. The uniform norm bounds come from the estimates for the model metrics.

Finally, $b_H$ is supported in $3r^e<|x|<4r^e$, where $b_V=0$, and $\varphi_V=0$ everywhere. The horizontal and vertical representatives are therefore pointwise orthogonal for every Hermitian bundle metric. For $v_Z=v_H+v_V$, the first estimate in~\eqref{eq:tangent-trial-comparison} now follows from~\eqref{eq:horizontal-trial-bound} and~\eqref{eq:vertical-trial-bound} because horizontal and vertical vectors are orthogonal for $g_{\sfm}$. By~\eqref{eq:approximate-to-harmonic-rates}, $h=h^{\app}e^s$ with $\|s\|_\infty\leq CT^{-\delta/2}$. The induced squared norm on endomorphisms changes by at most $e^{2\|s\|_\infty}$, proving the second estimate. The constants are uniform over $K$, the Picard fibers and all arguments of $\lambda$ by Propositions~\ref{prop:approximate-metric} and~\ref{prop:approximate-to-harmonic-metric}.
\end{proof}

\subsection{Boundary pairings when the graded residues vanish}
\label{sec:zero-residue-asymptotics}

Throughout this subsection $\sigma=0$ and the weights lie in a nonempty compact set $K\subset\Acal$. For the cyclic models of Section~\ref{subsec:strongly-parabolic-models} and the approximate metric on fixed discs in Section~\ref{sec:approximate-metrics}, we integrate by parts in the horizontal and vertical pairings. Cancellation of the resulting boundary terms leaves an exponentially small error.

On one rank-$e$ block, the scalar factor in~\eqref{eq:zero-residue-model-metric} cancels from endomorphism norms. Writing $H=\operatorname{diag}(e^{2u_j})$, we use the Hermitian pairing, linear in its first argument,
\begin{equation}\label{eq:zero-residue-local-pairing}
 \langle(a,\varphi),(b,\psi)\rangle_H
   =\frac{2\ii}{e}\int
       \tr(\varphi\wedge\psi^{\dagger_H}-a\wedge b^{\dagger_H}).
\end{equation}
The factor $1/e$ accounts for the degree of the ramified cover. We write
\begin{equation*}
 D''s=(\bar\partial_Es,[\Phi,s]),\qquad
 D's=\partial_hs+[\Phi^{\dagger_h},s],
\end{equation*}
with the usual graded extensions to forms. On an annulus $r_-<|w|<r_+$, integration by parts gives
\begin{equation}\label{eq:zero-residue-Green}
 \begin{aligned}
 \langle D''s,(a,\varphi)\rangle_H
   &=\langle s,(D'')^*(a,\varphi)\rangle_H
          +\mathfrak d_{r_+}(s,(a,\varphi))
          -\mathfrak d_{r_-}(s,(a,\varphi)),\\
 \mathfrak d_r(s,(a,\varphi))
   &=-\frac{2\ii}{e}\int_{|w|=r}\tr(s a^{\dagger_H}).
 \end{aligned}
\end{equation}
The circles in the last integral have counterclockwise orientation. On zero-forms the pairing is $(2/e)\int\tr(s t^{\dagger_H})\,\dd A$, where $\dd A$ is Euclidean area.

For an argument $\vartheta$ of $\lambda$, we have $\lambda=T e^{\ii\vartheta}$ and write $\dot\lambda=e^{\ii\vartheta}\xi$ with $\xi\in\C$. The horizontal representative on this block is
\begin{equation}\label{eq:zero-residue-horizontal-model}
 \begin{gathered}
 F_\xi=(0,e e^{\ii\vartheta}\xi P\,\dd w),\qquad
 s_\xi=\frac{\xi}{T}\rho\operatorname{diag}(u_j'),\\
 H_\xi=F_\xi-D''s_\xi,\qquad
 d_j=u_{j+1}-u_j.
 \end{gathered}
\end{equation}
Its Higgs component is
\begin{equation*}
 (H_\xi)^{1,0}
       =e e^{\ii\vartheta}\xi
          \sum_j(1+\rho d_j')E_{j+1,j}\,\dd w,
\end{equation*}
and its connection component is $-\bar\partial s_\xi$. These expressions are complex linear in $\xi$ and are $D''$-closed.

To verify the Coulomb condition for the unweighted $L^2$ pairing, we write
\begin{equation*}
 N_j=e^{2d_{j-1}}-e^{2d_j},\qquad q_j=\rho u_j'.
\end{equation*}
Differentiating~\eqref{eq:zero-residue-Toda} under radial scaling gives
\begin{equation}\label{eq:zero-residue-horizontal-adjoint}
 \begin{gathered}
 q_j''+\rho^{-1}q_j'=2\rho N_j'+4N_j,\qquad
 [P^{\dagger_H},[P,\operatorname{diag}(q_j)]]_{jj}
            =\frac{\rho}{2}N_j',\\
 (D'')^*D''\operatorname{diag}(q_j)
       =-2e^2T^2\operatorname{diag}(N_j).
 \end{gathered}
\end{equation}
For the last identity, the operator $(D'')^*D''$ on diagonal gauges is
$-2\partial_w\partial_{\bar w}+2e^2T^2[P^{\dagger_H},[P,\cdot]]$.
Also $(D'')^*F_\xi=-2e^2T\xi\operatorname{diag}(N_j)$, so
$(D'')^*H_\xi=0$.

Formula~\eqref{eq:zero-residue-model-zero} bounds the connection and Higgs coefficients of $H_\xi$ by
\begin{equation}\label{eq:zero-residue-horizontal-zero}
 C|\xi|\rho^{2\kappa-1},\qquad
 C|\xi|\rho^{\kappa-1},
\end{equation}
respectively, in the $\dd w,\dd\bar w$ coframe and the model endomorphism norm. For $\xi_2\in\C$, these bounds give
$|\mathfrak d_r(s_{\xi_2},H_\xi)|
 \leq C|\xi\xi_2|T^{-2}\rho^{2\kappa}$,
where $\rho=eTr$; this boundary pairing tends to zero as $r\to0$.

For a fixed outer radius $r$ we write $S(\rho)=\sum_j e^{2d_j(\rho)}$. Angular integration in~\eqref{eq:zero-residue-local-pairing} gives the exact identity
\begin{equation}\label{eq:zero-residue-horizontal-cancellation}
 \begin{aligned}
 &\langle H_\xi,F_{\xi_2}\rangle_{H,\,|w|<r}
       -\langle F_\xi,F_{\xi_2}\rangle_{\dec,\,|w|<r}\\
 &\quad=\frac{8\pi}{e}\frac{\xi\overline{\xi_2}}{T^2}
        \int_0^{eTr}\rho
               \left(S-e+\frac{\rho}{2}S'\right)\dd\rho,\\
 &\rho\left(S-e+\frac{\rho}{2}S'\right)
          =\frac12\frac{\dd}{\dd\rho}\bigl(\rho^2(S-e)\bigr).
 \end{aligned}
\end{equation}
The pairing marked $\dec$ uses the flat metric obtained by orthogonal direct image, for which $u_j=0$. Since each term $\rho^2e^{2d_j}$ is $O(\rho^{2e\delta_j})$ near zero, the boundary contribution at $\rho=0$ vanishes, while at the fixed outer radius it is exponentially small by~\eqref{eq:zero-residue-model-infinity}. Applying~\eqref{eq:zero-residue-Green} in the second argument replaces $F_{\xi_2}$ by $H_{\xi_2}$. The inner term vanishes by~\eqref{eq:zero-residue-horizontal-zero} and the outer term is exponentially small. Thus
\begin{equation}\label{eq:zero-residue-horizontal-Gram}
 \left|\langle H_\xi,H_{\xi_2}\rangle_{H,\,|w|<r}
       -\langle F_\xi,F_{\xi_2}\rangle_{\dec,\,|w|<r}\right|
       \leq C e^{-cT}|\xi\xi_2|.
\end{equation}

For a vertical tangent vector, the harmonic spectral $(0,1)$-form is $\upsilon=b\,\dd\bar w$ on this disc, since it is a constant multiple of $\overline{\eta_0}$. Its primitive vanishing at the center is $\beta_{\mathrm{loc}}=b\bar w$. With $M$ denoting multiplication on the spectral line bundle, we set
\begin{equation}\label{eq:zero-residue-vertical-model}
 S_\upsilon=M_{\overline{\beta_{\mathrm{loc}}}}^{\dagger_H}
       =b\bar w\sum_j e^{2d_j}E_{j,j+1},\qquad
 V_\upsilon=D''S_\upsilon.
\end{equation}
The adjoint of $S_\upsilon$ is holomorphic and commutes with $\Phi$, hence $D'S_\upsilon=0$. By the K\"ahler identities~\cite[Lemma~3.1(a), p.~876]{Simpson1988} for this harmonic bundle, we have
$(D'')^*D''S_\upsilon=(D')^*D'S_\upsilon=0$.
Integrability gives $D''V_\upsilon=0$. At zero, formula~\eqref{eq:zero-residue-model-zero} yields
\begin{equation}\label{eq:zero-residue-vertical-zero}
 \begin{aligned}
 |S_\upsilon|_H&\leq C|b|T^{-1}\rho^\kappa,\\
 |\bar\partial S_\upsilon|_H&\leq C|b|\rho^{\kappa-1},\qquad
 |[\Phi,S_\upsilon]|_H\leq C|b|\rho^{2\kappa-1}.
 \end{aligned}
\end{equation}
For these norms we express the one-forms in the Euclidean coframe. The bounds imply
$|\mathfrak d_r(S_\upsilon,V_\upsilon)|
 \leq C|b|^2T^{-2}\rho^{2\kappa}$,
which tends to zero as $r\to0$. Therefore
\begin{equation}\label{eq:zero-residue-vertical-boundary}
 \|V_\upsilon\|_{H,\,|w|<r}^2
         =\mathfrak d_r(S_\upsilon,V_\upsilon).
\end{equation}
For the flat metric obtained by orthogonal direct image, the corresponding $S_\upsilon$ and $V_\upsilon$ are $M_{\beta_{\mathrm{loc}}}$ and $(M_\upsilon,0)$, respectively. Their boundary data differ from those in~\eqref{eq:zero-residue-vertical-boundary} by exponentially small terms. Subtraction and polarization give
\begin{equation}\label{eq:zero-residue-vertical-Gram}
 \left|\langle V_\upsilon,V_{\upsilon_2}\rangle_{H,\,|w|<r}
       -\langle(M_\upsilon,0),(M_{\upsilon_2},0)
                                  \rangle_{\dec,\,|w|<r}\right|
       \leq C e^{-cT}|b b_2|.
\end{equation}

\subsection{Global representatives in the weighted Dolbeault complex}

On the $\lambda$-cover, a horizontal family has fixed spectral line bundle and fixed underlying bundle, with $\Phi=\lambda\Phi_1$. If $\nu=\dot\lambda\eta_0$, its closed representative $F_\nu=(0,M_\nu)$ restricts to $F_\xi$ above. Summing over all punctures and their ramified blocks, we define
\begin{equation}\label{eq:zero-residue-global-horizontal}
 H'=F_\nu-D''\sum\chi_p s_\xi.
\end{equation}
For the vertical representative we first define
\begin{equation}\label{eq:zero-residue-global-vertical}
 \begin{aligned}
 \upsilon'&=\upsilon-
      \sum\bar\partial\bigl((\pi^*\chi_p)\beta_{\mathrm{loc}}\bigr),\\
 V'&=(M_{\upsilon'},0)+D''\sum\chi_p S_\upsilon.
 \end{aligned}
\end{equation}
The form $\upsilon'$ vanishes on the smaller discs. Its direct image is a smooth Dolbeault variation preserving the flags and the vanishing graded residues. On a cutoff annulus,
\begin{equation*}
 V'=(M_\upsilon,0)
      +D''\sum\chi_p(S_\upsilon-M_{\beta_{\mathrm{loc}}}).
\end{equation*}
These expressions show that $H'$ and $V'$ are closed and represent the horizontal and vertical spectral classes of Proposition~\ref{prop:semiflat-period-formulas}.

All expressions are trace-free. The matrices $\operatorname{diag}(u_j')$ have trace zero, and the cyclic multiplication matrices in~\eqref{eq:zero-residue-cyclic-field} and~\eqref{eq:zero-residue-vertical-model} have zero diagonal. The trace of the harmonic spectral one-form on $\Pp^1$ also vanishes. Multiplication by the base cutoffs preserves these properties. Each local construction is invariant under the deck transformations of the ramified cover and therefore descends to its block on the base.

By~\eqref{eq:zero-residue-model-infinity}, the corrections to the direct image representatives are exponentially small on the annuli. On the smaller discs it follows from~\eqref{eq:zero-residue-horizontal-Gram} and~\eqref{eq:zero-residue-vertical-Gram} that the corresponding Gram matrices differ by exponentially small terms. For $h^{\app}$ the horizontal and vertical representatives are pointwise orthogonal. In each diagonal cyclic frame the connection and Higgs components of $H'$ have grades $0$ and $+1$, while those of $V'$ have grades $-1$ and $0$. Here $E_{ab}$ has grade $a-b$ modulo $e$, so corresponding components have distinct grades even when $e=2$. Outside the discs their nonzero components have different form types. The exact spectral norm formula of Proposition~\ref{prop:semiflat-period-formulas} now gives, for the tangent vector $Z$ represented by $H'+V'$,
\begin{equation}\label{eq:zero-residue-trial-upper}
 \|H'+V'\|_{h^{\app}}^2
        \leq(1+C e^{-cT})g_{\sfm}(Z,Z).
\end{equation}
Here the global squared norm is~\eqref{eq:tangent-trial-norm}. The squared spectral norm on horizontal vectors is a fixed positive multiple of $|\dot\lambda|^2$. On vertical vectors it is a fixed positive quadratic form on the tangent space of the character torus of $Y_1$. Thus the absolute local errors give the relative error in~\eqref{eq:zero-residue-trial-upper}, uniformly over the Picard fibers and the unit circle.

To apply $L^2$ harmonic projection, we place these representatives and their gauge corrections in the weighted Dolbeault spaces defined using $h$. In terms of $\tau_x=-\log|x|$ the radial variable satisfies, for each fixed $T$,
\begin{equation*}
 \rho^a=(eT)^a e^{-a\tau_x/e}.
\end{equation*}
In passing from $\dd w,\dd\bar w$ to the cylindrical coframe, the one-form coefficients acquire a factor $|w|$ up to the fixed constant $e$. Using~\eqref{eq:zero-residue-model-zero}, \eqref{eq:zero-residue-horizontal-zero} and~\eqref{eq:zero-residue-vertical-zero} and the corresponding estimates for cylindrical derivatives, we obtain
\begin{equation}\label{eq:zero-residue-weighted-coefficients}
 \begin{gathered}
 s_\xi=\frac{\xi}{T}\operatorname{diag}(\gamma_j)
                       +O(\rho^{2\kappa}),\qquad
 S_\upsilon=O(\rho^\kappa),\\
 H_\xi,V_\upsilon=O(\rho^\kappa).
 \end{gathered}
\end{equation}
In the last estimate the one-form coefficients are taken in the cylindrical coframe. The constant term in $s_\xi$ belongs to the Levi algebra, with the full multiplicity of every graded piece at a partial flag. We write
\begin{equation*}
 \varepsilon_K=\min_{\alpha\in K,\ p\in D}
                         \varepsilon_p>0,
\end{equation*}
with $\varepsilon_p$ as in~\eqref{eq:harmonic-weight-frequency}. We can therefore take a common Sobolev weight
\begin{equation}\label{eq:zero-residue-Sobolev-weight}
 0<d<\min\left\{\varepsilon_K,
                         \frac{2\kappa}{\max_p e_p}\right\}.
\end{equation}
Since the squared weight is $e^{d\tau_x}$, the coefficient estimates~\eqref{eq:zero-residue-weighted-coefficients} show that the gauge corrections belong to $\mathcal R_d^0$ and the representatives to $\mathcal V_d$. For each object these weighted norms are equivalent to those defined using $h$ by Proposition~\ref{prop:harmonic-weighted-form}, because the harmonic metric has positive Hermitian limits on the graded quotients and the remainders have bounded cylindrical derivatives. For membership in these spaces the constants of equivalence may depend on the object.

The uncorrected horizontal representative $F_\nu$ has the same cylindrical decay, and $(M_{\upsilon'},0)$ vanishes near the punctures. Thus
\begin{equation}\label{eq:zero-residue-quotient-equality}
 H'+V'=F_\nu+(M_{\upsilon'},0)
            +D''\left(-\sum\chi_p s_\xi
                            +\sum\chi_p S_\upsilon\right)
\end{equation}
is an equality in the weighted Dolbeault quotient of Proposition~\ref{prop:weighted-Dolbeault-comparison}. We can therefore project the representatives orthogonally in $L^2$ as in Proposition~\ref{prop:ordinary-L2-projection}.

\subsection{Metric comparison from tangent representatives}

\begin{lemma}\label{lem:metric-comparison-from-trials}
At a point of the Picard family over sufficiently large $|t|$, if
every tangent vector $Z$ has a closed representative
$v_Z\in\mathcal V_d$ of its algebraic cohomology class satisfying
\begin{equation*}
 \|v_Z\|_h^2\leq(1+q)g_{\sfm}(Z,Z)
\end{equation*}
for some $q\geq0$, then
\begin{equation*}
 (1+q)^{-1}g_{\sfm}\leq g_H\leq(1+q)g_{\sfm},
 \qquad
 \|g_{\sfm}^{-1}g_H-\Id\|_{\mathrm{op},g_{\sfm}}\leq q.
\end{equation*}
\end{lemma}

\begin{proof}
By Proposition~\ref{prop:ordinary-L2-projection} and~\eqref{eq:tangent-minimum-norm}, the harmonic representative satisfies $\|v_Z^h\|_h\leq\|v_Z\|_h$. Hence $g_H\leq(1+q)g_{\sfm}$.

By Propositions~\ref{prop:semiflat-period-formulas} and~\ref{prop:regular-end-hyperkahler}, the two metrics are hyperk\"ahler for the same complex structure and holomorphic symplectic form on the stable end. With the normalization~\eqref{eq:hk-volume-normalization}, both volume forms equal $\Omega\wedge\overline\Omega/4$. We denote the two complex eigenvalues of $g_{\sfm}^{-1}g_H$ by $0<\mu_1\leq\mu_2$. The volume identity and the upper bound give
\begin{equation*}
 \mu_1\mu_2=\det_{\C}(g_{\sfm}^{-1}g_H)=1,
 \qquad \mu_2\leq1+q\,.
\end{equation*}
Thus $(1+q)^{-1}\leq\mu_1\leq\mu_2\leq1+q$, proving the metric inequalities. Moreover,
\begin{equation*}
 \|g_{\sfm}^{-1}g_H-\Id\|_{\mathrm{op},g_{\sfm}}
 \leq\max\left\{q,1-\frac1{1+q}\right\}=q\,.
\end{equation*}
\end{proof}

For a compact set $K\subset\Acal\times\mathscr B$ the comparison uses the exponents $\eta,p_0,\delta$ as in~\eqref{eq:approximate-exponent-range}--\eqref{eq:approximate-defect-rate}.

\begin{proposition}\label{prop:harmonic-tangent-comparison}
There are constants $C,T_0$ depending on the compact set $K$ and the exponents above such that, on every end of the Hitchin fibration with $T\geq T_0$,
\begin{equation}\label{eq:harmonic-tangent-comparison}
 (1+\eta_T)^{-1}g_{\sfm}\leq g_H\leq(1+\eta_T)g_{\sfm},
 \qquad \eta_T=C(T^{-1}+T^{-\delta/2}).
\end{equation}
In particular,
\begin{equation}\label{eq:harmonic-tangent-operator-bound}
 \bigl\|g_{\sfm}^{-1}g_H-\Id\bigr\|_{\mathrm{op},g_{\sfm}}
                        \leq\eta_T.
\end{equation}
The bounds are uniform over $K$, the Picard fibers and all arguments of $\lambda$.
\end{proposition}

\begin{proof}
We apply Lemma~\ref{lem:metric-comparison-from-trials} to the representatives of Proposition~\ref{prop:tangent-trial-comparison}, with $q=C(T^{-1}+T^{-\delta/2})$. The bounds for these representatives hold uniformly over the compact parameter set, the Picard fibers and all arguments of the scaling parameter.
\end{proof}

For $\sigma=0$ we take the weights in an arbitrary nonempty compact set $K\subset\Acal$.

\begin{corollary}\label{cor:zero-residue-C0-comparison}
When all graded residues vanish, there are constants $C,c,T_0>0$ such that
\begin{equation}\label{eq:zero-residue-C0-comparison}
 (1+C e^{-cT})^{-1}g_{\sfm}\leq g_H
                  \leq(1+C e^{-cT})g_{\sfm},
 \qquad T\geq T_0.
\end{equation}
The constants are uniform for weights in $K$, all points of the Picard fibers and
all arguments of $\lambda$.
\end{corollary}

\begin{proof}
By Lemma~\ref{lem:zero-residue-metric-correction}, replacing the approximate metric by the harmonic metric changes the endomorphism-valued one-form norms by a factor $1+O(e^{-cT})$. The representatives in~\eqref{eq:zero-residue-quotient-equality} therefore satisfy the hypothesis of Lemma~\ref{lem:metric-comparison-from-trials}, by~\eqref{eq:zero-residue-trial-upper}, with $q=C e^{-cT}$. These representatives lie in the weighted Dolbeault space by the coefficient estimates near the punctures.
\end{proof}

\subsection{Higher derivatives}
\label{sec:metric-derivatives}

We fix a nonempty compact set $K\subset\Acal\times\mathscr B$ of weights and scalar graded residues. For each parameter in $K$ we use the Hitchin metric of Proposition~\ref{prop:regular-end-hyperkahler} on the end of the fibration, with the smooth identifications of the ends and the radius $R$ of~\eqref{eq:semiflat-ALG-model}. Derivatives are taken along each moduli space, with the weights and residues fixed. Strictly semistable objects may occur outside the end.

\begin{lemma}\label{lem:ricci-flat-interior-comparison}
For a family of pairs of smooth Ricci-flat K\"ahler metrics $G$ and $S$ on the Euclidean ball $B_2\subset\C^2$ with common ellipticity bounds on $B_2$ and uniform bounds for $S$ in every $C^j$ norm there, the following estimate holds for every $k\geq0$,
\begin{equation}\label{eq:metric-interior-difference-estimate}
 \|G-S\|_{C^k(B_1)}
       \leq C_k\|G-S\|_{C^0(B_2)}.
\end{equation}
The constants depend only on the ellipticity bounds and finitely many coefficient bounds for $S$. The norms use the fixed Euclidean coordinates.
\end{lemma}

\begin{proof}

We apply Sherman--Weinkove's interior estimate~\cite[Corollary~1.2, p.~492]{ShermanWeinkove} to the stationary K\"ahler--Ricci flow $G(s)=G$, $0\leq s\leq1$, with a fixed Euclidean background on fixed balls compactly contained in $B_2$. At $s=1$, the common ellipticity bounds give uniform $C^j$ bounds for $G$ on a fixed smaller ball for every $j$. Uniform ellipticity gives the same bounds for $G^{-1}$. The constants are uniform for the family.

For the difference of the metrics we use the equation satisfied by the coefficients of a Ricci-flat K\"ahler metric
\begin{equation}\label{eq:ricci-flat-coefficient-equation}
 g^{k\bar l}\partial_k\partial_{\bar l}g_{i\bar j}
   =g^{k\bar l}g^{p\bar q}
       (\partial_i g_{k\bar q})(\partial_{\bar j}g_{p\bar l}),
\end{equation}
which follows from $\partial_i\partial_{\bar j}\log\det g=0$ and the K\"ahler symmetry $\partial_i\partial_{\bar j}g_{k\bar l}=\partial_k\partial_{\bar l}g_{i\bar j}$. For $H=G-S$ we subtract the two equations, retaining $G^{k\bar l}\partial_k\partial_{\bar l}$ as the principal part and using
\begin{equation*}
 G^{-1}-S^{-1}=-G^{-1}HS^{-1}.
\end{equation*}
Every remaining difference is linear in $H$ or its first derivatives. The change in the principal coefficients multiplies second derivatives of $S$. Thus the four real components $H_A$ satisfy
\begin{equation}\label{eq:metric-difference-system}
 a^{uv}\partial_u\partial_vH_A
       +B_{AB}^{u}\partial_uH_B+C_{AB}H_B=0.
\end{equation}
Here $u,v$ index the four real coordinates and repeated indices are summed. The principal part is the same uniformly elliptic scalar operator on each component. The derivative bounds for $G,S$ and their inverses give uniform bounds of every fixed order for $a,B,C$.

We apply the scalar interior Schauder estimate~\cite[Theorem~6.2, p.~90]{GilbargTrudinger} to each component of~\eqref{eq:metric-difference-system} and absorb the lower-order terms using the interior interpolation inequality~\cite[Lemma~6.32, p.~130]{GilbargTrudinger}. Differentiating the system on fixed nested balls between $B_1$ and $B_2$ gives~\eqref{eq:metric-interior-difference-estimate}. The constants depend on finitely many coefficient bounds and the fixed separation of these balls.

\end{proof}

\begin{proposition}\label{prop:metric-derivative-comparison}
For every $0<\nu<1/4$, there exist a radius $R_0=R_0(K,\nu)>0$ and constants $C_{k,K,\nu}$ for all integers $k\geq0$ such that
\begin{equation}\label{eq:metric-derivative-comparison}
 \begin{aligned}
 &\bigl|(\nabla^{g_{\sfm}})^k(g_H-g_{\sfm})\bigr|_{g_{\sfm}}
   +\bigl|(\nabla^{g_{\sfm}})^k(\omega_I-\omega_{\sfm})\bigr|_{g_{\sfm}}\\
 &\hspace{45mm}\leq C_{k,K,\nu}R^{-\nu-k}
 \end{aligned}
\end{equation}
for $R\geq R_0$, uniformly over $K$, all points of the Picard fibers and all arguments of $\lambda$. The radius $R_0$ is independent of $k$. Here $\nabla^{g_{\sfm}}$ is the Levi--Civita connection of the semiflat metric. The identity for the holomorphic symplectic forms is
\begin{equation}\label{eq:metric-exact-symplectic-form}
 \omega_J+\ii\omega_K=\Omega,
\end{equation}
where $\Omega$ is the form in~\eqref{eq:semiflat-period-formulas}.
\end{proposition}

\begin{proof}
On each end of the Hitchin fibration, the metric $g_H$ of Proposition~\ref{prop:regular-end-hyperkahler} is smooth and hyperk\"ahler and has squared norm~\eqref{eq:tangent-trial-norm} on $L^2$ harmonic representatives. With this normalization, the trace pairing~\eqref{eq:semiflat-trace-form} yields the identity~\eqref{eq:metric-exact-symplectic-form} for the holomorphic symplectic form.

For fixed $0<\nu<1/4$, the bound $a_K<1$ in~\eqref{eq:uniform-weight-margins} allows us to choose common exponents $\eta,p_0$ in~\eqref{eq:approximate-exponent-range}. With $p_0$ sufficiently close to one, we have
\begin{equation*}
 \frac{\delta}{2}
       =\frac14-3\left(1-\frac1{p_0}\right)>\nu.
\end{equation*}
Above a common threshold in $T$, we obtain from Proposition~\ref{prop:harmonic-tangent-comparison}
\begin{equation}\label{eq:metric-C0-rate}
 |g_H-g_{\sfm}|_{g_{\sfm}}\leq C_{K,\nu}T^{-\nu}.
\end{equation}
The constants are uniform over $K$, the Picard fibers and all arguments of $\lambda$. Since $R$ is a fixed positive multiple of $T=|\lambda|$, this also gives the asserted zeroth-order bound in $R$.

We construct holomorphic covering charts whose radii grow with $T$. On a simply connected disc in the $\lambda$-plane the holomorphic fiber coordinate from Proposition~\ref{prop:semiflat-period-formulas} is
\begin{equation*}
 w=\theta_B-\tau(\lambda)\theta_A
\end{equation*}
on the universal cover of the Picard fibers. Its periods are $2\pi$ and $2\pi\tau(\lambda)$. The relative Picard family is the quotient by these holomorphic translations, so the map from the disc times $\C$ to the end is a holomorphic local diffeomorphism. We may therefore pull back the metric equations to discs of arbitrary radius in $w$.

The period functions in Lemma~\ref{lem:exterior-period-expansion} are holomorphic on a ball about zero in the space of graded residues, with vanishing first derivatives at zero. We use two smaller concentric balls whose closures are nested inside it. The projection of $K$ to the space of graded residues is bounded, so $\sigma/\lambda$ lies in the inner ball for all parameters once $T$ exceeds one common threshold. By Cauchy's estimates on these two fixed balls,~\eqref{eq:semiflat-period-derivative-rates} holds after any fixed number of derivatives in $\lambda$. The constants may depend on that number. Above a larger common threshold, $a'$ is nonzero and $\operatorname{Im}\tau$ lies in a fixed compact subinterval of $(0,\infty)$.

At a point over $\lambda_*$ we write $L=|\lambda_*|$ and choose a lift $w_*$ in a bounded fundamental parallelogram. The preceding period bounds make these parallelograms uniformly bounded over $K$. For fixed $0<c<1/10$ the covering chart is given by the map
\begin{equation}\label{eq:metric-growing-cover}
 \lambda=\lambda_*+cLz_1,\qquad
 w=w_*+cLz_2,\qquad |z_1|,|z_2|<1.
\end{equation}
For $L$ above a common threshold, a neighborhood of the closed polydisc maps into the end, and $|\lambda|$ is comparable to $L$ throughout it. The rescaled metrics $G$ and $S$ are the pullbacks of $g_H$ and $g_{\sfm}$, respectively, divided by $(cL)^2$. We write $b(\lambda)=\operatorname{Im}\tau(\lambda)$ and $b_\infty=\operatorname{Im}\tau_\infty$, with primes denoting $\lambda$ derivatives. The semiflat formula becomes
\begin{equation}\label{eq:metric-scaled-semiflat}
 S=4b|a'|^2|\dd z_1|^2
       +b^{-1}|\dd z_2+\theta_A\tau'\,\dd z_1|^2,
 \qquad \theta_A=-\frac{\operatorname{Im}w}{b}.
\end{equation}
Equation~\eqref{eq:semiflat-period-derivative-rates} gives, in these fixed coordinates,
\begin{equation*}
 a'-A_0=O(L^{-2}),\qquad
 \tau-\tau_\infty=O(L^{-2}),\qquad
 \tau'=O(L^{-3}),
\end{equation*}
with the same orders after any fixed number of real $z_1,z_2$ derivatives, uniformly over $K$. For each fixed $j$, the $C^j$ norm of $\theta_A$ in the $z$ coordinates is $O_j(L)$. Each derivative of a holomorphic coefficient contributes a factor $cL$ to its $\lambda$ derivative, and $w$ is affine in $z_2$. Consequently, for every $j\geq0$,
\begin{equation}\label{eq:metric-scaled-semiflat-bounds}
 \|S-S_\infty\|_{C^j}\leq C_{j,K}L^{-2},\qquad
 S_\infty=4b_\infty|A_0|^2|\dd z_1|^2
                         +b_\infty^{-1}|\dd z_2|^2.
\end{equation}
The norms here use real coordinates on the fixed polydisc. Thus $S$ and its inverse have uniform bounds in every $C^j$ norm. Their common ellipticity bounds require only the case $j=0$. Equation~\eqref{eq:metric-C0-rate} gives the same ellipticity bounds for $G$ and
\begin{equation}\label{eq:metric-scaled-C0}
 \|G-S\|_{C^0}\leq C_{K,\nu}L^{-\nu}.
\end{equation}

The metrics $G$ and $S$ are Ricci-flat K\"ahler, and \eqref{eq:metric-scaled-semiflat-bounds} gives the required bounds for $S$. Applying Lemma~\ref{lem:ricci-flat-interior-comparison} on fixed nested balls compactly contained in the polydisc and writing $H=G-S$, we obtain

\begin{equation}\label{eq:metric-scaled-difference-bounds}
 \|H\|_{C^k(B_{\mathrm{in}})}
       \leq C_{k,K}\|H\|_{C^0(B_{\mathrm{out}})}
       \leq C_{k,K,\nu}L^{-\nu}.
\end{equation}
The constants depend on finitely many coefficient bounds and the fixed separation of the balls. All derivative orders use the same lower threshold for $L$.

Since the connection coefficients of $S$ and their derivatives are bounded, the coordinate estimate~\eqref{eq:metric-scaled-difference-bounds} also bounds $|\nabla_S^kH|_S$ at the center. Scaling back multiplies this norm by $(cL)^{-k}$ because the difference tensor carries the factor $(cL)^2$, while its $k$th covariant derivative has $k+2$ covariant indices. Thus
\begin{equation*}
 \bigl|(\nabla^{g_{\sfm}})^k(g_H-g_{\sfm})\bigr|_{g_{\sfm}}
       \leq C_{k,K,\nu}L^{-\nu-k}.
\end{equation*}
The bounds are uniform over $K$, the center and its chosen Picard lift, so they hold throughout every end in the family. Finally $I$ is parallel for $g_{\sfm}$ and
$\omega_I-\omega_{\sfm}=(g_H-g_{\sfm})(I\,\cdot,\cdot)$.
The same estimate holds for this two-form. The relation between $R$ and $L$ proves~\eqref{eq:metric-derivative-comparison}. The threshold depends only on the balls chosen in the space of graded residues, the ellipticity bounds and the initial $C^0$ comparison, so the same threshold works for all derivative orders.
\end{proof}

\begin{proof}[Proof of Theorem~\ref{thm:main-asymptotics}\textup{(i)}]
The metric, radius, lattice and finite quotient are fixed by~\eqref{eq:semiflat-ALG-coordinates}--\eqref{eq:semiflat-ALG-model}. Equation~\eqref{eq:semiflat-ALG-rate} and the formula for the Levi--Civita connection give
\begin{equation*}
 \bigl|(\nabla^{g_0})^j
          (\nabla^{g_{\sfm}}-\nabla^{g_0})\bigr|_{g_0}
                  \leq C_{j,K}R^{-3-j}.
\end{equation*}
The period estimates make these bounds uniform over $K$ on one common end for every $j$. The bound for the difference of the two connections shows that~\eqref{eq:metric-derivative-comparison} holds with the same rate for covariant derivatives with respect to the flat connection. Adding~\eqref{eq:semiflat-ALG-rate} gives the metric and $\omega_I$ estimates in~\eqref{eq:main-ALG-rate}. The estimates for $\omega_J$ and $\omega_K$ follow from~\eqref{eq:metric-exact-symplectic-form} and the bound for $\Omega-\Omega_0$ in~\eqref{eq:semiflat-ALG-rate}. The finite deck action preserves the flat model, so the estimates descend to the quotient. The cone angle and limiting fiber area are $2\pi/n$ and $4\pi^2$, as in Proposition~\ref{prop:semiflat-ALG-end}.
\end{proof}

\begin{corollary}\label{cor:actual-ALG-end}
For every nonempty compact set $K\subset\Acal\times\mathscr B$ and every $0<\nu<1/4$, there exist $R_0=R_0(K,\nu)>0$ and constants $C_{k,K,\nu}$ for all integers $k\geq0$ such that
\begin{equation}\label{eq:actual-curvature-decay}
 \bigl|(\nabla^{g_H})^k\Rm(g_H)\bigr|_{g_H}
       \leq C_{k,K,\nu}R^{-2-\nu-k},
 \qquad R\geq R_0.
\end{equation}
The estimates are uniform over $K$, the Picard fibers and all arguments of $\lambda$, and $R_0$ is independent of $k$.
\end{corollary}

\begin{proof}
In local Cartesian coordinates for $g_0$, equation~\eqref{eq:main-ALG-rate} bounds the $j$th derivatives of the metric difference by $C_{j,K,\nu}R^{-\nu-j}$. The metric and its inverse have common bounds. The coordinate curvature formula and differentiation of the inverse metric give~\eqref{eq:actual-curvature-decay}. Conversion to covariant derivatives uses $\nabla^{g_H}-\nabla^{g_0}=O(R^{-\nu-1})$ and its differentiated estimates. These operations change only the constants, so the same $R_0$ suffices for every $k$.
\end{proof}

\begin{proof}[Proof of Theorem~\ref{thm:main-asymptotics}\textup{(ii)}]
For the compact set $K\subset\Acal$ in the statement we use the Hitchin metric of Proposition~\ref{prop:regular-end-hyperkahler} on the end.
By Proposition~\ref{prop:semiflat-ALG-end}, when all graded residues vanish,
\begin{equation*}
 g_{\sfm}=g_0,\qquad \omega_{\sfm}=\omega_0,
 \qquad \Omega=\Omega_0,
 \qquad R=2\sqrt{\operatorname{Im}\tau_\infty}|A_0|T.
\end{equation*}
The complex coordinate $\theta_B-\tau_\infty\theta_A$ on the universal cover of the Picard fiber is holomorphic, and $\tau_\infty$ is constant. Thus $(U,V)$ are common holomorphic coordinates for $g_H$ and $g_0$ upstairs. In these coordinates the estimate of Corollary~\ref{cor:zero-residue-C0-comparison} reads
\begin{equation*}
 |g_H-g_0|_{g_0}\leq C e^{-cR}.
\end{equation*}

The Euclidean balls of radius two in the $(U,V)$ cover lie in the end whenever $R$ at their centers exceeds a common threshold. The radial coordinate $R$ differs from its value at the center by at most two, giving a uniform exponential $C^0$ bound and ellipticity constant. The equations hold on these pullbacks even when the covering maps overlap in the torus direction.

The metric $g_H$ is Ricci-flat K\"ahler by Proposition~\ref{prop:regular-end-hyperkahler}, and $g_0$ is flat K\"ahler in the common holomorphic coordinates $(U,V)$. We can therefore apply Lemma~\ref{lem:ricci-flat-interior-comparison} with $G=g_H$ and $S=g_0$ on these balls. Since the flat metric $S$ has constant coefficients, we obtain
\begin{equation}\label{eq:zero-residue-interior-derivatives}
 \|g_H-g_0\|_{C^k(B_1)}
        \leq C_k\|g_H-g_0\|_{C^0(B_2)}
        \leq C_k e^{-cR}.
\end{equation}
Here $R$ is the center's radial coordinate, and the norms are Euclidean. The fixed radii and common ellipticity and $C^0$ bounds give constants uniform over $K$, with the same exponential rate and lower threshold for $R$ at every derivative order.

The complex structure $I$ is constant in these coordinates, and
$\omega_I-\omega_0=(g_H-g_0)(I\,\cdot,\cdot)$.
Hence the two-form difference satisfies the same estimates. The symplectic form in Proposition~\ref{prop:regular-end-hyperkahler} equals $\Omega_0$ here, proving~\eqref{eq:zero-residue-symplectic}. All tensors and estimates are invariant under~\eqref{eq:semiflat-ALG-deck-action}, so they descend to the asserted quotient end. This proves the theorem with common $c$ and $R_0$.
\end{proof}

\section{Elliptic compactification and ALG geometry}
\label{sec:elliptic-ALG-geometry}

We use the elliptic compactification to determine the topology of the moduli surface and to complete the proof of its ALG geometry.

\subsection{The elliptic compactification}
\label{sec:elliptic-compactification}

The relative spectral correspondence describes every object over the open set $U_{\sigma}$ of Proposition~\ref{prop:normalized-spectral-family}. To determine the whole surface, we first vary the scalar graded residues at fixed rational weights for which every semistable object is stable. A smooth trivialization reduces connectedness to the case where all graded residues vanish.

\subsubsection{The family over the residue space}

We work with rational weights $\widehat{\alpha}\in\Acal$ such that
\begin{equation}\label{eq:residue-family-weight-choice}
 d+\sum_{p,j}k_{p,j}\widehat\alpha_{p,j}\ne0
\end{equation}
for every numerical type with
$1\leq r<n$, $0\leq k_{p,j}\leq m_p$,
$\sum_jk_{p,j}=r$ at each puncture, and $d\in\Z$ with $|d|\leq3r$.
The equalities excluded by~\eqref{eq:residue-family-weight-choice} define finitely many proper hyperplanes, as in the proof of Proposition~\ref{prop:joint-walls}. In particular, such rational weights exist, and semistability equals stability for them at every $\sigma\in\mathscr B$.

\begin{proposition}\label{prop:residue-family}
The moduli spaces $\M_{\widehat{\alpha},\sigma}$ form a quasi-projective family of coarse moduli spaces
\begin{equation}\label{eq:residue-family-maps}
 \rho:\mathcal M_{\mathscr B}\longrightarrow\mathscr B,
 \qquad
 (\rho,h):\mathcal M_{\mathscr B}
                   \longrightarrow\mathscr B\times\C_t.
\end{equation}
The map $\rho$ is smooth and surjective of relative complex dimension two, and $(\rho,h)$ is projective. For every bounded set of residues, the latter map is smooth over all sufficiently large $|t|$.

For an open ball $V\subset\mathscr B$ with compact closure and center $\sigma_0$, and for sufficiently large $r_0>0$, there is a smooth trivialization
\begin{equation}\label{eq:residue-family-trivialization}
 \rho^{-1}(V)\simeq
       \M_{\widehat{\alpha},\sigma_0}\times V
\end{equation}
which preserves $h$ on $|h|\geq r_0$ and carries the compact manifolds with boundary $\{|h|\leq r_0\}$ in the fibers to one another.
\end{proposition}

\begin{proof}
The geometric invariant theory (GIT) presentation in Proposition~\ref{prop:algebraic-moduli} extends over $\mathscr B$ by imposing~\eqref{eq:universal-scalar-residues} with variable $\sigma$ on the universal graded quotients. These invariant closed equations describe the family over nonreduced bases as well. Exactness of invariants identifies the coarse quotient $\mathcal M_{\mathscr B}$ with a closed subscheme of the product of the parabolic Higgs moduli scheme and $\mathscr B$. The morphism recording the characteristic coefficients is projective. The same argument for the determinant over $\mathscr B$ shows that specifying a trivialization of the determinant line bundle gives the same coarse family.

To prove relative smoothness, we consider a small extension $A'\to A$ of local Artin $\C$-algebras with kernel $I$ annihilated by the maximal ideal. Given a stable object over $A$ and a residue lift $\sigma'\in\mathscr B(A')$, we choose trace-free local lifts of the Higgs field with these scalar graded residues in frames adapted to the flags and compatible with the determinant trivialization. The flag-preserving frame groups are smooth, so the transitions also lift. On each overlap the two lifts of the Higgs field differ by a term with zero graded residue, since the scalars agree and commute with the adapted transitions. This also holds when different graded pieces have the same scalar.

With $\sigma'$ fixed, the lifting argument in Proposition~\ref{prop:local-deformations} therefore uses the same complex~\eqref{eq:deformation-complex-weak}. The obstruction lies in $\mathbb H^2(\mathcal C^\bullet)\otimes I=0$, and the correcting coboundaries preserve the determinant and $\sigma'$. The isomorphism classes of lifts form a torsor for $\mathbb H^1(\mathcal C^\bullet)\otimes I$, where $\dim_\C\mathbb H^1(\mathcal C^\bullet)=2$. Since $\mathbb H^0(\mathcal C^\bullet)=0$, there are no infinitesimal automorphisms.

Successive lifting gives the relative formal ring $\widehat{\mathcal O}_{\mathscr B,\sigma}[[u_1,u_2]]$. Scalar automorphisms extend over every lift and fix the residue variables. Their action on this ring is therefore trivial, and the local quotient argument in Proposition~\ref{prop:algebraic-moduli} gives an isomorphism of completed local $\mathscr B$-algebras
\begin{equation}\label{eq:relative-coarse-completion}
 \widehat{\mathcal O}_{\mathcal M_{\mathscr B},x}
       \simeq
 \widehat{\mathcal O}_{\mathscr B,\sigma}[[u_1,u_2]].
\end{equation}
By the formal smoothness criterion for a morphism of finite presentation, $\rho$ is smooth of relative dimension two. Every residue parameter occurs by Proposition~\ref{prop:filtered-spectral-correspondence}, so the morphism is surjective.

The coefficient pencil in~\eqref{eq:coefficient-line} defines a closed immersion
\begin{equation*}
 \begin{aligned}
 \mathscr B\times\C_t&\longrightarrow\mathscr B\times\mathcal H_0,\\
 (\sigma,t)&\longmapsto
       (\sigma,F_{0,\sigma}-t r_n).
 \end{aligned}
\end{equation*}
Indeed, one coefficient recovers $t$ linearly, and the remaining coefficients satisfy polynomial graph equations. The coefficient morphism factors through this immersion: the equations hold at every closed point by Proposition~\ref{prop:coefficient-line}, and $\mathcal M_{\mathscr B}$ is reduced by smoothness. The factored map $(\rho,h)$ is the base change of the projective coefficient morphism, hence is projective.

For bounded $\sigma$, the parameter $\boldsymbol\gamma=\sigma/\lambda$ lies in the neighborhood of Proposition~\ref{prop:normalized-spectral-family} when $t=\lambda^n$ is sufficiently large. By Proposition~\ref{prop:filtered-spectral-correspondence}, the relative $\Pic^n$ represents the moduli sheaf of this family, including after base change to a nonreduced parameter space. Tensoring with $\pi^*\mathcal O_X(-1)$ gives the relative $\Pic^0$, a smooth proper family of elliptic curves. Descent from the finite \'etale $\lambda$-cover therefore proves smoothness over this region for $(\rho,h)$.

For sufficiently large $r_0$, this smoothness holds on $|h|>r_0-1$ over a neighborhood of $\overline V$. Then $|h|=r_0$ is a regular boundary in every fiber of $\rho$.

To trivialize the fibers while preserving the Hitchin map near infinity, we use a smooth right inverse of $\dd(\rho,h)$ on $|h|>r_0-1$. It lifts a base tangent vector $v$ to a vector mapping to $(v,0)$. A cutoff depending on $|h|$ interpolates between this lift and a smooth right inverse of $\dd\rho$ on the whole family. The resulting lift $W(v)$ satisfies
\begin{equation}\label{eq:residue-transport-lift}
 \dd\rho\,W(v)=v,\qquad
 \dd h\,W(v)=0\quad\text{on }|h|\geq r_0.
\end{equation}
We integrate $W$ along the radial path from $\sigma_0$ to a point of $V$. Along trajectories outside the boundary, $h$ is constant, and trajectories starting in $|h|\leq r_0$ remain there because the vector field is tangent to its boundary. Each trajectory is therefore confined to the inverse image under $(\rho,h)$ of a compact path times a bounded closed disc. This inverse image is compact by joint properness, so the flow exists throughout the parameter interval. Smooth dependence on the chosen point of the ball and integration along the reverse path give~\eqref{eq:residue-family-trivialization} with all the required properties.
\end{proof}

\subsubsection{Connectedness and compactification}

\begin{proposition}\label{prop:elliptic-compactification}
If semistability equals stability for
$(\alpha,\sigma)$, the surface
$\M_{\alpha,\sigma}$ is connected, irreducible and rational, and its projective surjective Hitchin morphism has connected fibers.

There is a smooth projective rational elliptic surface with section
$f:S\to\Pp^1_t$, relatively minimal over $\Pp^1_t$, and an algebraic isomorphism over $\C_t$
\begin{equation}\label{eq:whole-elliptic-complement}
 \M_{\alpha,\sigma}
       \simeq S\setminus F_\infty,\qquad
 F_\infty=f^{-1}(\infty).
\end{equation}
\end{proposition}

\begin{proof}
We begin with weights $\widehat{\alpha}$ satisfying~\eqref{eq:residue-family-weight-choice}. When all graded residues vanish, the action
$S_c(E,\Phi)=(E,c\Phi)$, $c\in\C^*$, satisfies
\begin{equation}\label{eq:zero-residue-scaling}
 h\circ S_c=c^n h.
\end{equation}
By Proposition~\ref{prop:filtered-spectral-correspondence}, a nonzero regular fiber is a connected curve of genus one. Scaling transports this whole fiber to every nonzero Hitchin value. The action also makes $h$ a submersion over $\C^*$. If $h^{-1}(\C^*)$ were disconnected, the images of a separation would separate $\C^*$, since the fibers are connected and the map is open.

Since the irreducible components of a smooth surface of finite type are disjoint, open and closed, exactly one meets $h^{-1}(\C^*)$. Any further component would be a smooth projective surface $C$ contained in $h^{-1}(0)$, by projectivity of $h$. The action of the connected group $\C^*$ preserves each component. In the trace pairing~\eqref{eq:semiflat-trace-form} for the holomorphic symplectic form of Proposition~\ref{prop:hyperkahler-structure}, scaling sends a tangent pair $(b,\varphi)$ to $(b,c\varphi)$. It follows that
\begin{equation*}
 S_c^*\Omega=c\Omega.
\end{equation*}
Its restriction to $C$ is nondegenerate. Change of variables under $S_c$ gives
\begin{equation*}
 \begin{aligned}
 0<I_C&:=\int_C\Omega\wedge\overline\Omega<\infty,\\
 I_C&=\int_C S_c^*(\Omega\wedge\overline\Omega)=|c|^2I_C,
 \end{aligned}
\end{equation*}
which is impossible for $|c|\ne1$. Thus the surface $\M_{\widehat{\alpha},\boldsymbol0}$ is connected. We apply Proposition~\ref{prop:residue-family} over a ball containing the segment from zero to any residue vector to obtain connectedness for every
$\M_{\widehat{\alpha},\sigma}$.

For the original weights we can choose rational $\widehat{\alpha}$ as in the proof of Proposition~\ref{prop:algebraic-moduli}, preserving the entire fixed-residue moduli functor and avoiding the numerical hyperplanes in~\eqref{eq:residue-family-weight-choice}. The resulting isomorphism of coarse moduli schemes proves connectedness of $\M_{\alpha,\sigma}$, which is then irreducible by smoothness. Its Hitchin morphism is projective and surjective by Propositions~\ref{prop:algebraic-moduli} and~\ref{prop:filtered-spectral-correspondence}. Since the generic fiber is geometrically connected, the Stein factor is finite and birational over the normal curve $\C_t$. It is isomorphic to $\C_t$, so all Hitchin fibers are connected.

We next use the full flag at $p=0$ to obtain an origin on the normalized spectral curve over $U_{\sigma}$. Here $m_0=1$, so~\eqref{eq:spectral-finite-modification} becomes
\begin{equation*}
 P(v)=xu,\qquad u+c_1(v,x,t)=0.
\end{equation*}
For any fixed ordered label $j$, these equations give an algebraic section $Q$ by
\begin{equation}\label{eq:full-flag-spectral-section}
 x=0,\qquad v=\sigma_{0,j},\qquad
 u=-c_1(\sigma_{0,j},0,t).
\end{equation}
The coordinate $v$ is the eigenvalue in the logarithmic frame; the corresponding global spectral coordinate at $z=0$ is $y=-\sigma_{0,j}$. For $t\in U_{\sigma}$, the function $u$ is a unit. The point in~\eqref{eq:full-flag-spectral-section} still defines a section when residue values coincide.

We set $\mathcal Y=\mathcal Y_{\sigma}|_{U_{\sigma}}$. On
$\mathcal Y\times_{U_{\sigma}}\mathcal Y$, with the first factor as parameter, the line bundle
\begin{equation}\label{eq:relative-Abel-line-bundle}
 \mathcal O\bigl(\Delta+(n-1)\operatorname{pr}_2^*Q\bigr)
\end{equation}
restricts over $y_0\in Y_t$ to $\mathcal O_{Y_t}(y_0+(n-1)Q_t)$. Here $\Delta$ is the relative diagonal, and $Q$ denotes the divisor given by the section. Its direct image with the flags of~\eqref{eq:newton-spectral-lattices} gives an algebraic map
$\mathcal Y\to h^{-1}(U_{\sigma})$.
On each fiber this is the isomorphism $Y_t\to\Pic^n(Y_t)$, $y\mapsto\mathcal O_{Y_t}(y+(n-1)Q_t)$, followed by the spectral correspondence of Proposition~\ref{prop:filtered-spectral-correspondence}. The map is proper, since its source is proper over $U_{\sigma}$ and its target is separated, and its fibers are finite. It is thus finite and birational. The target is normal, so the map is an isomorphism.

The total spectral surface defined by~\eqref{eq:coefficient-line} has function field $\C(z,y)$, since
\begin{equation}\label{eq:spectral-surface-function-field}
 t=\frac{F_{0,\sigma}(z,y)}{r_n(z)}.
\end{equation}
Indeed, $F_{0,\sigma}$ is monic in $y$ and relatively prime to $r_n(z)$; hence $F_{0,\sigma}-t r_n$ is primitive and irreducible as a polynomial in $t$. Normalization preserves its function field. This isomorphism therefore gives
$\C(\M_{\alpha,\sigma})=\C(z,y)$, proving rationality. It also gives a rational point of the generic Hitchin elliptic curve, which extends to a section over $\C_t$ by properness.

The relatively minimal regular model of this elliptic curve with origin gives a smooth projective elliptic surface $S\to\Pp^1_t$; we use the existence and uniqueness of the Kodaira--N\'eron model in~\cite[Sections~3.4--3.5, pp.~9--10]{SchuettShiodaEllipticSurfaces}. Its function field is $\C(z,y)$, so $S$ is rational.

The surface $\M_{\alpha,\sigma}$ is itself relatively minimal over $\C_t$. Its canonical bundle is holomorphically trivial by Proposition~\ref{prop:hyperkahler-structure}, so every complete smooth rational curve $C$ in it satisfies $C^2=-2$ by adjunction. Thus no fiber contains a smooth rational $(-1)$-curve. The morphism $h$ is flat because the local rings of this integral surface are torsion-free over the local discrete valuation rings of the base. The surface is thus a regular proper model over $\C_t$, and the algebraic isomorphism~\eqref{eq:whole-elliptic-complement} follows by uniqueness of the relatively minimal regular model.
\end{proof}

\subsubsection{The fiber at infinity}

\begin{proposition}\label{prop:infinity-fiber-type}
The fiber $F_\infty$ in Proposition~\ref{prop:elliptic-compactification} has Kodaira type $\mathrm{IV}$, $\mathrm{III}$ or $\mathrm{II}$ for $n=3$, $4$ or $6$, respectively. The canonical bundle of $S$ satisfies
\begin{equation}\label{eq:elliptic-canonical-bundle}
 K_S\simeq\mathcal O_S(-F_\infty).
\end{equation}
\end{proposition}

\begin{proof}
In the coordinates $w=\lambda^{-1}$ and $s=t^{-1}=w^n$, we use
$\zeta=\exp(2\pi\ii/n)$. The normalized family with residues
$\boldsymbol\gamma=w\sigma$ extends smoothly to $w=0$. Inverting the generator in~\eqref{eq:spectral-cyclic-descent} gives
\begin{equation}\label{eq:infinity-spectral-action}
 g:(w,y)\longmapsto(\zeta w,\zeta y).
\end{equation}
The section~\eqref{eq:full-flag-spectral-section} extends in the local coordinate as $v=w\sigma_{0,j}$ and is preserved by this action. Both divisors in~\eqref{eq:relative-Abel-line-bundle} are equivariant. Thus the isomorphism defined by~\eqref{eq:relative-Abel-line-bundle} intertwines this action with the descent of the generic Hitchin elliptic curve with its chosen origin.

On the central curve $y^n=r_n(z)$, the differential
$\eta_0=y\,\dd z/[z(z-1)]$ is the nowhere-zero holomorphic differential
$n\nu_1$ used in Section~\ref{subsec:semiflat-end}, and
$g^*\eta_0=\zeta\eta_0$. We choose a holomorphic frame $\widetilde\omega$ of the line bundle of fiberwise holomorphic one-forms whose value at zero is $\eta_0$. The average
\begin{equation*}
 \omega=\frac1n\sum_{k=0}^{n-1}
                      \zeta^{-k}(g^k)^*\widetilde\omega
\end{equation*}
equals $\eta_0$ on the central fiber. It is nonzero on a smaller disc and satisfies
\begin{equation}\label{eq:infinity-Hodge-character}
 g^*\omega=\zeta\omega.
\end{equation}

In short Weierstrass coordinates $p_1,p_2$ normalized by this origin and differential, the equation is
\begin{equation}\label{eq:infinity-normalized-Weierstrass}
 p_2^2=p_1^3+A(w)p_1+B(w),
 \qquad
 \omega=\frac{\dd p_1}{2p_2}.
\end{equation}
By Weierstrass uniformization with the local period lattice, the coefficients $A,B$ are holomorphic~\cite[Sections~2.1--2.3, pp.~3--4]{SchuettShiodaEllipticSurfaces}. In the local coordinate $z_0=\int_Q^\bullet\omega$ at the origin, the normalization is
\begin{equation*}
 p_1=z_0^{-2}+O(z_0^2),\qquad
 p_2=-z_0^{-3}+O(z_0).
\end{equation*}
These expansions determine the normalized functions uniquely. Since~\eqref{eq:infinity-Hodge-character} sends $z_0$ to $\zeta z_0$, they give
\begin{equation}\label{eq:infinity-Weierstrass-characters}
 \begin{gathered}
 g^*p_1=\zeta^{-2}p_1,\qquad
 g^*p_2=\zeta^{-3}p_2,\\
 A(\zeta w)=\zeta^{-4}A(w),\qquad
 B(\zeta w)=\zeta^{-6}B(w).
 \end{gathered}
\end{equation}
The central curve is smooth, so $4A(0)^3+27B(0)^2\ne0$. The characters then imply
\begin{equation*}
 \begin{cases}
 A(0)=0,\quad B(0)\ne0,&n=3,6,\\
 A(0)\ne0,\quad B(0)=0,&n=4.
 \end{cases}
\end{equation*}

On $w\ne0$, the invariant coordinates $X_1=w^2p_1$ and $Y_1=w^3p_2$ descend to
\begin{equation}\label{eq:infinity-descended-Weierstrass}
 \begin{gathered}
 Y_1^2=X_1^3+a(s)X_1+b(s),\\
 a(w^n)=w^4A(w),\qquad b(w^n)=w^6B(w).
 \end{gathered}
\end{equation}
By the congruences for the characters in~\eqref{eq:infinity-Weierstrass-characters}, $a$ and $b$ are holomorphic in $s$. Together with the central values above and~\eqref{eq:infinity-descended-Weierstrass}, they give the following orders for $a$, $b$ and $\Delta=-16(4a^3+27b^2)$:
\begin{equation}\label{eq:infinity-Weierstrass-valuations}
 \begin{array}{c|ccc|c}
 n&\operatorname{ord}_s a&\operatorname{ord}_s b
       &\operatorname{ord}_s\Delta&F_\infty\\ \hline
 3&\geq2&2&4&\mathrm{IV}\\
 4&1&\geq2&3&\mathrm{III}\\
 6&\geq1&1&2&\mathrm{II}.
 \end{array}
\end{equation}
In each case, the two summands of the discriminant have different orders, so its leading coefficient is nonzero. These equations are minimal, since either
$\operatorname{ord}_s a<4$ or $\operatorname{ord}_s b<6$.
The last column follows from the valuation table and minimality criterion in~\cite[Section~4.5, Table~1, p.~14, and Section~4.8, p.~16]{SchuettShiodaEllipticSurfaces}.

Finally, $S$ is rational, so $\chi(\mathcal O_S)=1$. By the canonical bundle formula~\cite[Theorem~6.8, p.~26]{SchuettShiodaEllipticSurfaces},
$K_S\simeq f^*\mathcal O_{\Pp^1}(-1)$, which is~\eqref{eq:elliptic-canonical-bundle}.
\end{proof}

\subsection{Completeness and curvature}
\label{subsec:metric-completeness}

\begin{proposition}\label{prop:metric-completeness}
For $(\alpha,\sigma)$ with semistability equal to stability, we equip $\M=\M_{\alpha,\sigma}$ with the Hitchin metric of Proposition~\ref{prop:hyperkahler-structure}. Then $(\M,g_H)$ is complete. For every $x\in\M$, there are positive constants $c_x,C_x,s_x$ such that
\begin{equation}\label{eq:quadratic-volume-growth}
 c_xs^2\leq\operatorname{Vol}_{g_H}(B_{g_H}(x,s))\leq C_xs^2
 \qquad(s\geq s_x).
\end{equation}
Its curvature is square integrable,
\begin{equation}\label{eq:finite-curvature-energy}
 \int_{\M}|\operatorname{Rm}(g_H)|_{g_H}^{2}
                     \,\dd\operatorname{vol}_{g_H}<\infty.
\end{equation}
\end{proposition}

\begin{proof}
The surface $\M$ is connected by Proposition~\ref{prop:elliptic-compactification}, and its Hitchin morphism is proper by Proposition~\ref{prop:algebraic-moduli}. By Propositions~\ref{prop:normalized-spectral-family} and~\ref{prop:filtered-spectral-correspondence}, the inverse image of the region where $|t|$ is sufficiently large is the end of the Hitchin fibration. Its complement is compact by properness.

On the end,~\eqref{eq:main-ALG-rate} gives uniform equivalence of $g_H$ and $g_0$, so $|\dd R|_{g_H}$ is bounded. A smooth extension of $R$ over the compact core still has bounded gradient. A metric Cauchy sequence in $\M$ has bounded $R$, hence bounded $|t|$ by~\eqref{eq:semiflat-ALG-model}. Properness puts the sequence in a compact subset, where smoothness of $g_H$ gives convergence. This proves completeness.

For $x\in\M$ we choose a sufficiently large level set of the radial coordinate. It is a compact connected torus bundle over a circle, so its points can be joined to $x$ by paths of uniformly bounded length. Outside it, radial paths in the model coordinates have $g_H$ length bounded by a constant times their radial increment. The bound on $|\dd R|_{g_H}$ gives the reverse distance inequality, up to an additive constant. Thus $d_{g_H}(x,\cdot)$ and $R$ are comparable up to additive constants on the end. The volume forms of $g_H$ and $g_0$ are uniformly comparable. In the flat model, integration over the angular and torus variables gives a constant multiple of $R\,\dd R$ because the cone angle and torus area are fixed. Comparison of radial sublevel sets with metric balls proves~\eqref{eq:quadratic-volume-growth}.

By~\eqref{eq:actual-curvature-decay}, the integral of the squared curvature norm on the end is bounded by a constant times
\begin{equation*}
 \int_{R_0}^{\infty}R^{-4-2\nu}R\,\dd R<\infty.
\end{equation*}
On the compact complement the squared curvature norm is integrable by smoothness, proving~\eqref{eq:finite-curvature-energy}.
\end{proof}

\begin{proof}[Proof of Theorem~\ref{thm:main-geometry}]
The algebraic and spectral assertions follow from Propositions~\ref{prop:algebraic-moduli}, \ref{prop:filtered-spectral-correspondence}, \ref{prop:elliptic-compactification} and~\ref{prop:infinity-fiber-type}. The metric is constructed in Proposition~\ref{prop:hyperkahler-structure}, and its completeness is proved in Proposition~\ref{prop:metric-completeness}.
\end{proof}

\subsection{ALG coordinates with decay of order two}
\label{sec:ALG-order-two}

For $(\alpha,\sigma)$ such that semistability equals stability, we write $\M=\M_{\alpha,\sigma}$. We denote the spectral coordinate map at infinity from Proposition~\ref{prop:semiflat-ALG-end} by
\begin{equation*}
 \Phi_{\mathrm{sp}}:\mathcal C(R_0)\longrightarrow\M,
 \qquad
 \mathcal C(R_0)=
 \bigl(\{|U|>R_0\}\times\C/\Lambda_\infty\bigr)
 \big/\bigl((U,V)\sim(\zeta U,\zeta^{-1}V)\bigr).
\end{equation*}
We retain the flat model~\eqref{eq:semiflat-ALG-model} and the ordered forms of Theorem~\ref{thm:main-asymptotics}. The model parameters are
\begin{equation*}
 (\beta,\tau,L)
   =\left(\frac1n,\tau_\infty,
          \frac{2\pi}{\sqrt{\operatorname{Im}\tau_\infty}}\right).
\end{equation*}

\begin{proof}[Proof of Theorem~\ref{thm:main-global-ALG}]
By Propositions~\ref{prop:elliptic-compactification} and~\ref{prop:metric-completeness}, the moduli space is connected and complete and its curvature is square integrable. We first prove nonflatness.

The elliptic compactification in Proposition~\ref{prop:elliptic-compactification} has a section. We write $\M^\circ$ for the complement of the finitely many singular Hitchin fibers and $B^\circ\subset\C$ for the set of regular values. The homotopy exact sequence of the torus bundle $\M^\circ\to B^\circ$, together with its section, shows that $\pi_1(\M^\circ)$ is generated by the fiber group $\Z^2$ and the section lifts of base loops. Inclusion $\M^\circ\hookrightarrow\M$ induces a surjection on fundamental groups, since a loop can be perturbed off the finitely many complex curves and their singular points. Each lifted base loop contracts in the section over $\C$. The conjugation relations for these lifts therefore impose $v=M_\gamma v$ in $\pi_1(\M)$, where $M_\gamma$ is the fiber monodromy along the base loop $\gamma$. Thus $\pi_1(\M)$ is a quotient of the fiber group on which every monodromy acts trivially.

The monodromy at infinity has eigenvalues $\zeta$ and $\zeta^{-1}$, by~\eqref{eq:reference-cycle-action} and~\eqref{eq:angle-deck-action}. Neither eigenvalue is one. Therefore $\pi_1(\M)$ is a quotient of the finite group
\begin{equation*}
 \Z^2/(M_n-\Id)\Z^2.
\end{equation*}
If $g_H$ were flat, its simply connected complete cover would be Euclidean $\R^4$. A finite group of Euclidean isometries fixes the average of any orbit, so a free finite deck group is trivial. This would give $\M\simeq\R^4$, so the closed K\"ahler form $\omega_I$ would be exact. A regular Hitchin fiber $F$ is a compact complex curve by Proposition~\ref{prop:filtered-spectral-correspondence}, so $\int_F\omega_I>0$. By exactness and Stokes' theorem, $\int_F\omega_I=0$, a contradiction. Thus $g_H$ is nonflat.

The metric and all three ordered forms approach the flat model at a positive polynomial rate by Theorem~\ref{thm:main-asymptotics}\textup{(i)}, and~\eqref{eq:actual-curvature-decay} gives faster than quadratic curvature decay. In the coordinates of Chen--Chen~\cite[author version, Theorems~1.1 and~4.6, pp.~1--2 and~14]{ChenChenALG}, the decay order is
\begin{equation*}
 \min_{\ell\in\Z,\ \ell<2\beta}\frac{2\beta-\ell}{\beta}.
\end{equation*}
The formulation for the ordered forms with the fixed model is recorded in~\cite[Section~5, following Lemma~5.2]{CV}. For $\beta=1/n$, the largest permitted integer is $\ell=0$, so the improved order is two. These results use the same flat tensors and model parameters as~\eqref{eq:semiflat-ALG-model}. The same estimates hold for covariant derivatives of every order on the end.

We verify that the coordinate changes extend over $\M$ and act trivially on integral homology. In the proof of Chen--Chen's Theorem~4.5, and its ALG application in Theorem~4.6, the changes are finite compositions of time-one flows of vector fields $Z$ satisfying
\begin{equation}\label{eq:ALG-coordinate-field}
 |(\nabla^{g_0})^k Z|_{g_0}=O(R^{1-\delta-k})
 \quad(k\geq0)
\end{equation}
for some $\delta>0$ at each step~\cite[author version, Section~2 and pp.~12--14]{ChenChenALG}. We push each field forward by $\Phi_{\mathrm{sp}}$, multiply it by a smooth radial cutoff that vanishes near the inner boundary of the end and equals one sufficiently far out, and extend it by zero over the core. For the growth estimate we use a complete auxiliary metric equal to $(\Phi_{\mathrm{sp}})_*g_0$ outside a compact set. The resulting field has at most linear growth in the distance from a fixed point in this metric. With this growth bound on a complete manifold, the flow exists for all finite times by the continuation theorem for ordinary differential equations.

On a trajectory in the unmodified end, the zeroth-order bound in~\eqref{eq:ALG-coordinate-field} gives
\begin{equation*}
 \left|\frac{\dd}{\dd s}R^\delta\right|\leq C\delta.
\end{equation*}
A sufficiently large initial radius keeps the whole trajectory for $0\leq s\leq1$ in the region where the cutoff is one. By uniqueness, the global extension agrees with the prescribed end flow there. Since there are only finitely many coordinate changes, one initial radius suffices for all these extensions. Their global time-one maps compose to $f$. Concatenating the isotopies and smoothly reparametrizing time at the joins gives an isotopy from the identity to $f$, so $f$ acts as the identity on integral homology. For the chosen parameter, the chart
\begin{equation}\label{eq:order-two-end-chart}
 \Psi=f\circ\Phi_{\mathrm{sp}}
\end{equation}
satisfies~\eqref{eq:main-order-two-rate} by the cited theorem. The model constants and the area are those of Proposition~\ref{prop:semiflat-ALG-end}.
\end{proof}

\subsection{The topology of the elliptic complement}
\begin{lemma}\label{lem:integral-homology-complement}
For $n=3,4,6$, respectively, the group $H_2(\M,\Z)$ is free of rank $7$, $8$, or $9$. Its intersection pairing has primitive radical $\Z[F]$, where $F$ is a smooth Hitchin fiber. The quotient by this radical is isometric to $E_6(-1)$, $E_7(-1)$, or $E_8(-1)$, respectively. In particular, it is negative definite of rank $6$, $7$, or $8$ and has determinant of absolute value $3$, $2$, or $1$.
\end{lemma}

\begin{proof}
By Propositions~\ref{prop:elliptic-compactification} and~\ref{prop:infinity-fiber-type}, $\M=S\setminus F_\infty$ for a rational elliptic surface with section, and the fiber at infinity has type $\mathrm{IV}$, $\mathrm{III}$, or $\mathrm{II}$. Writing $D_\infty$ for its support, we find that topologically $D_\infty$ is a wedge of three spheres, a wedge of two spheres, or a sphere, respectively. For types $\mathrm{III}$ and $\mathrm{IV}$ the rational components meet at one point; for type $\mathrm{II}$ the normalization of the cuspidal rational curve is a homeomorphism. In particular $H^1(D_\infty,\Z)=0$, and $H^2(D_\infty,\Z)$ has the components as its dual basis.

The support is compact and locally contractible in the closed oriented four-manifold $S$. By duality~\cite[Theorem~3.44, p.~254]{HatcherAlgebraicTopology},
\begin{equation*}
 H_i(S,\M;\Z)\simeq H^{4-i}(D_\infty,\Z).
\end{equation*}
Under the relative cap product duality, the map from $H_2(S,\Z)$ is evaluation against the boundary components. The exact sequence of the pair therefore gives an injection
\begin{equation}\label{eq:integral-homology-injection}
 H_2(\M,\Z)\hookrightarrow H_2(S,\Z)
\end{equation}
whose image is the integral orthogonal complement of those components.

We write $O$ for the section and $F$ for the fiber class in $S$. Their intersection matrix is
\begin{equation*}
 \begin{pmatrix}-1&1\\1&0\end{pmatrix}.
\end{equation*}
This unimodular plane splits integrally from $H_2(S,\Z)$, and its orthogonal complement is $E_8(-1)$~\cite[Section~8.8]{SchuettShiodaEllipticSurfaces}. We use positive definite root lattices in the notation $E_\ell$, so $E_\ell(-1)$ has the negative form. The components of $F_\infty$ disjoint from $O$ generate $R=A_2(-1)$, $A_1(-1)$, or $0$. All components have multiplicity one, so their total span is $\Z F+R$. Equation~\eqref{eq:integral-homology-injection} consequently gives
\begin{equation}\label{eq:integral-homology-splitting}
 H_2(\M,\Z)=\Z[F]\oplus R^{\perp}_{E_8(-1)}.
\end{equation}
This also proves freeness and primitivity of the radical. In the three cases the determinant of $R$ is $3$, $2$, or $1$. Each is square-free and divisible by the square of the index of the lattice in its saturation. It follows that $R$ is primitive. Primitive nondegenerate complementary sublattices in an even unimodular lattice have discriminant groups of the same order~\cite[Section~12.16, Proposition~12.14]{SchuettShiodaEllipticSurfaces}. Applying this to $R\subset E_8(-1)$ gives the ranks and determinants of the orthogonal complements.
The embeddings of $A_1$ and $A_2$ into $E_8$ are unique up to lattice automorphisms~\cite[Theorem~11.7]{SchuettShiodaEllipticSurfaces}. Their complements can therefore be computed in the standard coordinate realization of the positive definite lattice $E_8$. With $(e_i)$ the orthonormal basis of $\R^8$, we use the simple roots
\begin{equation*}
 \begin{gathered}
 \alpha_1=\tfrac12(e_1+e_8-e_2-e_3-e_4-e_5-e_6-e_7),\\
 \alpha_2=e_1+e_2,\qquad
 \alpha_j=e_{j-1}-e_{j-2}\quad(3\leq j\leq8).
 \end{gathered}
\end{equation*}
The roots $\alpha_1,\ldots,\alpha_7$ generate an $E_7$ sublattice perpendicular to $e_7+e_8$, which generates $A_1$. The roots $\alpha_1,\ldots,\alpha_6$ generate an $E_6$ sublattice perpendicular to the $A_2$ generated by $e_6+e_8$ and $e_7-e_6$. The $E_7$ and $E_6$ sublattices have determinants $2$ and $3$, respectively~\cite[Lemma~6.7]{SchuettShiodaEllipticSurfaces}, equal to the determinants of the corresponding orthogonal complements computed above. Each inclusion thus has index one. Restoring the negative sign in~\eqref{eq:integral-homology-splitting} proves the lattice identifications.

\end{proof}

\section{Hitchin metrics at boundary weights}
\label{sec:boundary}

When some parabolic weights coincide, we retain the eigenspaces of the residue on the corresponding graded piece. If the first and last labelled weights differ by one, they coincide modulo integers, so we use the periodic lattice chain. At an admissible parameter the residue eigenspaces recover each omitted step. We use this description to compare the algebraic moduli problems at nearby weights.

\subsection{Periodic lattices and residue eigenspaces}

\begin{lemma}\label{lem:joint-separation}
For $(\alpha,\sigma)\in\mathscr U$ and distinct labels $i,j$ at a puncture, $\alpha_{p,i}-\alpha_{p,j}\in\Z$ implies $\sigma_{p,i}\ne\sigma_{p,j}$.
\end{lemma}
\begin{proof}
The class with $r=0$, $k_{p,i}=1$, $k_{p,j}=-1$, all other entries zero, and $d=-(\alpha_{p,i}-\alpha_{p,j})$ belongs to $\Lambda$. This integral class has square $-2$. Its real period is zero, and equality of the two residues would make its complex period zero as well.
\end{proof}

The residue projectors reconstruct the labelled filtration in families, as follows.

\begin{proposition}\label{prop:canonical-refinement}
For an admissible parameter, forgetting the refinements at equal weights gives an equivalence between labelled presentations~\eqref{eq:periodic-chain} and the filtered Higgs bundles with the specified graded residue conjugacy classes. The equivalence holds for families over arbitrary bases, preserves Higgs morphisms and invariant saturated subobjects, and recovers the original bundle $E$ and its determinant line bundle.
\end{proposition}
\begin{proof}
We consider consecutive equal weights with labels $a,\ldots,b$ in the periodic chain. On $V=E_p^a/E_p^{b+1}$ the induced residue $R_V$ satisfies
\begin{equation*}
 \prod_{j=a}^b(R_V-\sigma_{p,j}\Id)=0\,.
\end{equation*}
Indeed, the successive factors carry each lattice term into the next. Their product annihilates the quotient because the factors commute. By Lemma~\ref{lem:joint-separation}, the corresponding residue eigenvalues are distinct, so the operators
\begin{equation}\label{eq:residue-projectors}
 \Pi_j=\prod_{\substack{a\leq h\leq b\\h\ne j}}
              \frac{R_V-\sigma_{p,h}\Id}{\sigma_{p,j}-\sigma_{p,h}}
\end{equation}
are orthogonal idempotents with sum $\Id$. Their images have rank $m_p$. For $a<q\leq b$, the missing lattice $E_p^q$ is the inverse image in $E_p^a$ of $\bigoplus_{j=q}^b\operatorname{im}\Pi_j$. This reconstructs every omitted step and the prescribed scalar actions on its successive quotients.

The same calculation applies when the equal weights occur on both sides of the periodic cut. For example, if $\alpha_{p,e_p-1}=\alpha_{p,0}+1$, the common graded piece includes $E_p^{e_p-1}/E_p^{e_p+1}$, with residue labels $\sigma_{p,e_p-1}$ and $\sigma_{p,0}$. Its projectors recover $E_p^{e_p}$ and hence $E=E_p^{e_p}(p)$. The construction for a longer sequence of equal weights uses the same formula with periodic indices.

These identities are polynomial over the base. Since the image of an idempotent on a locally free module is locally free of rank determined on geometric fibers, the construction commutes with base change to any base, including a nonreduced one. A filtered Higgs morphism commutes with every induced residue and hence with its projectors, so it preserves all the reconstructed lattices, in particular $E$. The same projectors restrict to invariant subobjects, giving their labelled multiplicities.
\end{proof}

At a puncture, a cut outside the finite set of weight values in $\R/\Z$ gives a periodic lattice whose weights $\beta_a$ lie in an interval of width less than one. In local frames this replaces selected vectors $e_a$ by $z^{s_a}e_a$, with $s_a\in\Z$, and changes their weights to $\alpha_a+s_a$. The degree decreases by $\sum_as_a$ while the weight sum increases by the same number. The cancellation also holds for every invariant saturated subobject. Thus the parabolic degree and the filtered determinant line bundle are independent of this logarithmic extension. We use such extensions for analysis and the recovered centered presentation for the numerical classes in~\eqref{eq:numeric-def}.

\subsection{Stability and the algebraic surface}

A strictly semistable object would have a stable factor of smaller rank and parabolic degree zero. We compute the Euler characteristic of its deformation complex to show that its numerical class is a root. The periods of that root would vanish, since the parabolic degree is zero and the traces of the residues sum to zero by the residue theorem. This contradicts admissibility.

\begin{lemma}\label{lem:proper-rank-rigidity}
For either $\alpha\in\Acal$ or $(\alpha,\sigma)\in\mathscr U$, every stable filtered Higgs bundle $V$ of parabolic degree zero and rank $0<r<n$ with scalar graded residues $\sigma_{p,j}\Id$ satisfies the following identity for its graded dimensions $0\leq k_{p,j}\leq m_p$, where $\sum_jk_{p,j}=r$.
\begin{equation}\label{eq:proper-rank-q}
 q(k):=\frac12\left(\sum_{p,j}k_{p,j}^2-r^2\right)=1.
\end{equation}
Its deformation complex for fixed graded residues,
\begin{equation}\label{eq:proper-rank-complex}
 \mathcal C_V^\bullet=
 [\ParEnd V\xrightarrow{[\Phi_V,\cdot]}\SParEnd V\otimes K_X(D)],
\end{equation}
has hypercohomology dimensions $(1,0,1)$. In particular, $V$ is infinitesimally rigid with these graded residues.
\end{lemma}
\begin{proof}
At an admissible boundary parameter, the residue projectors in Proposition~\ref{prop:canonical-refinement} restrict to $V$ and recover its labelled refinement; zero-dimensional grades may be omitted. The complex~\eqref{eq:proper-rank-complex} therefore computes the filtered deformation problem in both cases of the statement. A nonzero endomorphism of a stable filtered Higgs bundle is invertible, by the parabolic degree inequalities for its quotient and saturated image. Subtracting an eigenvalue shows that every endomorphism is scalar, so $\dim\mathbb H^0(\mathcal C_V^\bullet)=1$. The trace pairing gives
\begin{equation*}
 (\ParEnd V)^\vee\otimes K_X
       =\SParEnd V\otimes K_X(D),
\end{equation*}
so by Serre duality $\dim\mathbb H^2(\mathcal C_V^\bullet)=1$ as well.

The bundle $\End V$ has degree zero and rank $r^2$. At $p$, the flag condition removes $(r^2-\sum_jk_{p,j}^2)/2$ matrix entries. By Riemann--Roch on $\Pp^1$,
\begin{equation*}
 \chi(\ParEnd V)
   =r^2-\frac12\sum_p\left(r^2-\sum_jk_{p,j}^2\right)=q(k).
\end{equation*}
By trace duality the second sheaf has Euler characteristic $-q(k)$. Consequently
\begin{equation}\label{eq:factor-euler}
 \chi(\mathcal C_V^\bullet)=2q(k)
       =2-\dim\mathbb H^1(\mathcal C_V^\bullet).
\end{equation}
Since $\sum_p1/e_p=1$, we have $\sum_{p,j}k_{p,j}^2\geq r^2$ by the Cauchy--Schwarz inequality. Equality would require $k_{p,j}=r/e_p$ for every $p,j$. At the full flag, $e_0=n$, so this is impossible because $k_{0,j}$ is integral and $0<r<n$. Moreover
\begin{equation*}
 2q(k)\equiv\sum_{p,j}k_{p,j}-r^2
             =3r-r^2\equiv0\pmod2.
\end{equation*}
Thus $q(k)$ is a positive integer, while~\eqref{eq:factor-euler} gives $q(k)\leq1$. We obtain $q(k)=1$ and $\mathbb H^1(\mathcal C_V^\bullet)=0$ for the endomorphism complex, which imposes neither a trivial determinant line bundle nor a trace-free Higgs field on $V$.
\end{proof}

We now apply this calculation to a stable factor of a strictly semistable object. Its numerical class would be a root with vanishing periods. Once strict semistability is excluded, a sufficiently small perturbation of the weights into $\Acal$ preserves the algebraic moduli surface.

\begin{proposition}\label{prop:boundary-algebraic}
For $(\alpha,\sigma)\in\mathscr U$, semistability equals stability. There are arbitrarily close weights $\alpha'\in\Acal$ such that the stable and semistable moduli functors for $(\alpha,\sigma)$ and $(\alpha',\sigma)$ coincide. Their common coarse moduli space $\M_{\alpha,\sigma}$ is a nonempty smooth connected rational surface, projective over its affine Hitchin base. Its rational elliptic compactification has a section and a fiber of type $\mathrm{IV}$, $\mathrm{III}$, or $\mathrm{II}$ at infinity, respectively.
\end{proposition}
\begin{proof}
If a semistable object $E$ were not stable, it would have a proper saturated invariant subobject of parabolic degree zero. Among these subobjects, we choose $V$ of minimal positive rank. If $W\subset V$ is saturated, the exact sequence
\begin{equation*}
 0\longrightarrow V/W\longrightarrow E/W\longrightarrow E/V\longrightarrow0
\end{equation*}
has locally free outer terms, so $W$ is also saturated in $E$. Semistability of $E$ and the minimal rank of $V$ therefore make $V$ stable. It has $0<r<n$, with $0\leq k_{p,j}\leq m_p$ and $\sum_jk_{p,j}=r$ at each puncture. Its degree and residue trace satisfy
\begin{equation}\label{eq:factor-walls}
 d+\sum k_{p,j}\alpha_{p,j}=0,\qquad \sum k_{p,j}\sigma_{p,j}=0\,.
\end{equation}
The second equality is the residue theorem for $\tr\Phi_V$. By Lemma~\ref{lem:proper-rank-rigidity}, the class $[(r,d,k)]$ has square $-2$. Equations~\eqref{eq:factor-walls} contradict admissibility.

Only finitely many slope expressions enter the stability inequalities in a bounded relative neighborhood of $\alpha$ in $\overline{\Acal}$, as in the proof of Proposition~\ref{prop:joint-walls}. We can therefore find arbitrarily close rational weights $\alpha'\in\Acal$ preserving the sign of every such expression that is nonzero at $\alpha$. Since semistability equals stability at $\alpha$, the same implications between stable and semistable objects give equality of the moduli functors, fiberwise in families.

The assertions about the surface, projectivity, and elliptic compactification follow from Theorem~\ref{thm:main-geometry}. Its hypotheses hold because the common moduli problem has no strictly semistable objects. The equivalence preserves the Higgs field and hence the Hitchin map and its normalized spectral description.
\end{proof}

The equivalence in Proposition~\ref{prop:boundary-algebraic} preserves local families and their deformation functors. By Proposition~\ref{prop:local-deformations}, the deformation space is smooth of dimension two, and the determinant-preserving stabilizer is~$\mu_n$. To apply Propositions~\ref{prop:normalized-spectral-family} and~\ref{prop:filtered-spectral-correspondence}, we observe that the characteristic coefficients and the lattices~\eqref{eq:newton-spectral-lattices} involve only the scalar residues and the recovered labelled flag. The smooth Hitchin fibers are therefore $\Pic^n(Y_{t,\sigma})$, with origin given by $\pi^*\mathcal O_X(1)$.

\subsection{Harmonic metrics and deformations}
\label{sec:analysis}

The analytic boundary condition is determined jointly by the real weight and the complex residue. Within each subspace of a given weight, commutation with the residue distinguishes its different eigenspaces. We use the common centralizer of the weight and residue matrices to obtain the local estimates and construct the quotient. To compute periods we also need smooth dependence on the weights.

\subsubsection{The local normal operator}

The model separates the weight and residue contributions to the linearized equation. We use a logarithmic extension with weights $\beta_a$ in an interval of width less than one. In a holomorphic frame adapted to its filtration, the semisimple residue is diagonal on each graded piece. Repeating each label $m_p$ times gives
\begin{equation}\label{eq:local-model}
 \mathsf A=\diag(\beta_a),\qquad \mathsf B=\diag(\sigma_a),\qquad
 h_0=\diag(|z|^{2\beta_a}),\qquad \tau=-\log|z|\,.
\end{equation}
In the frame $|z|^{-\beta_a}e_a$, unitary for $h_0$, the Chern connection and Higgs field are
\begin{equation}\label{eq:local-fields}
 D_0=\dd+\ii\mathsf A\,\dd\theta,\qquad
 \Phi=(\mathsf B+q)(-\dd\tau+\ii\,\dd\theta)\,.
\end{equation}
If $a=\max\beta_a-\min\beta_a<1$, we choose $\gamma>0$ no larger than $1-a$, $1$, and all positive differences between distinct weights. When there is only one weight, we may take $\gamma=1$. Flag preservation and the chosen graded diagonalization give
\begin{equation}\label{eq:local-remainder}
 |\partial_\tau^i\partial_\theta^j q|\leq C_{ij}e^{-\gamma\tau}\,.
\end{equation}
The constants in this subsection may depend on the particular local object.

The trace-free joint centralizer is
\begin{equation}\label{eq:joint-centralizer}
 \mathfrak l_p=\{C\in\End_0\C^n:[\mathsf A,C]=[\mathsf B,C]=0\}\,.
\end{equation}
It consists of the original $\mathfrak{gl}(m_p)$ blocks subject to one trace condition. Indeed, distinct labels with the same real weight have distinct complex residues. The normal operator of the linearization of twice the harmonic metric equation is
\begin{equation}\label{eq:joint-normal}
 P=-\partial_\tau^2-(\partial_\theta+\ii\ad\mathsf A)^2
                       +4\ad(\mathsf B^*)\ad(\mathsf B)\,.
\end{equation}
For a matrix unit $E_{ab}$ and a Fourier index $k\in\Z$, the operator on the component $e^{\ii k\theta}E_{ab}$ is $-\partial_\tau^2+\mu_{abk}^2$, where
\begin{equation}\label{eq:joint-frequencies}
 \mu_{abk}^2=(k+\beta_a-\beta_b)^2+4|\sigma_a-\sigma_b|^2\,.
\end{equation}
The kernel of the angular part of $P$ is exactly the constant joint centralizer. Every positive $\mu_{abk}$ is bounded below by a number $\mu_*>0$: only finitely many $k$ need be checked, and the remaining frequencies grow like $|k|$.

\begin{proposition}\label{prop:local-boundary-decay}
Every stable filtered object under consideration has a unique adapted harmonic metric with the fixed determinant metric. In the local frame~\eqref{eq:local-fields}, its matrix $H$ has a positive Hermitian limit $H_p$ commuting with $\mathsf A,\mathsf B$. For $0<\nu<\min(\gamma,\mu_*)$ and all $i,j\geq0$,
\begin{equation}\label{eq:boundary-local-decay}
 |\partial_\tau^i\partial_\theta^j(H-H_p)|\leq C_{ij,\nu}e^{-\nu\tau}\,.
\end{equation}
\end{proposition}
\begin{proof}
The filtered Higgs bundle is regular, stable and of parabolic degree zero. On a curve the condition on the second Chern number is vacuous, and an adapted harmonic metric exists by~\cite[Theorem~2.23]{MochizukiGood}. Its uniqueness up to a positive scalar follows from~\cite[Proposition~2.22]{MochizukiGood}. The filtered determinant line bundle is trivial, so a constant rescaling gives the specified determinant metric. After decomposing by the distinct eigenvalues of the semisimple graded residue, we argue as in the proof of Lemma~\ref{lem:tame-power-comparison}. On each summand the nilpotent part is zero, so the nilpotent weight filtration is concentrated in degree zero. The diagonal metric in this estimate is $h_0$, without logarithmic factors. Hence the two metrics are mutually bounded.

We use the Fourier calculation of Lemma~\ref{lem:cylinder-linear-dirichlet} and the interior estimates of Lemma~\ref{lem:harmonic-local-compactness}(i) for the commuting normal fields~\eqref{eq:local-fields}. In the constructions of Lemmas~\ref{lem:harmonic-half-cylinder} and~\ref{lem:local-harmonic-exponential-decay}, the constant space is now $\mathfrak l_p\cap\Herm_0(\C^n)$, and the positive frequency bound is $\mu_*$. For the boundary problem and the limits of translates, we use decay of the remainder and all its derivatives at rate~$\gamma$, the lower bound~$\mu_*$ for positive normal frequencies, and the fact that every constant in~$\mathfrak l_p$ commutes with both normal fields. Thus, on the spaces $H^j_\nu$ used in Lemma~\ref{lem:harmonic-half-cylinder}, the Fourier inverse gives an isomorphism
\begin{equation}\label{eq:half-cylinder-inverse}
 (P,\operatorname{ev}_0):
 (\mathfrak l_p\cap\Herm_0(\C^n))\oplus H^4_\nu
 \longrightarrow H^2_\nu\oplus H^{7/2}(S^1;\Herm_0(\C^n)),
 \qquad 0<\nu<\mu_*.
\end{equation}
The inverse on the positive frequencies has coercivity constant $\mu^2-\nu^2$; the radial inverse on the constant summand is~\eqref{eq:harmonic-constant-green} in the proof of Lemma~\ref{lem:cylinder-linear-dirichlet}. For $\nu<\min(\gamma,\mu_*)$, the decaying coefficient terms satisfy the same weighted bounds. By the implicit function theorem, the harmonic metric has a constant term in the joint centralizer.

The translates form a relatively compact family by the mutual bound with~$h_0$ and Lemma~\ref{lem:harmonic-local-compactness}(i). As in the proof of Lemma~\ref{lem:harmonic-local-compactness}(ii), we apply the trace identity~\eqref{eq:harmonic-trace-identity} to a limit on the whole cylinder. The limit is $D_0$-parallel and commutes with~$\mathsf B$. Since the weight interval has width less than one, periodicity makes this limit a constant in the joint centralizer. The comparison with the half-cylinder solution gives its uniqueness and the decay~\eqref{eq:boundary-local-decay}. Interior estimates give the same decay after any number of differentiations in the cylindrical coordinates.
\end{proof}

\subsubsection{The deformation complex and quotient}

By the local estimates, filtered deformations admit representatives that decay on the cylindrical ends. The constant endomorphisms adjoined to the domain of the gauge operator are precisely those commuting with both normal fields, so the same construction of the weighted spaces applies when weights coincide.

We use the weighted norm~\eqref{eq:harmonic-weighted-norm} and choose
\begin{equation*}
 0<\delta<\frac12\min_p(\gamma_p,\mu_{*,p}).
\end{equation*}
The spaces $\mathcal R^j_{2\delta}$ and $\mathcal V^j_{2\delta}$ are those of~\eqref{eq:hk-higher-spaces}, with the joint centralizer~\eqref{eq:joint-centralizer} as the constant summand and the model connection~$D_0$ as reference. Thus their decaying terms have norm equivalent to $\|e^{\delta\tau}u\|_{W^{j+2,2}}$ and $\|e^{\delta\tau}v\|_{W^{j+1,2}}$, respectively. By Proposition~\ref{prop:local-boundary-decay}, the reference and harmonic metrics define equivalent norms. The operators~\eqref{eq:tangent-differentials} give the complex
\begin{equation}\label{eq:boundary-deformation-complex}
 \mathcal R^j_{2\delta}\xrightarrow{\dd_1}\mathcal V^j_{2\delta}
 \xrightarrow{\dd_2}
 W^{j,2}_{2\delta}(X^*,\End_0E\otimes K\otimes\overline K).
\end{equation}

\begin{proposition}\label{prop:boundary-hodge}
The first cohomology $\ker \dd_2/\operatorname{im}\dd_1$ of~\eqref{eq:boundary-deformation-complex} is naturally isomorphic to the first hypercohomology of the algebraic deformation complex of the filtered object. At a stable solution,
\begin{equation}\label{eq:gauge-inverse}
 L_h=\dd_1^*\dd_1:\mathcal R^0_{2\delta}\longrightarrow L^2_{2\delta}(\End_0E)
\end{equation}
is an isomorphism. Every algebraic tangent class has a unique $L^2$ harmonic representative in $\mathcal V^0_{2\delta}$.
\end{proposition}
\begin{proof}
We first remove the directions tangent to the conjugacy orbits of the graded residues from the deformation complex. For the flag determined by the distinct weights, we write $\mathcal A$ for its parabolic endomorphism sheaf and $\mathcal B_{\mathcal O}$ for the sheaf of logarithmic Higgs variations whose graded residues are tangent to the prescribed semisimple conjugacy orbits. The subsheaf $\mathcal A_{\mathsf B}\subset\mathcal A$ consists of endomorphisms whose graded values commute with the residue. The sections of $\mathcal B_0\subset\mathcal B_{\mathcal O}$ have zero graded residue. The quotient of these two deformation complexes is a sum of skyscraper complexes
\begin{equation}\label{eq:orbit-acyclic}
 \mathfrak l_{\mathrm{coarse}}/\mathfrak z_{\mathfrak l_{\mathrm{coarse}}}(\mathsf B)
       \xrightarrow{\ad\mathsf B}\operatorname{im}\ad\mathsf B\,.
\end{equation}
Here $\mathfrak l_{\mathrm{coarse}}$ is the block-diagonal Lie algebra of the coarse weight decomposition, and $\mathfrak z_{\mathfrak l_{\mathrm{coarse}}}(\mathsf B)$ is the centralizer of the residue in that algebra. The map is an isomorphism because the residue is semisimple. Thus $[\mathcal A_{\mathsf B}\to\mathcal B_0]$ computes the same hypercohomology as the deformation complex with graded residues varying in their conjugacy orbits and, by Proposition~\ref{prop:canonical-refinement}, the deformation complex for the labelled flags. We impose the weighted $L^2$ condition on this quasi-isomorphic complex.

We split the local holomorphic logarithmic matrix by contour spectral projections around its distinct residue eigenvalues. The projections preserve the filtered lattices, which split into the corresponding images. Within a summand associated with one residue eigenvalue, different original labels have distinct real weights modulo integers by Lemma~\ref{lem:joint-separation}. The local calculations of Lemma~\ref{lem:weighted-local-dolbeault} apply on its diagonal Hom complexes, with $d=2\delta$ and the multiplicity blocks as constant summands. On these summands the residue differences vanish. The positive angular frequencies are consequently determined by the weights and are at least $\mu_{*,p}$. Thus $2\delta<\mu_{*,p}$ places the Sobolev weight below the first positive frequency. Between two summands with distinct residue eigenvalues, the Sylvester operator $u\mapsto R_a(z)u-uR_b(z)$ is holomorphically invertible on a smaller disc. Its inverse preserves the filtered Hom lattices by contour functional calculus. The Higgs differential therefore contracts these two-term complexes. In the model unitary frame the inverse and its cylindrical derivatives are bounded, since the leading coefficient of the Sylvester operator is a nonzero residue difference and its remaining coefficients decay. For a closed weighted cocycle $(b,\varphi)$ in a Hom summand between distinct residue eigenspaces, applying this inverse to the logarithmic Higgs coefficient gives $u\in W^{1,2}_{2\delta}$ with $[\Phi,u]=\varphi$. The cocycle equation and the holomorphicity of the Sylvester inverse give $\bar\partial_Eu=b$. Weighted local elliptic regularity then places $u$ in $W^{2,2}_{2\delta}$, so the contraction belongs to the degree-zero domain in~\eqref{eq:hk-higher-spaces}. These are local quasi-isomorphisms between the weighted Dolbeault resolution and the algebraic two-term complex. The sheaves of weighted forms admit partitions of unity, so this resolution is fine. Its global first cohomology therefore computes the first hypercohomology group of the algebraic complex.

The normal operator of $L_h$ is one half of~\eqref{eq:joint-normal}. The map~\eqref{eq:gauge-inverse} is Fredholm by the Fourier inverses at decay exponent $\delta$, the interior elliptic parametrices and the decay of the coefficient remainders. We compute its index in two steps, since for strictly increasing parabolic weights the differences between them determine the interval of Sobolev weights in the index formula.

We first replace $L_h$ by an operator equal to its translation-invariant normal operator sufficiently far along each end. By Proposition~\ref{prop:local-boundary-decay}, the coefficient remainders decay faster than $e^{-2\delta\tau}$. Cutting off these lower-order remainders is a compact perturbation from $\mathcal R^0_{2\delta}$ to $L^2_{2\delta}$, also on the adjoined constant summands. We perturb the real weights to a nearby point of $\Acal$, with the same periodic cuts and residue matrices. The joint centralizer remains the sum of the original multiplicity blocks subject to the trace condition, since labels whose equal real weights split have distinct complex residues. After shortening the path, all positive frequencies in~\eqref{eq:joint-frequencies} remain larger than $\delta$. Interpolating the coefficients on the compact part therefore gives a norm-continuous Fredholm homotopy on the extended domain at decay exponent $\delta$. At the perturbed parameter we take the linearization at the corresponding harmonic metric and apply the same cutoff to its decaying remainders.

The smallest positive angular frequency determined by the perturbed weights is $\epsilon_*>0$. We may choose $0<\delta'\leq\delta$ with $2\delta'<\epsilon_*$ so small that the coefficient remainders have the required decay. The operator at the perturbed parameter has index zero on $\mathcal R^0_{2\delta'}\to L^2_{2\delta'}$ by Proposition~\ref{prop:ordinary-L2-projection}. Truncating its remainders preserves this index. For the truncated operator, the weights between $\delta'$ and $\delta$ cross no indicial value: its positive normal frequencies exceed $\delta$, and zero lies below $\delta'$. More explicitly, multiplication by $e^{s\tau}$ is an isomorphism from the decaying spaces to unweighted Sobolev spaces for $s\in[\delta',\delta]$, while the adjoined constant summands are unchanged. The conjugated operators form a norm-continuous Fredholm family, so their indices are equal. The truncated operator at the perturbed parameter thus has index zero at $\delta$. The homotopy of the real weights and the compact perturbation on the ends give index zero for $L_h$ at the originally chosen weight.

Cutoff integration by parts gives $\ker L_h=\ker \dd_1$. The end constants commute with both normal fields and the remainders decay, so every cutoff boundary term tends to zero. A kernel element is a bounded filtered holomorphic Higgs endomorphism. It extends to the refined filtration by Proposition~\ref{prop:canonical-refinement}, is scalar by stability, and vanishes because its trace is zero. Since the index is zero, the operator is invertible. For a closed representative $v$, the solution of $L_hu=\dd_1^*v$ gives the orthogonal harmonic projection $v-\dd_1u$. The local comparison and the same cutoff argument give uniqueness.
\end{proof}

For $j\geq2$ we use the solution space $\mathcal Z_{2\delta}$ and gauge group $\mathcal G_{2\delta}$ from Section~\ref{sec:hyperkahler-structure}, with the joint centralizer at each puncture. The limiting unitary matrices act on the labelled residue eigenspaces. By elliptic regularity the smooth quotient is independent of~$j$.

\begin{proposition}\label{prop:boundary-hyperkahler}
The quotient of the stable solution space by the gauge group is a smooth hyperk\"ahler manifold with algebraic complex structure~$I$. Its Hitchin metric and holomorphic symplectic form have the normalizations~\eqref{eq:tangent-trial-norm} and~\eqref{eq:semiflat-trace-form}. The solution space is locally a smooth principal $\mathcal G_{2\delta}/\mu_n$-bundle over the quotient. Harmonic metrics and harmonic representatives depend smoothly on local moduli parameters in the weighted spaces above. Their remainders decay exponentially, as do all their mixed derivatives in the moduli parameters and cylindrical coordinates.
\end{proposition}
\begin{proof}
We use the inverse of the gauge Laplacian from Proposition~\ref{prop:boundary-hodge} both to construct the slice and to solve the metric equation. Sobolev multiplication in the decaying spaces with constants in the joint centralizer makes the gauge action, inversion, exponentiation and moment maps smooth. The linear splitting~\eqref{eq:hk-four-operator-inverse}, with~\eqref{eq:gauge-inverse} in place of~\eqref{eq:hk-higher-inverse}, gives the local slice by the implicit function theorem.

We apply Lemma~\ref{lem:gauge-local-properness} with $d=2\delta$ to the models~\eqref{eq:local-fields}. Its estimate uses both the connection and residue intertwining equations, so it also controls directions with equal weights and distinct residue eigenvalues. The limiting stabilizer is $\mu_n$ by the algebraic deformation comparison and Proposition~\ref{prop:local-deformations}. We obtain local principal charts for the effective gauge group in
\begin{equation}\label{eq:boundary-principal-bundle}
 \mathcal Z_{2\delta}\longrightarrow\M_{\alpha,\sigma}.
\end{equation}

The quaternionic operators are given by~\eqref{eq:hk-quaternionic-operators}. The adjoint identities and moment-map calculation in the proof of Lemma~\ref{lem:local-coulomb-slice} apply on these spaces: every puncture boundary term vanishes by exponential decay, and the adjoined constants commute with both normal matrices. Consequently the harmonic tangent space is preserved by $I,J,K$, the descended forms are closed, and their pairings are~\eqref{eq:tangent-trial-norm} and~\eqref{eq:semiflat-trace-form}. The deformation comparison gives the algebraic complex structure~$I$.

To obtain smooth dependence in local moduli charts, we differentiate the harmonic metric equation and use~\eqref{eq:gauge-inverse} in higher Sobolev orders together with uniqueness. The harmonic projection depends smoothly on the parameters by the identity for the difference of two inverses. For mixed derivatives we argue as in the proof of Lemma~\ref{lem:harmonic-family-dependence}, using~\eqref{eq:gauge-inverse} and the local representatives from Proposition~\ref{prop:boundary-hodge}. The normal estimates of Proposition~\ref{prop:local-boundary-decay} then bound the remainders and their mixed derivatives in the cylindrical and moduli variables by exponentially decaying functions.
\end{proof}

\subsubsection{Smooth dependence of the Hitchin metric on the weights}

Near an admissible boundary parameter, the algebraic moduli space stays unchanged as the weights enter the interior of the alcove. Smooth dependence of the harmonic representatives on this common surface will allow the period formulas to pass to the boundary.

\begin{proposition}\label{prop:weight-continuity}
For $(\alpha,\sigma)\in\mathscr U$, we work along a sufficiently short segment $\alpha(s)$ of centered weights with $\alpha(0)=\alpha$ and $\alpha(s)\in\Acal$ for $s>0$, on which the moduli functors coincide as in Proposition~\ref{prop:boundary-algebraic}. On every compact subset of their common algebraic surface, the Hitchin metrics and the ordered hyperk\"ahler forms depend smoothly on $s$, up to $s=0$.
\end{proposition}
\begin{proof}
On a relatively compact algebraic moduli chart we choose a local family and logarithmic lattices whose cuts avoid the weights modulo $\Z$ at $s=0$. The same integer cuts work for sufficiently small $s$. We split the logarithmic Higgs matrices by contour projectors around the distinct residue eigenvalues. This splitting varies holomorphically in the coordinates on the moduli chart. Inside one summand corresponding to a residue eigenvalue all distinct labelled weights are separated at $s=0$, by Lemma~\ref{lem:joint-separation}, and remain separated for small $s$. Between different summands the Higgs coefficients vanish. In these frames the remainders~\eqref{eq:local-remainder} therefore have a common positive decay exponent.

The constant joint centralizer is unchanged, and the positive frequencies in~\eqref{eq:joint-frequencies} have a common lower bound. We choose a common decay exponent $\delta$ satisfying the bounds used in Proposition~\ref{prop:boundary-hodge}, and use the corresponding spaces with subscript $2\delta$. The model unitary frames compare these spaces along the segment. The harmonic metric equation is a smooth map in these spaces. Derivatives of $|z|^{\alpha(s)}$ introduce powers of $\log|z|$ only in terms with a strict decay margin, so derivatives of every finite order satisfy the required weighted bounds. The derivative of the equation in the metric variable at $s=0$ is the isomorphism~\eqref{eq:gauge-inverse}. By the implicit function theorem with parameters, the normalized solutions form a smooth family. They agree with the harmonic metrics on the given bundles by uniqueness in Proposition~\ref{prop:local-boundary-decay}.

For tangent representatives we use the same residue splitting. Within each summand associated with one residue eigenvalue, the same local comparison applies with a common positive separation of the real weights. On the Hom spaces between distinct residue summands the invertible logarithmic Higgs differential contracts the complex, as in the proof of Proposition~\ref{prop:boundary-hodge}. Thus one obtains smooth algebraic tangent representatives in the common weighted space. At each parameter the local comparison sends their weighted cohomology classes isomorphically to the algebraic tangent space; after shrinking the chart, these representatives continue to span a smoothly varying space of the same dimension.

The positive joint frequencies of the gauge Laplacian have a uniform lower bound. The identity for the difference of two inverses then gives smooth dependence of the harmonic projection. The positive decay margin provides an integrable dominating function for all finite parameter derivatives of the $L^2$ pairings. Hence~\eqref{eq:tangent-trial-norm} and~\eqref{eq:semiflat-trace-form}, and their quaternionic real forms, vary smoothly on the chosen chart. A finite covering by such charts gives smooth dependence on the given compact subset.
\end{proof}

\subsection{The metric at infinity}
\label{sec:asymptotics}

We extend the metric estimates for $\alpha\in\Acal$ to admissible parameters on the boundaries of the real weight alcoves. The normalized spectral family and the flat metrics on the spectral line bundles in the exterior region are independent of the real weights. The additional estimates concern the polynomial models near zero and the horizontal deformation in the frame of Section~\ref{subsec:filtered-polynomial-models}. Throughout this section, constants may depend on $(\alpha,\sigma)$; the bounds hold uniformly for points of the Picard fiber and for all arguments of the large spectral parameter. For variation of the weights over compact subsets of the moduli space, we use Proposition~\ref{prop:weight-continuity}.

We use the semiflat metric~\eqref{eq:semiflat-period-formulas} and the spectral end coordinates~\eqref{eq:semiflat-ALG-coordinates}. The normalized spectral curves, the origin $\pi^*\mathcal O_X(1)$ of their Picard fibers, and these formulas depend only on the residues.

\subsubsection{Harmonic metrics on the polynomial models}

The approximate metric is obtained from the harmonic metrics of the polynomial models near the punctures and the orthogonal direct images of the flat metrics on the spectral line bundles outside the puncture discs. We first verify that the polynomial models have the required behavior at zero for weights in the closed alcoves; we use Section~\ref{sec:local-models} for their asymptotics at the irregular singularity at infinity.

At one puncture, with its index suppressed and $e=e_p$, each of the $m_p$ local blocks has spectral equation
\begin{equation}\label{eq:polynomial-model}
 P(w)=c\xi,\qquad P(w)=\prod_{j=0}^{e-1}(w-\sigma_j),\qquad c\ne0\,.
\end{equation}
In its Newton basis the logarithmic Higgs coefficient is
\begin{equation}\label{eq:boundary-Newton}
 A_cf_j=\sigma_jf_j+f_{j+1}\quad(j<e-1),\qquad
 A_cf_{e-1}=\sigma_{e-1}f_{e-1}+c\xi f_0\,.
\end{equation}
The weight of $f_j$ is $\alpha_j$ at zero and $(e-1-2j)/(2e)$ at infinity.

\begin{lemma}\label{lem:boundary-polynomial}
The filtered Higgs bundle~\eqref{eq:boundary-Newton} has a unique harmonic metric $h_c$ inducing the prescribed metric on the determinant line bundle. Near zero it satisfies the convergence and derivative estimates of Proposition~\ref{prop:local-boundary-decay}. At infinity it satisfies the conclusions of Proposition~\ref{prop:polynomial-model-tail}. In the ramified coordinate $r=\xi^{-1/e}$, its ratio $G(r)$ to the flat metric on the eigenline bundles extends smoothly to $r=0$ with $G(0)=\Id$. Its off-diagonal entries vanish to every order at zero, and every Taylor polynomial of $\log G_{ii}$ is harmonic.
\end{lemma}
\begin{proof}
The closed weight inequalities imply
\begin{equation*}
 \lfloor a+\alpha_j\rfloor\leq\lfloor a+\alpha_{j+1}\rfloor,
 \qquad
 \lfloor a+\alpha_{e-1}\rfloor-\lfloor a+\alpha_0\rfloor\leq1\,.
\end{equation*}
They are exactly the conditions that the subdiagonal entries and the upper-right entry $c\xi$ in~\eqref{eq:boundary-Newton} preserve every filtered lattice at zero. The determinant weight sum is zero at both ends. The polynomial $P(w)-c\xi$ is irreducible over $\C(\xi)$: the quotient of $\C[\xi,w]$ by the ideal it generates is $\C[w]$, and it is primitive as a polynomial over $\C[\xi]$. Hence the generic Higgs module has no proper invariant subspace, so the filtered Higgs bundle is stable. Its irregular filtered structure at infinity is the same good structure as in Section~\ref{subsec:filtered-polynomial-models}. Existence of $h_c$ follows from the filtered Kobayashi--Hitchin correspondence~\cite[Theorem~2.23]{MochizukiGood}.

By Lemma~\ref{lem:joint-separation}, the graded residues at zero are semisimple even when the common weight occurs on both sides of the periodic cut. We can therefore apply Proposition~\ref{prop:local-boundary-decay} at zero.

The good filtered structure at infinity is unchanged, so the decoupling estimates in the proof of Proposition~\ref{prop:polynomial-model-tail} apply. The Laplacians of the diagonal logarithms extend smoothly to $r=0$ and vanish there to every order. The scalar Poisson argument in that proof then gives smooth extensions of the diagonal logarithms, with every Taylor polynomial harmonic. The off-diagonal entries vanish to every order by the same decoupling estimates. The deck group acts transitively on the eigenline bundles, so their limiting diagonal constants are equal. Their product is one by determinant normalization, so their common positive value is one. The chain rule under $\zeta=r^{-1}$ gives the derivative estimates in~\eqref{eq:polynomial-model-tail}. For each global parameter there are only finitely many model blocks, so their constants can be chosen simultaneously.
\end{proof}

\subsubsection{Matching the graded residues}

After rescaling by $t=\lambda^n$ and $\xi=\lambda^ex$, with $T=|\lambda|$, we compare the spectral field near a puncture with its polynomial model. In the local spectral coordinates of Proposition~\ref{prop:local-spectral-matching}, the logarithmic Higgs field in~\eqref{eq:exact-inner-matrix} has the form
\begin{equation}\label{eq:boundary-exact-inner}
 \Psi_\varepsilon=
 \left(A_c+\xi\sum_jd_j(\xi,\varepsilon)E_{j,e-1}\right)\frac{\dd\xi}{\xi},
 \qquad \varepsilon=\lambda^{-e},\qquad d_j=O(\varepsilon)\,.
\end{equation}
The coefficients are holomorphic on a fixed inner disc and satisfy the same bounds after any fixed finite number of derivatives. The constants are uniform over the finite collection of charts on the Picard fibers and on the circle of arguments of the scaling parameter used to cover the end. A change of periodic cut can turn some regular terms into residue terms. We first remove their graded residue contributions.

\begin{lemma}\label{lem:graded-alignment}
In a fixed holomorphic frame adapted to a periodic lattice of weight width less than one, there is a filtration-preserving gauge $G_\varepsilon=\Id+O(\varepsilon)$ of determinant one such that the graded residues of $G_\varepsilon^{-1}\Psi_\varepsilon G_\varepsilon$ and $\Psi_0$ agree. The gauge may be chosen independent of $\xi$ in this frame and holomorphic in $\varepsilon$. Its norm and inverse norm in the power metric $h_0$ are bounded, and the difference from $\Id$ is $O(\varepsilon)$ in that metric. On a smaller parameter disc, the $\varepsilon$-derivatives of every fixed order of $G_\varepsilon$ and its inverse are bounded in $h_0$.
\end{lemma}
\begin{proof}
On each coarse graded piece $V$, the exact residue $R_\varepsilon$ and the model residue $R_0$ are triangular in the recovered labelled refinement. Their diagonal blocks are the prescribed $\sigma_j\Id$, with distinct eigenvalues for different labels in this grade. They are therefore semisimple, have the same eigenvalues and multiplicities, and differ by $O(\varepsilon)$. The exact residue is holomorphic in $\varepsilon$ by~\eqref{eq:actual-inner-polynomial}. Its spectral projectors $\Pi_{j,\varepsilon}$ are polynomials in $R_\varepsilon$ with fixed denominators, as in~\eqref{eq:residue-projectors}. Thus
\begin{equation}\label{eq:alignment-matrix}
 C_\varepsilon=\sum_j\Pi_{j,\varepsilon}\Pi_{j,0}
\end{equation}
is holomorphic in $\varepsilon$, satisfies $C_0=\Id$, and obeys $R_\varepsilon C_\varepsilon=C_\varepsilon R_0$ by the projector identities. For small $\varepsilon$ it is invertible and equals $\Id+O(\varepsilon)$. On $V$ we divide $C_\varepsilon$ by the branch of $(\det C_\varepsilon)^{1/\dim V}$ taking value one at $\varepsilon=0$; this preserves the intertwining identity. The block diagonal lift of these normalized maps, taken constant in $\xi$ in the chosen periodic frame, defines $G_\varepsilon$ and has determinant one. This lift and its inverse are holomorphic in $\varepsilon$, so their parameter derivatives of every fixed order are bounded on a smaller disc. Each matrix acts inside an equal-weight block, where $h_0$ is a scalar multiple of a fixed metric, so the matrix estimates also hold in $h_0$. The construction depends only on the induced residues on the coarse graded pieces.
\end{proof}

After the graded residues agree, every remaining coefficient decays at a positive rate in the unitary frame of the model. The resulting error in Hitchin's equation is integrable, even when real weights coincide.

\begin{proposition}\label{prop:boundary-inner-defect}
There is a number $0<\gamma\leq1$, depending on the weights and periodic cuts, for which the error in Hitchin's equation for the transported model metric is $O(|\varepsilon|)$ in Euclidean $L^{p_0}$ whenever
\begin{equation}\label{eq:boundary-p-range}
 1<p_0<\min\left\{\frac{12}{11},\frac{2}{2-\gamma}\right\}\,.
\end{equation}
The resulting global approximate metric satisfies
\begin{equation}\label{eq:boundary-global-defect}
 \|F_{h^{\app}}+[\Phi,\Phi^{\dagger_{h^{\app}}}]\|_{L^{p_0}(X)}
       \leq CT^{-\vartheta},\qquad
 \vartheta=\frac12-6\left(1-\frac1{p_0}\right)>0\,.
\end{equation}
It induces the unit metric on the trivialized determinant line bundle and the given filtration on the bundle.
\end{proposition}
\begin{proof}
After the change of frame in Lemma~\ref{lem:graded-alignment}, the difference of the logarithmic coefficients has zero coarse graded residue. We choose $\gamma$ as in~\eqref{eq:local-remainder}, simultaneously at the three punctures, and reduce it to at most one. In the frame unitary for $h_0$, a residue entry strictly lowering the coarse filtration acquires a positive power at least $\gamma$; a regular entry acquires a power at least $1-a\geq\gamma$, where $a<1$ is the weight width. The mutual bound for $h_c$ then gives
\begin{equation}\label{eq:aligned-inner-field-bound}
 |\widetilde\Psi_\varepsilon-\Psi_0|_{h_c,|\dd\xi|^2}
       \leq C|\varepsilon||\xi|^{-1+\gamma},\qquad
 |\Psi_0|_{h_c,|\dd\xi|^2}\leq C|\xi|^{-1}\,.
\end{equation}
Since $h_c$ solves the model equation, the two cross commutators and the quadratic difference term satisfy
\begin{equation}\label{eq:aligned-inner-defect-bound}
 |F_{h_c}+[\widetilde\Psi_\varepsilon,
                  \widetilde\Psi_\varepsilon^{\dagger_{h_c}}]|
       \leq C|\varepsilon||\xi|^{-2+\gamma}\,.
\end{equation}
The condition $p_0<2/(2-\gamma)$ is exactly the integrability condition, because
\begin{equation*}
 \int_0^1r^{1-p_0(2-\gamma)}\,\dd r
        =\frac1{2-p_0(2-\gamma)}<\infty\,.
\end{equation*}
Near zero we transport $h_c$ back by $G_\varepsilon$. On a fixed inner annulus we interpolate to the identity by applying a cutoff to the logarithm of $G_\varepsilon$. All metric, curvature, and Higgs errors on that annulus are $O(|\varepsilon|)$, since the annulus stays away from zero. Outside this annulus we use the holomorphic isomorphisms between the eigenline bundles constructed in Lemma~\ref{lem:transported-inner-metric}. The off-diagonal model expansion in Lemma~\ref{lem:boundary-polynomial} gives $O(|\varepsilon|\rho^{-N})$ for every $N$ there, with $\rho=|\xi|^{1/e}$. The same construction as for $\alpha\in\Acal$ thus gives an error in Hitchin's equation whose Euclidean $L^{p_0}$ norm is at most $C|\varepsilon|$ on the growing inner region.

On the exterior region we continue to use the orthogonal direct image of the flat metric on the spectral line bundle. At $|x|\asymp T^{-e/2}$, or $\rho\asymp T^{1/2}$, its logarithmic discrepancy with the inner metric has diagonal part $O(T^{-1/2})$ and off-diagonal part $O(T^{-N})$ for every $N$, with all derivatives in $z=T^{e/2}x$. For this comparison on the eigenline bundles, the proof of Lemma~\ref{lem:inner-outer-metric-matching} uses only the expansion of the model metric at infinity from Lemma~\ref{lem:boundary-polynomial}. We use the same normalization of the relative positive constants on the blocks and of the scalar fixing the determinant. The new local gauge has determinant one and already equals the identity on this annulus, so neither normalization changes.

For the interpolation~\eqref{eq:approximate-metric-cutoff}, the diagonal part of the curvature on the annulus contributes $O(T^{e(1-1/p_0)-1/2})$ to the $L^{p_0}$ norm by Lemma~\ref{lem:inner-outer-metric-matching}. Under $x=\lambda^{-e}\xi$, the inner two-form estimate acquires the factor $T^{2e-2e/p_0}$, so its contribution is
\begin{equation}\label{eq:boundary-inner-rescaling}
 CT^{-e}T^{2e-2e/p_0}=CT^{-e(2/p_0-1)}.
\end{equation}
For $2\leq e\leq6$ and~\eqref{eq:boundary-p-range}, both terms are bounded by $CT^{-\vartheta}$ with $\vartheta$ as in~\eqref{eq:boundary-global-defect}. On the exterior region the orthogonal direct image of the flat metric on the spectral line bundle solves Hitchin's equation, and summation over the punctures gives that estimate. The block and determinant normalizations are those of~\eqref{eq:inner-block-constants} and~\eqref{eq:inner-determinant-factor}.
\end{proof}

To compare the approximate and harmonic metrics, we combine the bound on the error in Hitchin's equation with a coercivity estimate. Connectedness of the spectral cover makes the coercivity constant uniform along the end.

\begin{proposition}\label{prop:boundary-harmonic-comparison}
For the harmonic metric of an object on the end,
\begin{equation}\label{eq:boundary-harmonic-comparison}
 h=h^{\app}e^s,\qquad \|s\|_{L^\infty}\leq CT^{-\vartheta/2}\,.
\end{equation}
For fixed weights and graded residues, the constant is uniform on the Picard fiber and independent of the argument of $\lambda$.
\end{proposition}
\begin{proof}
We check the two analytic hypotheses of Proposition~\ref{prop:small-defect-metric-comparison}. First, on each fixed compact subset of $X\setminus D$, the approximate metric is the exterior flat metric for large $T$. Its connection and $T^{-1}\Phi$ converge smoothly along subsequences for which the Picard characters and the arguments of the scaling parameters converge. The limiting Higgs field is normal with distinct eigenvalues there; its normalized spectral covering is the connected curve $Y_{1,0}$ with its puncture preimages removed. By Lemma~\ref{lem:energy-bound-under-convergence}, there is a uniform positive coercivity constant for the trace-free metric variation. The coercivity argument uses the Kato inequality, compactness on the base, and the connected limiting covering, and imposes no strict inequality on the weights at the removed points.

Second, for each individual object the harmonic metric and the transported inner metrics are mutually bounded at the punctures by Proposition~\ref{prop:local-boundary-decay}, since they induce the same filtration. On the compact complement they are smooth positive metrics. Thus their global ratio is bounded for that object. We can now apply Proposition~\ref{prop:small-defect-metric-comparison}, whose estimate is independent of the finite initial bound. Combining~\eqref{eq:boundary-global-defect} with uniform coercivity gives~\eqref{eq:boundary-harmonic-comparison}. Finitely many charts with compact closures cover the Picard fibers and the circle of arguments of the scaling parameter, so one constant suffices.
\end{proof}

\subsubsection{Horizontal tangent representatives}

The comparison of Hermitian metrics controls the norms of the chosen representatives of the tangent classes. To estimate the Hitchin metric itself, the chosen horizontal representative must also belong to the weighted complex. Differentiating the frame in which the graded residues are constant removes the otherwise nondecaying graded residue variation.

\begin{proposition}\label{prop:boundary-tangent-comparison}
On the end, the Hitchin metric and the spectral semiflat metric satisfy
\begin{equation}\label{eq:boundary-C0-metric}
 |g_H-g_{\sfm}|_{g_{\sfm}}
       \leq C(T^{-1}+T^{-\vartheta/2})\,.
\end{equation}
They have the same complex structure and holomorphic symplectic form.
\end{proposition}
\begin{proof}
The vertical representative is supported away from the puncture and ramification discs. Its construction and spectral Hodge norm are unchanged from Proposition~\ref{prop:tangent-trial-comparison}. For a horizontal derivative $\partial_\lambda$, we first construct an algebraic cocycle near the punctures to which we can apply Proposition~\ref{prop:boundary-hodge}.

The exact logarithmic coefficient in~\eqref{eq:boundary-exact-inner} is $A=A(\xi,\varepsilon)$. We use the gauge $G_\varepsilon$ from Lemma~\ref{lem:graded-alignment} in the chosen periodic frame. A matrix $B$, constant in $\xi$ and independent of $\lambda$, acts within the coarse weight spaces and diagonalizes the model graded residue. The full change from the Newton frame and the conjugated coefficient are
\begin{equation}\label{eq:horizontal-moving-frame}
 U_\xi=\diag(\xi^{s_a}),\qquad S_\lambda=U_\xi G_\varepsilon B,
 \qquad A'=S_\lambda^{-1}AS_\lambda\,.
\end{equation}
The integers $s_a$ implement the change of extension, and $G_\varepsilon$ and $B$ act in the periodic frame. The resulting local frames vary holomorphically in the family. With a dot denoting differentiation at fixed $x$, we set $\Gamma=\dot S_\lambda S_\lambda^{-1}$. Direct differentiation gives
\begin{equation}\label{eq:horizontal-corrected-derivative}
 \dot A'=S_\lambda^{-1}(\dot A+[A,\Gamma])S_\lambda,
 \qquad \dot A=\frac e\lambda
           (\xi\partial_\xi-\varepsilon\partial_\varepsilon)A\,.
\end{equation}
The coefficient $A'$ is holomorphic in $(\xi,\varepsilon)$ in the periodic extension, with bounded derivatives on smaller fixed discs. Its graded residue is the fixed diagonal model residue. Applying $(e/\lambda)(\xi\partial_\xi-\varepsilon\partial_\varepsilon)$ therefore gives $\dot A'$ with zero graded residue and $O(T^{-1})$ bounds for every fixed number of $\xi$-derivatives. On that disc, the horizontal Higgs variation is $\varphi_H=\dot A'\,\dd\xi/\xi$ in the new frame. In the following norm, $h_c$ is expressed in the periodic frame of the model. Both $G_\varepsilon B$ and its inverse are bounded in $h_0$, and $h_c$ is mutually bounded with $h_0$, so the estimate in the frame $S_\lambda$ also holds in the model frame. The weight calculation in~\eqref{eq:aligned-inner-field-bound} gives
\begin{equation}\label{eq:boundary-horizontal-core}
 |\varphi_H|_{h_c,|\dd\xi|^2}\leq CT^{-1}|\xi|^{-1+\gamma},
 \qquad
 \|\varphi_H\|_{L^2(|\xi|<1)}^2
       \leq CT^{-2}\int_0^1r^{2\gamma-1}\,\dd r
       =\frac{C}{2\gamma}T^{-2}\,.
\end{equation}

We use a cutoff $\chi$ supported in a fixed larger disc and equal to one on the original disc, and hold it fixed when differentiating at the given parameter. With $\Gamma^0=\Gamma-\tr(\Gamma)\Id/n$, the corrected cocycle in the original local description is
\begin{equation}\label{eq:boundary-horizontal-cocycle}
 b_H=\bar\partial(\chi\Gamma^0),\qquad
 \varphi_H=\dot A\,\frac{\dd x}{x}+[\Phi,\chi\Gamma^0]\,.
\end{equation}
Here $\tr\Gamma=(e/\lambda)\sum_a s_a$, since $\det G_\varepsilon=1$ and $B$ is independent of $\lambda$. This scalar is independent of $x$, so its subtraction removes the scalar variation of the determinant frame and leaves the commutator unchanged. Thus~\eqref{eq:boundary-horizontal-cocycle} is the trace-free \v{C}ech--Dolbeault cocycle for the same algebraic deformation. Near zero it is the variation in~\eqref{eq:horizontal-corrected-derivative}, whose residue vanishes on every graded piece. Outside the larger disc it is the original representative. With $\tau_\xi=-\log|\xi|$, the cylindrical norm in~\eqref{eq:boundary-horizontal-core} is $O(T^{-1}e^{-\gamma\tau_\xi})$, with the same decay after differentiation. The weighted square integral is bounded by a constant times
\begin{equation*}
 T^{-2}\int^\infty e^{-2(\gamma-\delta)\tau_\xi}\,\dd\tau_\xi<\infty.
\end{equation*}
The cylindrical coordinate on the original punctured disc is $\tau_x=-\log|x|$, and $\tau_\xi=\tau_x-e\log T$. For each parameter, the two weighted norms therefore differ by a finite scalar factor. The cutoff terms are supported away from the puncture, so the corrected cocycle belongs to the weighted deformation complex. The uniform estimate for the corrected representative uses its unweighted $L^2$ norm from~\eqref{eq:boundary-horizontal-core}.

On the fixed inner annulus,
\begin{equation*}
 \begin{aligned}
 \dot U_\xi U_\xi^{-1}&=\frac e\lambda\diag(s_a),\\
 \dot{(G_\varepsilon B)}(G_\varepsilon B)^{-1}
       &=\dot G_\varepsilon G_\varepsilon^{-1}\\
       &=-\frac{e\varepsilon}{\lambda}
          (\partial_\varepsilon G_\varepsilon)G_\varepsilon^{-1}
         =O(|\varepsilon|/T)\,.
 \end{aligned}
\end{equation*}
Here $U_\xi$ and its inverse have bounded derivatives. The additional $(0,1)$ and Higgs terms therefore have squared norm $O(T^{-2})$. Beyond this annulus the calculation in the proof of Proposition~\ref{prop:tangent-trial-comparison} uses only the normalized model at infinity and yields
\begin{equation*}
 |\dot A|_{h_c}\leq CT^{-1}\rho(1+T^{-e}\rho^e),
 \qquad 1\leq\rho\leq CT^{1/2}\,.
\end{equation*}
After integration over the inner region, the squared norm is bounded by
\begin{equation}\label{eq:boundary-horizontal-energy}
 CT^{-2}\left(1+\int_1^{CT^{1/2}}\rho^2\frac{\dd\rho}{\rho}\right)
       \leq CT^{-1}\,.
\end{equation}
On the outer unramified region the representative has precisely its spectral Hodge norm. Every added $(0,1)$ term lies where the vertical representative vanishes. The remaining horizontal Higgs part and vertical $(0,1)$ part are orthogonal in~\eqref{eq:tangent-trial-norm}, so the two representatives are pointwise orthogonal for every Hermitian metric.

The trace pairing under the deformation comparison is the algebraic Serre pairing. In the subcomplex $[\mathcal A_{\mathsf B}\to\mathcal B_0]$ used in Proposition~\ref{prop:boundary-hodge}, weighted decay makes the puncture boundary terms vanish. The quotient~\eqref{eq:orbit-acyclic} is acyclic, so this pairing agrees with that of the recovered labelled complex. Hence~\eqref{eq:semiflat-trace-form} equals the spectral holomorphic symplectic form $\Omega$ in~\eqref{eq:semiflat-period-formulas}.

Orthogonal $L^2$ projection decreases the norm. The bounds on the corrected representatives, the comparison~\eqref{eq:boundary-harmonic-comparison}, and the symplectic linear-algebra argument of Lemma~\ref{lem:metric-comparison-from-trials} now give~\eqref{eq:boundary-C0-metric}.
\end{proof}

\subsubsection{ALG estimates and completeness}

The $C^0$ comparison and the common complex symplectic form give derivative estimates on expanding covers of the end. These estimates and properness of the Hitchin map imply completeness. In the coordinates of Chen--Chen, the resulting ALG metric has decay of order two.

\begin{proof}[Proof of Theorem~\ref{thm:boundary-ALG}]
The flat model and the decay of order two of the semiflat tensors in Proposition~\ref{prop:semiflat-ALG-end} are independent of~$\alpha$. The radial coordinate~$R$ of the flat model is a positive constant multiple of~$T$. Given $0<\nu<1/4$, we choose $p_0$ sufficiently close to one in~\eqref{eq:boundary-p-range} that $\vartheta/2>\nu$. The required $C^0$ comparison follows from Proposition~\ref{prop:boundary-tangent-comparison}.

Both metrics are Ricci-flat and K\"ahler for the same complex structure, and their holomorphic symplectic forms agree. We apply Lemma~\ref{lem:ricci-flat-interior-comparison} on the expanding holomorphic covers~\eqref{eq:metric-growing-cover}. The spectral bounds~\eqref{eq:metric-scaled-semiflat-bounds} depend only on the residues, so they remain valid here. After rescaling, the $C^0$ bound gives uniform ellipticity and bounds the metric difference by $O(R^{-\nu})$. On the smaller covers, interior estimates bound derivatives of each order. Rescaling back contributes $R^{-1}$ per derivative, exactly as in the proof of Proposition~\ref{prop:metric-derivative-comparison}, and gives~\eqref{eq:boundary-ALG-all-orders}. The estimates for $\omega_2,\omega_3$ follow from the common holomorphic symplectic form, and the estimate for $\omega_1$ follows from the common complex structure. A finite cover of the Picard fibers and the circle of arguments by charts with compact closures gives uniform constants and a radius independent of the derivative order.

The Hitchin map is proper and the elliptic compactification has a section by Proposition~\ref{prop:boundary-algebraic}. Together with~\eqref{eq:boundary-ALG-all-orders}, these are the hypotheses used in Proposition~\ref{prop:metric-completeness} and Section~\ref{sec:ALG-order-two}. As in those proofs, the metric is complete and nonflat and its curvature is square integrable. The change to coordinates with decay of order two extends to a global diffeomorphism isotopic to the identity. In particular, the correction acts trivially on integral homology.
\end{proof}

\section{Periods and modularity}
\label{sec:periods}
\label{sec:global}

We use the integral period calculation of Dimakis--Fredrickson--Yae for $\xi\in\mathscr U^\circ$ and extend it to all of $\mathscr U$ by smooth variation on compact subsets of the moduli space. Compactly supported changes of marking realize every normalized period and preserve the ALG coordinates needed for the Torelli theorem.

\subsection{The integral lattice and period formula}

We express the integral basis of~\cite[Sections~6.1 and~6.3]{DFY} in the numerical lattice~\eqref{eq:numeric-def}. For each puncture, $e_{p,j}$ denotes the unit vector in the $k_{p,j}$ coordinate. The affine diagram has central vertex $*$ and vertices $(p,h)$, where $0\leq h\leq e_p-2$. We set $a=(0,0)$ and define
\begin{equation}\label{eq:numerical-affine-basis}
 \begin{aligned}
 b_*&=\left[\left(1,0,\sum_p e_{p,0}\right)\right],\\
 b_{p,h}&=\left[(0,-\delta_{(p,h),a},
             e_{p,e_p-1-h}-e_{p,e_p-2-h})\right].
 \end{aligned}
\end{equation}
The multiplicities of these vertices are $\rho_*=n$ and $\rho_{p,h}=m_p(h+1)$; in particular, $\rho_a=1$.

\begin{lemma}\label{lem:numerical-affine-basis}
The classes in~\eqref{eq:numerical-affine-basis} form an integral basis of $\Lambda$. Their pairing is the negative affine $E_\ell$ Cartan matrix, and
\begin{equation}\label{eq:numerical-affine-fiber}
 F=\sum_v\rho_v b_v\,.
\end{equation}
The radical is $\Z F$, and $\Lambda/\Z F$ is isometric to $E_\ell(-1)$.
\end{lemma}
\begin{proof}
Each of these classes has square $-2$. Consecutive classes on one arm have pairing $1$, and $b_*$ has pairing $1$ with $b_{p,e_p-2}$ at each puncture. All other pairings vanish. Their intersection graphs are the affine diagrams in~\cite[Figure~4]{DFY}. Summing the $k$ coordinates with multiplicities $\rho_v$ gives $m_p$ in every coordinate at $p$, while the rank and degree are $n$ and $-1$. Thus the sum is $\delta_0+F$, which proves~\eqref{eq:numerical-affine-fiber} in the quotient.

To prove integral generation, we subtract a multiple of $b_*$ to set the rank equal to zero. Each remaining array of $k$ coordinates has sum zero and is therefore an integral linear combination of differences along the corresponding arm. The class with only its degree coordinate remaining is a multiple of $F$, which belongs to the span by~\eqref{eq:numerical-affine-fiber}. The vector $\delta_0$ is primitive because $m_0=1$, so $\Lambda$ is free. Its rank is $\ell+1$, equal to the number of generators, so they form a basis. Since $\rho_a=1$, the classes $F$ and $b_v$, $v\ne a$, are also an integral basis. Their pairing is the orthogonal sum of the zero form on $\Z F$ and the negative finite $E_\ell$ Cartan matrix.
\end{proof}

In the descending flag indices of~\cite{DFY}, the adjacent differences are $\beta_{p,h}=\alpha_{p,e_p-1-h}-\alpha_{p,e_p-2-h}$ and $\tau_{p,h}=\sigma_{p,e_p-1-h}-\sigma_{p,e_p-2-h}$, with central values $\beta_*=\sum_p\alpha_{p,0}$ and $\tau_*=\sum_p\sigma_{p,0}$. In these conventions the real K\"ahler form of Dimakis--Fredrickson--Yae equals $\omega_I$ and their complex symplectic form equals $-\Omega$, as follows from~\cite[equations~(36) and~(66)]{DFY} and~\eqref{eq:semiflat-trace-form}.

We choose a reference parameter $\xi_*=(\widehat\alpha,0)\in\mathscr U^\circ$ with strictly increasing weights and vanishing graded residues for which the compact fiber of Kodaira type $\mathrm{IV}^*$, $\mathrm{III}^*$ or $\mathrm{II}^*$ is as in~\cite[Section~6.1]{DFY}, and put $\M_*:=\M_{\widehat\alpha,0}$. For the integral basis $(\eta_v)$ on this surface from~\cite[Section~6.3]{DFY}, the assignment $b_v\mapsto\eta_v$ gives an integral isometry
\begin{equation}\label{eq:reference-integral-marking}
 \iota_*:\Lambda\longrightarrow H_2(\M_*,\Z)
\end{equation}
carrying $F$ to the class of a smooth Hitchin fiber. We use this marking throughout the remainder of the proof.

\begin{proposition}\label{prop:strict-integral-periods}
For every $\xi=(\alpha,\sigma)\in\mathscr U^\circ$, parallel transport from $\xi_*$ gives an integral isometry $\iota_\xi:\Lambda\to H_2(\M_\xi,\Z)$, and
\begin{equation}\label{eq:strict-numerical-periods}
 \int_{\iota_\xi(C)}\omega_I=P_I(C),\qquad
 \int_{\iota_\xi(C)}\Omega=P_{\C}(C),\qquad C\in\Lambda\,.
\end{equation}
\end{proposition}
\begin{proof}
The root condition places $\xi$ in the generic locus of~\cite[Definition~5.3 and Remark~6.12]{DFY}. Its coprimality hypothesis holds because $m_0=1$. We use parallel transport in integral homology as in their Corollary~5.4. For the basis $(\eta_v)$ of their Sections~6.1 and~6.3, the pairing and fiber relation agree with those in Lemma~\ref{lem:numerical-affine-basis}. With our convention for the forms, the formula in Theorem~6.9 of~\cite{DFY} becomes
\begin{equation*}
 \int_{\eta_v}\omega_I=A_F\delta_{va}-A_F\beta_v,
 \qquad
 \int_{\eta_v}\Omega=2A_F\tau_v\,.
\end{equation*}
Substitution of~\eqref{eq:numerical-affine-basis} into~\eqref{eq:formal-periods} gives exactly these values. Equality on the integral basis proves~\eqref{eq:strict-numerical-periods}.
\end{proof}

\subsection{Admissible parameters and boundary periods}

The affine formula has an inverse on the centered parameter space before the order and width inequalities are imposed. We use it to recover the parameters of an arbitrary normalized period triple.

\begin{lemma}\label{lem:formal-period-inverse}
Formula~\eqref{eq:formal-periods} is an affine isomorphism from the centered parameter space onto the affine space $P(F)=(A_F,0,0)$.
\end{lemma}
\begin{proof}
For $1\leq l<e_p$, we set
\begin{equation*}
 \Gamma_{p,l}=[(0,0,e_{p,l-1}-e_{p,l})],\qquad
 u_{p,l}=A_F^{-1}P_I(\Gamma_{p,l}),\qquad
 v_{p,l}=(2A_F)^{-1}P_{\C}(\Gamma_{p,l})\,.
\end{equation*}
The differences $\alpha_{p,l}-\alpha_{p,l-1}=u_{p,l}$ and $\sigma_{p,l-1}-\sigma_{p,l}=v_{p,l}$ determine the individual coordinates uniquely once the centering conditions are imposed. Explicitly,
\begin{equation}\label{eq:formal-period-inverse}
 \begin{aligned}
 \alpha_{p,0}&=-\frac1{e_p}\sum_{l=1}^{e_p-1}(e_p-l)u_{p,l},&
 \alpha_{p,j}&=\alpha_{p,0}+\sum_{l=1}^j u_{p,l},\\
 \sigma_{p,0}&=\frac1{e_p}\sum_{l=1}^{e_p-1}(e_p-l)v_{p,l},&
 \sigma_{p,j}&=\sigma_{p,0}-\sum_{l=1}^j v_{p,l}\,.
 \end{aligned}
\end{equation}
The classes $F$ and $\Gamma_{p,l}$ span $\Lambda\otimes\Q$: subtraction of a rational multiple of $\delta_0$ sets $r=0$, after which each remaining array of $k$ coordinates is a rational linear combination of adjacent differences. Their evaluations therefore determine every real or complex functional on $\Lambda$. The formulas in~\eqref{eq:formal-period-inverse} give the inverse.
\end{proof}

Every root of $\Lambda$ has the form $\gamma+kF$, where $\gamma$ belongs to the finite root system in the splitting of Lemma~\ref{lem:numerical-affine-basis}. On the bounded real alcoves, only finitely many integers $k$ can satisfy $P_I(\gamma+kF)=0$. The excluded root loci are therefore the intersections of the parameter domain with a finite union of real codimension-three affine subspaces. In particular, $\mathscr U$ is relatively open in the closed parameter domain, and $\mathscr U^\circ$ is path connected.

\begin{proposition}\label{prop:hyperkahler-periods}
For every $\xi=(\alpha,\sigma)\in\mathscr U$, there is an integral isometry $\iota_\xi:\Lambda\to H_2(\M_\xi,\Z)$ carrying $F$ to a smooth Hitchin fiber and satisfying~\eqref{eq:strict-numerical-periods}.
\end{proposition}
\begin{proof}
For $\xi\in\mathscr U^\circ$ the claim follows from Proposition~\ref{prop:strict-integral-periods}. At a boundary parameter, we use Proposition~\ref{prop:boundary-algebraic} to choose a short segment $\alpha(s)$, $0\leq s\leq s_0$, with $\alpha(0)=\alpha$, $\alpha(s)\in\Acal$ for $s>0$, and the same algebraic moduli functor throughout. Relative openness allows the segment to remain in $\mathscr U$. We regard all its fibers as one smooth algebraic surface $M$.

We choose compact smooth singular cycles $\Gamma_v$ on $M$ representing the transported basis at $s=s_0$. The common algebraic surface gives a product trivialization for $0<s\leq s_0$, so these cycles represent the same homology transport as in~\cite[Corollary~5.4]{DFY}. On a compact set containing their supports, all three forms converge smoothly by Proposition~\ref{prop:weight-continuity}. Hence, for $j=I,J,K$,
\begin{equation*}
 \int_{\Gamma_v}\omega_{j,0}
 =\lim_{s\downarrow0}\int_{\Gamma_v}\omega_{j,s}
 =\lim_{s\downarrow0}P_{j,\alpha(s),\sigma}(b_v)
 =P_{j,\alpha,\sigma}(b_v)\,.
\end{equation*}
The last equality uses the affine formulas. The cycles remain an integral basis with the same intersection form, and their weighted sum is the fiber class. They define $\iota_\xi$, and linearity proves the formulas on every class.
\end{proof}

Under Poincar\'e duality the map $H^2_c(\M_*,\R)\to H^2(\M_*,\R)$ corresponds to the intersection pairing. By Lemma~\ref{lem:numerical-affine-basis}, we obtain
\begin{equation}\label{eq:compact-period-subspace}
 \operatorname{im}\bigl(H^2_c(\M_*,\R)\to H^2(\M_*,\R)\bigr)
 =\{u:u[F]=0\}\,.
\end{equation}
Thus subtraction of the reference triple converts the absolute normalized domain~\eqref{eq:full-period-domain} into the relative period domain used in~\cite[Theorem~1.7]{ChenViaclovskyZhang}.

\subsection{Markings and the ALG coordinates}

We now realize the homology markings by diffeomorphisms that also match the ALG coordinates. Starting from the spectral chart $\Phi_\xi$ in~\eqref{eq:semiflat-ALG-coordinates}, we choose a chart $\Psi_\xi=f_\xi\circ\Phi_\xi$ with decay of order two as in Theorem~\ref{thm:boundary-ALG}. Each $f_\xi$ is a global diffeomorphism isotopic to the identity. We use the resulting chart $\Psi_*$ on $\M_*$.

\begin{proposition}\label{prop:end-compatible-marking}
For every $\xi\in\mathscr U$, there is a diffeomorphism $D_\xi:\M_*\to\M_\xi$ satisfying
\begin{equation}\label{eq:compatible-ALG-marking}
 D_\xi\circ\Psi_*=\Psi_\xi
\end{equation}
on a sufficiently far end and $(D_\xi)_*\iota_*=\iota_\xi$ for a marking in Proposition~\ref{prop:hyperkahler-periods}. Consequently its pulled-back Hitchin structure has decay of order two in $\Psi_*$ and the periods~\eqref{eq:formal-periods}.
\end{proposition}
\begin{proof}
For $\xi\in\mathscr U^\circ$, we choose a smooth path $\xi_s$ in $\mathscr U^\circ$ from $\xi_*$ to $\xi$. The smooth family and the local independence of the algebraic moduli problem from the real weights follow from~\cite[Proposition~4.7, Section~4.3 and Corollary~5.4]{DFY}. Locally one may choose rational weights avoiding the numerical slope equalities and use Proposition~\ref{prop:residue-family}. On each such neighborhood the joint map given by the residues and the Hitchin morphism is projective, so the closed core over a compact subinterval is compact. Finitely many of these subintervals cover the path.

The spectral charts vary smoothly along the path on one sufficiently far end. Indeed, for large $|\lambda|$, all scaled residues $\sigma_s/\lambda$ lie in the parameter ball of Proposition~\ref{prop:normalized-spectral-family}. Gauss--Manin transport there fixes the cycles, and $\pi^*\mathcal O_X(1)$ fixes the origin of the Picard fibers. The coordinates in~\eqref{eq:semiflat-ALG-coordinates} are therefore compatible along the path. Since $U$ is a constant multiple of $\lambda$ and $t=\lambda^n$, differentiation while keeping $(U,V)$ unchanged lifts the parameter velocity and preserves $t$.

We extend this lift over the family by a cutoff and a smooth right inverse of the parameter projection. The lift agrees with the prescribed one near and outside a large level $|t|=r_0$, so it is tangent to that level. Compactness of the core over the path gives a flow for $0\leq s\leq1$, while the flow on the end preserves $(U,V)$. Its final map $D^0:\M_*\to\M_\xi$ satisfies $D^0\circ\Phi_*=\Phi_\xi$. This flow trivializes the family over the path, so its action on homology is the parallel transport used in Proposition~\ref{prop:strict-integral-periods}. Thus
\begin{equation*}
 D_\xi=f_\xi\circ D^0\circ f_*^{-1}
\end{equation*}
has the same homology action and satisfies~\eqref{eq:compatible-ALG-marking}.

For a boundary parameter, we choose a nearby $\xi'=(\alpha',\sigma)\in\mathscr U^\circ$ as in Proposition~\ref{prop:boundary-algebraic}, using the choice in the proof of Proposition~\ref{prop:hyperkahler-periods}. The two moduli spaces are the same algebraic surface and have the same spectral chart, since the normalized spectral curve and its Picard origin depend on the residues. We set
\begin{equation*}
 D_\xi=f_\xi\circ f_{\xi'}^{-1}\circ D_{\xi'}\,.
\end{equation*}
Both correction maps act trivially on homology. This map therefore realizes the marking constructed in the proof of Proposition~\ref{prop:hyperkahler-periods} and agrees near infinity with the change between the two ALG charts in which the decay has order two.
\end{proof}

\subsection{Compact spheres and period realization}

The reference surface with vanishing graded residues has a compact fiber of type $\mathrm{IV}^*$, $\mathrm{III}^*$, or $\mathrm{II}^*$ whose component classes form an integral basis of $H_2(\M_*,\Z)$ with intersection matrix the negative affine Cartan matrix~\cite[Section~6.1]{DFY}. Dehn twists about these spheres change the marking and preserve the ALG chart.

\begin{proposition}\label{prop:compact-root-reflections}
Every reflection
\begin{equation}\label{eq:root-reflection}
 s_\gamma(C)=C+(C\cdot\gamma)\gamma,
 \qquad \gamma^2=-2\,,
\end{equation}
is induced by an orientation-preserving diffeomorphism of $\M_*$ supported in a compact set.
\end{proposition}
\begin{proof}
The oriented normal bundle of a component sphere $S$ has Euler number $-2$. Hence an orientation-preserving diffeomorphism identifies a tubular neighborhood of the sphere with a neighborhood of the zero section in $T^*S^2$. The smooth generalized Dehn twist~\cite[Section~6]{SeidelDehnTwists}, extended by the identity, gives an orientation-preserving diffeomorphism $\tau_S$ with compact support and $(\tau_S)_*[S]=-[S]$. By the relative homology sequence of the tubular neighborhood, $(\tau_S)_*C-C=b[S]$. Since the diffeomorphism preserves the intersection form,
\begin{equation*}
 -C\cdot S=(\tau_S)_*C\cdot S=C\cdot S-2b\,.
\end{equation*}
Thus $b=C\cdot S$, which proves~\eqref{eq:root-reflection} for the component roots.

Under $\iota_*^{-1}$, the component classes give a splitting $\Lambda=\Z F\oplus E_\ell(-1)$, with the distinguished component equal to $F-\theta$, where $\theta$ is the highest root. The other component reflections generate the finite Weyl group. Conjugating the reflection in $F-\theta$ gives $s_{\delta-F}$ for every finite root $\delta$. Direct calculation gives
\begin{equation}\label{eq:affine-transvection}
 s_\delta s_{\delta-F}(C)=C-(C\cdot\delta)F\,.
\end{equation}
For every finite root $\gamma_0$, there is a finite root $\delta$ with $\gamma_0\cdot\delta=\pm1$. Indeed, after a Weyl group element carries $\gamma_0$ to a simple root, we may take an adjacent simple root and apply the inverse Weyl group element. Powers of~\eqref{eq:affine-transvection} therefore carry $\gamma_0$ to $\gamma_0+kF$ for every integer $k$. Conjugation realizes the reflection in each such class. These are all the roots, and the compositions remain compactly supported.
\end{proof}

The component twists generate the group~$G$ introduced in Section~\ref{sec:main-modularity}. For an integral array $u\in\Z^{e_p}$ with $\sum_j u_j=0$, we define
\begin{equation}\label{eq:numerical-translation}
 T_u[(r,d,k)]=[(r,d+\sum_jk_{p,j}u_j,k)]\,.
\end{equation}
This is an integral isometry fixing $F$, and the period formula gives
\begin{equation}\label{eq:translation-period-action}
 P_{\alpha,\sigma}\circ T_u=P_{\alpha+u,\sigma}\,.
\end{equation}
For $u=e_i-e_j$ and $\gamma=[(0,0,u)]$, its reflection permutes the $i,j$ labels, and
\begin{equation}\label{eq:translation-reflection-factorization}
 T_{e_i-e_j}=s_\gamma s_{\gamma+F}\,.
\end{equation}
Indeed, the right-hand side sends $C$ to $C+(C\cdot\gamma)F$, and $C\cdot\gamma=-\sum_jk_{p,j}u_j$. By Proposition~\ref{prop:compact-root-reflections} every such reflection is induced by a compactly supported diffeomorphism. Thus elements of $G$ realize all zero-sum integral translations and permutations of the labels.

\begin{lemma}\label{lem:closed-alcove-reduction}
Every centered real array can be carried into $\overline{\Acal}$ by zero-sum integral translations and permutations of the labels at each puncture. Their period action preserves admissibility.
\end{lemma}
\begin{proof}
At a given puncture, we choose a point of minimum Euclidean norm in the discrete set $\alpha+\{u\in\Z^e:\sum_j u_j=0\}$. If its smallest and largest coordinates $a,b$ satisfy $b-a>1$, increasing $a$ by one and decreasing $b$ by one changes the squared norm by
\begin{equation*}
 (a+1)^2+(b-1)^2-a^2-b^2
 =2\bigl(1-(b-a)\bigr)<0\,,
\end{equation*}
a contradiction. Sorting the coordinates gives a point of the closed alcove. We repeat this at each puncture and apply the same permutations to the residue labels. The corresponding integral isometries permute the roots and hence preserve the simultaneous nonvanishing condition.
\end{proof}

\begin{proof}[Proof of Theorem~\ref{thm:full-period-realization}]
For $P\in\mathcal D$, the inverse in Lemma~\ref{lem:formal-period-inverse} determines centered formal parameters $(\alpha,\sigma)$. By Lemma~\ref{lem:closed-alcove-reduction} there is an isometry $h$ induced by an element of $G$ such that $P\circ h=P_{\alpha',\sigma'}$ with $(\alpha',\sigma')\in\mathscr U$. The complete Hitchin metric of Theorem~\ref{thm:boundary-ALG} pulls back, by Propositions~\ref{prop:hyperkahler-periods} and~\ref{prop:end-compatible-marking}, to a metric on $\M_*$ with these periods and decay of order two in $\Psi_*$. Pulling back once more by an element of $G$ inducing $h^{-1}$ gives periods $P$. Compact support preserves $\Psi_*$ on a sufficiently far end. Conversely, the period formulas and admissibility place every Hitchin period triple in $\mathcal D$, and every element of $G$ preserves this domain.
\end{proof}

\subsection{Modularity}

We finish by matching a given exceptional ALG instanton with one of the Hitchin metrics. After identifying the two ends and pulling both structures back to the reference surface, we can compare their periods there.

\begin{proof}[Proof of Theorem~\ref{thm:modularity}]
For $(Y,g,\boldsymbol\omega)$ with cone parameter $1/3$, $1/4$, or $1/6$, the bound $\beta\leq1/2$ ensures that there are ALG coordinates with decay of order two~\cite[Remark~1.4]{ChenViaclovskyZhang}. We orient the primitive fiber class and rotate the ordered triple so that its fiber period is $(A,0,0)$ with $A>0$. Multiplication of the metric and forms by $A_F/A$ normalizes the area to $A_F$. The finite torus monodromy gives $\tau_\infty=e^{2\pi\ii/3},\ii,e^{2\pi\ii/3}$, respectively, up to a change of lattice basis, and the area determines $L=2\pi/\sqrt{\operatorname{Im}\tau_\infty}$.

By Theorems~D and~E of~\cite{CV}, there is a diffeomorphism from $\M_*$ to $Y$ for which the pulled-back normalized structure has decay of order two in $\Psi_*$ with the same ordered model forms. The period triple of this structure $P$ is nonzero on every integral class of square $-2$ by~\cite[Theorem~1.7]{ChenViaclovskyZhang}, so $P\in\mathcal D$.

By Theorem~\ref{thm:full-period-realization}, a Hitchin structure on $\M_*$ has the same periods and decay of order two in the same chart. Both metrics are complete and nonflat, with square-integrable curvature. Their ordered de Rham classes agree because their integrals agree on every integral two-cycle. By the ALG Torelli theorem~\cite[Theorem~1.5]{ChenViaclovskyZhang}, there is a diffeomorphism carrying one metric and its ordered triple to the other. Composing with the two markings and undoing the normalization proves the modularity statement.
\end{proof}

\appendix

\section{The characteristic coefficients}\label{app:coefficients}

We prove Proposition~\ref{prop:coefficient-line} in the notation of Section~\ref{sec:hitchin-base}.

\begin{proof}
We expand $a_k(z)=\sum_{j=0}^k a_{k,j}z^j$ and solve the conditions separately in the three ranks.

For $n=3$, we write
\begin{equation*}
 P_0=y^3+A_2y+A_3,\qquad P_1=y^3+B_2y+B_3,\qquad
 P_\infty=v^3+C_2v+C_3.
\end{equation*}
All multiplicities are one, so the constant terms, the values at $z=1$, and the highest coefficients of $a_2$ and $a_3$ are prescribed. These conditions determine $a_2$ and leave one free multiple of $z(z-1)$ in $a_3$. With $t=-a_{3,2}$, the solutions are
\begin{equation}\label{eq:rank-three-coefficients}
 \begin{aligned}
 a_2(z)&=A_2+(B_2-A_2-C_2)z+C_2z^2,\\
 a_3(z)&=A_3+(B_3-A_3-C_3)z+C_3z^3-tz(z-1).
 \end{aligned}
\end{equation}
For $n=4$, the residue polynomials are
\begin{equation*}
 \begin{aligned}
 P_0&=y^4+A_2y^2+A_3y+A_4,\\
 P_1&=y^4+B_2y^2+B_3y+B_4,\qquad P_\infty=v^2+b,
 \end{aligned}
\end{equation*}
and we set $L=B_2-A_2-2b$. With $t=-a_{4,2}$, the solutions are
\begin{equation}\label{eq:rank-four-coefficients}
 \begin{aligned}
 a_2(z)&=A_2+Lz+2bz^2,\\
 a_3(z)&=A_3+(B_3-A_3)z,\\
 a_4(z)&=A_4+(B_4-A_4-bL-b^2)z\\
       &\quad+bLz^3+b^2z^4-tz(z-1).
 \end{aligned}
\end{equation}
The value $F^\infty(0,v)=(v^2+b)^2$ and divisibility of the coefficient of $u$ by $v^2+b$ give, respectively,
\begin{equation*}
 \begin{aligned}
 a_{2,2}&=2b,& a_{3,3}&=0,& a_{4,4}&=b^2,\\
 a_{3,2}&=0,& a_{4,3}&=b\,a_{2,1}.
 \end{aligned}
\end{equation*}
The values at $0$ and $1$ then determine $a_2$ and $a_3$. The difference between two possible choices of $a_4$ has degree at most two and vanishes at both points, so the only freedom is a multiple of $z(z-1)$, as in~\eqref{eq:rank-four-coefficients}.

For $n=6$, we write
\begin{equation*}
 \begin{aligned}
 P_0&=y^6+A_2y^4+A_3y^3+A_4y^2+A_5y+A_6,\\
 P_1&=y^3+cy+d,\qquad P_\infty=v^2+b.
 \end{aligned}
\end{equation*}
The coefficients are given recursively by
\begin{equation}\label{eq:rank-six-constants}
 \begin{aligned}
 L&=2c-A_2-3b,\qquad M=2d-A_3,\\
 N&=A_4-c^2+c(6b+L)-9b^2-4bL,\\
 N_5&=A_5-2cd+d(6b+L)+(c-2b)M,\\
 J&=bN-b^2A_2,\\
 K&=A_6-d^2+dM-5b^3-4b^2L-3J.
 \end{aligned}
\end{equation}
Then all solutions are
\begin{equation}\label{eq:rank-six-coefficients}
 \begin{aligned}
 a_2(z)&=A_2+Lz+3bz^2,\\
 a_3(z)&=A_3+Mz,\\
 a_4(z)&=A_4+(c^2-A_4-N-2bL-3b^2)z\\
       &\quad+Nz^2+2bLz^3+3b^2z^4,\\
 a_5(z)&=A_5+(2cd-A_5-N_5-bM)z+N_5z^2+bMz^3,\\
 a_6(z)&=A_6+(d^2-A_6-K-J-b^2L-b^3)z+Kz^2\\
       &\quad+Jz^4+b^2Lz^5+b^3z^6-tz(z-1)^2.
 \end{aligned}
\end{equation}
The value at $0$ gives $a_k(0)=A_k$ for $2\leq k\leq6$, while $F(1,y)=(y^3+cy+d)^2$ prescribes
\begin{equation*}
 (a_2(1),a_3(1),a_4(1),a_5(1),a_6(1))
       =(2c,2d,c^2,2cd,d^2).
\end{equation*}
Divisibility of $\partial_zF(1,y)$ by $y^3+cy+d$ gives
\begin{equation}\label{eq:rank-six-at-one}
 \begin{aligned}
 a_4'(1)&=c\,a_2'(1),\\
 a_5'(1)&=d\,a_2'(1)+c\,a_3'(1),\\
 a_6'(1)&=d\,a_3'(1).
 \end{aligned}
\end{equation}
At infinity, the value $F^\infty(0,v)=(v^2+b)^3$ gives
\begin{equation*}
 (a_{2,2},a_{3,3},a_{4,4},a_{5,5},a_{6,6})
       =(3b,0,3b^2,0,b^3).
\end{equation*}
The coefficients of $u$ and $u^2$ in $F^\infty$ must be divisible by $(v^2+b)^2$ and $v^2+b$, respectively. Polynomial division gives precisely
\begin{equation}\label{eq:rank-six-at-infinity}
 \begin{aligned}
 a_{3,2}&=0,& a_{5,4}&=0,\\
 a_{4,3}&=2b\,a_{2,1},& a_{6,5}&=b^2a_{2,1},\\
 a_{5,3}&=b\,a_{3,1},&
 a_{6,4}&=b\,a_{4,2}-b^2a_{2,0}.
 \end{aligned}
\end{equation}
The value conditions and~\eqref{eq:rank-six-at-infinity} determine $a_2$ and $a_3$. For $a_4$ and $a_5$, the values at $0$ and $1$ express the remaining coefficients in terms of $a_{4,2}$ and $a_{5,2}$. In the first two equations in~\eqref{eq:rank-six-at-one}, these unknowns each occur with coefficient one, giving $a_{4,2}=N$ and $a_{5,2}=N_5$. The last equation in~\eqref{eq:rank-six-at-infinity} gives $a_{6,4}=J$.

With $a_2,\ldots,a_5$ determined, the difference $\Delta a_6$ between two possible choices of $a_6$ has degree at most three and satisfies $\Delta a_6(0)=\Delta a_6(1)=\Delta a_6'(1)=0$. Hence $\Delta a_6$ is a multiple of $z(z-1)^2$. Taking $a_{6,3}=-t$ gives~\eqref{eq:rank-six-coefficients}. Substitution verifies all the value and divisibility conditions, which are exactly the conditions of Lemma~\ref{lem:local-coefficients}.

Each $A_k,B_k,C_k$ is homogeneous of degree $k$ in the residues; $b,c$ have degree two and $d$ has degree three. The formulas therefore make $a_k(z;0,\sigma)$ homogeneous of degree $k$ in all three ranks. They involve no division by a residue or a difference of residues, so they hold for all prescribed scalar values. The local realizations follow from Lemma~\ref{lem:local-coefficients}.
\end{proof}

\section{The algebraic moduli construction}

\subsection{Deformations with fixed graded residues}\label{proof:prop:local-deformations}
\begin{proof}
To prove Proposition~\ref{prop:local-deformations}, we use the trace pairing to identify the dual sheaf
\begin{equation*}
 (\ParEnd_0E)^*\otimes K_X
   =\SParEnd_0E\otimes K_X(D).
\end{equation*}
Indeed, the annihilator of the flag-preserving algebra under trace is the algebra strictly lowering the flag; the factor $\mathcal O_X(D)$ accounts for the simple poles in the dual lattice. The trace pairing carries the dual commutator to the original differential with the indicated sign change. This is the parabolic deformation duality used in~\cite[Sections 2.1 and 3.2.1]{LogaresMartens}.

The common shift $\alpha'_{p,j}=\alpha_{p,j}-\alpha_{p,0}$ puts the weights in $[0,1)$ and changes the parabolic degree of each subbundle $F$ by $-\rank(F)\sum_p\alpha_{p,0}$. It therefore preserves all slope differences, while leaving the flags, Higgs field and parabolic endomorphisms unchanged. By~\cite[Remark~3.5, p.~432]{MondelloParabolic}, every holomorphic endomorphism preserving the flags and commuting with $\Phi$ is scalar. The trace-free condition gives $\mathbb H^0(\mathcal C^\bullet)=0$, and by Serre duality $\mathbb H^2(\mathcal C^\bullet)=0$. By Riemann--Roch,
\begin{equation*}
 \dim_\C\mathbb H^1(\mathcal C^\bullet)
 =-2(n^2-1)+\sum_{p\in D}(n^2-e_pm_p^2)=2.
\end{equation*}

We next consider higher-order deformations with the same fixed graded residues. In local frames adapted to the flags and determinant, Higgs fields with the fixed scalar graded residues form an affine space modeled on the sections of $\SParEnd_0E\otimes K_X(D)$. The flag-preserving frame groups are smooth. Over a small extension of local Artin $\C$-algebras with kernel $J$, local lifts therefore exist. The differences in the cocycle condition for their transition maps and in the compatibility of their Higgs fields form a total degree-two \v Cech cocycle for $\mathcal C^\bullet\otimes J$. Changing the local lifts changes this cocycle by a total coboundary. The obstruction class vanishes by $\mathbb H^2=0$, so every deformation lifts. The isomorphism classes of lifts form a torsor for $\mathbb H^1(\mathcal C^\bullet)\otimes J$, and automorphisms reducing to the identity are controlled by $\mathbb H^0(\mathcal C^\bullet)\otimes J=0$. Once a basis of $\mathbb H^1(\mathcal C^\bullet)$ is fixed, successive lifting constructs a formal deformation over $\C[[s_1,s_2]]$. At each small extension its tangent map induces an isomorphism between the two torsors of lifts; induction on the length of the Artin algebra proves that it represents the deformation functor with fixed determinant and graded residues.

Every determinant-preserving automorphism is a scalar $\zeta\Id$ with $\zeta^n=1$. The same scalar extends over every deformation, so its action on the formal deformation space is trivial.
\end{proof}

\subsection{The moduli scheme and its local structure}\label{proof:prop:algebraic-moduli}
\begin{proof}
To prove Proposition~\ref{prop:algebraic-moduli}, we return to the neighborhood constructed in the proof of Proposition~\ref{prop:joint-walls}. In this neighborhood there are rational weights $\alpha'$ outside the finitely many numerical weight hyperplanes. These weights preserve the moduli problem with the given graded residues, fiberwise in families. They are also generic for parabolic Higgs bundles with unrestricted residues. The weights $\alpha'_{p,j}-\alpha'_{p,0}$ lie in $[0,1)$; this translation changes every parabolic slope by the same constant and preserves stability. The construction recalled in~\cite[Sections 2.1--2.3.1]{LogaresMartens} therefore gives the degree-zero moduli scheme $\mathcal P$ of semistable parabolic Higgs bundles with these rational weights and flag types. Its Hitchin morphism to $\bigoplus_{i=1}^nH^0(X,K_X(D)^{\otimes i})$ is projective by Yokogawa's result, as recalled in~\cite[Section 2.3.1]{LogaresMartens}.

We take $q:R\to\mathcal P$ to be the good quotient of the GIT parameter scheme by the reductive group $G$ acting on the frames. On $R$, the universal flags have locally free graded quotients $Q_{p,j}$, and the universal residues induce endomorphisms $R_{p,j}$ of $Q_{p,j}$. On these quotients we impose the equations
\begin{equation}\label{eq:universal-scalar-residues}
 R_{p,j}-\sigma_{p,j}\Id_{Q_{p,j}}=0
 \quad\text{for all }p,j,
 \qquad \operatorname{tr}\Phi=0.
\end{equation}
These endomorphism identities define a $G$-invariant closed subscheme $Z\subset R$ even when prescribed scalars coincide. Over any base, including a nonreduced one, a framed family is parametrized by $Z$ exactly when the identities hold on its graded quotients. Every semistable family admits such a frame locally on its base. Exactness of invariants for $G$ shows that $Z\mathbin{/\!/}G$ is a closed subscheme of $\mathcal P$. Since the chosen rational weights make every semistable object stable, this quotient is geometric and is the coarse moduli scheme of objects with the given weights and graded residues. As the conditions agree in families, this coarse moduli scheme is isomorphic to $\M_{\alpha,\sigma}$. Over a discrete valuation ring the equations persist in semistable specializations. Indeed, the graded quotients are locally free, so an endomorphism section vanishing on the generic fiber vanishes on the whole family. Restricting the projective Hitchin morphism to this closed subscheme, and then to the subspace $\mathcal H_0$ of its target, proves projectivity. The closed immersion into the quasi-projective scheme $\mathcal P$ proves quasi-projectivity.

Every degree-zero determinant line bundle on $\Pp^1$ is isomorphic to $\mathcal O_X$. For a family over an algebraic or analytic base $S$, we write $p:X\times S\to S$ for the projection. Since the determinant line bundle on each fiber has $h^0=1$ and $h^1=0$, its direct image $p_*(\det E)$ is a line bundle by cohomology and base change, and evaluation gives an isomorphism $p^*p_*(\det E)\simeq\det E$. Thus determinant trivializations exist locally on $S$. Two such trivializations differ by a unit, which has an $n$th root \'{e}tale-locally in the algebraic setting and locally in the analytic setting. Multiplication by this scalar on $E$ absorbs the difference. The analytic assertion includes nonreduced bases: a unit has a local holomorphic $n$th root near a nonzero value. Forgetting the determinant trivialization therefore gives the same coarse moduli scheme and the same sheaf of isomorphism classes over analytic bases. The argument also applies over local Artin $\C$-algebras.

The local construction applies at any stable pair $x$ with $\alpha\in\Acal$ and fixed $\sigma$. Near $x$, the framed parameter schemes for bundles, flags and logarithmic Higgs fields are of finite type. Only finitely many numerical slope expressions $d+\sum_{p,j}k_{p,j}\alpha_{p,j}$ can change sign in a small weight neighborhood. A rational weight in the same face of these rational hyperplanes preserves their signs and zeros, so the rational GIT construction gives the same stable open locus. Imposing~\eqref{eq:universal-scalar-residues} and a determinant trivialization gives a finite-type neighborhood $\mathfrak S$ of $x$ in the stable moduli stack. It is quasi-separated with affine diagonal, since isomorphisms of the bundles form affine schemes and preservation of the flags, Higgs fields and determinant trivialization imposes closed equations.

The stabilizer at $x$ is the smooth linearly reductive group $\mu_n$. By Theorem~1.1 of~\cite{AlperHallRydh}, there is an \'etale quotient chart $[\operatorname{Spec}A/\mu_n]\to\mathfrak S$ with a fixed point $a$ over $x$. By the deformation calculation in Proposition~\ref{prop:local-deformations} and~\cite[Remark~4.17]{AlperHallRydh}, the completed local ring is isomorphic to $\widehat{\mathcal O}_{\operatorname{Spec}A,a}\simeq\C[[s_1,s_2]]$, with trivial $\mu_n$ action. After shrinking, $U=\operatorname{Spec}A$ is smooth of dimension two and the action is trivial: the differences of each group element on a finite set of algebra generators vanish in the local ring because they vanish in its completion. The atlas map determines a family with the prescribed determinant trivialization and scalar graded residues,
\begin{equation}\label{eq:stable-local-family}
 (U,a)\longrightarrow(\mathfrak S,x),\qquad
 \widehat{\mathcal O}_{U,a}\simeq\C[[s_1,s_2]]\,.
\end{equation}
Its Kodaira--Spencer map $T_aU\to\mathbb H^1(\mathcal C^\bullet)$ is an isomorphism.

For the GIT stack with only the isomorphism class of the determinant line bundle fixed, the local existence of determinant trivializations gives the same miniversal ring, and its scalar stabilizer $\C^*$ acts trivially. Applying Theorem~4.12 and Remark~4.17 of~\cite{AlperHallRydh}, we obtain
\begin{equation*}
 \widehat{\mathcal O}_{\M_{\alpha,\sigma},x}
 \simeq\C[[s_1,s_2]]^{\C^*}
 =\C[[s_1,s_2]]\,.
\end{equation*}
Every closed point is thus smooth of dimension two, which proves the smoothness and dimension of $\M_{\alpha,\sigma}$.
\end{proof}

\section{Decay on a punctured disc}

\subsection{The half-cylinder boundary problem}\label{proof:lem:harmonic-half-cylinder}

We use the local data $A,B$ in~\eqref{eq:harmonic-local-data}; in particular, $B$ is scalar on every equal-weight subspace. In $V=\Herm_0(\C^n)$ we write $K\subset V$ for the constant block diagonal matrices on these subspaces. The spaces $H^j_\nu$ are as defined before Lemma~\ref{lem:harmonic-half-cylinder}. For the model $q=0$, the linearization of twice the contracted harmonic metric equation is
\begin{equation}\label{eq:harmonic-cylinder-linearization}
 P=-\partial_x^2-(\partial_\theta+\ii\ad A)^2
       +4\ad(B^*)\ad(B).
\end{equation}

\begin{lemma}\label{lem:cylinder-linear-dirichlet}
For $0<\nu<\varepsilon_p$, the map
\begin{equation}\label{eq:harmonic-cylinder-inverse}
 (P,\operatorname{ev}_0):K\oplus H^4_\nu
       \longrightarrow H^2_\nu\oplus H^{7/2}(S^1;V)
\end{equation}
is a bounded isomorphism. If $P(c+v)=f$ and $(c+v)|_{x=0}=b$, then
\begin{equation*}
 |c|+\|v\|_{H^4_\nu}
 \leq C\bigl(\|f\|_{H^2_\nu}+\|b\|_{H^{7/2}}\bigr)\,.
\end{equation*}
The constant depends only on the fixed matrices $A,B$ and on $\nu$.
\end{lemma}

\begin{proof}
After complexification, the action of $P$ on the $e^{\ii k\theta}E_{ab}$ component is $-\partial_x^2+\mu_{abk}^2$, where
\begin{equation}\label{eq:harmonic-cylinder-frequencies}
 \mu_{abk}^2=(k+\alpha_a-\alpha_b)^2
                      +4|\sigma_a-\sigma_b|^2.
\end{equation}
The zero eigenvalue of the angular operator occurs precisely when $k=0$ and $a,b$ belong to the same weight block. Every positive frequency satisfies $\mu_{abk}\geq\varepsilon_p$.

The inverse is explicit. For a component in $K$, the right-hand side $f$ and the boundary value $b$ give
\begin{equation}\label{eq:harmonic-constant-green}
 v(x)=-\int_x^\infty(y-x)f(y)\,\dd y,
 \qquad c=b+\int_0^\infty y f(y)\,\dd y,
 \qquad u=c+v.
\end{equation}
The kernels obtained after multiplication by $e^{\nu x}$ are integrable, and the equation controls the remaining derivatives. By the Cauchy--Schwarz inequality, $|c|$ is also bounded by a constant times $|b|+\|e^{\nu x}f\|_2$. For $\mu>0$ the solution is
\begin{equation}\label{eq:harmonic-positive-green}
 u(x)=b e^{-\mu x}+\int_0^\infty
 \frac{e^{-\mu|x-y|}-e^{-\mu(x+y)}}{2\mu}f(y)\,\dd y.
\end{equation}
To see that the estimates are uniform over the angular Fourier components, we use $w=e^{\nu x}u$ for the solution with zero boundary value. Then
\begin{equation*}
 -w''+2\nu w'+(\mu^2-\nu^2)w=e^{\nu x}f.
\end{equation*}
Its real $L^2$ pairing with $w$ controls $\|w'\|_2^2+(\mu^2-\nu^2)\|w\|_2^2$. The equation then bounds $\mu^2\|w\|_2$, $\mu\|w'\|_2$, and $\|w''\|_2$ by a fixed multiple of $\|e^{\nu x}f\|_2$, since $\mu\geq\varepsilon_p>\nu$. Applying these estimates to angular derivatives and then using the equation for two further radial derivatives gives the $H^4_\nu$ bound from $H^2_\nu$. For the homogeneous boundary term, direct integration gives
\begin{equation*}
 \sum_{i+j\leq4}\int_0^\infty e^{2\nu x}
 \bigl|\partial_x^i\partial_\theta^j
       (b e^{-\mu x+\ii k\theta})\bigr|^2\,\dd x
 \leq C(1+k^2)^{7/2}|b|^2\,.
\end{equation*}
Here $\mu$ is comparable to $1+|k|$ on the components with $\mu_{abk}>0$, with constants depending on the fixed data. Summing proves~\eqref{eq:harmonic-cylinder-inverse}; the formulas also prove uniqueness and preserve the Hermitian condition.
\end{proof}

\begin{proof}[Proof of Lemma~\ref{lem:harmonic-half-cylinder}]
For the model $q=0$, every constant $c\in K$ gives a harmonic metric $h_0e^c$. We translate the end by $x=\tau-\tau_0$. The remainder $q_{\tau_0}$ in its Higgs coefficient satisfies
\begin{equation*}
 |\partial_x^i\partial_\theta^j q_{\tau_0}|
 \leq C_{i+j}e^{-\varepsilon_p\tau_0}e^{-\varepsilon_p x}.
\end{equation*}
For $s=c+v\in K\oplus H^4_\nu$ near zero, we use the metric $h_s=h_0e^s$ to define
\begin{equation}\label{eq:harmonic-local-residual}
 \mathcal F_{\tau_0}(s)=
 2e^{s/2}\ii\Lambda_{\mathrm{cyl}}\bigl(F_{h_s}+[\Phi,\Phi^{\dagger_{h_s}}]\bigr)e^{-s/2},
\end{equation}
where $\Lambda_{\mathrm{cyl}}$ is contraction for the cylindrical metric. In the fixed model frame the value of this expression is a trace-free Hermitian matrix. Near zero we have the smooth map
\begin{equation*}
 K\oplus H^4_\nu \longrightarrow H^2_\nu\oplus H^{7/2}(S^1;V),
 \qquad
 s\longmapsto\bigl(\mathcal F_{\tau_0}(s),s|_{x=0}\bigr).
\end{equation*}
A constant $c\in K$ commutes with $A$, $B$ and $B^*$ and solves the model equation. Every remaining term contains $v$, $q_{\tau_0}$ or one of their derivatives. Sobolev multiplication in dimension two and smooth matrix exponentiation therefore give the weighted mapping property and bound the second derivative on a fixed small ball, uniformly for large $\tau_0$. Moreover, $\mathcal F_{\tau_0}(0)\to0$ in $H^2_\nu$, and the derivative at zero of $(\mathcal F_{\tau_0},\operatorname{ev}_0)$ converges in operator norm to $(P,\operatorname{ev}_0)$. Using the inverse in Lemma~\ref{lem:cylinder-linear-dirichlet} and the quantitative implicit function theorem, we can solve
\begin{equation*}
 (\mathcal F_{\tau_0}(s),s|_{x=0})=(0,b)
\end{equation*}
for all sufficiently large $\tau_0$ and all $b$ in one fixed neighborhood of zero. We obtain~\eqref{eq:harmonic-half-cylinder-metric}. Since $s=c+v$ is trace-free Hermitian, the metric is positive and has determinant one.

Since $v\in H^4_\nu$, it and its first two derivatives are $O(e^{-\nu x})$. Differentiating the harmonic metric equation on nested cylindrical strips and applying interior elliptic estimates gives this decay for every derivative of $v$. At each step the derivatives of $v$ and of $q_{\tau_0}$ on the right-hand side have already been estimated. Moreover, the coefficients remain uniformly elliptic, and the constant $c$ solves the model equation. Smoothness of matrix exponentiation gives the same decay estimates for $\widetilde H-e^c$. A constant block diagonal change of frame preserves $A$ and $B$ and conjugates $q_{\tau_0}$. On a compact set of such changes, all the coefficient bounds and smallness thresholds used above are uniform.
\end{proof}

\subsection{Local estimates and translated limits}\label{proof:lem:local-harmonic-exponential-decay}

Interior estimates first control translated metrics on finite strips. The trace identity then determines their possible limits on the entire cylinder.

\begin{lemma}\label{lem:harmonic-local-compactness}
\begin{enumerate}[label=\textup{(\roman*)}]
\item For coordinate discs $U'\Subset U$ and a holomorphic Higgs field $\phi\,\dd w$ on $U$, we write $M$ for the matrix of a harmonic metric in a holomorphic frame. If $C^{-1}\Id\leq M\leq C\Id$ and the derivatives of $\phi$ in the disc coordinate satisfy fixed bounds of every order on $U$, then for each $k\geq0$ the norm $\|M\|_{C^k(U')}$ is bounded by a constant depending only on these positive metric bounds, the required coefficient bounds, and the distance from $\overline{U'}$ to $\partial U$.
\item For a harmonic metric $h$ of determinant one with local data~\eqref{eq:harmonic-local-data} satisfying $C^{-1}h_0\leq h\leq Ch_0$, we write $H$ for its matrix in the frame $(f_a)$ of~\eqref{eq:harmonic-power-model}. Every sequence $\tau_j\to\infty$ has a subsequence for which $H(\tau_j+x,\theta)$ converges in $C^\infty_{\mathrm{loc}}(\R\times S^1)$ to a constant positive matrix of determinant one that is block diagonal with blocks indexed by the distinct weights. For every $T>0$ and $k\geq0$, the matrices $H(\tau+\cdot,\cdot)$ approach the set of such constants in $C^k([-T,T]\times S^1)$ as $\tau\to\infty$.
\end{enumerate}
\end{lemma}

\begin{proof}
We use the local compactness argument in the proof of~\cite[Proposition~2.6, pp.~18--19]{LiMochizukiToda}, keeping track of the Higgs coefficient.
For part~(i), we set $a=M^{-1}\partial_wM$. Then
\begin{equation}\label{eq:harmonic-compactness-equation}
 \partial_{\bar w}a=[\phi,M^{-1}\phi^\dagger M],\qquad
 \partial_{\bar w}\partial_w\tr M
 =\tr\bigl((\partial_wM)^\dagger M^{-1}\partial_wM\bigr)
   +\tr(M\partial_{\bar w}a).
\end{equation}
Since $C^{-1}\Id\leq M\leq C\Id$ and the coefficients of $\phi$ are uniformly bounded, the commutator is uniformly bounded as well. Hence
\begin{equation*}
 (\partial_{\operatorname{Re}w}^2+\partial_{\operatorname{Im}w}^2)\tr M
 \geq c|\dd M|^2-C'.
\end{equation*}
Integrating against a cutoff gives uniform $L^2$ bounds for $\dd M$ and $a$ on a smaller disc. The interior $\bar\partial$ estimate gives $a\in W^{1,2}$. By Sobolev embedding in two dimensions, $a$ is then bounded in every finite $L^p$, so a second application of the interior estimate yields $a\in W^{1,p}$. For $p>2$, we obtain $a\in L^\infty$; the identities $\partial_wM=Ma$ and $\partial_{\bar w}M=a^\dagger M$ give first $M\in W^{1,\infty}$ and then $M\in W^{2,p}$. After differentiating, we bound the products by Sobolev multiplication on smaller discs. This gives bounds for derivatives of every order. The constants depend only on the positive bounds for $M$, the bounds for the derivatives of $\phi$ in the disc coordinate, and the separation between the boundaries of the nested discs. A finite cover gives the same estimates on each smaller cylindrical strip.

For part~(ii), on a fixed simply connected disc in the flat cylinder we choose a $D_0$-parallel holomorphic $h_0$-unitary frame. Its transition from $(f_a)$ has uniformly bounded derivatives, and~\eqref{eq:harmonic-taylor-decay} gives the coefficient bounds required in part~(i).
For any sequence $\tau_j\to\infty$, a diagonal subsequence of $H(\tau_j+x,\theta)$ therefore converges smoothly on compact subsets of $\R\times S^1$. Its limit $H_*$ has determinant one, remains between $C^{-1}\Id$ and $C\Id$, and is harmonic for the exact model $D_0$, $B(-\dd x+\ii\,\dd\theta)$, because the translated $q$ tends smoothly to zero. We apply~\eqref{eq:harmonic-trace-identity} to the model metric $\Id$ and $H_*$. The circle average
\begin{equation*}
 F(x)=\frac1{2\pi}\int_{S^1}\tr H_*(x,\theta)\,\dd\theta
\end{equation*}
is constant because it is bounded and convex on $\R$. The nonnegative right-hand side of~\eqref{eq:harmonic-trace-identity} has zero circle average and vanishes pointwise. Its two terms give $D_0H_*=0$ and $[B,H_*]=0$, so
\begin{equation*}
 H_*(x,\theta)=e^{-\ii A\theta}C_0e^{\ii A\theta}
\end{equation*}
for a constant positive matrix $C_0$. Periodicity implies $(C_0)_{ab}=0$ unless $\alpha_a-\alpha_b\in\Z$. The weight interval has width less than one, so $C_0$ is block diagonal on the equal-weight subspaces and $H_*=C_0$.

For any fixed $T$ and $k$, a sequence of translates staying a positive $C^k([-T,T]\times S^1)$ distance from the set of constant positive block diagonal matrices of determinant one would have a subsequential limit of the form just proved. This contradiction proves the last assertion.
\end{proof}

\begin{proof}[Proof of Lemma~\ref{lem:local-harmonic-exponential-decay}]
By Lemma~\ref{lem:harmonic-local-compactness}(ii), the translates approach the constant positive block diagonal matrices of determinant one in $C^4([-1,1]\times S^1)$. We first use this conclusion to normalize one boundary value and construct a single limiting matrix.
With $\Pi$ the projection onto the block diagonal matrices, we define
\begin{equation}\label{eq:harmonic-boundary-normalization}
 D_{\tau_0}=
 \bigl(\det\Pi H(\tau_0,0)\bigr)^{-1/n}\Pi H(\tau_0,0).
\end{equation}
These positive matrices of determinant one and their inverses stay in a fixed compact set. The normalization in~\eqref{eq:harmonic-boundary-normalization} fixes every positive block diagonal matrix of determinant one, so
\begin{equation*}
 D_{\tau_0}^{-1/2}H(\tau_0,\cdot)D_{\tau_0}^{-1/2}
 \longrightarrow\Id\quad\text{in }C^4(S^1).
\end{equation*}

We fix $0<\nu_0<\varepsilon_p$ and choose $\tau_0$ sufficiently large. After translating by $x=\tau-\tau_0$ and changing frame by $D_{\tau_0}^{-1/2}$, the coefficient of the Higgs field is
\begin{equation*}
 B+\widehat q_{\tau_0},\qquad
 \widehat q_{\tau_0}
 =D_{\tau_0}^{1/2}q(\tau_0+x,\theta)D_{\tau_0}^{-1/2}.
\end{equation*}
Its remainder satisfies $\|e^{\nu_0x}\widehat q_{\tau_0}\|_{H^j(\End\C^n)}=O(e^{-\varepsilon_p\tau_0})$ for every fixed $j$. The logarithm of the normalized boundary metric tends to zero in $H^{7/2}(S^1)$. By Lemma~\ref{lem:harmonic-half-cylinder}, there is a harmonic metric with that boundary value and matrix $e^{c+v}$, where $c\in K$ and $v\in H^4_{\nu_0}$. In the original frame the metric just constructed and $h$ have the same boundary value, determinant one, and mutual bounds with $h_0$. By Lemma~\ref{lem:harmonic-trace-comparison}, they agree throughout the punctured disc. In particular,
\begin{equation*}
 H\longrightarrow H_\infty
 =D_{\tau_0}^{1/2}e^cD_{\tau_0}^{1/2}.
\end{equation*}
Since $\tau_0$ is fixed, the limit is a single matrix. By the estimates in Lemma~\ref{lem:harmonic-half-cylinder},
\begin{equation}\label{eq:harmonic-initial-convergence}
 \sup_{[s,s+1]\times S^1}
 |\partial_\tau^i\partial_\theta^j(H-H_\infty)|\longrightarrow0
 \quad\text{as }s\longrightarrow\infty,
 \qquad i+j\leq4.
\end{equation}

For any $0<\nu<\varepsilon_p$, we normalize by the constant frame change $H_\infty^{-1/2}$ and apply Lemma~\ref{lem:harmonic-half-cylinder} on a sufficiently far half-cylinder. Smallness of the boundary value follows from~\eqref{eq:harmonic-initial-convergence}, and by Lemma~\ref{lem:harmonic-trace-comparison} the resulting solution $e^{c+v}$ equals the normalized $H$. Its limit is $\Id$, so $c=0$. By the derivative estimates in Lemma~\ref{lem:harmonic-half-cylinder}, we have~\eqref{eq:harmonic-local-exponential-decay} for this $\nu$. The convergence of $H$ determines the limiting matrix uniquely.
\end{proof}

\section{Weighted deformations and the hyperk\"ahler quotient}

\subsection{The weighted Dolbeault comparison}\label{proof:prop:weighted-Dolbeault-comparison}

We work near a puncture $p$ in a holomorphic frame compatible with the flag, with the diagonal metric~\eqref{eq:harmonic-power-model} and $0<d<\varepsilon_*$ as in~\eqref{eq:tangent-weight-window}. With $s_{ab}=\alpha_a-\alpha_b\in(-1,1)$, the holomorphic matrix unit $E_{ab}$ has norm $r^{s_{ab}}$. The squared weighted norms use $e^{d\tau}\,\dd\tau\,\dd\theta$, as in~\eqref{eq:harmonic-weighted-norm}.

\begin{lemma}\label{lem:weighted-local-dolbeault}
After restriction to a smaller punctured disc, the following statements hold.
\begin{enumerate}[label=\textup{(\roman*)}]
\item In the frame unitary for $h_0$, every smooth section of $\ParEnd_0E$ is a constant in $\mathfrak l_p$ plus a remainder in $W^{m,2}_d$ for every $m\geq0$. Smooth sections of $\SParEnd_0E\otimes K_X(D)$ and smooth $(0,1)$-forms with values in $\ParEnd_0E$ belong to $W^{m,2}_d$ for every $m\geq0$.
\item For every integer $m\geq0$, the local operator $\bar\partial_E$ admits a bounded right inverse from $W^{m,2}_d(\End E\otimes\overline K)$ to $W^{m+1,2}_d(\End E)$. We obtain a right inverse on a smaller punctured disc by a fixed cutoff. A fixed bounded extension across the finite boundary gives a right inverse on the half-cylinder itself without a boundary condition. Both constructions preserve the trace-free subspaces. The bounds depend on $m,d$, the local weights, and the chosen cutoff or extension.
\item A holomorphic endomorphism is bounded in $h_0$ near the puncture if and only if it extends in the parabolic endomorphism lattice. A holomorphic logarithmic one-form belongs to $L^2_d$ near the puncture if and only if it extends in the strongly parabolic logarithmic lattice. These assertions also hold in the trace-free subspaces.
\end{enumerate}
\end{lemma}

\begin{proof}
If a local endomorphism has entry $f_{ab}$, and a logarithmic one-form has entry $g_{ab}\,\dd z/z$, their holomorphic lattice conditions are
\begin{equation}\label{eq:tangent-lattice-orders}
 \begin{array}{c|cc}
  &\text{any holomorphic germ}
  &\text{divisible by }z\\ \hline
  \ParEnd_0E & s_{ab}\geq0 & s_{ab}<0\\
  \SParEnd_0E\otimes K_X(D) & s_{ab}>0 & s_{ab}\leq0.
 \end{array}
\end{equation}
The positive decay exponents obtained from these lattice conditions are bounded below by $\varepsilon_*$. Smooth sections of $\ParEnd_0E$ consequently have constant limits in $\mathfrak l_p$, with remainders decaying at least as $e^{-\varepsilon_*\tau}$. Smooth sections of $\SParEnd_0E\otimes K_X(D)$ have the same decay. A smooth $(0,1)$-form with values in $\ParEnd_0E$ acquires an additional factor $r$ from $\dd\bar z$ and also decays at this rate. The same estimates hold after every cylindrical derivative. Since $d<\varepsilon_*$, integration against $e^{d\tau}\,\dd\tau\,\dd\theta$ proves part~(i), including the inclusions into~\eqref{eq:tangent-weighted-spaces}.

For part~(ii), we solve $\bar\partial_E\eta=b$ entrywise.
The coefficient of one matrix entry in the frame unitary for $h_0$ is $u=r^{s_{ab}}\eta_{ab}$. In the logarithmic frame $\dd\bar z/\bar z$, the coefficient of $\bar\partial_E\eta$ is
\begin{equation}\label{eq:tangent-local-dbar}
 r^{s_{ab}}\bar z\,\partial_{\bar z}(r^{-s_{ab}}u)
   =-\frac12(\partial_\tau-\ii\partial_\theta+s_{ab})u.
\end{equation}
After conjugation by $e^{d\tau/2}$, the operator on the full cylinder has Fourier symbol
\begin{equation*}
 -\frac12(\ii\xi+k+s_{ab}-d/2),
 \qquad \xi\in\R,\quad k\in\Z.
\end{equation*}
Every $k+s_{ab}$ is zero or has absolute value at least $\varepsilon_*$. Hence this symbol has no zero, and its modulus is bounded below by a positive constant times $(1+\xi^2+k^2)^{1/2}$. The reciprocal symbol is bounded by a constant times $(1+\xi^2+k^2)^{-1/2}$, so its Fourier multiplier maps $W^{m,2}$ to $W^{m+1,2}$ for every integer $m\geq0$. After a fixed bounded extension of the conjugated source across the finite boundary, Fourier inversion and conjugation back give a bounded right inverse on the half-cylinder. One can also extend the source by zero after multiplying by a cutoff that equals one on a smaller punctured disc. This gives a local right inverse there. These constructions act on all endomorphisms. The diagonal entries use the same extension and Fourier inverse, so a trace-free right-hand side has a trace-free solution.

For part~(iii), boundedness of a holomorphic endomorphism $f$ in $h_0$ gives $|f_{ab}|\leq Cr^{-s_{ab}}$. A nonzero Laurent term $z^kE_{ab}$ is bounded precisely when $k+s_{ab}\geq0$, and Cauchy's coefficient estimate imposes this condition on every nonzero coefficient of $f_{ab}$. Since $-1<s_{ab}<1$, all nonzero terms have nonnegative degree. If $s_{ab}<0$, the constant coefficient also vanishes. Thus $f$ extends holomorphically and preserves the flag. Conversely, the first lattice in~\eqref{eq:tangent-lattice-orders} makes every such extension bounded in $h_0$. In particular, the entire constant block is allowed when $s_{ab}=0$.

For a holomorphic logarithmic one-form, the cylindrical norm of $\dd z/z$ is constant. A nonzero Laurent term $z^kE_{ab}\,\dd z/z$ is therefore in $L^2_d$ precisely when
\begin{equation*}
 \int_{\tau_0}^\infty e^{(d-2(k+s_{ab}))\tau}\,\dd\tau<\infty\,.
\end{equation*}
Equivalently,
\begin{equation}\label{eq:tangent-Laurent-integrability}
 k+s_{ab}>d/2.
\end{equation}
By angular orthogonality, every nonzero Laurent coefficient satisfies this condition. It follows that $k<0$ is impossible and that $k=0$ is allowed precisely when $s_{ab}>0$, since every positive $k+s_{ab}$ is at least $\varepsilon_*>d$. Conversely, these lattice conditions give the weighted integrability on a smaller disc. This proves the second extension assertion.
\end{proof}

\begin{proof}[Proof of Proposition~\ref{prop:weighted-Dolbeault-comparison}]
For the diagonal metric, the smooth Dolbeault complex embeds in the weighted complex by Lemma~\ref{lem:weighted-local-dolbeault}(i). The equivalence of norms following~\eqref{eq:tangent-weight-window} gives the same inclusion for $h$.

For $(b,\varphi)\in\ker \dd_2$, we cover $X$ by coordinate discs, one about each puncture, with all overlaps disjoint from $D$. By Lemma~\ref{lem:weighted-local-dolbeault}(ii) on the puncture discs and the local Dolbeault lemma elsewhere, the equation $\bar\partial_E\eta_i=b$ has a solution on each disc. Then $\varphi_i:=\varphi-[\Phi,\eta_i]$ is holomorphic. On a puncture disc it remains in $W^{1,2}_d$ and extends in the second lattice by Lemma~\ref{lem:weighted-local-dolbeault}(iii). On an overlap, $\eta_i-\eta_j$ is holomorphic. We choose a smooth partition of unity $\rho_i$, with each puncture function equal to one on a smaller disc, and put $\eta=\sum_i\rho_i\eta_i$. Subtracting $\dd_1\eta$ gives
\begin{equation*}
 b'=-\sum_i\bar\partial\rho_i\,\eta_i,
 \qquad \varphi'=\sum_i\rho_i\varphi_i.
\end{equation*}
The first expression is smooth because we may subtract any one $\eta_j$ locally from all its terms, using $\sum_i\bar\partial\rho_i=0$. It vanishes near $D$. The second is smooth in the logarithmic lattice and holomorphic near $D$. Thus $(b',\varphi')$ is a smooth Dolbeault cocycle, proving surjectivity and the existence of a representative of this form.

For injectivity we consider a smooth Dolbeault cocycle equal to $\dd_1\eta$ with $\eta\in\mathcal R^0_d$. On a puncture disc, there is a smooth section $u$ satisfying $\bar\partial_Eu=b$, by the Dolbeault lemma in the locally free sheaf $\ParEnd_0E$. Both $u$ and $\eta$ are bounded in the diagonal metric; for $\eta$ this follows from weighted Sobolev embedding and its constant Levi limits. By Lemma~\ref{lem:weighted-local-dolbeault}(iii), the holomorphic difference $\eta-u$ extends as a section of $\ParEnd_0E$. Away from $D$ the section is smooth by interior elliptic regularity. Hence $\eta$ is a global smooth section of the parabolic endomorphism sheaf, and the cocycle is already exact in its smooth Dolbeault resolution. The Dolbeault resolution computes $\mathbb H^1(X,\mathcal C^\bullet)$, proving~\eqref{eq:weighted-Dolbeault-comparison}.
\end{proof}

\subsection{The gauge Laplacian and harmonic projection}\label{proof:prop:ordinary-L2-projection}
\begin{proof}
In the proof of Proposition~\ref{prop:ordinary-L2-projection} we first establish boundedness on the weighted Sobolev domain. The leading matrices in~\eqref{eq:harmonic-unitary-end} commute with every element of $\mathfrak l_p$ and with its adjoint. This commutation proves the bound for~\eqref{eq:tangent-gauge-isomorphism} on the constant summands as well. The indicial exponents $\mu_{abk}$ are those computed in Lemma~\ref{lem:cylinder-linear-dirichlet}, equation~\eqref{eq:harmonic-cylinder-frequencies}: on the Fourier component $e^{\ii k\theta}E_{ab}$, the normal operator of $L_h$ is
\begin{equation}\label{eq:tangent-normal-operator}
 \frac12\left(-\partial_\tau^2+\mu_{abk}^2\right).
\end{equation}
Here the factor $1/2$ agrees with the $L^2$ adjoint of $\bar\partial_E$. The indicial exponents $\mu_{abk}$ vanish exactly when $k=0$ and $a,b$ belong to the same weight block, since the graded residue on that block is scalar. The corresponding trace-free constant matrices form $\mathfrak l_p$. Every positive indicial root is at least $\varepsilon_*$, also when residues on different blocks coincide.

For $u\in[0,1]$, replacing $\Phi$ by $u\Phi$ defines an operator $L_{h,u}$. The frequency $\mu_{abk}(u)$ of its normal operator is given by~\eqref{eq:harmonic-cylinder-frequencies} with each $\sigma_a$ replaced by $u\sigma_a$. The same lower bound for the nonzero frequencies holds throughout this interval. Conjugating the normal operator by $e^{d\tau/2}$ gives the symbol
\begin{equation*}
 \frac12\left(- (\ii\xi-d/2)^2
       +\mu_{abk}(u)^2\right).
\end{equation*}
Its modulus is bounded below by a positive multiple of $1+\xi^2+k^2$, uniformly in $u\in[0,1]$ for these fixed data. Indeed, a zero requires $\xi=0$ and an indicial root equal to $d/2$, contrary to~\eqref{eq:tangent-weight-window}. The bound holds directly at high Fourier frequencies and by compactness at the remaining frequencies. The inverses on full cylinders give parametrices on the ends. Patching them to interior elliptic parametrices gives Fredholm operators on $W^{2,2}_d$: cutoff commutators have compact support and lower order, and the coefficient errors are compact because they have order at most one and decay with their derivatives. Appending the finite-dimensional constant summands preserves Fredholmness. The operators $L_{h,u}$ form a norm-continuous Fredholm family on the same domain.

At $u=0$, the index is zero by~\cite[Proposition~3.24]{CollierFredricksonWentworth}. The domain of the operator at this parameter agrees with~\eqref{eq:tangent-weighted-spaces}, since the boundary map in~\cite[Proposition~3.3]{CollierFredricksonWentworth} has kernel $W^{2,2}_d$ and image the full Levi algebra. To use the weight interval $[0,1)$ of that index formula, we subtract the minimum weight at each puncture. This leaves the induced model connection on $\End_0E$ unchanged. Since the zero Higgs field is strongly parabolic, the partial flags here satisfy the hypotheses of the index formula. Constancy of the index proves that $L_h$ has index zero.

For integration by parts on the extended domain, we use a cutoff $\kappa_R$ equal to one on the compact part and on $\tau\leq R$, and zero on $\tau\geq R+1$, with uniformly bounded derivative. Every $\zeta\in\mathcal R^0_d$ is bounded, and therefore
\begin{equation}\label{eq:tangent-cutoff-bound}
 \|\zeta\,\dd\kappa_R\|_{L^2_{-d}}
                    \leq C_\zeta e^{-d R/2}.
\end{equation}
For $\eta\in\mathcal R^0_d$, boundedness of $\eta$ and the inclusion $\dd_1\eta\in L^2_d$ allow the cutoff argument in the proof of~\cite[Proposition~3.24, p.~26]{CollierFredricksonWentworth}, which gives $\ker L_h=\ker \dd_1$. An element of this kernel is holomorphic on $X^*$, bounded in the diagonal metric, and commutes with $\Phi$. By Lemma~\ref{lem:weighted-local-dolbeault}(iii), it extends as a trace-free parabolic Higgs endomorphism on $X$. It follows from Proposition~\ref{prop:local-deformations} that $\eta=0$. Since the index is zero,~\eqref{eq:tangent-gauge-isomorphism} follows by the bounded inverse theorem.

Since $\dd_1^*:\mathcal V_d\to L^2_d$ is bounded, we may take $\eta=-L_h^{-1}\dd_1^*v$. We obtain~\eqref{eq:tangent-coclosed-projection} from $\dd_2\dd_1=0$ and~\eqref{eq:tangent-projection-bound} from the bound on the inverse. Uniqueness follows from the trivial kernel of $L_h$. Since the weight $d$ is positive, a section in the weighted space is $L^2$ integrable. Finally, we integrate the pairing of $v^h$ with $\dd_1(\kappa_R\zeta)$ for $\zeta\in\mathcal R^0_d$. The equation $\dd_1^*v^h=0$ makes this pairing zero;~\eqref{eq:tangent-cutoff-bound} removes the cutoff term as $R\to\infty$. Hence $v^h$ is orthogonal to every $\dd_1\zeta$.
\end{proof}

\subsection{Local quotients and families of harmonic metrics}\label{proof:prop:smooth-harmonic-quotient}

For $\alpha\in\Acal$, Dimakis--Fredrickson--Yae constructed the local smooth quotient and proved smooth dependence of harmonic representatives on stable algebraic families~\cite[Propositions~4.6--4.7, Claim~4.6.1 and Section~4.3]{DFY}. Their weighted norm~\cite[equation~(4)]{DFY} is equivalent to~\eqref{eq:harmonic-weighted-norm} when their decay exponent is $d/2$; the Sobolev orders for configurations and gauge transformations here are $j+1$ and $j+2$. We record the linear operators, normalizations and derivative estimates used below, and prove a common properness estimate for the models in Sections~\ref{sec:hyperkahler-structure} and~\ref{sec:boundary}.

\begin{lemma}\label{lem:local-coulomb-slice}
Our base point is a stable solution $x=(D_h,\Phi)$ with weights $\alpha\in\Acal$ and scalar graded residues $\sigma$. We take $d$ as in~\eqref{eq:tangent-weight-window}. For $j\geq2$, the Hitchin solutions near $x$ satisfying the real Coulomb condition form a smooth slice with tangent space $\mathcal H_{h,d}$. Nearby configurations admit a unique sufficiently small unitary gauge transformation into Coulomb gauge. The quaternionic splitting and moment maps for the pairing~\eqref{eq:tangent-trial-norm} are given by~\eqref{eq:hk-four-operator-inverse} and~\eqref{eq:hk-moment-maps}.
\end{lemma}
\begin{proof}
By Proposition~\ref{prop:ordinary-L2-projection} and elliptic regularity, we have a bounded isomorphism
\begin{equation}\label{eq:hk-higher-inverse}
 L_h=\dd_1^*\dd_1:\mathcal R^j_d
                   \longrightarrow W^{j,2}_d(\End_0E).
\end{equation}
After subtracting the cutoff Levi constants from the solution in $\mathcal R^0_d$, we apply the second-order elliptic estimate on uniformly sized cylindrical charts. Conjugation by $e^{d\tau/2}$ and summation give
\begin{equation*}
 \|\eta\|_{\mathcal R^j_d}
 \leq C_j\bigl(\|L_h\eta\|_{W^{j,2}_d}
                       +\|\eta\|_{\mathcal R^0_d}\bigr).
\end{equation*}
The coefficients and all their derivatives are bounded by Proposition~\ref{prop:harmonic-weighted-form}, and $L_h(\chi_pc_p)$ decays in every order. The inverse estimate in the lower Sobolev order now gives~\eqref{eq:hk-higher-inverse}.

At a solution, $L_h$ commutes with the Hermitian adjoint, and hence restricts to real isomorphisms on $\mathcal U^j_d$ and $\mathcal S^j_d$. This is the linear identity used in~\cite[Theorem~3.33, pp.~30--31]{CollierFredricksonWentworth}. For the slice calculation, we express the surface metric as $\rho^2|\dd z|^2$, put $\Phi=\phi\,\dd z$, and take $v=(B\,\dd\bar z,P\,\dd z)$. Then
\begin{equation}\label{eq:hk-adjoint-identities}
 \begin{aligned}
 \dd_1^*v&=2\rho^{-2}\bigl(-\nabla_zB+[\phi^\dagger,P]\bigr),\\
 \dd_1^*(Jv)&=2\ii C^\dagger,\qquad
 C=\rho^{-2}\bigl(-\nabla_{\bar z}P+[\phi,B]\bigr),
 \end{aligned}
\end{equation}
The maps in~\eqref{eq:hk-quaternionic-operators} are isometries for the real polarization of~\eqref{eq:tangent-trial-norm} and satisfy the quaternion identities.

We write $\dd_0=\dd_1|_{\mathcal U^j_d}$ and
$\dd_0^*v=\operatorname{skew}(\dd_1^*v)$, where
$\operatorname{skew}(T)=(T-T^\dagger)/2$. We set
$Q_0=\Id$, $Q_1=I$, $Q_2=J$ and $Q_3=K$. Since $\dd_2\dd_1=0$, equations~\eqref{eq:hk-adjoint-identities} give
\begin{equation*}
 \dd_0^*\dd_0=L_h,\qquad
 \dd_0^*I\dd_0=\dd_0^*J\dd_0=\dd_0^*K\dd_0=0.
\end{equation*}
Consequently the operator and its right inverse
\begin{equation}\label{eq:hk-four-operator-inverse}
 \begin{aligned}
 \mathscr Dv&=(\dd_0^*Q_a^{-1}v)_{a=0}^3,\\
 \mathcal Q(f_0,f_1,f_2,f_3)
   &=\sum_{a=0}^3Q_a \dd_0L_h^{-1}f_a
 \end{aligned}
\end{equation}
are bounded between $\mathcal V^j_d$ and
$\bigl(W^{j,2}_d(\mathfrak{su}(E,h))\bigr)^4$, with
$\mathscr D\mathcal Q=\Id$. The equations with $a=0,1$ are equivalent to $\dd_1^*v=0$; those with $a=2,3$ are equivalent to $\dd_2v=0$ by~\eqref{eq:hk-adjoint-identities}. The kernel is
$\mathcal H_{h,d}$ from~\eqref{eq:tangent-harmonic-space}. Its elements belong to $\mathcal V^j_d$ for every $j\geq0$ by elliptic regularity. Thus
\begin{equation*}
 \mathcal V^j_d=\mathcal H_{h,d}
                           \oplus\operatorname{im}\mathcal Q.
\end{equation*}

We apply this splitting to the equations and the real Coulomb condition. The constant forms on the configuration space
$\omega_a(v,w)=g(Q_av,w)$, $a=1,2,3$, where $g$ is the real polarization of~\eqref{eq:tangent-trial-norm}, have moment maps $\mu_a$ satisfying
\begin{equation}\label{eq:hk-moment-derivative}
 \dd\mu_a|_x(v)=\dd_0^*Q_a^{-1}v.
\end{equation}
Here functions are paired by
$2\operatorname{Re}\int_{X^*}\tr(\xi\eta^\dagger)\,\dd\operatorname{vol}_X$,
where $\dd\operatorname{vol}_X$ is the volume form of the fixed surface metric.
For a general nearby pair, we set
\begin{equation*}
 R=\rho^{-2}\bigl(F_{z\bar z}+[\phi,\phi^\dagger]\bigr),
 \qquad C=-\rho^{-2}\nabla_{\bar z}\phi.
\end{equation*}
The moment maps are
\begin{equation}\label{eq:hk-moment-maps}
 \mu_1=\ii R,\qquad
 \mu_2=-\operatorname{skew}(2\ii C^\dagger),\qquad
 \mu_3=\operatorname{skew}(2C^\dagger).
\end{equation}
The common zero set of these moment maps consists of the solutions of Hitchin's equations. By~\eqref{eq:hk-adjoint-identities}, their pairings with $\xi$ have derivatives $\omega_a(\dd_0\xi,v)$. The integration by parts is valid for every bounded Levi limit of $\xi$ by~\eqref{eq:tangent-cutoff-bound} and~\eqref{eq:tangent-trace-boundary}. The model matrices $A,B$ in~\eqref{eq:harmonic-local-data} are scalar on each grade and commute with all these limits, so each residual belongs to $W^{j,2}_d$.

The smooth map
\begin{equation*}
 v\longmapsto
 \bigl(\dd_0^*v,\mu_1(x+v),\mu_2(x+v),\mu_3(x+v)\bigr)
\end{equation*}
has derivative $\mathscr D$ at zero. By the implicit function theorem, its zero set is a smooth graph over a neighborhood of zero in $\mathcal H_{h,d}$. The function whose graph is this zero set has derivative zero at the origin. This graph is a real four-dimensional slice in the solution space. The equation
$\dd_0^*(\exp(\xi)^*(x+v)-x)=0$ has derivative $L_h$ in $\xi$. The real isomorphism~\eqref{eq:hk-higher-inverse} then gives a unique gauge transformation near the identity and a smooth projection to the slice.
\end{proof}

The next estimate treats unitary gauges on all matrix blocks, including scalar matrices. It will also apply when equal weights have distinct residue eigenvalues.

\begin{lemma}\label{lem:gauge-local-properness}
The model connection and residue are $D_0=\dd+\ii A\,\dd\theta$ and $B$, where $A=\diag(\beta_a)$ and $B=\diag(\sigma_a)$ commute and $\max\beta_a-\min\beta_a<1$. For $j\geq2$ and $d>0$, we consider any two sequences of configurations with these models whose remainders converge in $W^{j+1,2}_d$. If unitary gauge transformations $g_\ell$ of determinant one carry the first sequence to the second, with limiting matrices commuting with $A,B$ and remainders in $W^{j+2,2}_d$, then $(g_\ell)$ has a subsequence whose limiting matrices converge at every puncture and whose remainders converge in $W^{j+2,2}_d$.
\end{lemma}
\begin{proof}
For an $\End E$-valued section $u$ with limiting constant $c$ commuting with $A,B$, Fourier decomposition gives
\begin{equation}\label{eq:joint-gauge-tail}
 \|u-\chi c\|_{W^{j+1,2}_d}
 \leq C\bigl(\|D_0u\|_{W^{j,2}_d}+\|[B,u]\|_{W^{j,2}_d}
              +\|u\|_{W^{j+1,2}(\mathcal A_0)}+|c|\bigr),
\end{equation}
where $\chi$ equals one sufficiently far along the end and $\mathcal A_0$ is a fixed finite annulus. On the angular matrix components outside the joint kernel, the sum of the squared connection and residue symbols is
\begin{equation*}
 (k+\beta_a-\beta_b)^2+|\sigma_a-\sigma_b|^2.
\end{equation*}
It is bounded below by a positive multiple of $1+k^2$: only finitely many $k$ require a separate check, and every parallel matrix component is constant because the weights lie in an interval of width less than one. The radial derivative controls the remaining derivatives. On the joint kernel, $u-c=-\int_\tau^\infty\partial_su\,\dd s$; conjugation by $e^{d\tau/2}$ gives an integral kernel $e^{-d(s-\tau)/2}$, which is integrable for $s\geq\tau$. These estimates, a cutoff on $\mathcal A_0$, and differentiation prove~\eqref{eq:joint-gauge-tail} in every derivative order. The scalar component belongs to this same kernel.

Unitarity bounds $g_\ell$ pointwise. We extend $D_0$ over the compact part and denote the connection and Higgs remainders by $a_\ell,b_\ell$ and $a'_\ell,b'_\ell$. The intertwining equations are
\begin{equation}\label{eq:hk-gauge-equation}
 D_0g_\ell=g_\ell a'_\ell-a_\ell g_\ell,
 \qquad
 [B,g_\ell]=g_\ell b'_\ell-b_\ell g_\ell.
\end{equation}
On compact subsets, the first equation gives bounded $W^{j+2,2}$ norms and a convergent subsequence in that norm, by compact embedding followed by the equation. The limiting unitary matrices lie in compact groups, so they also have a convergent subsequence.

On sufficiently far ends the coefficient remainders are uniformly small. The translated version of~\eqref{eq:joint-gauge-tail} absorbs their contributions and bounds the gauges in $W^{j+1,2}_d$ modulo their limits. Since the coefficients lie in $W^{j+1,2}_d$ and $j\geq2$, the Sobolev product estimates and~\eqref{eq:hk-gauge-equation} give one additional derivative. Applying~\eqref{eq:joint-gauge-tail} at this order bounds the $W^{j+2,2}_d$ remainders. The same estimates applied to differences give convergence in the gauge topology: the coefficients converge in their weighted norms, the gauges converge on a compact annulus, and their limiting constants converge.
\end{proof}

\begin{lemma}\label{lem:harmonic-family-dependence}
Over an analytic polydisc $S$ centered at zero, we consider a holomorphic family of stable parabolic Higgs bundles with a trivialized determinant line bundle, fixed weights $\alpha\in\Acal$ and fixed scalar graded residues $\sigma$. For a weight $d$ satisfying~\eqref{eq:tangent-weight-window}, we may shrink $S$ and choose smooth bundle isomorphisms with the central fiber that are holomorphic near the punctures and preserve the flags and the determinant trivialization. If $h$ is the harmonic metric at zero, we write the harmonic metrics as $h_u=he^{s(u)}$. Then $u\mapsto s(u)$ is smooth into $\mathcal S^j_d$ for every $j\geq0$. The induced map to a local solution slice has derivative at zero equal to the Kodaira--Spencer map followed by $L^2$ harmonic projection for $h$. If the Kodaira--Spencer map at zero is an isomorphism, this map is a local diffeomorphism and carries the complex structure of $S$ to $I$.
\end{lemma}
\begin{proof}
We choose holomorphic frames on puncture discs times $S$ compatible with the flags and the fixed determinant trivialization. The resulting bundle isomorphisms with the central fiber extend smoothly over the compact remainder of $X$. For $u\in S$, this gives $\bar\partial_u=\bar\partial_0$ near $D$, while $\Phi_u-\Phi_0$ has zero graded residue and depends smoothly on $u$ in the weighted Sobolev spaces of every order by Lemma~\ref{lem:weighted-local-dolbeault}(i). With $h$ the harmonic metric at $u=0$, we define the harmonic metric operator on $s\in\mathcal S^j_d$ by
\begin{equation}\label{eq:hk-metric-equation}
 \mathcal N(u,s)=e^{s/2}\ii\Lambda_{\mathrm{cyl}}
 \bigl(F_{\bar\partial_u,he^s}
       +[\Phi_u,\Phi_u^{\dagger_{he^s}}]\bigr)e^{-s/2}.
\end{equation}
The conjugation makes $\mathcal N(u,s)$ take values in $W^{j,2}_d(\Herm_0(E,h))$, and the algebra~\eqref{eq:hk-ambient-algebra} gives smoothness. Every limiting Levi block commutes with $A,B$ in~\eqref{eq:harmonic-local-data}, so each remaining term contains a derivative along the punctured curve or a remainder that decays there. As in~\cite[Theorem~3.33]{CollierFredricksonWentworth},
$D_s\mathcal N(0,0)=L_h$ on $\mathcal S^j_d$. By the real implicit function theorem, a solution $s(u)$ exists, and uniqueness implies that $he^{s(u)}$ is the harmonic metric of the given family.

The positive isometry $e^{s(u)/2}$ and the real slice projection give a smooth map from $S$ to the solution slice. Its derivative at zero sends a Dolbeault deformation to its $L^2$ harmonic representative for $h$: the variation of the metric contributes an infinitesimal complex gauge transformation, and the slice equations impose $\dd_2v=\dd_1^*v=0$. If the Kodaira--Spencer map at zero is an isomorphism, the derivative is an isomorphism by Proposition~\ref{prop:weighted-Dolbeault-comparison}. The map is then a local diffeomorphism by the finite-dimensional inverse function theorem. Repeating the tangent calculation at each parameter shows that the map carries multiplication by $\ii$ on deformation cohomology to $I$ in~\eqref{eq:hk-quaternionic-operators}.

Finally, at each higher Sobolev order we obtain the same metrics by the implicit function theorem and uniqueness. On a fixed smaller parameter chart, we differentiate~\eqref{eq:hk-metric-equation} in the parameters and apply the estimates giving~\eqref{eq:hk-higher-inverse}. These estimates show that every parameter derivative of finite order, after subtraction of its smooth Levi limit, belongs to weighted Sobolev spaces of every order. By Sobolev embedding these remainders decay at any fixed rate below $d/2$, uniformly on compact subsets of the parameter chart. The gauge-fixing equations give the same conclusion for transition gauges. All derivatives here are taken with $(\alpha,\sigma)$ fixed.
\end{proof}

\begin{proof}[Proof of Proposition~\ref{prop:smooth-harmonic-quotient}]
For gauges relating points in the solution slices of Lemma~\ref{lem:local-coulomb-slice}, we have the properness estimate of Lemma~\ref{lem:gauge-local-properness} because each graded residue is scalar and the limiting groups are $K_p$. A stabilizer extends to a parabolic Higgs automorphism by Lemma~\ref{lem:weighted-local-dolbeault}; stability and the determinant condition make it $\mu_n$. Local properness and the uniqueness of the gauge near the identity therefore give locally trivial principal $\mathcal G_d/\mu_n$ charts.

We compare these charts with the algebraic moduli space. The Fourier inverse of Lemma~\ref{lem:weighted-local-dolbeault}(ii) acts in every Sobolev order on all endomorphisms, including their scalar part. A bounded extension across the finite boundary of a half-cylinder gives a right inverse $T_0$ on that half-cylinder. On a sufficiently far half-cylinder, the perturbation in $\bar\partial_0+b$ is small in the translated weighted norm. We solve
\begin{equation*}
 u=-T_0\bigl(b(\Id+u)\bigr)
\end{equation*}
in $W^{j+2,2}_d(\End E)$. The contraction yields a complex gauge $\Id+u$ tending to $\Id$ and satisfying
\begin{equation}\label{eq:hk-holomorphic-frame-gauge}
 (\bar\partial_0+b)\circ(\Id+u)
       =(\Id+u)\circ\bar\partial_0.
\end{equation}
It therefore carries a model holomorphic frame to a holomorphic frame for $\bar\partial_0+b$. Since $\tr b=0$, the determinant satisfies $\bar\partial\det(\Id+u)=-(\tr b)\det(\Id+u)=0$. We divide the change-of-frame matrix by the holomorphic $n$th root of its determinant that tends to one. The resulting matrix has determinant one. The Laurent test~\eqref{eq:tangent-Laurent-integrability} then extends the Higgs field with its prescribed flag and scalar graded residue. The construction is smooth in the decaying perturbation. Bounded holomorphic isomorphisms extend by Lemma~\ref{lem:weighted-local-dolbeault}(iii). Conversely, an isomorphism of stable filtered objects is unitary for their harmonic metrics by uniqueness in Proposition~\ref{prop:harmonic-metric-existence}. This gauge belongs to $\mathcal G_d$ by Proposition~\ref{prop:harmonic-weighted-form}. This proves the asserted bijection on points.

The local atlas~\eqref{eq:stable-local-family} carries a family through the stable pair $x$, with determinant trivialization and scalar graded residues fixed and with Kodaira--Spencer map an isomorphism. We restrict this family to an analytic polydisc $S$ centered at $x$.

For this stable family, the hypotheses of Lemma~\ref{lem:harmonic-family-dependence} hold. Since its Kodaira--Spencer map is an isomorphism, the algebraic and analytic charts are smoothly equivalent and their complex structures agree. The estimates at the end of that lemma's proof give smooth block limits, all finite mixed derivatives and the same decay for transition gauges, proving the second assertion.
\end{proof}


\end{document}